\input xy
\xyoption{all}

\documentclass{conm-p-l}
   \usepackage{amsthm}
   \usepackage{amsmath}
   \usepackage{amssymb}
   \usepackage{amsfonts}
   \usepackage{amscd}
   \usepackage[mathscr]{eucal}
   \usepackage{verbatim}
   \usepackage{graphics}

\newtheorem*{main}{Main~Theorem}

\newtheorem{athm}{Theorem}[subsection]
\newtheorem{alem}[athm]{Lemma}
\newtheorem{aprop}[athm]{Proposition}
\newtheorem{acor}[athm]{Corollary}

\theoremstyle{definition}
  \newtheorem{arem}[athm]{Remark}
  
 \newtheorem{adefi}[athm]{Definition}
 \newtheorem{aexam}[athm]{Example}
 \newtheorem{aexams}[athm]{Examples}

\def\bilap#1{\hbox to 0pt{\hss#1\hss}}
 \def\Rarrow#1{\bilap{\hbox to#1{\rightarrowfill}}}
 \def\Larrow#1{\bilap{\hbox to#1{\leftarrowfill}}}
\def\Equals#1{\bilap
                  {\hbox{\rule[3.5pt]{#1} {.5pt}}
                   \kern-#1
                   \hbox{\rule[1pt]{#1}{.5pt}}%
                 }}

\def\UnderElement#1#2#3#4{\vbox to 0pt{
\hbox{$
\llap{$\scriptstyle#1$}
\left#2\vbox to #3{}\right.
\rlap{$\scriptstyle#4$}
     $}
\vss}}
\newcommand{\under}[2]
{\vbox to 0pt{\vskip-#1 ex\hbox{$\scriptstyle #2$}\vss}}

\newcommand{\EQAL}[1]%
{\,\begin{picture}(#1,0)%
\put(0,3){\line(1,0){#1}}%
\put(0,1){\line(1,0){#1}}%
\end{picture}\,}%

\newcommand{\vlto}[1]%
{\,\begin{picture}(#1,3)%
\put(0,2){\vector(1,0){#1}}%
\end{picture}\,}%

\newcommand{\vllarrow}[1]%
{\,\begin{picture}(#1,3)%
\put(#1,2){\vector(-1,0){#1}}%
\end{picture}\,}%

\newcommand{\dirlm}[1]%
  {
     {\lim\hskip-1.58em\lower.65ex
       \hbox{$
                {}_{\stackrel{\lower1ex\hbox
                                        {$\scriptstyle -\!\!\<\longrightarrow$}
                             }{ ^{#1} }
                   }
            $}
     }
\:}

\newcommand{\subdirlm}[1]%
  {
     {\lim\hskip-1.5em\lower.6ex
       \hbox{$
                   {}_{\stackrel{\lower1ex\hbox
                                           {$\scriptstyle\longrightarrow$}
                                }{ ^{#1} }
                      }
             $}
     }
\:}

\newcommand{\inlm}[1]%
   {
      {\lim\hskip-1.58em\lower.65ex
        \hbox{$
                 {}_{\stackrel{\lower1ex\hbox
                                        {$\scriptstyle \longleftarrow\!\!\<-$}
                              }{ ^{#1} }
                    }
             $}
      }
\:}

\def\Iso{\vbox to 0pt{\vss\hbox{$\widetilde{\phantom{nn}}$}\vskip-7pt}}
\def\>{\mspace {1mu}}
\def\<{\mspace{-1mu}}
\def\({{\textup(}}
\def\){{\textup)}}
\def\bigl#1{{\textup{\begin{large}#1\end{large}}}}
\def\bigr#1{{\textup{\begin{large}#1\end{large}}}}

\newcommand{\X}{{\mathscr X}}

\newcommand{\Y}{{\mathscr Y}}
\newcommand{\Z}{{\mathscr Z}}
\newcommand{\V}{{\mathscr V}}
\newcommand{\U}{{\mathscr U}}
\newcommand{\W}{{\mathscr W}}
\newcommand{\I}{{\mathscr I}}
\newcommand{\J}{{\mathscr J}}
\newcommand{\sL}{{\mathscr L}}
\newcommand{\sK}{{\mathscr K}}
\newcommand{\sB}{{\mathscr B}}

\newcommand{\Ufr}{{\mathfrak U}}
\newcommand{\Cfr}{{\mathfrak C}}

\newcommand{\A}{{\mathcal A}}

\newcommand{\C}{{\mathcal C}}
\newcommand{\cD}{{\mathcal D}}
\newcommand{\E}{{\mathcal E}}
\newcommand{\F}{{\mathcal F}}
\newcommand{\G}{{\mathcal G}}
\renewcommand{\H}{{\mathcal H}}
\newcommand{\Hr}{{\mathrm H}}

\newcommand{\M}{{\mathcal M}}
\newcommand{\N}{{\mathcal N}}
\newcommand{\Esf}{{\mathsf E}}
\newcommand{\Ssf}{{\mathsf S}}
\newcommand{\ssf}{{\mathsf s}}
\newcommand{\cL}{{\mathcal L}}
\newcommand{\cI}{{\mathcal I}}
\newcommand{\cJ}{{\mathcal J}}
\newcommand{\cK}{{\mathcal K}}

\newcommand{\cCb}{{\mathcal C}^{\bullet}}
\newcommand{\eCb}{{\mathscr C}^{\bullet}}

\newcommand{\cDb}{{\mathcal D}^{\bullet}}
\newcommand{\cEb}{{\mathcal E}^{\bullet}}
\newcommand{\cFb}{{\mathcal F}^{\>\bullet}}
\newcommand{\cGb}{{\mathcal G}^{\bullet}}
\newcommand{\cIb}{{\mathcal I}^{\bullet}}

\newcommand{\cKb}{{\mathcal K}^{\bullet}}
\newcommand{\cLb}{{\mathcal L}^{\bullet}}
\newcommand{\cMb}{{\mathcal M}^{\bullet}}
\newcommand{\cNb}{{\mathcal N}^{\bullet}}
\newcommand{\cPb}{{\mathcal P}^{\bullet}}
\newcommand{\cRb}{{\mathcal R}^{\bullet}}
\newcommand{\Ab}{A^{\bullet}}
\newcommand{\Bb}{B^{\bullet}}
\newcommand{\Cb}{C^{\bullet}}
\newcommand{\Ccb}{{\check{C}}^{\bullet}}
\newcommand{\Eb}{E^{\bullet}}
\newcommand{\Fb}{F^{\bullet}}
\newcommand{\Fbp}{F^{\prime\bullet}}
\newcommand{\Gb}{G^{\bullet}}
\newcommand{\Ib}{I^{\bullet}}

\newcommand{\Ibsim}{I^{\bullet \sim}}
\newcommand{\Jb}{J^{\bullet}}
\newcommand{\Kbi}{K^{\bullet}_{\infty}}
\newcommand{\Kb}{K^{\bullet}}
\newcommand{\Lb}{L^{\bullet}}
\newcommand{\Mb}{M^{\bullet}}
\newcommand{\Nb}{N^{\bullet}}

\newcommand{\cFbp}{\F^{\prime\bullet}}

\newcommand{\cKbi}{{\mathcal K}^{\bullet}_{\infty}}
\newcommand{\cKi}{{\mathcal K}_{\infty}}

\newcommand{\cO}{{\mathcal O}}
\newcommand{\OX}{\cO_{\X}}
\newcommand{\OY}{\cO_{\Y}}
\newcommand{\OU}{\cO_{\U}}
\newcommand{\OV}{\cO_{\V}}
\newcommand{\OW}{\cO_{\W}}
\newcommand{\OZ}{\cO_{\Z}}
\newcommand{\OXx}{\cO_{\X\<,\>x}}
\newcommand{\OYy}{\cO_{\Y\<,\>y}}
\newcommand{\OUx}{\cO_{\U\<,\>x}}
\newcommand{\OVy}{\cO_{\V\<,\>y}}
\newcommand{\OZz}{\cO_{\Z\<,\>z}}
\newcommand{\OWw}{\cO_{\W\<,\>w}}
\newcommand{\XOX}{(\X,\cO_{\X})}
\newcommand{\YOY}{(\Y,\cO_{\Y})}
\newcommand{\Ox}{\cO_{x}}
\newcommand{\Oxp}{\cO_{x'}}

\newcommand{\bbFc}{{\mathbb F}_{\mkern-1.5mu\rm c}}
\newcommand{\Coz}{{\mathrm {Coz}}}
\newcommand{\Cou}{{\mathrm {Cou}}}
\newcommand{\Mod}{{\mathrm {Mod}}}
\newcommand{\Spec}{{\mathrm {Spec}}}
\newcommand{\Spf}{{\mathrm {Spf}}}
\newcommand{\Supp}{{\mathrm {Supp}}}
\newcommand{\Der}{{\mathrm {Der}}}

\newcommand{\Forget}{{\texttt{F\small{gt}}}}

\newcommand{\hit}{{\mathrm {ht.}}}
\newcommand{\cont}{{\mathrm {c}}}
\newcommand{\cR}{{\mathcal R}}
\newcommand{\Rc}{{\cRb_{\textup{c}}}}
\newcommand{\Ohm}{{\widehat{\Omega}}}

\newcommand{\bC}{{\mathbf C}}
\newcommand{\D}{{\mathbf D}}
\newcommand{\K}{{\mathbf K}}

\newcommand{\Dqc}{\D_{\mkern-1.5mu\mathrm {qc}}}

\newcommand{\Dqct}{\D_{\mkern-1.5mu\mathrm{qct}}}

\newcommand{\Dc}{\D_{\mkern-1.5mu\mathrm c}}

\newcommand{\qc}{{\mathrm{qc}}}

\newcommand{\Bf}{{\bf f}}
\newcommand{\bg}{{\bf g}}
\newcommand{\bh}{{\bf h}}
\newcommand{\bt}{{\bf t}}
\newcommand{\bs}{{\bf s}}
\newcommand{\bu}{{\bf u}}

\newcommand{\afr}{\mathfrak a}
\newcommand{\bfr}{\mathfrak b}
\newcommand{\cfr}{\mathfrak c}

\newcommand{\mfr}{\mathfrak m}
\newcommand{\nfr}{\mathfrak n}
\newcommand{\pfr}{\mathfrak p}
\newcommand{\qfr}{\mathfrak q}

\newcommand{\R}{{\mathbf R}}

\newcommand{\bL}{{\mathbf L}}
\newcommand{\Hom}{{\mathrm {Hom}}}
\newcommand{\Homc}{{\mathrm {Hom}}^{\cont}}
\newcommand{\Homb}{{\mathrm {Hom}}^{\bullet}}
\newcommand{\Ac}{\A_{\mathrm c}}
\newcommand{\Aqc}{\A_{\qc}}
\newcommand{\Avc}{\A_{\vec {\mathrm c}}}
\newcommand{\Aqct}{\A_{\mathrm {qct}}\<}
\newcommand{\Act}{\A_{\mathrm {ct}}\<}

\newcommand{\At}{\A_{\mathrm t}\<}
\newcommand{\Dt}{\D_{\mathrm t}\<}

\newcommand{\sHom}{\H{om}}
\newcommand{\sHomb}{\H{om}^{\bullet}}
\newcommand{\iGp}[1]{{\varGamma_{\<\!#1}'}}
\newcommand{\iG}[1]{{\varGamma_{\<\!#1}^{\phantom\prime}}}

\newcommand{\iGo}[1]{\varGamma_{\!\overline{\!\{#1\}\!}\,}}
\newcommand{\set}{\!:=}
\newcommand{\lra}{\longrightarrow}
\newcommand{\Lra}{\Longrightarrow}

\newcommand{\onto}{\twoheadrightarrow}
\newcommand{\into}{\hookrightarrow}
\newcommand{\xto}[1]{\xrightarrow{#1}}
\newcommand{\xgets}[1]{\xleftarrow{#1}}

\newcommand{\bw}{\bigwedge}
\newcommand{\wh}[1]{{\widehat{#1}}}
\newcommand{\ov}[1]{{\overline{#1}}}
\newcommand{\uv}[1]{{\underline{#1}}}
\newcommand{\Zb}[1]{\overline{Z}^{#1}}
\newcommand{\fXY}{f \colon \X \to \Y}

\newcommand{\iso}%
{{\mkern8mu\longrightarrow \mkern-25.5mu{}^\sim\mkern17mu}}
\newcommand{\osi}%
{{\mkern8mu\longleftarrow \mkern-24.5mu{}^\sim\mkern16mu}}
\def\Otimes{\underset
  {\vbox to 0pt {\vskip-1ex\hbox{$\scriptscriptstyle=$}\vss}}
    {\otimes}\vadjust{\kern.4pt}}

\newcommand{\BL}{{\boldsymbol\Lambda}}

\newcommand{\smcirc}%
  {{\raise.15ex\hbox to.7em{$\hss \scriptstyle\circ\hss$}}} 

\newcommand{\DsssW}{\Delta_{\mspace{-2mu}\scriptscriptstyle \W}}
\newcommand{\DsssX}{\Delta_{\mspace{-2mu}\scriptscriptstyle \X}}
\newcommand{\DsssY}{\Delta_{\mspace{-2mu}\scriptscriptstyle \Y}}
\newcommand{\DsssZ}{\Delta_{\mspace{-2mu}\scriptscriptstyle \Z}}

\makeindex

\title[Pseudofunctorial behavior of Cousin complexes]
{\phantom{i}\vspace{-35pt}\phantom{i} Pseudofunctorial behavior of Cousin complexes\\on formal schemes}

\author[ J.\:Lipman]{Joseph Lipman}
\address{Dept.\ of Mathematics, Purdue University\\
              W. Lafayette IN 47907, USA}
\email {lipman@math.purdue.edu}

\author[ S.\:Nayak]{Suresh Nayak}
\address{Chennai Mathematical Institute \\
   92, G.\,N.\,Chetty Road\\
   Chennai-600017, INDIA}
\email {snayak@cmi.ac.in}

\author[P.\;Sastry]{Pramathanath Sastry}
\address{Dept.\ of Mathematics\\
University of Toronto, 
Toronto, Ont., CANADA M5S 3G3}
\email {pramath@math.toronto.edu}

\thanks{J.\,Lipman partially supported by National Security Agency, 
S.\,Nayak by
National Board of Higher Mathematics,
P.\,Sastry by GANITA laboratory, University of Toronto.
Close collaboration among all three authors
enabled by Purdue University, the Mathematisches 
Forschungsinstitut Oberwolfach, and the Banff International Research
Station.\vspace{-7pt}}

\begin{document}

 \def\copyrightyear{2005}
\def\copyrightholder{American Mathematical Society}

\setcounter{page}{3}

\begin{abstract}
On a suitable category
of formal schemes equipped with codimension functions
we construct a \emph{canonical pseudofunctor} $(-)^\sharp$
taking values in the corresponding categories of \emph{Cousin
complexes.} Cousin complexes on such a formal scheme $\X$ 
functorially represent\vspace{.4pt} derived-category objects $\F$ by the local\- cohomologies 
\smash{$H_x^{{\rm codim}(x)}\F\ (x\in\X\>)$} together with ``residue maps" from the cohomology at $x$ to that at each immediate specialization of $x$; this representation is faithful when restricted to $\F$ which are \emph{Cohen-Macaulay} (CM), i.e., $H^i_x
\F=0$ whenever $i\ne{\rm codim}(x)$. Formal schemes provide a framework for 
treating local and global duality as aspects of a 
single theory. One motivation has been to gain a better understanding of the close relation between local properties of residues and global variance properties of  dualizing complexes (which are CM). Our construction, depending heavily on local phenomena,
 is inspired by, but generalizes and makes
more concrete, that of the classical pseudofunctor~$(-)^\Delta$ taking
values in residual complexes, on which the proof of  Grothendieck's (global) Duality \index{Grothendieck, Alexandre!Grothendieck Duality}Theorem in Hartshorne's ``Residues and Duality'' is based. Indeed, it
is shown in the following paper by Sastry that $(-)^\sharp$
is a good ``concrete approximation'' to the fundamental duality 
pseudo\-functor~$(-)^!$. The pseudofunctor~$(-)^\sharp$ takes residual
complexes to residual complexes, so contains a canonical
representative of $(-)^\Delta$; and it generalizes as well 
several other functorial
(but not pseudofunctorial) constructions of residual complexes which
appeared in the 1990s.  \vspace{-24pt}
\end{abstract}

\enlargethispage*{15pt}

\maketitle


\tableofcontents
\vspace{-27pt}

\renewcommand{\theenumi}{\roman{enumi}}

\bibliographystyle{amsplain}

\newpage

\section{Introduction and main results}
\label{sec:intromain}

\subsection{Introduction}
\label{subsec:intro}

At the heart of the foundations of 
Grothendieck Duality\index{Grothendieck, Alexandre!Grothendieck Duality} lies the duality pseudofunctor $(-)^!$ 
described in the Preface. As indicated in~\S0.6,
under suitable  hypotheses on the map~$f\colon \X\to \Y$
of noetherian formal schemes, 
the functor $f^!\colon\Dqct^+(\Y)\to\Dqct^+(\X\>)$
can be realized in terms of \emph{dualizing complexes}
on~$\X$ and~$\Y$. Anyway, 
the study of dualizing complexes has its own importance,
for example as a natural 
generalization of the oft-appearing notion of \emph{dualizing sheaf.}
Concrete models for dualizing 
complexes---the residual complexes---are 
found in the category of Cousin complexes. This category,
which, among other virtues,
is an abelian subcategory of the usual category of complexes, 
provides fertile ground for a concrete pseudofunctorial
(or ``variance'') theory modeled after that of~$(-)^!$.

Our purpose here is to develop such a canonical
pseudofunctorial construction of Cousin complexes 
over a suitably general category of formal schemes.
(The notion of ``pseudofunctor'' is recalled at the
beginning of \S\ref{sec:pseudo}.)
The construction is motivated by well-known concrete 
realizations of the duality pseudofunctor $(-)^!$,
and indeed,  is shown in \cite{Sa2} to provide a ``concrete
approximation'' to~$(-)^!\<.$\looseness=-1

Before stating the main result (in \S\ref{subsec:outline}), 
we highlight some of its salient features, 
and relate it to some results in the literature on Cousin 
complexes.

First, the underlying category $\mathbb F$ on which 
we work is that of morphisms of noetherian
formal schemes with additional mild hypotheses specified
below in~\S\ref{subsec:mainre}. In particular $\mathbb F$ contains
many ordinary schemes, which can be regarded as
formal schemes whose structure sheaf (of topological rings) has the discrete topology.  
Also, $\mathbb F$ contains the 
opposite category of the category $\Cfr$\index{ $\Cfr$@$\Cfr$ (complete local rings subcategory)} of those
local homomorphisms of complete noetherian local 
rings which induce a finitely generated extension at
the residue fields. A key advantage of working in a 
category of formal schemes is that it offers a framework for 
treating local and global duality as aspects of a 
single theory, see, e.g.,~\cite[\S 2]{AJL2}. This paper continues
efforts to generalize all of Grothendieck's\index{Grothendieck, Alexandre!Grothendieck Duality!and formal schemes} duality theory to the 
context of formal schemes.

Second, while originally inspired by a study of the
classical construction in \cite[Chap.\,6]{RD} of a pseudofunctor
on residual complexes over schemes (see also \cite[Chap.\,3]{BCo}), 
we work more generally with Cousin complexes, the only restriction
being that the underlying modules be quasi-coherent and torsion.
(Over ordinary schemes the ``torsion'' condition is vacuous.)
The pseudofunctor we construct does however take
residual complexes to residual complexes (Proposition
\ref{prop:residual6}). So our construction generalizes the one
in \cite{RD}. 
 
Our construction is based on
the canonical pseudofunctor 
of Huang (\cite{Hu}),\index{Huang, I-Chiau!pseudofunctor} which is defined over~$\Cfr$. This pseudofunctor expands readily to
one with values in the category of graded objects underlying 
Cousin complexes, that is, Cousin complexes with vanishing\- differentials.
Most of our effort lies in working out what to do with
nontrivial differentials. 

Many details are thus already absorbed into 
the local theory of residues, through its basic role in Huang's work. 
In fact much of what we do in this paper comes down
ultimately to the relation between local operations involving
residues and global operations on Cousin complexes.

Finally, we note that several canonical 
constructions of residual complexes came out during the
1990s, see \cite{Hu2}, 
\cite{Hubl}, \cite{Ye}, \cite{Sa}
(some of which also use \cite{Hu}).
These constructions---functorial, 
but not pseudo\-functorial---lead by various methods to 
special cases of our result.
\enlargethispage*{0pt}

\subsection{Terminology and remarks on basic issues}
\label{subsec:mainre}

The main theorem of this paper is stated in 
\S\ref{subsec:outline}. To prepare the way, we first describe  
various underlying notions, referring to later sections for precise 
definitions; and  give a preliminary discussion of some  
of the basic issues involved.

Consider the category $\mathbb F$\index{ $\Forget$1@$\mathbb F$ (formal schemes subcategory)} whose objects are all (in some universe)
the \emph{noetherian\- universally 
catenary formal schemes admitting a codimension function} 
and~whose morphisms $\X'\to\X$ are all those formal-scheme maps which are 
\emph{essentially of 
pseudo-finite type} (\S\ref{subsec:formal}).
We will usually work 
with the refined category~$\bbFc$%
\index{ $\Forget$1@$\mathbb F$ (formal schemes subcategory)!$\bbFc$ (refinement of $\mathbb F$)} whose objects are pairs
$(\X,\Delta)$ with $\X$ in 
~$\mathbb F$ and $\Delta$ a codimension function on~$\X$,%
\index{ $\Delta$|see{codimension function}}
and whose morphisms 
$(\X',\Delta') \to (\X,\Delta)$ are those $\mathbb F$-maps
$f \colon \X' \to \X$ such that 
for any $x'\in \X'$ and $x \set f(x')$,  
$\Delta(x) - \Delta'(x')$ is the transcendence
degree of the residue field extension $k(x')/k(x)$ 
($k(x)$ being the residue field of the local ring~$\OXx\dots$,\vspace{.6pt}
see \ref{exam:formal1}).\looseness=-1

Let $(\X,\Delta)\in \bbFc$. 
A \emph{Cousin complex on} $(\X,\Delta)$,\index{Cousin complex} or a \emph{$\Delta$-Cousin complex on\/ $\X,$}
is an  $\OX$-complex~$\cMb$
such that for each $ n\in\mathbb Z$, $\M^n$ is the direct sum of
a family  of $\OX$-submodules $(i_xM_x)_{x\in\X\<,\,\Delta(x)=n}$, 
where\vspace{-1pt} $M_x$  is an $\OXx$-module and $i_xM_x$%
\index{ $\I$@$i_x$ (produces constant sheaf on $\,\ov{\<\<\<\{x\}\<\<\<}\,$)} is the extension by 0 of
the constant sheaf~$\,\ov{\! M}_x$ on the closure
$\,\ov{\!\{x\}\!}\,$ such that 
for all nonempty\vspace{1pt} open subsets $V$ of~$\,\ov{\!\{x\}\!}$,
$\,\ov{\! M}_x(V)$ is the 
$\cO_{\>\ov{\!\{x\}\!}}(V)$-module~$M_x$\vspace{-2pt} 
(see \S \ref{subsec:cousin}). For such an $\cMb$, and
any $x\in\X$, $i_xM_x$ is uniquely determined:
$i_xM_x=\iG{\ov{\!\{x\}\!}}\,\M^{\Delta(x)}$, 
where $\iG{\ov{\!\{x\}\!}}\,$ is the
subfunctor\vspace{-1.5pt} of the identity functor\vspace{.6pt} taking any $\OX$-module to its
sheaf of sections supported in~$\,\ov{\!\{x\}\!}\,$.  Thus 
$M_x$,\vspace{1pt} which is $(i_xM_x)(U)$ for any open neighborhood $U$ of $x$, is
determined by~$\cMb$ and~$x$, and so we denote it by $\cMb(x)$.\vspace{1pt}

Let $\Coz_{\Delta}(\X\>)$\index{  $\Coz_{\Delta}$ (Cousin $\Aqct$ complexes category)} be the full subcategory of the category of
$\OX$-complexes with objects
those $\Delta$-Cousin $\OX$-complexes~$\cMb$
whose underlying graded modules  are 
\emph{quasi-coherent torsion $\OX$-modules} 
(\S \ref{subsec:mod}, \S \ref{subsec:tor}). By \ref{lem:mod6}
this last condition on~$\cMb$
means simply that for any $x\in\X$, each element of~$\cMb(x)$  is annihilated
by some power of the maximal ideal~$m_x$ of~$\OXx$, or, as we will say, 
$\cMb(x)$ is a \emph{zero-dimensional\/ $\OXx$-module}.\index{zero-dimensional module} Thus we can, and will,
view $\cMb(x)$ as a zero-dimensional module over the $m_x$-adic completion\vspace{1pt}
$\wh{\OXx}$. 
The category of such modules will be denoted by $(\cO_x)_\#\>$.\vspace{1.5pt}
\index{ N@$(\cO_x)_\#$  (zero-dimensional $\OXx$-modules category)}

Let $\Coz^0_{\Delta}(\X\>)$ be the full subcategory%
\index{  $\Coz_{\Delta}$ (Cousin $\Aqct$ complexes category)!$\Coz^0_{\Delta}$ $\Coz_{\Delta}$-complexes with 0 differential)}  
of~$\Coz_{\Delta}(\X\>)$ with objects those complexes
whose differentials are all zero.
For any $\OX$-module $\F$ the natural map 
$$
\Hom_{\>\cO_{\<\<\X}}(\F\<,\>i_xM_x)\to \Hom_{\>\OXx}(\F_{\<x} ,M_x)
$$
 is easily seen
to be bijective (i.e., $i_x$ is right-adjoint to the 
functor $\F\mapsto \F_{\<x})$. Hence 
$$
[\,\Hom_{\>\cO_{\<\<\X}}(i_yM_y\<,i_xM_x)\ne0\,]
\implies [\,(i_yM_y)_x\ne 0\,]\implies [\,x\in\,\ov{\!\{y\}\!}\,\,].
$$
Consequently, a morphism $\varphi\colon\cMb\to\cNb$ in $\Coz^0_{\Delta}(\X\>)$ 
is the same\vspace{.6pt} thing as
a family of $\>\wh{\OXx}$-homomorphisms 
$\bigr(\varphi(x)\colon \cMb(x)\to\cNb(x)\bigl)_{x\in\X}\>.$ 
It follows that $\Coz^0_{\Delta}(\X\>)$ is naturally equivalent to the
disjoint union of the family of categories 
$\bigl((\cO_x)_\#\bigr)_{x\in\X}$. (It also follows that the category~ $\Coz_{\Delta}(\X\>)$ is \emph{abelian.})\vspace{1pt}
 
The \emph{forgetful functor} 
$\Forget_{\Delta}(\X\>) \colon \Coz_{\Delta}(\X\>) \to \Coz^0_{\Delta}(\X\>)$%
\index{ $\Forget$5@$\Forget$ (forgetful functor)}
sends a complex 
to its underlying graded module, and a map of complexes to itself.

{\small{
\begin{aexam}
The scheme $\X= \Spec(\mathbb Z)$ is an affine formal scheme 
with $(0)$ as an ideal of definition. Every $\OX$-module
is a torsion module. Let $\Delta$ be the codimension
function sending the generic point to 0 and all other points to 1.

The natural surjection 
$\>\mathbb Q\twoheadrightarrow \mathbb Q/\mathbb Z\cong
\oplus_p\, ( {\mathbb Q}/{\mathbb Z}_{p\mathbb Z})$
(with $p$ ranging over the positive primes, so that $p\mathbb Z$ ranges
over all nonzero prime ideals) may be viewed as a 
$\mathbb Z$-complex concentrated in degrees 0 and 1.
Application of the sheafification functor ${}^\sim$ produces a
complex in $\Coz_{\Delta}(\X\>)$
\[
\mathbb Q\>^{\sim} \lra \,\bigoplus_p
\; ({\mathbb Q}/{\mathbb Z}_{p\mathbb Z})^{\sim}.
\]
This is a residual complex, with homology $\OX$ concentrated 
in degree~0.

Analogous statements hold for any Dedekind domain and its
fraction field. For a generalization to any irreducible
regular scheme, with the codimension function which sends the generic
point to 0, see \cite[p.\,304, Example]{RD}.
\end{aexam}
}}

\smallskip
By a $\Coz$-valued (resp.~$\Coz^0$-valued) pseudofunctor on~$\bbFc$ 
we shall always mean a \emph{contravariant pseudofunctor} 
(see \S \ref{sec:pseudo}) on the category $\bbFc$ which assigns to any
$(\X,\Delta)$ in $\bbFc$ the category
~$\Coz_{\Delta}(\X\>)$ (resp.~$\Coz^0_{\Delta}(\X\>)$).

We aim to construct a canonical $\Coz$-valued
pseudofunctor on~$\bbFc$, one which, when restricted to pseudo-proper
maps will be shown in \cite{Sa2}
to be right adjoint to the direct image functor of Cousin complexes, 
and further, to be a ``concrete approximation''
to the basic duality pseudofunctor $(-)^!$. Let us discuss briefly
some of the issues involved.

To begin with, the following examples serve as inspiration
for our construction. These are generalized Cousin versions of 
well-known concrete realizations of $(-)^!$ on ordinary schemes.

\begin{aexam}
\label{exam:mainre1}
Let $f \colon (\X, \Delta_1) \to (\Y, \Delta)$ be an $\bbFc$-map
and $\cMb\in\Coz_{\Delta}(\Y)$.\vspace{1.5pt}

(i) If $f$ is smooth, with 
relative dimension~$d$ (a locally constant function, see 
\ref{def:morph1a}, \ref{def:diff4a}), then the complex
$E_{\Delta_1}\R\iGp{\X}(f^*\cMb \otimes_{\X} \omega_{\<\<f}[d\>])$ is 
in~$\Coz_{\Delta_1}(\X\>)$ (see \ref{def:diff6} and  \S\ref{subsec:gl2lo}),
where: 
\newline---$E_{\Delta_1}$ is the Cousin functor corresponding to the
filtration of $\X$ induced by $\Delta_1$ (see \S\ref{subsec:cousin});
\newline---$\iGp{\X}$ is the subfunctor of the identity functor on
$\OX$-modules 
associating to any such module
the largest submodule each section of which, over any open set~$\U$, is
annihilated by some open $\cO_\U$-ideal (see \S\ref{subsec:tor}), 
and $\R\iGp{\X}$ is the right-derived functor of $\iGp{\X}$;
\newline---$\omega_{\<\<f}$ is the $d$-th exterior power of the%
\index{ WV@$\omega_{\<\<f}$ (relative dualizing sheaf)}%
\index{ O@$\R$ (right-derived functor of)}
complete differential module of $\>\X$ over $\Y$ 
(see \S\ref{subsec:differ}, \ref{prop:diff4}, \ref{def:diff6}); and
$[d\>]$ is the usual translation operator on complexes, where
both instances of ``$d\>\>$'' should be
replaced separately on each connected
component of~$\X$ by the value of $d$ on that component.

Thus we get a functor 
$\mathbb E_{\<f}\colon\Coz_{\Delta}(\Y)\to\Coz_{\Delta_1}\<(\X\>).$
In \S5 we expand this functor to a pseudofunctor on smooth maps.
This expansion is canonical in that it 
involves no choices other than the sign convention fixed in 
\S\ref{subsec:conv}\kern.5pt\eqref{conv5} to handle the relation between
$\otimes$ and translation of complexes.  
For open immersions $f\<$, $\mathbb E_f$ is then
canonically pseudofunctorially isomorphic 
to the usual restriction functor~$f^*$.\vspace{1.5pt}

(ii) If $f$ is a closed immersion (\cite[p.\,442]{EGAI}),
the complex $f^{-1}\sHom_{\Y}(f_*\OX, \cMb)$ is in~$\Coz_{\Delta_1}(\X\>)$ 
(see \ref{prop:clim1}). 

 \enlargethispage{-3pt}
 
Thus we get a functor 
$f^\flat\colon\Coz_{\Delta}(\Y)\to\Coz_{\Delta_1}\<(\X\>).$
In \S6 we canonically expand this functor to a 
pseudofunctor on closed immersions. For open-and-closed
immersions $f\<$, $f^\flat$ is
canonically pseudofunctorially isomorphic 
to the usual restriction functor~$f^*\<$.
\end{aexam}

\pagebreak

What we want is to glue these examples together, i.e.,
to get a Coz-valued pseudofunctor $(-)^{\sharp}$ on $\bbFc$
whose res\-triction to
smooth maps (resp.~closed immersions) 
is isomorphic to $\mathbb E_{(-)}$
(resp.~$(-)^\flat$).  The result should be \emph{canonical,} i.e.,
unique up to unique isomorphism, as
explicated in the Main Theorem in \S\ref{subsec:outline}.

For this we need  
concrete descriptions of a number of ``glueing'' isomorphisms
involving $\mathbb E_{(-)}$ and~$(-)^\flat$,
associated to various compositions of smooth maps and closed
immersions. 
For instance, suppose $f$ and $\cMb$ are as in~\ref{exam:mainre1}, and
assume further that $f$ admits factorizations $f = h_1i_1$ and $f =h_2i_2$ 
where $h_1$, $h_2$ are 
smooth maps in $\bbFc$ and $i_1$, $i_2$ are closed immersions. 
(A map $g$ will be called \emph{factorizable}\index{factorizable map} 
if~$\>g = hi$ with $h$ a smooth map in $\bbFc$
and $i$ a closed immersion.) Then, at the least, there has to be a 
canonical functorial isomorphism 
$i_1^\flat \mathbb E_{h_1}\iso i_2^\flat \mathbb E_{h_2}$. 

Moreover, for a general map~$g$ in $\bbFc$, one that is not
factorizable, it is not immediately obvious how to 
define $g^{\sharp}$, canonically or otherwise.
This issue is related to the previous one in that if one
somehow obtains a good definition of $f^{\sharp}$ for 
any factorizable $f$, then for general $g$\vspace{.8pt} 
one could use \emph{local}  
factorizations (\ref{cor:morph1b}) to 
canonically define $(g\big|_{\U_{\<\lambda}})^{\sharp}$ \vspace{-.4pt} for 
a suitable open cover $\{\U_{\lambda}\}$ of~$\X$,\vspace{.6pt}
and, having noted that $\mathbb E_f$ for $f$ an
open immersion ``is''  the restriction functor $f^*\<\<,$ 
anticipate further that  canonicity would enable pasting 
the various~ $(g\big|_{\U_{\<\lambda}}\<\<)^{\sharp}$
together to form a global~$g^\sharp$.

Thus the issue of canonicity is basic. One 
must specify with care all the maps that go into the
construction, 
and relations between them,
and in particular, 
pay close attention\vspace{3pt} to signs. 

\enlargethispage*{0pt}

We now briefly describe our construction of $(-)^{\sharp}$.
One can think of a Cousin complex as comprising two parts, 
its underlying graded object and its
differential; and accordingly
the problem of canonicity is addressed in two stages.
We start  with a 
pseudofunctor $(-)^{\natural}$ on graded objects
only---in other words, $(-)^{\natural}$ is a $\Coz^0$-valued 
pseudofunctor on~$\bbFc$.%
\index{    ${\boldsymbol{{-}^\natural}}$ ($\Coz^0$-valued pseudofunctor on~$\bbFc$)}  The idea is to
use Huang's work\index{Huang, I-Chiau!pseudofunctor} (\cite[Chap.\,6]{Hu}) wherein he 
constructs a \emph{canonical} pseudofunctor $(-)_{\#}$ which 
takes values in categories $R_\#$ of zero-dimensional\vspace{.4pt} 
modules over local
rings $(R,\mathfrak m)$ in the previously-mentioned category~$\Cfr$, 
that is, modules such that each element is
annihilated by some power of the maximal ideal $\mathfrak m$.
A map~$f\colon\X\to\Y$ such as in~\ref{exam:mainre1}
provides, for each 
$x \in \X$, a map $\wh{f_{\<\<x}^{}}$ in~$\Cfr$, namely,
the naturally induced map $\wh{\cO_{\Y\<,\>f(x)}} \to \wh{\OXx}$. 
We can therefore define~$f^{\natural}$ by using 
Huang's $(-)_{\#}$, pointwise; and furthermore define~$(-)^{\natural}$ as a 
pseudofunctor in a similar manner.
 (Actually, we will need to impose a  subsequently important sign modification on $(-)_{\#}$, see proof of \ref{thm:huang1}.)

Then we upgrade the pseudofunctor $(-)^{\natural}$  to a
pseudofunctor $(-)^{\sharp}$ at the level of complexes. 
This means setting up a canonical differential on 
every $f^{\natural}\cMb$. 
For that purpose we show that in each case, (i) or (ii), 
of~\ref{exam:mainre1},
the complex obtained in~$\Coz_{\Delta_1}(\X\>)$ is, at the 
graded level, canonically isomorphic to $f^{\natural}\cMb$. 
(See \S\S\ref{subsec:fin-sm} and~\ref{subsec:fin-clim}.)
So we can transfer the differential from $\mathbb E_f\cMb$
(resp.~$f^\flat$) to $f^\natural\cMb$ when $f$ is smooth 
(resp.~a closed immersion). The important step, the one lying at the 
heart of this paper, is to show that
\emph{if\/~$f$ is a factorizable map, then  the 
differentials for\/ $f^{\natural}\cMb$ so obtained via different 
factorizations of\/~$f$ are the same.} This is accomplished in
Proposition \ref{prop:fin-fact6}. 
So we have a canonical differential for $f^\natural\cMb$ whenever $f$ is
factorizable. As mentioned before, a definition of the differential in the 
general case, indeed a canonical one, then follows easily, 
see \S\ref{subsec:fin-gen}.
Upgrading of the remaining data of $(-)^{\natural}$ 
is straightforward.

We can summarize as follows. Let us say that a $\Coz$ or $\Coz^0$ functor
is \emph{canonical} if it satisfies the appropriate
analogs of the conditions on pages 319--321 of \cite{RD}, with ``finite
morphism" replaced by ``closed immersion". (Some of these conditions spell out 
what we have said above about pasting the basic examples together; and the
others say that this pasting should be compatible with certain base-change 
and residue isomorphisms, as explained in detail in \S\S6 and 7 below.)
In addition, the pseudofunctor should commute with restriction to open
subsets.
Any two canonical pseudofunctors are canonically isomorphic.  What is
accomplished in \S\S5-7 is, in essence, showing that \emph{the $\Coz^0\<\<$-valued
pseudofunctor derived as above from Huang's zero-dimensional pseudofunctor%
\index{Huang, I-Chiau!pseudofunctor}
is, modulo some sign modifications, canonical;} 
and what is done in \S8 shows, in essence, 
that \emph{any canonical $\Coz^0\<\<$-valued
pseudofunctor can be upgraded to a canonical\/ $\Coz$-valued pseudofunctor.}

\subsection{Outline of construction of $(-)^{\sharp}$}
\label{subsec:outline}
In more detail, our construction of~$(-)^{\sharp}$ is realized through the
steps~(A)-(D) below. The main theorem is stated after~(C).
Step~(D) occupies the bulk of this paper.
 
\medskip
{\textbf{(A)}} We start with the punctual case, i.e., we (temporarily) 
restrict attention
to the full subcategory $\bbFc^{\!\smcirc}\subset\bbFc$ with objects those $(\X,\Delta)$ 
such that the underlying space of~$\X$
is a single point---so that $\Delta$ can be identified with a single
integer. Associating to each such $(\X,\Delta)$ the
complete local ring $A_\X\set\cO_\X(\X\>)$ leads to an antiequivalence
from $\bbFc^{\!\smcirc}$ to the category of pairs $(A,n)$ with $A\in\Cfr$ 
(see~\S\ref{subsec:intro}) and \mbox{$n\in\mathbb Z$,} a morphism $(A,n)\to
(B,m)$ being a $\Cfr$-map $A\to B$ of residual transcendence degree $n-m$.
The category $\Coz_\Delta(\X\>)$ is
isomorphic to the category
of zero-dimensional $A_\X$-modules. 
In effect, then, the punctual case is covered by
\cite{Hu}, wherein Huang constructs\index{Huang, I-Chiau!pseudofunctor} 
a \emph{canonical} covariant pseudofunctor $(-)_{\#}$ 
on~$\Cfr$ taking values in zero-dimensional modules.

For our purposes a variant
$(-)_{\sharp}$ of Huang's pseudofunctor $(-)_{\#}$ is more convenient. 
The changes made are just to multiply the comparison isomorphisms 
of~$(-)_{\#}$ by a suitable $\pm$ factor, and to
reverse the order in which modules appear in certain tensor products.
The description of~$(-)_{\sharp}$ begins with the following
data, cf.~Example\-~\ref{exam:mainre1}: 
\begin{enumerate}
\item For any formally 
smooth $\Cfr$-map $\phi \colon R \to S$ with 
$r\set\text{dimension of }S/\mathfrak m^{}_{\<\<R}S$
and \mbox{$t\set{}\;$residual transcendence
degree,} and $\omega_{S/R}\set\widehat{\Omega_{S/R}^{r+t}}\>$,%
\index{ WW@$\omega_{S/R}$ (relative dualizing module)}
the $(r+t)$-th exterior power of the complete
module of $R$-differentials of~$S$,
and for any zero-dimensional $R$-module $M\<$, a 
specific isomorphism 
$$
\phi_{\sharp}M \iso H^r_{\mathfrak m^{}_{\<\<S}}(M
\otimes_R\omega_{S/R}).
$$
  
\item For any surjective $\Cfr$-map $\phi \colon R \to S$ and
zero-dimensional $R$-module~$M\<$, a 
specific isomorphism 
$$
\phi_{\sharp}M \iso \Hom_R(S,M).
$$   
\end{enumerate}
The full characterization of $(-)_{\sharp}$  in 
\ref{thm:huang1} also
involves concrete descriptions of some special comparison maps 
via the specific isomorphisms in~(i) and~(ii).

\medskip
{\textbf{(B)}} Using $(-)_{\sharp}$, we construct a $\Coz^0$-valued 
pseudofunctor $(-)^{\natural}$ on~$\bbFc$ as follows.
For any $\bbFc$-map $(\X, \Delta_1) \xto{\; f \;} (\Y, \Delta)$ 
and points $x \in \X$, $y\set f(x) \in \Y$, 
the natural induced map of complete local rings 
$\wh{f_{\<\<x}^{}} \colon \wh{\OYy} \to \wh{\OXx}$ (completions 
being along the respective maximal ideals) is in~$\Cfr$. For  
$\cMb \in \Coz^0_{\Delta}(\Y)$
the $\OYy$-module~$\cMb(y)$, being a 
zero-dimensional $\OYy$-module, is also 
naturally an $\wh{\OYy}$-module, from which we get the zero-dimensional 
$\OXx$-module
$\wh{f_{\<\<x}^{}}_{\sharp}(\cMb(y))$. We now let  $f^{\natural}\cMb$
be the unique object in~$\Coz^0_{\Delta_1}(\X\>)$ given by
\begin{equation}
\label{eq:out1}
\bigl(f^{\natural}\cMb\bigr)(x) = \wh{f_{\<\<x}^{}}_{\sharp}\bigl(\cMb(y)\bigr)
\qquad\bigr(x\in\X,\;y=f(x)\bigl). 
\end{equation}
Then $f^{\natural}$ is a $\Coz^0$-valued functor. Furthermore, 
pseudofunctoriality of $\>(-)_{\sharp}\>$ induces a 
pseudofunctor $(-)^{\natural}$ in an obvious manner.

\medskip
{\textbf{(C)}} The upgrading of the $\Coz^0$-valued 
pseudofunctor $(-)^{\natural}$ to a $\Coz$-valued 
one~$(-)^{\sharp}$ is carried out over two 
subcategories of~$\bbFc$, namely, the subcategory 
of smooth~maps and the subcategory of closed 
immersions. (Dealing with the smooth subcategory is not at all
straightforward, see Proposition \ref{prop:itloco1}, 
but dealing with closed immersions is,
see Proposition \ref{prop:clim3}.)

Let us elaborate. Notation remains as in~(B),
but now $\cMb\in \Coz_{\Delta}(\Y)$. 

Suppose $f$ is smooth. Assume $f$ has constant relative 
dimension~$d$. (This can be arranged by restricting to 
connected components if necessary). 
Then we specify a natural graded isomorphism (see \eqref{eq:fin-sm1})
\begin{equation}
\label{eq:out2}
f^{\natural}\cMb \iso 
E_{\Delta_1}\R\iGp{\X}\bigl(f^*\cMb \otimes_{\Y} \omega_{\<f}[d\>]\bigr).
\end{equation} 
The inputs into the definition of~\eqref{eq:out2}
can be organized into two parts. The first involves a
specific  isomorphism (\S \ref{subsec:gl2lo}), with $r =
\dim(\OXx/m_y\OXx)$,
$M\set\cMb(y)$:
\[
\bigl(E_{\Delta_1}\R\iGp{\X}(f^*\cMb \otimes_{\Y} \omega_f[d\>])\bigr)(x)
\iso H^r_{m_x}(M \otimes_y (\omega_{\<f})_x\bigr).
\]
The second involves an isomorphism
$$
H^r_{m_x}\bigl(M \otimes_y (\omega_{\<f})_x\bigr) \iso 
\wh{f_{\<\<x}^{}}_{\sharp}M  = (f^{\natural}\cMb)(x)
$$
given up to a sign by completion and the 
isomorphism of case~(i) in step~(A) above. 

If $f$ is a closed immersion, then we specify a natural
graded isomorphism
\begin{equation}
\label{eq:out3}
f^{\natural}\cMb \iso f^{-1}\sHom_{\Y}(f_*\OX, \cMb).
\end{equation}
As with \eqref{eq:out2}, the inputs used in 
defining \eqref{eq:out3} can be organized into two parts,
one consisting of a specific natural isomorphism (with $M = \cMb(y)$)
\[
(f^{-1}\sHom_{\Y}(f_*\OX, \cMb))(x) \iso
\Hom_{\OYy}(\OXx, M) 
\]
and the other of an isomorphism
$$
\Hom_{\OYy}(\OXx, M) \iso  
\wh{f_{\<\<x}^{}}_{\sharp}M = (f^{\natural}\cMb)(x)
$$
given by completion and the 
isomorphism of case~(ii) in step~(A) above. 

The graded isomorphisms in \eqref{eq:out2} and \eqref{eq:out3} 
are canonical, so we have a canonical choice for a 
differential on~$f^{\natural}\cMb$ whenever $f$ is a smooth
$\bbFc$-map or a closed immersion, namely the unique one such that 
\eqref{eq:out2} or \eqref{eq:out3} (as the case may~be) 
\mbox{underlies} an isomorphism
of complexes. (When $f$ is a smooth closed immersion
the two possible differentials coincide, see \ref{rem:fin-clim2}.)
Thus we 
have a canonical choice for a Coz-valued functor $f^{\sharp}$
for all $f$ in either of the two subcategories under consideration.

\medskip

We are now in a position to state our main result. The notation remains
as in steps~(A)--(C) above. 
Refer to the beginning of~\S\ref{sec:pseudo} for the
notations $C^\#$ (for comparison isomorphism) and $\delta^\#$ 
(for unit isomorphism)
used in the definition of ``pseudofunctor.''
Recall that for any $(\X,\Delta) \in \bbFc\>$, the forgetful functor
$\Forget_{\Delta}(\X\>) \colon \Coz_{\Delta}(\X\>) \to
\Coz^0_{\Delta}(\X\>)$ 
operates by forgetting the differential on
a Cousin complex. For convenience, we shall abuse notation and 
use only the symbol~$\Forget$, the rest being clear from
the context. 
 \medskip

\begin{main}\index{Main Theorem on existence and uniqueness of Cousin-complex pseudofunctor}
There exists a unique $\Coz$-valued pseudofunctor $(-)^{\sharp}$ on\/ $\bbFc$
satisfying the following conditions.
\begin{enumerate}
\item The forgetful functor \emph{\Forget} makes $(-)^{\sharp}$
into the above pseudofunctor~$(-)^{\natural}$.
In other words:

\begin{enumerate}
\item For any $\bbFc$-morphism\/ $f,$ 
we have\/ 
$\emph{\Forget}\smcirc f^{\sharp} = f^{\natural}\smcirc\emph{\Forget}$.
\item For any\/ $\bbFc$-morphisms\/ $f\<,\>g$ such that the
composition~$gf$ 
is defined it holds that, via\/ \emph{(a),}
$\emph{\Forget}\bigl(C^{\sharp}_{f,g}(-)\bigr) = 
C^{\natural}_{f,g}\bigl(\emph{\Forget}(-)\bigr)$.\vspace{-1pt}
\item 
For any $(\X, \Delta) \in \bbFc$ it holds that, via\/ \emph{(a),} 
$\emph{\Forget}\bigl(\delta_\X^{\sharp}(-)\bigr) 
= \delta_\X^{\natural}\bigl(\emph{\Forget}(-)\bigr)$.
\end{enumerate}\vspace{1.5pt}
\item 
If the\/ $\bbFc$-morphism\/ $f\colon (\X, \Delta_1) \to (\Y, \Delta)$ 
is smooth\/ 
$($resp.~a closed immersion\/$)$ 
then for\/ $\cMb\in\Coz_{\Delta}(\Y)$
the differential on\/~$f^{\sharp}\cMb$ is the 
unique one such that 
the graded isomorphism in~\eqref{eq:out2}
$($resp.~\eqref{eq:out3}$)$ underlies
an isomorphism of complexes $(\<$see \textup{(i)(a)}$)$.
\vspace{1pt}
\item
If\/ $(\X,\Delta)\in\bbFc,$ $1^{}_{\<\X}$ is
the identity map of\/ $\X,$ and\/ 
$u\colon\U\hookrightarrow \X$ is the inclusion map of an
open subset of\/ $\X,$ then\/ for any\/ $\cMb\in\Coz_\Delta(\X\>),$
$u^\sharp\cMb$ is the restriction to\/ $\U$ of the complex\/
$1_\X^\sharp\cMb$.
\end{enumerate}
\end{main}

\emph{Remark.} The functor from $\OX$-modules to $\cO_\U$-modules
given by ``restriction to~$\U$,'' being left-adjoint to the direct
image functor $u_*$, can be identified with $u^*\<$. When this is done,
the condition in (iii) becomes
\begin{equation*}
u^\sharp=u^*\<\smcirc 1_{\<\X}^\sharp.\tag*{$(\textup{iii})'$}
\end{equation*}
And then if the functor $(-)^\sharp$ is replaced by its normalization 
(see \S\ref{sec:pseudo}),  the condition becomes simply that
$u^\sharp=u^*\<$.\vspace{1.5pt}

The uniqueness part of the theorem is easy to prove.
Indeed if~$(-)^{\sharp}$ and~$(-)^{\sharp '}$ are pseudofunctors
satisfying conditions (i) and~(ii) of the theorem, then
for $f,\X,\Y,\cMb$ as in~(ii) of the theorem,\vspace{-.6pt} by (i)(a) we have 
$f^{\sharp}\cMb = f^{\sharp '}\cMb$ as graded torsion modules. 
It suffices to check locally on~$\X$ that the differentials
on these graded objects agree, i.e., that for $u$ and 
$1^{}_{\<\X}$ as in (iii),
$u^*\<f^{\sharp}\cMb = u^*\<f^{\sharp '}\cMb$ as complexes.
A straightforward check of definitions shows that the functorial
composition
$$
(fu)^\sharp\underset{\under{1.3}{C^\sharp_{u,\>f}}}\iso
u^\sharp f^\sharp\>\overset{\:\textup{(iii)}'} {=\mspace{-6mu}=} 
u^*\<1_{\<\X}^\sharp f^\sharp
\xrightarrow[\under{1.3}{\!\!u^*C^\sharp_{1_{\<\<\X},\>f}\!\!\!}]{\Iso}
u^*\<f^\sharp\vspace{1pt}
$$
is the identity. \big(In other words, for every $\cMb$ and $x\in\U$, 
the composition induces the identity 
map of~$[(fu)^\sharp\cMb](x)=\wh{f_{\<\<x}^{}}_{\sharp}\bigl(\cMb(f(x))\bigr)=
[u^*\<f^\sharp\cMb](x)$.\big) This can also be stated via
``restriction'' as\vspace{-2pt}
$$
(f|_\U^{})^\sharp=f^\sharp|_\U^{}.\vspace{1pt}
$$
The same holds for $(-)^{\sharp'}$.
Hence we can replace $f$ by $fu$,\vspace{.7pt} and so by
\ref{cor:morph1b} we may assume without loss of
generality that $f$ factors as $f = hi$ where~$h$ is a smooth $\bbFc$-map
and~$i$ is a closed immersion. By~(ii), 
$h^{\sharp} = h^{\sharp '}$ and~$i^{\sharp} = i^{\sharp '}$. By~(i)(b)
it follows that $f^{\sharp}\cMb = f^{\sharp '}\cMb$ as 
complexes too, q.e.d.\vspace{1pt}

\enlargethispage*{3pt}

Existence, which is the final step 
of our construction, is the difficult part of the Theorem.
For instance, though we have obtained, in step~(C), a definition 
for~$f^{\sharp}$
over the subcategory of smooth $\bbFc$-maps, it is not at all
obvious that $(-)^{\sharp}$ is a pseudofunctor on that 
subcategory. Even the seemingly simple condition~(iii) is not
trivial to verify, see Proposition~\ref{lem:fin-sm3}.

\pagebreak[3]

\medskip
{\textbf{(D)}} 
First suppose the ~$\bbFc$-map $f$ admits a factorization
$f = hi$ where $h$ is a separated smooth $\bbFc$-map and 
$i$ is a closed immersion. Then, using step~(C) twice, 
we get a differential on~$h^{\natural}i^{\natural}\cMb$.
The graded comparison\vspace{-1pt} 
$C^\natural_{i,h}\colon h^{\natural}i^{\natural}\cMb \iso 
f^{\natural}\cMb$ then induces a differential $d_{h,i}$
on~$f^{\natural}\cMb$.
A necessary condition for the existence of~$(-)^{\sharp}$ in
the Main Theorem is that this differential 
on~$f^{\natural}\cMb$ not depend on the choice of 
factorization $f = hi$. In other words, if there is 
another factorization $f = h'i'$ with~$h'$ smooth and~$i'$
a closed immersion then one must have
$d_{h,i} = d_{h'\!,\>i'}$.
 
In Proposition \ref{prop:fin-fact6} we prove that this necessary
condition holds. The key inputs are \ref{prop:itloco1},
\ref{prop:clim3}, \ref{prop:cartsq9} and~\ref{prop:ret1a}.
Since $\bbFc$-maps all have local factorizations as above, the canonical nature
of our construction makes it fairly easy to patch the resulting local
pseudofunctors together to a global one, see \S\ref{subsec:fin-gen}.

\smallskip
In~\S\ref{sec:residual} we discuss residual complexes.
We end in~\S\ref{sec:sexd} with some results 
on~$(-)^{\sharp}$ for \'etale maps and on 
certain Cohen-Macaulay complexes that we encounter while 
constructing~$(-)^{\sharp}$ over smooth maps.

\enlargethispage*{2pt}
\subsection{Conventions}
\label{subsec:conv}
We use the following notation and conventions in this paper.
The Bourbaki dangerous-bend symbol indicates
potential sources of ambiguity or conflict with some other
convention.

\begin{enumerate}
\item\label{conv1} 
Let $\A$ be an abelian category. Set\index{  $\A$ (abelian category)}

$\bC(\A)$ \set the category of\index{ $\bC(\A)$ (category of $\A$-complexes)}
$\A$-complexes,

$\K(\A)$ \set the corresponding homotopy category,\index{ $\K(\A)$ (homotopy category
of $\A$-complexes)} 

$\D(\A)$ \set the corresponding derived category. \index{ $\D(\A)$ (derived category
of $\A$-complexes)}

\noindent The differential~$d_C^\bullet$ in a complex $C^\bullet$%
\index{ $\Cfr$2@$d_C^\bullet$ (differential of a complex~$C^\bullet)$}
is always understood to increase degree: $d_C^n$ maps $C^n$ to~$C^{n+1}$ 
for all $n\in\mathbb Z$.
\vspace{1pt}

\item\label{conv2} 
Let $(X, \cO_{\<\<X})$ be a ringed space, and $x\in X$. Set

$\A(X)$ \set the abelian category of \mbox{$\cO_{\<\<X}$-modules,}%
\index{  $\A(X)$ (category of $\cO_{\<\<X}$-complexes on ringed space $(X, \cO_{\<\<X})$)}

$\Aqc (X)$%
\index{  $\A(X)$ (category of $\cO_{\<\<X}$-complexes on ringed space $(X, \cO_{\<\<X})$)!$\Ac$@$\Aqc (X)$ (quasi-coherent subcat.\ of $\A(X)$)} 
 (resp.~$\Ac (X)$,%
\index{  $\A(X)$ (category of $\cO_{\<\<X}$-complexes on ringed space $(X, \cO_{\<\<X})$)!$\A$@$\Ac (X)$ (coherent subcat. of~$\A(X)$)}
 resp.~$\Avc (X)$)%
\index{  $\A(X)$ (category of $\cO_{\<\<X}$-complexes on ringed space $(X, \cO_{\<\<X})$)!$\A_{?}$@$\Avc (X)$ 
(\smash{$\subdirlm{}\<\<\<\<$}-coherent\vspace{1pt} subcat.\ of $\A(X)$)} 
 $\set{}$the full subcategory of~$\A(X)$ whose objects are the quasi-coherent 
(resp.~coherent, resp.~\smash{$\dirlm{}\!\!$'s} 
of co\-herent) \mbox{$\cO_{\<\<X}$-modules,}

$\E\otimes_{\<\<X}\<\F \set $ the tensor product of the 
$\cO_{\<\<X}$-modules $\E$, $\F$,

$E\otimes_x\< F\set $ the tensor product of the $\cO_{\<\<X,x}$-modules
$E$, $F$,

$\bC(X) \set \bC(\A(X))$, \quad 
\index{ $\bC(\A)$ (category of $\A$-complexes)!$\bC(X)\set\bC(\A(X))$}
$\K(X) \set \K(\A(X))$, \quad 
\index{ $\K(\A)$ (homotopy category of $\A$-complexes)!$\K(X)\set\K(\A(X))$}
$\D(X) \set \D(\A(X))$.\vspace{1pt}
\index{ $\D(\A)$ (derived category of $\A$-complexes)!$\D(X)\set\D(\A(X))$}

\noindent For any full subcategory $\A_{?}(X)$ 
\index{  $\A(X)$ (category of $\cO_{\<\<X}$-complexes on ringed space $(X, \cO_{\<\<X})$)!$\A$@$\A_{?}(X)$ (full subcategory of~$\A(X)$)} 
of~$\A(X)$, let $\D_{?}(X)\subset\D(X)$ 
\index{ $\D(\A)$ (derived category of $\A$-complexes)!$\D(X)\set\D(\A(X))$!$\D_{?}(X)\subset\D(X)$ (homology in $\A_{?}(X)$)}
be the full subcategory of~$\>\D_{?}(X)$ whose objects are the complexes 
$\cFb$ with homology $H^m(\cFb)$ in~$\A_{?}(X)$, and 
$\D_{?}^+(X)\subset\D_{?}(X)$ 
\index{ $\D(\A)$ (derived category of $\A$-complexes)!$\D(X)\set\D(\A(X))$!$\D_{?}^+(X)\subset\D_{?}(X)$ (homologically\vspace{.4pt} bounded\\ below)}
(resp.~$\D_{?}^-(X)\subset\D_{?}(X))$
\index{ $\D(\A)$ (derived category of $\A$-complexes)!$\D(X)\set\D(\A(X))$!$\D_{?}^-(X)\subset\D_{?}(X)$ (homologically bounded\\ above)}
the full subcategory whose objects are the complexes 
$\cFb \in \D_{?}(X)$ such that $H^m(\cFb)$
vanishes for all $m \ll 0$ (resp.~$m \gg 0$).\vspace{1pt}

\item\label{conv3} 
\makebox[0pt]{\raisebox{-1ex}{\hspace{-6em}\;{\Huge $\mathsf{Z}$}}}
\kern-2.3pt For a formally smooth local homomorphism 
$\phi \colon (A, \mathfrak m_A) \to (B,\mathfrak  m_B)$ 
of noetherian complete local rings, the usual definition 
of the relative dimension of~$\phi$, viz., $\dim(B/\mathfrak m_AB)$,
is inconsistent with the definition  in~\ref{def:diff4a} of
relative dimension  for 
the induced formal scheme-map \mbox{$\Spf(B) \to \Spf(A)$}. 
The meaning of the term `relative dimension'
will therefore depend on whether the map under consideration
is in the category of local rings or that of formal schemes.\vspace{1pt}

\item\label{conv4} 
We adopt the usual sign convention for the differential~
$d$ in the tensor product of 
two complexes $A^\bullet$, $B^\bullet$
(over rings, ringed spaces, etc.), expressed symbolically by
$$
d^n|(A^p\otimes B^{n-p})=
d^p_A\otimes 1 + (-1)^p\otimes d^{n-p}_B.
$$

By way of convention for how to connect
$\otimes = \otimes_{\<\<X}$ with the translation functor on complexes, 
noting that
$A^{p+i} \otimes B^{q+j}$ occurs as a direct summand 
in degree~$\>p+q\>$ of both $\Ab[i] \otimes \Bb[\>j]$ 
and $(\Ab \otimes \Bb)[i+j]$, we choose the unique 
isomorphism of complexes
$$
\theta=\theta_{i,\>j}^{A,B} \colon \Ab[i] \otimes \Bb[\>j] \iso 
  (\Ab \otimes \Bb)[i+j]\qquad(i,j \in \mathbb Z)
$$
\pagebreak[2]
such that for any $p,q \in \mathbb Z$,
\[
\theta_{i,j}\big|(A^{p+i} \otimes B^{q+j}) = 
\text{multiplication by $(-1)^{pj}$ }.
\]

\item\label{conv5} In particular, \eqref{conv4} 
applies when $\Ab$ and  $\Bb$ are concentrated in
degree 0. Thus when  $A$ and $B$ are $\cO_{\<\<X}$-modules,  
our choice for a natural isomorphism of complexes 
$A[i] \otimes B[\>j] \iso (A \otimes B)[i+j]$
is given in degree $-i-j$ by $(-1)^{ij}$ times the identity map
of~$A \otimes B$. (In (iv), take $p=-i$, $q=-j$.)\vspace{1pt}

\item\label{conv6} Recall that a $\delta$-functor%
\index{ $\Cfr$02@$\delta$-functor} between 
two triangulated categories $\D_1,\D_2$, with translation functors 
$T_1,T_2$ respectively, is a pair $(F,\Theta)$ consisting of
an additive functor $F \colon \D_1 \to \D_2$
and a natural isomorphism $\Theta \colon FT_1 \iso T_2F$ 
such that for any triangle
$A \xto{\;u\;} B \xto{\;v} C \xto{\;w\;} T_1A \;$  in $\D_1$, 
the corresponding diagram
\[
FA \xto{\;Fu\;} FB \xto{\;Fv} FC \xto{\;\Theta \;\circ\; Fw\;} T_2FA
\] 
is a triangle in $\D_2$. Explicit reference to $\Theta$ is 
frequently suppressed once $\Theta$ has been specified.\vspace{1pt}

\item\label{conv7} With notation as in \eqref{conv4}, for a fixed
$\Bb \in \bC(X)$ where $\Bb$ consists of flat $\cO_X$-modules, 
the functor sending $\Ab$ to $\Ab \otimes \Bb$
induces a functor from $\D(X)$ to itself, thus
yielding a $\delta$-functor (see \eqref{conv6}) 
where $\Theta$ is given by the identity 
map $= \theta_{1,0}$ from \eqref{conv4}. Similarly, if we 
fix $\Ab$ as a complex of flat $\cO_X$-modules then the
functor sending $\Bb$ to $\Ab \otimes \Bb$, induces a functor
from $\D(X)$ to itself which also upgrades to a $\delta$-functor.
However in this case $\Theta = \theta_{0,1}$ is \emph{not} the 
identity map.\vspace{1pt}

\item\label{conv8} 
For a complex $\cFb\in\bC(\A)$ (see (i)), we have
$(\cFb[n])^{i} = \F^{i+n}$, and we can identify
the submodule of $i$-cocycles in $\cFb[n]$ with that of $i+n$-cocycles 
in $\cFb$, and similarly for  coboundaries.  
Accordingly, 
we make the identification
$H^i(\cFb[n]) = H^{i+n}(\cFb)$ without introducing
any signs.\vspace{1pt}

\item\label{conv9} Let $\A$ be an abelian category.
For an exact sequence of $\A$-complexes
\[
0 \xto{\quad} \Ab \xto{\quad} \Bb \xto{\quad} \Cb \xto{\quad} 0,
\]
the connecting homomorphism $H^i\Cb \to H^{i+1}\Ab$ is the usual
one described via ``chasing elements.'' This 
connecting map, modulo 
the identification $H^i(\Ab[1]) = H^{i+1}\Ab$,  
is $(-1)$ times the map obtained by 
applying~$H^i$ to the derived category map 
$\Cb \to \Ab[1]$ in the 
triangle associated to the above exact sequence.\vspace{1pt}

\item\label{conv10} For a complex $\cFb$ we use the following truncation
operators 
{\small{
\[
\begin{array}{rrcrcrcrcr}
\sigma_{\ge p}\cFb \set &\cdots \xto{\quad} &0 &\xto{\quad} &\F^{p} 
&\xto{\quad} &\F^{p+1} &\xto{\quad} &\F^{p+2} & \xto{\quad} \cdots  \\ 
\sigma_{\le p}\cFb \set &\cdots \xto{\quad} &\F^{p-1} &\xto{\quad} 
&\F^{p} &\xto{\quad} &0 &\xto{\quad} &0 & \xto{\quad} \cdots  \\ 
\tau_{\ge p}\cFb \set &\cdots  \xto{\quad} &0 &\xto{\quad} 
&\text{coker $d^{p-1}$} &\xto{\quad} 
&\F^{p+1} &\xto{\quad} &\F^{p+2} &\xto{\quad} \cdots  \\ 
\tau_{\le p}\cFb \set &\cdots \xto{\quad} &\F^{p-1} &\xto{\quad} 
&\ker d^{p} &\xto{\quad} &0 &\xto{\quad} &0 &\xto{\quad} 
\cdots 
\end{array}
\]
}}%
where $d^{\bullet}$ is the differential in $\cFb\<$;
and $\sigma_{> p}\set \sigma_{\ge p+1}$, 
$\sigma_{< p} \set \sigma_{\le p+1}$, etc.\vspace{1pt}%
\index{ OOO@$\sigma_{\ge}$, $\sigma_{>}$, $\sigma_{\le}$, $\sigma_{<}$ (truncation functors)}%
\index{ OOO@$\tau_{\ge}$, $\tau_{>}$, $\tau_{\le}$, $\tau_{<}$ (truncation\vspace{1pt} functors)}

These operators induce functors $\D(\A)\to\D(\A)$.
\end{enumerate}

\newpage

\section{Preliminaries on formal schemes}
\label{sec:prelim}
In this section we recall and develop some subsequently-needed basic notions pertaining
to formal schemes, their morphisms, and their (sheaves of) modules. For definitions and properties of formal schemes, 
see \cite[\S10]{EGAI}. \emph{Unless otherwise indicated, all formal schemes will be understood to be noetherian.}\vspace{1pt}

In \S\ref{subsec:formal} we define  \emph{essentially pseudo-finite-type maps} of formal schemes, a mild generalization
of the notion of pseudo-finite-type maps treated in 
\cite{AJL2} and \cite{Ye}; and we discuss the behavior of \emph{codimension functions.} 
In \S\ref{subsec:mod} we discuss \emph{quasi-coherent 
modules} over formal schemes, especially in relation to modules over noetherian adic rings. In \S\ref{subsec:tor} we discuss \emph{torsion modules} over formal schemes, an important class which includes arbitrary modules over an 
ordinary scheme. In \S\ref{subsec:morphism} we discuss \emph{smooth maps}---those which are formally smooth and essentially of pseudo-finite type. A key property (\ref{cor:morph1b}) of essentially pseudo-finite-type maps is that they factor locally as
(smooth)$\smcirc$(closed immersion). In the remaining two subsections we discuss  \emph{modules of continuous differentials} relative to essentially pseudo-finite-type maps, first  of noetherian adic rings, and then of formal schemes.
 
 \enlargethispage*{0pt}
 
\subsection{Codimension functions and maps of formal schemes} 
\label{subsec:formal} 

A topological ring $R$ is \emph{adic}\index{adic!ring} 
if there exists an $R$-ideal $I$ whose powers 
form a base of neighborhoods of $(0)$, and if $R$ is complete and Hausdorff. 
Any such  $I\>$ is called an {\it ideal of definition
of\/ $R$} (or a {\it defining ideal of\/ $R$}).
\index{defining ideal $(={}$ideal of definition)}%
\index{ideal of definition $(={}$defining ideal)}
If  $\XOX$ is a  formal scheme  then there exists  a coherent $\OX$-ideal~$\I$  whose sections over any affine 
open set $\U$ of $\X$ form a defining ideal for the 
adic ring $\Gamma(\U,\OX)$; any such $\I$ is called an {\it ideal of definition of\/ $\XOX$} (or simply of~$\X$)  (\cite[\S\S10.5, 10.10]{EGAI}). The
ringed space $ (\X, \OX/\I)$ is then a \emph{noetherian ordinary scheme,} that is, a noetherian formal scheme with discretely topologized structure sheaf, or equivalently, with $(0)$ as a defining ideal.

The notion in \cite[\S10.13]{EGAI} of  finite-type maps of formal schemes  
is not adequate for our purposes, and instead we consider a
slight generalization of the notion in~\cite{AJL2} and~\cite{Ye}
of pseudo-finite type maps. A homomorphism 
$f\colon A \to B$ of noetherian adic rings is  
{\em $(\<$essentially\/$)$ of pseudo-finite type\/} if it is continuous and if for one---hence
any---defining ideal $\bfr \subset B$
the composition $A \to B \to B/\bfr$ is (essentially) of finite type, i.e.,
$B/\bfr$ is (a localization of) a finite-type $A$-algebra. 
\index{pseudo-finite type!homomorphism}
\index{essentially of pseudo-finite type!homomorphism}

A map of not-necessarily-noetherian ordinary schemes $f_0 \colon X \to Y$ is 
{\em essentially of finite type\/} if every~$y \in Y$ 
has an affine open neighborhood $V = \Spec(A)$ such 
that $f_0^{-1}V$ can be covered by finitely many affine 
open $U_i = \Spec(C_i)$ such that the corresponding maps 
$A \to C_i$ are essentially of finite type. 
For any morphism of  formal schemes  
$ f \colon \XOX \to \YOY $, there exist ideals of definition
$ \I \subset \cO_{\X}$ and $\J \subset \cO_{\Y} $ 
satisfying $ \J \cO_{\X} \subset \I $ (\cite[10.6.10]{EGAI}); and  correspondingly
there is an induced map of ordinary schemes $f_0\colon (\X, \cO_{\X}/ \I) 
\to (\Y, \cO_{\Y}/ \J) $ (\cite[10.5.6]{EGAI}).
We say~$f$ is ({\em essentially\/}) {\it of pseudo-finite type}  
if $f_0$ is ({essentially})  of finite type.\vspace{1pt}
\index{pseudo-finite type!formal-scheme map}
\index{essentially of pseudo-finite type!formal-scheme map}

This property of maps is independent of the choice of defining 
ideals $\I,\J$. It behaves well with respect to composition and base change:  if $f\colon\X\to \Y$ and \mbox{$g\colon\Y\to\Z$} are formal-scheme maps,
and if both $f$ and $g$ are essentially of pseudo-finite type, then so is the composition~$gf$; and if $gf$ and $g$ are  essentially of pseudo-finite type then so is~$f$; moreover, if $f\colon\X\to \Y$ is essentially of pseudo-finite type and $\Y'\to\Y$ is any map of formal schemes, then $\X'\set\Y'\times_\Y \X$ is noetherian and the projection 
$ \X'\to\Y'$ is essentially of pseudo-finite type  (cf.~ proof of \ref{lem:form2}).

{\makebox[0pt]{\raisebox{-1ex}{\hspace{-6em}\;{\Huge $\mathsf{Z}$}}}}
We don't know for an arbitrary formal-scheme map 
$\U = \Spf(B) \to \V = \Spf(A)$ that is 
essentially of pseudo-finite type, whether the corresponding 
g $A \to B$ is essentially of pseudo-finite type. \vspace{1pt}

A  formal scheme $(\X, \cO_{\X})$ is called 
{\it \textup{(}universally\/\textup) catenary\/}%
\index{catenary (formal scheme)}%
\index{universally catenary (formal scheme)}
if there exists one defining ideal $ \I $ such that the 
scheme $ (\X, \cO_{\X}/ \I)$ is (universally) catenary.
This implies that for {\it any\/} defining ideal $ \J \subset \cO_{\X}$
the scheme $ (\X, \cO_{\X}/ \J)$ is (universally) 
catenary---an easy consequence of the fact that on a  
formal scheme, for any two defining ideals $\I_1$ and $\I_2$ there
exists an integer~$n$ such that $\I_1^n \subset \I_2$.

A {\it codimension function}%
\index{codimension function}
on the underlying topological space $\X$
of a formal scheme is a function $\Delta \colon \X \to \mathbb Z $ 
such that $\Delta (x^{\prime}) = 
\Delta (x) + 1 $ for every immediate specialization%
\index{specialization}\index{specialization!immediate}%
\index{immediate specialization}
$ x \leadsto x^{\prime} $ of points in $\X$.%
\footnote{A specialization $ x \leadsto x^{\prime}$---i.e., $x'$ is in the closure of $\{x\}$---is called \emph{immediate} if 
$x \ne x^{\prime}$ and there are no  specializations $ x \leadsto x^{\prime \prime} \leadsto x^{\prime} $ with $x\ne x''\ne x'$.}
If $\X$ admits a codimension function, then $\X$ has to be catenary.
When $\X$ is connected and noetherian, any two codimension functions 
on~$\X$ differ by a constant.
  
\begin{aexam}
\label{exam:formal0}
If $\X$ is catenary and irreducible, then for any integer~$n$ there is a unique
codimension function on~$\X$ assigning $n$ to the generic point.
If $\X$ is catenary and biequidimensional (\cite[14.3.3]{EGAOIV})
then we can assign one fixed integer to all the closed points of $\X$,
and this extends uniquely to a codimension function.
\end{aexam}

\begin{aexam}
\label{exam:formal1}
Let $ f\colon \X \to \Y $ be a map of 
formal schemes that is essentially of pseudo-finite type. 
If~$\Y$ is universally catenary then so is~$\X$.
Furthermore, if~$\Y$ admits a codimension function~$\Delta$ then 
\emph{the function $f^{\sharp} \Delta$ on $\X$ such that}
\index{ $\Forget$4@$f^{\sharp} \Delta$ (lifted codimension function)}
\[ 
\bigl(f^{\sharp} \Delta\bigr)(x) = \Delta(y) - \mbox{tr.deg.}_{k(y)}k(x) 
\qquad\bigl(x\in\X,\;y\set f(x)\bigr)
\]
(where $\mbox{tr.deg.}_{k(y)}k(x)$  is the 
transcendence degree of the residue field extension $k(x)/k(y)$)
\emph{is a codimension function on} $\X$, as follows from the dimension formula (\cite[5.6.2]{EGAOIV}).
\end{aexam}

In particular, from \ref{exam:formal1} we see that any 
formal scheme that is essentially of pseudo-finite type 
over a field $k$ admits a 
codimension function. Furthermore if $f \colon \X \to \Spec (k)$
is essentially of pseudo-finite type and if~$\Delta$ is the 
function on~$\Spec (k)$ sending the unique point to~0 then 
$f^{\sharp} \Delta(x)=0$ if and only if $x$ is a closed point 
of~$\X$.

\smallskip
Let $\mathbb F$ be the category whose
objects are the (noetherian) universally catenary formal schemes 
admitting a codimension function, and whose morphisms are the formal-scheme maps which are essentially of pseudo-finite type. 
Let $\bbFc$ be the category whose
objects are the pairs $(\X,\Delta)$ with $\X\in\mathbb F$ and
$\Delta$ a codimension function on $\X$,  and whose
morphisms $(\X^{\prime},\Delta^{\prime}) \to (\X,\Delta)$ are the  
$\mathbb F$-morphisms
$f \colon \X^{\prime} \to \X$  such that 
$\Delta^{\prime} = f^{\sharp}\Delta$. (Note that the formula 
in~\ref{exam:formal1} behaves well with respect to compositions
in~$\mathbb F$.) The next Lemma is about fiber products in these categories. (See \cite[10.7.3]{EGAI} for the construction of fiber products of formal schemes).

\medskip

\begin{alem}
\label{lem:form2}
\textup{(a)} If\/ $f\colon\X\to\Z,$ $g\colon\Y\to\Z$ are\/ $\mathbb F$-morphisms and\/  
$\W\set\X\times_\Z\Y,$ with projections\/ $\W\xto{p} \X, $ $\W\xto{q} \Y,$  then\/ $(\W,p,q)$ is an\/ $\mathbb F$-fiber product of\/ $f$ and\/~$g$.
\vspace{1pt}

\textup{(b)}  Moreover,\vspace{.6pt}
if\/ $\Delta$ is a codimension function on\/~$\Z,$ then\/ 
$p^\sharp f^\sharp\Delta=q^\sharp g^\sharp\Delta=\textup{(say)}\:\Delta^\times\<,
$
and\/ $\bigl((\W,\Delta^\times\<),\>p,q\bigr)$ is an\/ $\bbFc$-fiber product of\/ 
$f\colon(\X,f^\sharp\Delta)\to(\Z,\Delta)$ and\/ $g\colon(\Y,g^\sharp\Delta)\to(\Z,\Delta)$.
 \end{alem}
 
 \pagebreak[2]
 
\begin{proof}
Let $\I,\J,\cK$ be defining 
ideals for $\X,\Y,\Z$ respectively, so that $\W$ has the defining ideal $\cL \set \I\OW + \J\OW$.  If $\X_0,\Y_0,\Z_0,\W_0$
are the schemes obtained from $\X,\Y,\Z,\W$ by going modulo the
corresponding defining ideals $\I,\J,\cK,\cL$, then
$\W_0\set \X_0 \times_{\Z_0} \Y_0$.
The natural map $\W_0 \to \Y_0$ is essentially of finite type,
so the induced map $q \colon \W \to \Y$ is 
essentially of pseudo-finite type. Moreover $\W_0$ is noetherian,
hence so is~$\W$ (\cite[10.6.4]{EGAI}). 
Now  \ref{exam:formal1} shows that $\W\in \mathbb F$, and the rest is straightforward.  
\end{proof}

\medskip

\subsection{Quasi-coherent modules}
\label{subsec:mod}

We recall some basic facts about modules over formal schemes (assumed, as always, to be noetherian). Proofs can
be found in \cite[\S3]{Ye}, or \cite[\S3]{AJL2}.

Let $A$ be a noetherian adic ring and $ \U := \Spf (A) $ the 
corresponding affine formal scheme. Let $\Mod(A)$%
\index{ Mod (module category)}%
\index{ Mod (module category)!$\Mod_f$ (finitely generated modules)}
(resp.~$\Mod_f(A)$) denote the category of $A$-modules
(resp.~finitely generated $A$-modules). For any $A$-module $M$, 
consider the presheaf 
that assigns to any open set $\V$, the  
module $\Gamma(\V,\cO_{\U})\> {\otimes}_{\<A}\> M$. Let $M^{\sim_A}$ 
(or~$M^{\sim}$ when there is no cause for confusion) denote%
\index{   $\{\:\}^{\sim}$ ($\Avc(\Spf( A))$-sheaf associated to\\ 
the $A$-module~$\{\:\}\>$)\vspace{1pt}} 
the associated sheaf.\vspace{-.6pt} This defines a functor 
$\sim_A \colon \Mod(A) \to \A(\U)$ 
(see \eqref{conv2} of \S\ref{subsec:conv}). 

In fact, if $\kappa=\kappa_A\colon\U\to U\set\Spec(A)$ is the canonical map,
and $\widetilde{M}$ is the sheaf over $\Spec(A)$ corresponding to $M\<$, 
then there is a functorial isomorphism
$$
M^\sim\iso\kappa^*\widetilde{M}.
$$
In other words, $M^\sim$ represents the functor 
$\Hom_{\cO_{\<U}^{}}\<(\widetilde M, \kappa_*\G)$ of $\cO_\U$-modules~$\G$. Indeed\-, an $\cO_\U$-homomorphism of sheaves $M^\sim\to \G$ corresponds naturally to  an \mbox{$\cO_\U$-homomorphism} of
presheaves $\Gamma(\V,\cO_{\U})\> {\otimes}_{\<A}\> M\to\Gamma(\V,\G)$,
which is determined by the single $A$-homomorphism
$M\to\Gamma(\U,\G)$ obtained by taking $\V=\U$, and thus, in view of  \cite[p.\,213, Cor.\,(1.7.4)]{EGAI}, we have natural isomorphisms
$$
\Hom_{\>\cO_{\<\<\U}}(M^\sim\<,\>\G)\cong \Hom_A\bigl(M, \Gamma(\U,\G)\bigr)
\cong \Hom_{\cO_{\<U}^{}}\<(\widetilde M, \kappa_*\G).
$$
The first of these isomorphisms shows, moreover, that \emph{the
functor $\sim_A$ is left-adjoint to the functor} $\Gamma(\U, -)$ of
$\cO_\U$-modules.%
\footnote{These considerations hold for \emph{any} map of ringed
spaces $\varkappa\colon X\to\Spec(A)$. In that generality, 
\ref{prop:mod1}\kern.5pt\eqref{mod1i2} below (with $\vec c$ meaning ``$\subdirlm{}\!\!$ of finitely-presented") \vspace{-2.5pt}and \ref{lem:mod4a}  (with $f$ replaced
by~$\varkappa$) hold; and if $\varkappa$
is \emph{flat} then so do \ref{prop:mod1}\kern.5pt\eqref{mod1i1}, remark~(a)
following \ref{prop:mod1}, 
and \ref{lem:mod2a}. What distinguishes~$\kappa$ from the horde is \ref{prop:mod1}\kern.5pt\eqref{mod1i3}.}
 
\begin{aprop}[{\cite[p.\,874, Prop.\,3.2]{Ye},\;\cite[p.\,31, Prop.\,3.1.1]{AJL2}}]
\label{prop:mod1} 
With the preceding notation, and\/ $\sim\ \set\;\sim_A\>,$
\begin{enumerate}
\item\label{mod1i1}The functor\/ $\sim$ is exact and commutes with 
direct limits.\vspace{1pt}
\item\label{mod1i2}For any\/ $A$-module\/ $M,$ 
$M^{\sim} \in \Avc(\U)$ and is quasi-coherent.\vspace{1pt}
\item\label{mod1i3}The functor\/ $\sim$ is an equivalence of 
categories between\/ $\Mod(A)$ and\/ $\Avc(\U),$ with quasi-inverse\/ 
$\Gamma(\U,\<-) \colon \<\Avc(\U)\< \to \Mod(A).\!$ 
These quasi-inverse equivalences restrict to quasi-inverse
equivalences 
between\/ $\Mod_f(A)$ and\/ $\Ac(\U)$.\vspace{1pt}
\item\label{mod1i4} For any affine open\/ $\V\subset\U$
and any\/ $A$-module $M,$ the natural map is an isomorphism
\[ 
\Gamma(\V,\cO_{\U}) \otimes_A M \iso \Gamma(\V,M^{\sim}). 
\]
\end{enumerate}
\end{aprop}

\emph{Remarks.} (a) We can (and will) identify $A^\sim$ with $\cO_\U$; and for any $A$-ideal $I$ we  can (and will) identify $I^\sim$ with its natural image in $A^\sim\<$, so that $I^\sim$ will be regarded as an $\cO_\U$-ideal.  In view of \eqref{mod1i1}, we will also identify $(A/I)^\sim$ with 
$\cO_\U/I^\sim\<.$\vspace{1pt}

(b)  If $I$ is a defining ideal of the adic ring $A$ then 
$I^\sim$ is a defining ideal of $\U$.\vspace{1pt}

(c) For any  formal scheme~$\X$,  (ii) and~(iii) in \ref{prop:mod1}
 imply that
$\Avc(\X) \subset \Aqc(\X)$, as it is enough to check 
this locally. Equality holds if~$\X$ is an ordinary scheme.\vspace{1pt}

(d) When $M$ is finitely generated, $M^\sim$ is the same as $M^\Delta$
in \cite[\S10.10]{EGAI}, so the second part of (iii) is identical with
\emph{loc.\,cit.,} Thm. (10.10.2).\vspace{1pt}

(e) Note that (iv) follows from (iii), because with
$B\set\Gamma(\V,\cO_{\U})$ we have 
$(B\otimes_A M)^{\sim_{\<B}}=M^{\sim_{\<\<A}}|_\V^{}\>$.

\begin{aprop}[{\cite[p.\,875, Cor.\,3.4]{Ye}, or \cite[p.\,32, Cor.\,3.1.4]{AJL2}}]
\label{prop:mod2}
Let\/ $\X$ be a formal scheme and\/ $x \in \X$.
For any quasi-coherent\/ $\cO_{\X}$-module~$\M$ 
there exists an affine open neighborhood\/
$ \U := \Spf (A) $ of\/~$x$ such that the natural map is an isomorphism\/
$\Gamma (\U,\M)^{\sim_{\<\<A}} \iso \M {\big|}_{\U}$.
\end{aprop}  

\begin{alem}
\label{lem:mod2a}
Let $\U := \Spf(A)$ be an
affine formal scheme and let~$M,N$ be $A$-modules. Set\/ $\sim\ \set\;\sim_{\<\<A}.$
\begin{enumerate}
\item If\/ $N$ is finitely generated 
then the map\/ $\varphi$ corresponding  to the natural map\/ 
$\Hom_A(N,M) \to \Hom_{\>\>\cO_{\<\<\U}}\<\<(N^{\sim},\> M^{\sim})$
is an isomorphism
\[
\Hom_A(N,M)^{\sim} \iso \sHom_{\>\>\OU}(N^{\sim}\<, M^{\sim}).
\]
\item The map 
corresponding to the natural map 
$$
M \otimes_A N \cong M^{\sim}(\U) \otimes_{A} N^{\sim}(\U)
\lra (M^{\sim}\< \otimes_{\>\cO_{\<\<\U}} N^{\sim})(\U)
$$
is an isomorphism
$$
(M \otimes_A N)^{\sim} \iso M^{\sim} \otimes_{\>\cO_{\<\<\U}}\! N^{\sim}.
$$
\item
For any\/ $A$-ideals\/~$I,J$ the\/ $\cO_\U$-ideals \/ $(IJ)^\sim$ and\/ 
$I^\sim J^\sim$ coincide.\vspace{1pt}
\item For any\/ $m>0$ the natural map is an isomorphism
$$
\big(\bw^m_AM\big)^\sim\iso \bw^m_{\>\cO_{\<\U}}M^\sim.\\[-2pt]
$$
\end{enumerate}
\end{alem}  
\begin{proof}
(i). Let $A^i \to A^j \to N \to 0$ be a presentation of $N$. 
In the following natural commutative diagram 
{\small{
\[ 
\begin{CD}
0 @>>> \Hom_A(N,M)^{\sim} @>>> \Hom_A(A^j\<,M)^{\sim} @>>> 
\Hom_A(A^i\<,M)^{\sim}  \\
@. @V\varphi VV  @VV{\cong}V @VV{\cong}V \\
0 @>>> \sHom_{\>\cO_{\<\<\U}}(N^{\sim}\<,M^{\sim}) @>>> 
\sHom_{\>\cO_{\<\<\U}}(A^{j\sim}\<,M^{\sim}) @>>> 
\sHom_{\>\cO_{\<\<\U}}(A^{i\sim}\<,M^{\sim})  
\end{CD}
\]
}}%
the rows are exact by \ref{prop:mod1}\kern.5pt (i).
Since~$\sim$ commutes with direct sums, the functors 
$\Hom_A$ and  $\sHom_{\>\cO_{\<\<\U}}$ commute with finite direct sums, and since the
natural map $\Hom_A(A,M)^{\sim} \to \sHom_{\>\cO_{\<\<\U}}(A^{\sim}, M^{\sim})$
is an isomorphism,  the two 
vertical arrows on the right are isomorphisms, whence so is $\varphi$.
\vspace{1pt}

(ii). Proceed as in (i), using a presentation of~$N$;\vspace{.6pt} or simply use that with
$Q=M$ or $N\<$, one has ${Q}^\sim=\kappa^*(\widetilde Q)$.\vspace{1pt}

(iii). Since $\sim$ is exact, and so ``commutes with image," the assertion amounts to equality of
the images of the natural maps $\mu$ and $\hat\mu$ in
the following diagram:
$$
\begin{CD}
(I\otimes_A J)^\sim @>\mu>> A^\sim \\
@V\cong V\rm(ii)V @| \\
I^\sim\otimes_{\cO_\U^{}} J^\sim @>>\under{1.2}{\hat\mu}> \cO_\U
\end{CD}
$$
It is clear that the diagram  commutes, whence the
conclusion.\vspace{1pt}

(iv) The assertion results, upon passage to associated sheaves, from the natural 
presheaf isomorphisms (for affine open $\V\subset \U$):
$$
\Gamma(\V,\cO_{\U}) \otimes_A \>\bw^m\< M 
\cong\bw^m\<\big(\Gamma(\V,\cO_{\U}) \otimes_A M  \big)
\>\underset{ \rm\ref{prop:mod1}\kern.5pt(iv)}\iso\,
\bw^m\<\big(\Gamma(\V,M^\sim)\big).
$$
\end{proof}
\pagebreak[3]

\begin{alem}
\label{lem:mod4a}
Let\/ $\V := \Spf(B) \xto{\; f \;} \U := \Spf(A)$ be a map of 
formal schemes. Then for any\/ $A$-module $M\<,$
the  map\/ $\phi$ corresponding to the natural composition
\begin{align*}
M \otimes_A B \iso
\Gamma(\U,  M^{\sim_{\<\<A}}\<)\otimes_A B
&\to
\Gamma(\U, f_*f^* (M^{\sim_{\<\<A}}\<))\otimes_A B\\
&\to
\Gamma(\V,f^*(M^{\sim_{\<\<A}}\<))\otimes_A B
\to
\Gamma(\V,f^*(M^{\sim_{\<\<A}}\<))
\end{align*}
is an isomorphism
\[
(M \otimes_A B)^{\sim_B} \iso f^*(M^{\sim_{\<\<A}}\<).
\]
\end{alem}
\begin{proof}
Let $A^I \to A^J \to M \to 0$ be a presentation of $M$. 
In the following diagram of natural induced maps, the top row is
obtained by applying $(- \otimes_A B)^{\sim_B}$ to the presentation
of~$M$ while the bottom row is obtained by applying $f^*(-^{\sim_{\<\<A}})$.
\[ 
\begin{CD}
(B^I)^{\sim_B} @>>> (B^J)^{\sim_B} @>>> (M \otimes_A B)^{\sim_B} @>>> 0 \\
@V{\cong}VV  @V{\cong}VV @VV\phi V \\
f^{*}(A^{I})^{\sim_{\<\<A}}  @>>> f^{*}(A^{J})^{\sim_{\<\<A}} 
@>>> f^*M^{\sim_{\<\<A}} @>>> 0   
\end{CD}
\]
The rows are exact, by \ref{prop:mod1}\kern.5pt (i).  Since $\sim_{\<\<A}$, $\sim_B$ and $f^*$ commute with direct sums, and the natural map 
$(A \otimes_A B)^{\sim_B} \lra f^*(A^{\sim_{\<\<A}})$ 
is an isomorphism, the vertical maps on the left are isomorphisms and hence so is $\phi$.

Alternatively, since $f^*$ is left-adjoint to $f_*$ one can extract the assertion from the sequence of natural isomorphisms (with $\G\in\A(\V)$)
\begin{align*}
\Hom_{\>\cO_{\<\<\V}}\<\bigl((M\otimes_A B)^{\sim_B}\<\<,\>\G\bigr)
&\iso
\Hom_B\bigl(M\otimes_A B,\>\Gamma(\V,\G)\bigr)\\
&\iso 
\Hom_A\bigl(M,\>\Gamma(\U,f_*\G)\bigr)
\iso
\Hom_{\>\>\cO_{\<\<\U}}\<\<\bigl(M^{\sim_{\<\<A}}\<\<,\>f_*\G\bigr).\\[-20pt]
\end{align*}
\end{proof}

\smallskip

Following \cite[(1.9.1)]{Li}, we say that a subcategory $\A_{?}(\X)\subset\A(\X)$ is  \emph{plump}\index{plump subcategory} 
if it is full and if for every exact sequence
$\M_1 \to \M_2 \to \M \to \M_3 \to \M_4$ in $\A(\X)$  with
$\M_1, \M_2, \M_3$ and $\M_4$ in~$\A_{?}(\X)$, $\M$ is in $\A_{?}(\X)$ too.
Then $\A_{?}(\X)$ is an \emph{abelian} subcategory of $\A(\X)$.
Moreover, $\D_{?}(\X)$ is a \emph{triangulated} subcategory of~$\D(\X)$. 

\makebox[0pt]{\raisebox{-1ex}{\hspace{-6em}\;{\Huge $\mathsf{Z}$}}}
\kern-2.3pt(However, the natural functor  $\D(\A_{?}(\X))\to\D_{?}(\X)$ need not be an equivalence.)
\medskip

\begin{aprop}[{\cite[p.\,34, Prop.\,3.2.2]{AJL2}}]
\label{prop:mod3}
For any formal scheme\/ $\X,$\vspace{.5pt} the 
subcategories\/ $\Avc(\X)$ and\/~$\Aqc(\X)$ of\/~$\A(\X)$ are  
plump.
\end{aprop}  

Let $\U = \Spf(A)$ be an affine formal scheme and $x$ a point 
in~$\U$. Let $\pfr$ be the open prime ideal in the adic ring $A$ 
corresponding to $x$ and let $m_x$ denote the maximal ideal 
of the local ring $\OUx$. Let $\afr$ be a defining ideal in $A$
and $\I \set \afr^{\sim_A}$ the corresponding defining ideal in $\OU$.
Recall that $\OUx$ is a direct limit of the 
rings $A_{\{f\}}$ for $f \notin \pfr$
where $A_{\{f\}}$ is%
\index{   $A_{\{S^{-1}\}}$ (completed localization)!$A_{\{f\}}$ ($S=\{f,f^2,f^3,\dots\}\>$)}  
the completion of $A_f$ along $\afr_f$.
In particular there is a natural map $A_{\pfr} \to \OUx$.
For any open ideal $J$ in $A$, with $\J \set J^{\sim_A}$,
the natural induced map $A_{\pfr}/J_{\pfr} \to \OUx/\J_x$
is an isomorphism. Since any power of $J,\J$ is open,
there results a canonical isomorphism
$(A_{\pfr}, J_{\pfr})^{\,\wh{}} \iso (\OUx, \J_x)^{\,\wh{}}$.
In particular, we can associate canonically to 
$x \in \U$ the following faithfully flat inclusions of 
noetherian local rings
\[
A_{\pfr} \hookrightarrow \OUx \hookrightarrow B \hookrightarrow C,
\]
where $B = (A_{\pfr}, \afr A_{\pfr})^{\,\wh{}} \cong (\OUx, \I_x)^{\,\wh{}}$ 
and $C = (A_{\pfr}, \pfr A_{\pfr})^{\,\wh{}} \cong (\OUx, m_x)^{\,\wh{}}
\cong (B, m_B)^{\,\wh{}}$. 

\medskip

\subsection{Torsion modules}
\label{subsec:tor}

Let $A$ be a ring,  $\afr$ an $A$-ideal. For any $A$-module~$M\<$, let
$\iG{\afr}M$%
\index{ $\Hom^{\cont}$2@$\iG{}$ \kern-2pt (torsion functor)}%
\index{ $\Hom^{\cont}$2@$\iG{}$ \kern-2pt (torsion functor)!$\iG{\I}$@$\iG{\afr}M$ ($\afr$ an ideal)}
be the submodule consisting of those elements which are 
annihilated by a power of~$\afr$.  
Note that $\iG{\afr}M$ is naturally isomorphic to\vspace{-4pt}
$\dirlm{}_{\!\!n}\, \Hom_{A}(A/{\afr}^n\<, M)$. We say that $M$ is an 
{\em $\afr$-torsion module\/} if $\iG{\afr}M = M.$
If~$A$ is a topological ring whose
topology is defined by powers of $\afr$ then an $\afr$-torsion module is
referred to simply as a {\em torsion module}.
\index{torsion module}

\begin{alem}
\label{lem:mod4b}
Let\/ $A$ be a noetherian adic ring and\/ $M$ a torsion $A$-module. 
For any multiplicatively closed set $S \subset A$  let $A\{S^{-1}\}$ be%
\index{   $A_{\{S^{-1}\}}$ (completed localization)}  
the completion of the 
localization $A_S$ along a defining\/ $A$-ideal. Then the natural map 
is an isomorphism 
$$
M_S \iso M \otimes_A A\{S^{-1}\}.
$$
\end{alem}

\begin{proof}
Let $N$ be a finitely generated submodule of $M$. Then
the torsion module~$N$ is annihilated by some 
defining ideal, say $\afr$, in $A$.
Hence we have the following isomorphisms 
\[
N_S \iso N \otimes_A A_S/(\afr) \iso N \otimes_A A\{S^{-1}\}/(\afr) \osi
N \otimes_A A\{S^{-1}\}.
\]
Taking direct limit over all finitely-generated submodules of $M$
gives the result. 
\end{proof}

Let $(X, \cO_{\<\<X})$ be a ringed space. For any $\cO_{\<\<X}$-ideal 
$\I$ and any $\M \in \A(X)$, set
$$ 
 \iG{\I}\< \M := \dirlm{n} \sHom_{\cO_{\<\<X}}(\cO_{\<\<X}/{\I}^n, \M).\\[-3pt] 
$$
We\index{ $\Hom^{\cont}$2@$\iG{}$ \kern-2pt (torsion functor)!$\iG{\I}$ ($\I$ a sheaf of ideals)}
regard $ \iG{\I}$ as a subfunctor of the identity functor on 
$\cO_{\<\<X}$-modules. As such it is idempotent and left exact. 
Since taking direct limits commutes with restriction to open 
 $\U\subset X$,
therefore $(\iG{\I}\<\M)|_{\U}^{} \cong \iG{\I|_{\U}^{}}\!(\M|_{\<\U}^{})$.\vspace{1pt}

For a  formal scheme~$\X$ with ideal of definition~$\I$, 
we set $\iGp{\X} :=  \iG{\I}$,%
\index{ $\Hom^{\cont}$2@$\iG{}$ \kern-2pt (torsion functor)!$\iG{x}$@$\iGp{\X}$ ($\X$ a formal scheme)}
this definition being independent of the choice of the defining ideal~$\I$. 
We call $\M \in \A(X)$ a \emph{torsion $\OX$-module} if 
$\iGp{\X} \M = \M$. Thus for any $\N \in \A(\X)$,
$\iGp{\X} \N $ is the largest torsion
submodule of $\N\<$. Note that $\N$ is a torsion $\OX$-module
if and only if there is an open cover $\{\U_i\}$ of $\X$ for
which each $\N|_{\U_i}^{}$ is a torsion $\cO_{\U_i}$-module.

\pagebreak[3]

\begin{alem}
\label{lem:mod5a}
Let\/ $\X$ be a  formal scheme. 
For any injective\/ $\A(\X)$-module\/~$\sL$, $\iGp{\X}\sL$ is a flasque sheaf.
\end{alem}
\begin{proof}
For any defining ideal $\I$  of $\OX$, 
 $\iGp{\X}\sL$ is a direct limit\vspace{.6pt} over the
noetherian  space $\X$ of the flasque sheaves 
$\sHom(\OX/\I^n\<,\> \sL)$,   and so is  flasque.
\end{proof}

\pagebreak[3]

\begin{alem}
\label{lem:mod4c}
Let\/ $\U = \Spf(A)$ be an affine  formal scheme,
and let\/ $M$ be an $A$-module. Let\/ $I$ be an\/ $A$-ideal  and\/
$\I \set I^{\sim}$ the corresponding $\OU$-ideal.
Then there is a natural isomorphism
\[
\iG{\I}M^{\sim} \iso (\iG{I}M)^{\sim},
\]
whose composition with the natural injection\/ $j\colon (\iG{I}M)^\sim\to M^\sim$
is the inclusion\/ $i\colon \iG{\I}M^{\sim}\into M^\sim$. In particular, $M$ is a torsion\/ $A$-module 
$\iff$ $M^{\sim}$ is a torsion\/ $\OU$-module.
\end{alem}
\begin{proof}  By remark (a) following \ref{prop:mod1}, and by \ref{lem:mod2a}\kern.5pt (iii), one can identify $(A/I^n)^\sim$ and $\cO_\U/\I^n$.
One checks then that the composed isomorphism 
\begin{align*}
\iG{\I}M^{\sim}=
\dirlm{n}\sHom_{\>\cO_{\<\<\U}}(\OU/\I^n, M^{\sim})
&\underset{\ref{lem:mod2a}}\iso \;\dirlm{n}\Hom_A(A/I^n,M)^{\sim}\\
&\underset{\ref{prop:mod1}}\iso\: (\dirlm{n}\Hom_A(A/I^n,M))^{\sim}
=(\iG{I}M)^{\sim}
\end{align*}
has the required properties. In particular, $i$ is an isomorphism iff $j$ is.

Taking $I$ to be an ideal of definition in~$A$, so that $\I$ is an ideal of definition
in~$\cO_\U$ (remark (b) following \ref{prop:mod1}), one has that $i$ is an isomorphism iff $M^\sim$ is a torsion $\cO_\U$-module, and by
\ref{prop:mod1}\kern.5pt(iii), $j$ is an isomorphism iff $M$ is a torsion $A$-module. The last assertion results.
\end{proof}

\smallskip
For a formal scheme $\X$,  $\At(\X)$  denotes%
\index{  $\A(X)$ (category of $\cO_{\<\<X}$-complexes on ringed space $(X, \cO_{\<\<X})$)!$\Ac$@$\At(\X)$ (torsion subcat.\ of $\A(\X)$ on formal scheme $\X$)} 
the full subcategory of~$\A(\X)$
whose objects are the torsion modules; and $\Aqct(\X)$
$:=\Aqc(\X) \cap \At(\X)$ (respectively \mbox{$\Act(\X)$ $:=\Ac(\X) \cap \At(\X)$})%
\index{  $\A(X)$ (category of $\cO_{\<\<X}$-complexes on ringed space $(X, \cO_{\<\<X})$)!$\Ac$@$\At(\X)$ (torsion subcat.\ of $\A(\X)$ on formal scheme $\X$)!$\Act(\X)$ ($:=\Ac(\X) \cap \At(\X)$)} 
denotes the full subcategory of~$\A(\X)$
whose objects are the quasi-coherent torsion modules (resp.~coherent
torsion modules). 
\index{  $\A(X)$ (category of $\cO_{\<\<X}$-complexes on ringed space $(X, \cO_{\<\<X})$)!$\Ac$@$\At(\X)$ (torsion subcat.\ of $\A(\X)$ on formal scheme $\X$)!$\Aqct(\X)$ ($:=\Aqc(\X) \cap \At(\X)$)} 

If~$\X$ is an ordinary  scheme 
(i.e., $(0)$ is a defining ideal), then $\At(\X) = \A(\X)$ and 
$\Aqct(\X) = \Aqc(\X) = \Avc(\X)$.   

\medskip

\begin{aprop}
\label{prop:mod5}
With the preceding notation$\>:$
\begin{enumerate}
\item\label{mod5i1}For any\/ $\M \in \Aqc(\X),$ it holds that\/ 
$\iGp{\X} \M \in \Aqct(\X)$.
\item\label{mod5i2}The subcategories\/ $\At(\X)$ and\/ $\Aqct(\X)$ 
of\/ $\A(\X)$ are plump, and hence are abelian categories. 
\item\label{mod5i3}$\Aqct(\X) \subset \Avc(\X)$. So 
if\/ $\X$ is affine, say\/ $\X = \Spf(A),$ and\/ 
$\M \in \Aqct(\X)$ then by\/ \textup{\ref{prop:mod1}\kern.5pt(iii),} the natural map is an isomorphism\/ 
$\Gamma(\U,\M)^{\sim_{\<\<A}} \iso \M$.
\end{enumerate}
\end{aprop}  
\begin{proof}
See \cite[\S5.1]{AJL2}. 
\end{proof}

\smallskip

We consider next a class of torsion modules 
that includes the injectives in the 
categories $\Aqct(\X)$ and~$\At(\X)$. These torsion 
modules will be the building blocks of the objects of 
main concern to us, viz., Cousin complexes. 

For an abelian group~$G$ and  a point  $x$ of the formal scheme~$\X$, let $i_{\X,\>x} \>G$ (or simply $i_x\> G$%
\index{ $\I$@$i_x$ (produces constant sheaf on $\,\ov{\<\<\<\{x\}\<\<\<}\,$)} 
if no confusion results) be the 
sheaf whose sections over any open $\U\subset\X$ are the elements  
of $G$ if $x\in\U$ and~(0) otherwise
(restriction maps being the obvious ones). For an 
$\OXx$-module $M$ the sheaf $i_xM$ has a natural 
$\OX$-module structure. Our interest lies in 
the situation where~$M$ is a zero-dimensional $\OXx$-module.

\begin{alem}
\label{lem:mod6}
Let\/ $\X$ be a formal scheme, $x\in\X$,
$m_x$ the maximal ideal of\/ $\OXx$ and\/ $M$ an\/ $\OXx$-module. 
The following are equivalent.
\renewcommand{\labelenumi}{(\roman{enumi})}
\begin{enumerate}
\item\label{mod6i1} $M$ is a zero-dimensional\/ $\OXx$-module, 
i.e., $M$ is\/ $m_x$-torsion.\vspace{.6pt}
\item\label{mod6i2} $i_xM \in \Aqct(\X)$.\vspace{.6pt}
\item\label{mod6i3} For every affine open neighborhood\/ 
$\U=\Spf(A)$ of\/ $x,$ $i_xM|_{\U}\in \Aqct(\U)$. Moreover, $M$ being  an\/ $A$-module 
via the natural map\/ $A \to \OXx$, there is a natural isomorphism\/ 
$i_xM|_{\U} \iso M^{\sim}$.  \vspace{.6pt}
\item\label{mod6i4} There is an affine open neighborhood\/ 
$\U$ of\/ $x$ such that\/ $i_xM|_{\U}\in \Aqct(\U)$. 
\end{enumerate}
\end{alem}

\begin{proof}
We show that (i) $\Lra$ (iii)$\:\iff\:$(ii), and (iii) $\Lra$ (iv) $\Lra$ (i).

(i)$\:\Rightarrow\:$(iii). Let $\U = \Spf(A)$ be an affine open 
neighborhood of~$x$, and for any \mbox{$f \in A$,}  $\U_f\subset\U$  the 
open subset $\Spf(A_{\{f\}})$. We claim that 
\emph{if\/ $x \notin \U_f$ then\/ $\Gamma(\U_f, M^{\sim}) = 0,$ and if\/ 
$x \in \U_f$ then the natural map\/ $M \to \Gamma(\U_f, M^{\sim})$ is 
an isomorphism.} Indeed, since the open
prime ideal ~$\pfr$ in~$A$ corresponding to~$x$ is taken into 
~$m_x$ by the natural map $A \to \OXx$, therefore $M$ is~$\pfr$-torsion. Thus $M$ is a torsion $A$-module, 
and so  there are natural isomorphisms
$$
M_f \underset{\ref{lem:mod4b}}{\iso}\ M 
  \otimes_A A_{\{f\}}\>\underset{\ref{prop:mod1}\textup{(iv)}}{\iso}\, \Gamma(\U_f, M^{\sim});
$$ 
and  if $x \notin \U_f$ then $f \in \pfr$, so 
$M$ is $(f)$-torsion, i.e., $M_f =0$.
If $x \in \U_f$ then $f \notin \pfr$, so $f$ maps to a unit
in $\OXx$, and the natural map $M\to M_f$ is an isomorphism.
The claim results.

 Any open subset $\V \subset \U$ is a union of open subsets of the type $\U_f$, and it follows that if $\>x \notin\< \V\>$ then
$\>M^{\sim}(\V) = 0$, and that if $\>x \in \V\>$ then the natural map 
$M\< \to M^{\sim}(\V)$ is an isomorphism. 
Thus there is a natural isomorphism $i_xM|_{\U} \iso M^{\sim}\<$.
Since $M^\sim\in \Aqct(\U)$  (by \ref{prop:mod1}\kern.5pt(ii) and 
\ref{lem:mod4c}), therefore  $i_xM|_{\U}  \in \Aqct(\U)$.\vspace{1pt}

(iii) $\Rightarrow$ (ii). The property of being quasi-coherent and torsion
is local. If $\V$ is an open set not containing $x$ then $i_xM$ 
restricts to the zero sheaf and hence is in~$\Aqct(\V)$. The implication 
follows.\vspace{1pt}

(ii) $\Rightarrow$ (iii) follows easily from  \ref{prop:mod5}\kern.5pt(iii);  and (iii) $\Rightarrow$ (iv) is obvious.\vspace{1pt}

(iv) $\Rightarrow$ (i). With $\U = \Spf(A)$  as in (iv),  \ref{prop:mod5}\kern.5pt(iii) gives \mbox{$M^\sim\cong i_xM|_{\U}\in\Aqct(\U)$,}
so by \ref{lem:mod4c}, $M$ is a torsion $A$-module. Let $\pfr$ be
 the open~prime $A$-ideal corresponding to $x \in \U$. For any $f \in \pfr$, we have \mbox{$x\notin \U_f \set\Spf(A_{\{f\}})$}, whence, as above, there are natural isomorphisms
\[
M_f \underset{\ref{lem:mod4b}}{\iso}\ M 
  \otimes_A A_{\{f\}}\>\underset{\ref{prop:mod1}\textup{(iv)}}{\iso}\, \Gamma(\U_f, M^{\sim})\iso \Gamma(\U_f, i_xM)=0,
\]
so that $M$ is $(f)$-torsion. As $\pfr$ is finitely 
generated,   $M$ must be $\pfr$-torsion; and since 
$m_x=\dirlm{f\notin\pfr}\pfr_{\{\<f\<\}}^{}=\pfr\OXx$
(\cite[p.\,186, proof of (7.6.17)]{EGAI}),
$M$ is $m_x$-torsion.\looseness=-1
\end{proof} 

\smallskip
Let $\X$ be a noetherian  formal scheme.
For $x \in \X$, let $J(x)$ denote the injective 
hull of the residue field $k(x)$ over the local 
ring $\OXx$. It is easily checked that $i_xJ(x)$ 
is an injective object in $\A(\X)$ (\cite[page 123]{RD}).
By \ref{lem:mod6}, $i_xJ(x)\in \Aqct(\X)$.

Recall that in a locally noetherian category a direct sum 
of injectives is also an injective.

\medskip

\begin{aprop}
\label{prop:mod7}
Let $\X$ be a noetherian formal scheme.
\begin{enumerate}
\item\label{mod7i1}The categories\/ $\At(\X)$ and\/ $\Aqct(\X) $ are 
locally noetherian and have enough injectives.\vspace{1pt} 
\item\label{mod7i2}The indecomposable injectives of\/ $\Aqct(\X)$ are
those of the form\/ $i_xJ(x)$ defined above. In particular,
any injective object\/ $\cI$ of $\Aqct(\X) $ is a 
direct sum of injectives of the form\/ $i_xJ(x)$
and hence is injective in\/ $\At(\X)$ too.\vspace{1pt} 
\item\label{mod7i3}For any\/ $\M \in \At(\X) \cup \Ac(\X)$ and 
for any injective\/ $\cI$ in\/ $\Aqct(\X),$ the sheaf\/ 
$\sHom_{\OX}(\M,\cI)$ is flasque.\vspace{1pt}
\item\label{mod7i4}For any injective\/ $\cI$ in $\Aqct(\X),$ and\/ $x \in \X,$ the\/ $\OXx$-module\/ $\cI_x$ is injective.
\end{enumerate}
\end{aprop}  
\begin{proof}
For proofs of (i), (ii) and (iii), we 
refer to \cite[p.\,876]{Ye}. 

For (iv), 
we may assume,\vspace{1pt} by (ii), that $\cI=i_{x^{\prime}}J(x^{\prime})$ for suitable 
${x^{\prime}}$. If ${x^{\prime}}$ doesn't specialize to $x$
then $\cI_x = 0$. 
It remains to show that\vspace{1.5pt} \emph{if\/ ${x^{\prime}}$ specializes to\/ $x$
then\/ $J(x^{\prime})$ is an 
injective $\OXx$-module.} For this it suffices,\vspace{-.3pt} 
with $\cO_{\<y}\set\cO_{\X\<,y}\ (y\in\Y)$,   that
$\wh{\cO_{\<x^{\prime}}} \cong \wh{(\cO_{\<x})_{\wp}}$ 
for some prime $\wp$ in $\cO_{\<x}$,
since then the injective hull $J(x^{\prime})$ of the residue
field $k(x^{\prime})$ over ${\cO_{\<x^{\prime}}}$, being 
injective over the local ring $\wh{\cO_{\<x^{\prime}}}$ too, 
is  injective over $\wh{(\cO_{\<x})_{\wp}}$ 
and hence over $\cO_{\<x}\>$. \vspace{1.5pt}

Let $\U=\Spf(A)$ be an open neighborhood of $x$. Recalling that $\cO_{\<y}$ is noetherian for all $y\in\U$ \cite[p.\,403, (10.1.6)]{EGAI}, one deduces  from \cite[p.\,184, (7.6.9)]{EGAI} that 
if~$\afr$ is an open $A$-ideal  and 
$\qfr=\qfr_y$ is the  open prime $A$-ideal  corresponding to~$y$, then
\mbox{$\cO_{\X,y}/\afr\cO_{\<y} \cong A_{\qfr}/\afr A_{\qfr}$.} So if $\pfr\set\qfr_{x'}\subset\qfr_y$  then $\wp_y\set\pfr\cO_{\<y}$ is prime,\vspace{.8pt} 
and for all~ $i\ge 0$,
$(\cO_{\<y})_{\wp_y}/\wp_y^i (\cO_{\<y})_{\wp_y}$ is isomorphic to the localization of the $A$-module 
$A_{\qfr}/\pfr^i \<A_{\qfr}$ at the prime ideal $\pfr$, i.e., to $A_{\pfr}/\pfr^i \<A_{\pfr}$.
Hence for all $y$ in the closure\vspace{-1pt} of $x'\<$,  $\wh{(\cO_{\<y})_{\wp_y}}\cong\wh{A_\pfr}$. Taking $y=x'$ and $x$ respectively, 
one gets $\wh{\cO_{\X,x^{\prime}}} \cong \wh{(\cO_{\<x})_{\wp_x}}$, as desired.
\end{proof}

\subsection{Smooth maps}
\label{subsec:morphism}

A homomorphism of topological rings 
$\phi \colon A \to B$ is 
\mbox{\it formally smooth\/} if $\phi$ is continuous and if%
\index{formally smooth!ring homomorphism}
for every discrete topological $A$-algebra~$C$
and~every nilpotent ideal $I$ of $C$, any continuous $A$-homomorphism
$B \to C/I$ factors as $B \stackrel{v\,}{\to} C \onto C/I$
with $v$  a continuous $A$-homomorphism (\cite[19.3.1]{EGAOIV}).\vspace{.6pt}

A morphism of formal schemes $f \colon \X \to \Y$ is said to be 
{\em formally smooth\/} if for any morphism $Z \to \Y$%
\index{formally smooth!formal-scheme map}
where $Z = \Spec(C)$ is an affine
scheme, and for any closed subscheme $Z_0
\subset Z$ defined by a nilpotent ideal in $C$, every $\Y$-morphism 
$Z_0 \to \X$ extends to a $\Y$-morphism $Z \to \X$ 
(cf.~\cite[\S17.1]{EGAIV}).

\begin{aexams}
\label{exam:morph0}
We recall some of the elementary properties and standard examples 
of formally smooth maps. These will often be used without explicit mention.
For proofs cf.~ \cite[\S17.1]{EGAIV} and~\cite[\S19.3]{EGAOIV}.
\begin{enumerate}
\item An open immersion is a formally smooth map. 
A composition of formally smooth maps is formally smooth.
Formal smoothness is preserved under base change.
\item A map of affine formal schemes 
$\Spf(B) \to  \Spf(A)$ is formally smooth if and only if 
the corresponding homomorphism of topological rings
$A \to B$ is formally smooth.
\item For a discrete ring $A$, any polynomial algebra under the 
discrete topology is formally smooth over $A$. 
\item If $\phi\colon A \to B$ is a formally
smooth map of topological rings then for any multiplicative sets
$S\subset A$ and $T\subset B$ such that $\phi(S)\subset T$,
the induced map $S^{-1}A \to T^{-1}B$ is formally smooth. 
\item A map of topological rings 
$A \to B$ is formally smooth if and only if the induced map
of  respective completions $\wh{A} \to \wh{B}$ is formally smooth.
\item Let $\phi \colon A \to B$ be a 
formally smooth map of topological rings.
Then for any choice of coarser topologies on $A,B$ for which the
square of any open $B$-ideal is open and for which $\phi$
remains continuous, $\phi$ is formally smooth under the coarser 
topologies too.
\end{enumerate}
\end{aexams}

\begin{adefi}
\label{def:morph1a}
A continuous homomorphism $\phi \colon A \to B$ of 
noetherian adic rings (resp.~a map $\fXY$ of noetherian formal
schemes)
is \emph{smooth} if it is 
essentially of pseudo-finite type and formally smooth. 
\end{adefi}
\index{smooth!ring homomorphism}\index{smooth!formal-scheme map}
 
\smallskip

Recall that a homomorphism $\phi \colon R \to S$ of noetherian adic rings is 
{\it adic}\index{adic!homomorphism} if for one (hence for every) defining ideal~$I$
of~$R$, $\phi(I)$ is a defining ideal of $S$. 

\begin{alem}
\label{lem:morph1}
Let\/ $\phi \colon A \to B$ be an essentially pseudo-finite type 
homomorphism of noetherian adic rings. Then\/ $\phi$ factors as 
$A \xto{\;\sigma\;} C \xto{\;\pi\;} B$
where\/ $\sigma$ is a smooth map and\/ $\pi$ is a surjective adic 
homomorphism. More specifically, $C$ can be obtained as 
the completion of a localization of a polynomial algebra 
over\/ $A$---say\/ $(A[X_1,\ldots,X_n])_S$---along an ideal\/ $I$ 
that contracts to an open\/ $A$-ideal. 
\end{alem}
\begin{proof} 
Let $\afr,\bfr$ be defining ideals of $A,B,$ respectively, such 
that $\phi(\afr) \subset \bfr$. 
By hypo\-thesis there exists a finitely generated $A$-algebra
$D$ and a multiplicative set~$T$ in~$D$ such that 
$B/\bfr \cong D_T$. Then there is a surjection from a polynomial ring
\mbox{$P \set A[X_1,\ldots,X_r]$} to~$D$ which maps $\{X_i\}$
to a set of generators of~$D$.
Let~$S$ be the inverse image of~$\>T$ in~$P$.
The composition $P \onto D \to D_T \cong B/\bfr$ lifts to an 
$A$-homomorphism 
$\psi \colon P \to B$. Any 
element of~$\psi(S)$ maps to a unit in~$D_T \cong B/\bfr$. 
Since~$\bfr$, being a defining ideal of~$B$, is in the Jacobson 
radical, it follows that $\psi(S)$ consists of units of~$B$ and so we obtain an induced map 
$\psi_S \colon P_S \to B$. 

Let $\pi_0^{} \colon Q = P_S[X_{r+1}, \ldots, X_n] \to B$ 
be a map extending $\psi_S$ and taking the indeterminates
$X_{r+1}, \ldots, X_n$
to a generating set of the $B$-ideal~$\bfr$. 
Let $I$ be the $Q$-ideal generated 
by~$\afr$ and~$X_{r+1}, \ldots, X_n$. 
Then $\pi_0^{}$ is continuous for the $I$-adic topology on~$Q$ 
and hence extends to a map~$\pi$, clearly adic,  from the $I$-adic completion 
$C \set \wh{Q}$ to~$B$. Since $\widehat{I}B = \bfr$, 
therefore $B$ is an $\widehat{I}$-adically 
separated $C$-module and~$C$ surjects onto $B/(\widehat{I}B)$. 
It follows from \cite[Theorem 8.4]{Ma} that  $\pi$ is 
surjective.
\end{proof}

\smallskip

\begin{acor}
\label{cor:morph1b}
If\/ $f \colon \X \to \Y$ is an essentially pseudo-finite type 
map of formal schemes then any neighborhood\/ of a point $x \in \X$ contains
an affine open neighborhood\/ $\U$ such that\/ 
$f|_{\U}^{} = hi$ where for some formal scheme\/ $\Z,$ $i \colon \U \to \Z$ is a closed immersion  and\/ $h \colon \Z \to \Y$
is smooth. \looseness=-1
\end{acor}

\pagebreak[2]

\medskip

Let $A \xto{\;\phi\;}B$ be a continuous homomorphism\vspace{-1pt} of 
noetherian adic rings and let
$\U \set \Spf(B) \xto{\; f \;} \V \set \Spf(A)$ be 
the induced map of formal schemes.
Let\vspace{.5pt} $x \in \U$ and  $y\set f(x)$. 
Let $k(x)$ and~$k(y)$ denote the residue fields 
of the local rings $\OUx$ and~$\OVy$
respectively. We regard these local rings as being topologized 
by the powers of the stalks of defining 
ideals on~$\U$ and~$\V$ respectively. Thus if $\afr,\bfr$ are
defining ideals of $A,B$ 
then $\afr\OUx,\: \bfr\OVy$ are defining ideals
of $\OUx,\: \OVy,$ respectively.

\medskip

\begin{aprop}
\label{prop:morph2}
In the preceding situation we have\textup{:}\begin{enumerate}
\item\label{morph2i1}If\/ $\phi$ is essentially of pseudo-finite type 
then\/ $k(x)$ is a finitely generated field extension of\/ $k(y)$.\vspace{1pt}
\item\label{morph2i2}If\/ $\phi$ is formally smooth, then\/ $\phi$ is flat.
Furthermore, the induced map of topological rings\/ 
$\OVy \to \OUx$ is also formally smooth and flat. In particular,
$f \colon \U \to \V$ is a flat map of formal schemes.
\end{enumerate}
\end{aprop}
\begin{proof}
(i). Let $\bfr$ be a defining ideal in $B$ and 
$\afr \subset \phi^{-1}\bfr$ a defining ideal in~$A$. 
Then $A/\afr \to B/\bfr$ is essentially
of finite type. As in the proof of \ref{prop:mod7}\kern.5pt(iv), 
$\OVy/\afr\OVy$ and $\OUx/\bfr\OUx$
are localizations of $A/\afr$ and $B/\bfr$ respectively.
Therefore the induced map $k(x) \to k(y)$ is essentially 
of finite type, whence the result.\vspace{1pt}

(ii). 
Let $\nfr$ be a maximal $B$-ideal and $\mfr = \phi^{-1}\nfr$.
Let $\bfr$ be a defining ideal in $B$ and $\afr \subset \phi^{-1}\bfr$
a defining ideal in $A$. Then the induced map 
$A_{\mfr} \to B_{\nfr}$ is formally smooth under the 
$\afr$-adic and $\bfr$-adic topologies respectively. 
The map $A_{\mfr} \to B_{\nfr}$ is also continuous for the 
coarser $\mfr$-adic and $\nfr$-adic topologies 
respectively. Therefore $A_{\mfr} \to B_{\nfr}$
is formally smooth for the coarser topologies too. 
By~\cite[19.7.1]{EGAOIV}, $B_{\nfr}$ is flat over~$A_{\mfr}$, 
hence over~$A$. This being so for any~$\nfr$, $B$ is flat over~$A$.

For the final part, let $\pfr \subset B$ (resp.~$\qfr \subset A$)
be the open prime corresponding to the point $x$ (resp.~$y$).
As in the proof of \ref{prop:mod7}(iv), the completion 
of~$\OUx$ (resp.~$\OVy$) along the ideal $\bfr\OUx$ 
(resp.~$\afr\OVy$) is isomorphic to the completion 
of~$B_{\pfr}$ (resp.~$A_{\qfr}$) along the ideal $\bfr_{\pfr}$ (resp,~$\afr_\qfr$).
Therefore [$A \to B$  formally smooth] $\Lra$
[$A_{\qfr} \to B_{\pfr}$ formally smooth] $\Lra$ 
[$\wh{A_{\qfr}} \to \wh{B_{\pfr}}$ formally smooth] $\Lra$ 
[$\OVy \to \OUx$ formally smooth] 
(\cite[19.3.5, (iv)]{EGAOIV}, \cite[19.3.6]{EGAOIV}). 
Again, by~\cite[19.7.1]{EGAOIV}, $\OVy \to \OUx$ is flat.
\end{proof}

\subsection{Differentials on topological rings}
\label{subsec:differ}
We now describe various elementary properties of 
the module of relative differentials for a map of 
topological rings, and its 
behavior under smooth maps.
Let $A \to B$ be a continuous homomorphism of adic rings.
Let~$\bfr$ be any defining ideal in~$B$. Let
$\Omega^1_{B/A}$ be the relative $B$-module of differentials 
and~$\Omega^m_{B/A}$ its~$m$-th exterior power $(m\ge 0)$, $\bfr$-adically topologized
(cf.~ \cite[20.4.5]{EGAOIV}).
Set $B_i \set B/\bfr^{i+1}$, and 
\[
\widehat{\Omega^m_{B/A}}\set 
{\inlm{i}}(\Omega^{m}_{B/A}\otimes_B B_i),
\] 
the $\bfr$-adic completion of 
$\Omega^m_{B/A}$.  Also, set 
\[
\Omega^{\mathrm{sep}}_{B/A}\set 
\Omega^1_{B/A}/(\cap_i\bfr^i\Omega^1_{B/A}), 
\]
the universal separated module of differentials of $B/A$. These definitions do not 
depend on the choice of~$\bfr$. Indeed, the canonical $A$-derivation\vspace{-2.5pt} $d_{B/A}\colon B\to \Omega^1_{B/A}$ induces derivations \mbox{$d^{\>\mathrm{sep}}_{B/A}\colon B\to \Omega^{\mathrm{sep}}_{B/A}$}%
\index{ $\Cfr$2@$d^{\>\mathrm{sep}}$ (universal separated derivation)} 
and $\wh d_{B/A}\colon B\to \widehat{\Omega^1_{B/A}}$ which
are universal for continuous $A$-derivations\vspace{1pt} of $B$ into separated (resp.~complete separated) $B$-modules.
\index{ $\Cfr$2@$\hat d$ (universal derivation)! for continuous homomorphisms}

Recall that for any  topological $B$-module $N$, composition with $d_{B/A}$ gives  an isomorphism
$\Hom_B^{\cont}(\Omega^1_{B/\<A\>},N)\iso\Der_A^{\cont}(B,N)$,%
\index{ $\Hom^{\cont}$3@$\Hom^{\cont}$ (continuous homomorphisms)}%
\index{ $\Der^{\cont}$ (continuous derivations)}
where the superscript ``c'' signifies, respectively, \emph{continuous} homomorphisms and derivations 
(see \cite[20.4.8.2]{EGAOIV}). Let~$\bfr$ be a defining
ideal of~$B$. If the topology on~$N$ is~$\bfr$-adic and~$N$ 
is separated (e.g., if~$N$ is finitely generated) then we 
also have the following relations.
\[
\Der_A^{\cont}(B,N) = \Der_A(B,N), 
\]
\[
\Hom_B^{\cont}(\Omega^1_{B/A},N) \osi 
\Hom_B^{\cont}(\Omega^{\mathrm{sep}}_{B/A},N) = 
\Hom_B(\Omega^{\mathrm{sep}}_{B/A},N) \iso \Hom_B(\Omega^1_{B/A},N).
\]
All these relations are easily verified. For example, the first one 
holds because for any $\delta \in \Der_A(B,N)$ we have 
$\delta(\bfr^{i+1}B) \subset \bfr^i N$.

   
\begin{aprop}
\label{prop:diff1}
Let\/ $f\colon A \to B$ be an essentially pseudo-finite type 
map of noetherian adic rings. Let\/ $\bfr$ be a defining\/ 
~$B$-ideal and\/ $B_i \set B/\bfr^{i+1}$. Let\/ $m\ge0.$ Then\/$:$
\renewcommand{\labelenumi}{(\roman{enumi})}
\begin{enumerate}
\item The\/ $B$-module\vspace{-2.5pt} $\Omega^{\mathrm{sep}}_{B/A}$ is finitely generated and complete; and there is a natural isomorphism
$\Omega^{\mathrm{sep}}_{B/A} \iso \widehat{\Omega^1_{B/A}}.$
 
 \item There are natural isomorphisms\/ 
$$
\bigwedge^m_B \widehat{\Omega^1_{B/A}}\iso\widehat{\Omega^m_{B/A}} \iso
{\inlm{\under{.3}i}}\Omega^{m}_{B_i/A}. 
$$
In particular, $\widehat{\Omega^m_{B/A}}$ is a finitely generated\/
$B$-module. \vspace{2pt}

\item If\/ $f$ is  smooth then 
$\smash{\widehat{\Omega^m_{B/A}}}$ is a projective $B$-module.\vspace{2pt}

\item If\/ $\bfr$ is a defining ideal of\/ $B,$ $\bfr'\subset\bfr,$
and\/ $B'\set B$ with\/ $\bfr'$-adic topology\/ $($also
an adic ring\/$)$ is essentially of pseudo-finite type over\/ $A,$ then the natural map is an isomorphism\vspace{-3pt}
$$
\widehat{\Omega^1_{B'/A}}\iso\widehat{\Omega^1_{B/A}}\>.
$$
\end{enumerate}
\end{aprop}  
\begin{proof}
$B_0$ is essentially of finite type 
over $A$, so $\>\Omega^1_{B_0/\<A}\>$ is a finitely generated 
$B_0$-module and hence a finitely generated $B$-module. 
Then the instance $i=0$ of the natural exact sequences 
\begin{equation}
\label{eq:diff1a}
\bfr^{i+1}/\bfr^{i+2} \underset{\textup{via}\;d}\lra \,\Omega^1_{B/A}\otimes_B B_i \lra \Omega^1_{B_i/A} \lra 0 
\qquad(i\ge 0)\\[-2pt]
\end{equation}
shows that $\Omega^1_{B/A}\otimes_B B_0$ is a 
finitely generated $B$-module,\vspace{1pt} so that the $B$-module
$\Omega^{\mathrm{sep}}_{B/A} \otimes_B B_0 
\cong \Omega^1_{B/A} \otimes_B B_0\>$ is finitely 
generated.\vspace{-2.5pt} By~\cite[Thm.\,8.4]{Ma}  
it follows that the $B$-module $\Omega^{\mathrm{sep}}_{B/A}$ is 
finitely generated, hence complete;  and so $\Omega^{\mathrm{sep}}_{B/A}= \widehat{\Omega^1_{B/A}} $.\vspace{1pt}

Furthermore, the exact sequence\vspace{1pt} (\ref{eq:diff1a}) gives that the kernel $k_i$ of the 
natural surjection $\Omega^1_{B/A} \onto \Omega^1_{B_i/A}$ satisfies
$$
\bfr^{i+1}\Omega^1_{B/A}\subset \,k_i=\bfr^{i+1}\Omega^1_{B/A} + d(\bfr^{i+1}) \,\subset \bfr^i\Omega^1_{B/A},
$$
giving rise to a sequence 
$$
\cdots\lra \Omega^1_{B_{i+1}/A}\lra \Omega^1_{B/A}\otimes_B B_i\lra 
\Omega^1_{B_i/A}\lra
\Omega^1_{B/A}\otimes_B B_{i-1}\lra\cdots,
$$
from which follows, upon application of $\>\bw^{\!m}_B\>$, a natural sequence
$$
\cdots\lra\Omega^m_{B_{i+1}/A}\lra \Omega^m_{B/A}\otimes_B B_i\lra 
\Omega^m_{B_i/A}\lra
\Omega^m_{B/A}\otimes_B B_{i-1}\lra\cdots,
$$
and hence a natural isomorphism $\,\widehat{\Omega^m_{B/A}} \iso 
{\inlm{\under{.3}i}}\Omega^{m}_{B_i/A}$.

\pagebreak[3]

By part (i), $\bw^m\widehat{\Omega^1_{B/A}}$ is a finitely generated 
$B$-module and hence is complete. With $B_i \set B/\bfr^{i+1}$ there are 
natural isomorphisms
\[ 
(\bw^m \widehat{\Omega^1_{B/A}}) \otimes_B B_i \cong
\bw^m (\widehat{\Omega^1_{B/A}} \otimes_B B_i) \cong
\bw^m (\Omega^1_{B/A} \otimes_B B_i) \cong
(\bw^m \Omega^1_{B/A}) \otimes_B B_i.                 
\]
Taking inverse limit over $i$ we conclude that
$\bw^m\wh{{\Omega}^1_{B/A}} \cong \widehat{\Omega^m_{B/A}}$,
proving~(ii).
 
Now suppose $f$ is  smooth.\vspace{-1pt} By \cite[Cor.~20.4.10]{EGAOIV}, 
$\widehat{\Omega^1_{B/A}}\otimes_B B_i$ is a projective $B_i$-module for 
each $i$. Hence\vspace{1pt} by \cite[Thm.~22.1]{Ma}, $\widehat{\Omega^1_{B/A}}$ is a 
projective $B$-module. Taking exterior powers and using  (ii), we
deduce (iii).\vspace{1.5pt}

The natural map in (iv) is the obvious one\vspace{.6pt} from the completion of a
module to the completion of the same module with coarser topology.\vspace{1pt}
That it is an isomorphism follows from the fact that every finitely-generated $B$-module is $\bfr$- and $\bfr'$-complete and separated, whence 
the canonical derivations $d\colon B\to \widehat{\Omega^1_{B/A}}$ and
$d'\colon B\to \widehat{\Omega^1_{B'/A}}$ are both universal for
$A$-derivations of $B$ into finitely-generated $B$-modules.
\end{proof}

From now on,\vspace{.6pt} for any continuous homomorphism $A \to B$ that is essentially of 
pseudo-finite type we shall denote 
$\,\widehat{\Omega^m_{B/A}}\,$ by $\,\Ohm^m_{B/A}\>$.
\vspace{1pt}  
\index{ WT@$\Ohm^m$ (continuous order-$m$ relative differentials)!for continuous homomorphisms}

\begin{alem}
\label{lem:diff2}
Let\/ $A \xto{\;f\;} B \xto{\;g\;} C$ be essentially pseudo-finite-type
homomorphisms of noetherian adic rings. Then there is a natural exact sequence
of continuous $C$-module maps
\[ 
 \Ohm^1_{B/A}\otimes_B C \to \Ohm^1_{C/A} \to \Ohm^1_{C/B} \to 0. 
\]
If, moreover, $g$ is smooth
\emph{(\ref{def:morph1a})} then this sequence  is part of a 
split exact sequence
\[ 
0 \to \Ohm^1_{B/A} \otimes_B C \to \Ohm^1_{C/A} \to \Ohm^1_{C/B} \to 0. 
\]
\end{alem}
\begin{proof}
The map \vspace{4pt} $ \Ohm^1_{B/A}\otimes_B C \to \Ohm^1_{C/A}$ corresponds to the 
continuous\vspace{-3pt} $A$-derivation 
$B\to C\xto{\under{1.7}{\hat d_{C/A}}}\Ohm^1_{C/A}$. (Note here that 
by \ref{prop:diff1}(i),\vspace{1pt} $\Ohm^1_{B/A}$ is a 
finitely generated $B$-module,  so that
$\Ohm^1_{B/A} \:\wh{\otimes}_B\: C \cong \Ohm^1_{B/A} \otimes_B C$ 
(see \cite[p.\,189, 7.7.9]{EGAI}); and that  $\>\Ohm^1_{C/A}$, being finitely generated over $C$, is a complete separated $B$-module.)\vspace{-1pt} The map $\Ohm^1_{C/A} \to \Ohm^1_{C/B}$ corresponds to the continuous $A$-derivation ~$\hat d_{C/B}$.\vspace{1pt}
Exactness\- follows from that of the dual sequence---for any
separated complete $C$-module~$N$:
\[ 
\Der_A^{\cont}(B,N)\gets\Der_A^{\cont}(C,N)\gets\Der_B^{\cont}(C,N)\gets 0.
\]
(See also  \cite[p.\,152, 20.7.17.3]{EGAOIV}.)
For smooth $g$ one can then apply \cite[20.7.18]{EGAOIV}, in view of 
\emph{ibid.,} p.\,114, Definition 19.9.1. (In fact, arguing as in the proof of \emph{loc.\,cit.,} one 
sees that $\Der_A^{\cont}(C,N)\to\Der_A^{\cont}(B,N)$ is surjective.)
\end{proof}

\begin{alem}
\label{lem:diff3} 
Let\/ $A,$ $B,$ $C,$ be noetherian adic rings, let $f\colon A\to B$ be
a homomorphism essentially of pseudo-finite type, and let\/ 
$\pi\colon B\to C$ be a continuous surjection, with kernel\/ $I$.   
Suppose that the topological
quotient algebra\/ $B/I$ is smooth over\/ $A$. Then there is a 
split exact sequence
\[ 
0 \to I/I^2 \to \Ohm^1_{B/A} \otimes_B C \to \Ohm^1_{C/A} \to 0. 
\]
\end{alem}

\begin{proof} In view of \ref{prop:diff1}\kern.5pt(iv), one may assume that $C=B/I$.
Let $\bfr$ be a defining $B$-ideal. Consider the  diagram
\[
\begin{CD}
0 @>>> (\cap_i\:\bfr^i\Omega^1_{B/A})\otimes_B C @>>> 
\cap_i\:\bfr^i\Omega^1_{C/A} \\
@VVV @VVV @VVV \\
I/I^2 @>>> \Omega^1_{B/A}\otimes_B C @>>> \Omega^1_{C/A} \\
\end{CD}
\]
wherein the bottom row consists of usual natural maps, 
the top row is induced by the bottom row, and the vertical 
maps are the natural ones. The diagram clearly commutes
and so the cokernels of the vertical maps
lie in a sequence
\begin{equation*} 
0 \to I/I^2 \to \Omega^{\mathrm{sep}}_{B/A} \otimes_B C \to 
\Omega^{\mathrm{sep}}_{C/A} \to 0. 
\tag*{$\Eb :$}
\end{equation*}
By \ref{prop:diff1}(i) every module in $\Eb$ is
finitely generated over $C$, so for $\Eb$ to be  split exact
it suffices  that for 
any finitely generated $C$-module $N$ the induced sequence
$\Hom_C(\Eb,N)$ is exact. Rewriting $\Hom_C(\Eb,N)$ as 
\[ 
0 \to \Der_A(C, N) \to \Der_A(B, N) \stackrel{u\:}{\to} 
\Hom_C(I/I^2,N) \to 0 
\]
one sees that the only nontrivial thing to show is surjectivity
of~$u$. 

\pagebreak[3]
Since the kernel $I/I^2$ of the natural surjection
$\ov{\pi}\colon B/I^2 \to C$ is nilpotent and closed,
and $C$ is smooth over $A$, therefore 
the  identity map of $C$ lifts~to an \mbox{$A$-homomorphism} 
$l\colon C\to B/I^2$ \cite[p.\,84, 19.3.10]{EGAOIV}. 
Then
$$
\sigma \set 1-l\ov{\pi} \colon B/I^2 \to I/I^2
$$ is 
an $A$-derivation: for all $x,y \in B/I^2,$
$$
\sigma(xy) =  xy -l\ov{\pi}(xy) =  xy - l\ov{\pi}(x)l\ov{\pi}(y) =
 xy - (x-\sigma(x))(y-\sigma(y)) \notag \\
= x\sigma(y) + y\sigma(x). 
$$
So for any
$\delta\in \Hom_C(I/I^2,N)$, 
 $\delta\sigma\in\Der_A(B/I^2\<,\>N)=\Der_A(B,N)$ and
$u(\delta\sigma)=\delta$. This proves surjectivity of~$u$, 
and hence split-exactness of~ $\Eb\<$. 
\end{proof}

\begin{aexam}
\label{exam:diff3a}
Let $B$ be a noetherian adic ring, $S$ a multiplicative 
set in $B$ and $I$ an open ideal in $B$.
Set $C \set B\{S^{-1}\}$, the completion of~$B_S$ along~$I_S$.
Then the induced map $B \to C$ is smooth. 

Any continuous $B$-derivation from $C$ into a topological $C$-module  vanishes on~ $B_S$, hence, by continuity, on $C$.
Therefore,  $\Ohm^1_{C/B} = 0$.

Any continuous $A$-derivation from $B$ into a complete separated\vspace{1pt} $C$-module
extends uniquely to $C$. Hence
(or by \ref{lem:diff2})  there is a natural isomorphism
$$
\Ohm^1_{B/A}\:\wh\otimes_B\:B\{S^{-1}\}\iso \Ohm^1_{B\{S^{-1}\}/A}\>.
$$
\end{aexam}

\begin{aexam}
\label{exam:diff3c}
Let $A$ be a noetherian adic ring and 
$P = A[X_1,\ldots,X_r]$ a polynomial ring over $A$. Let $I\>$ be 
an ideal in $P$ such that $I \cap A$ is open in $A$, and let~$C$ 
be the completion of $P$ along~$I$. Then the natural map 
$A \to C$ is smooth. As any $A$-derivation of $P$ into an $I$-adically complete separated $C$-module  extends uniquely to $C$,
one sees that the $C$-module $\Ohm^1_{C/A}$ is free, 
with basis  $dX_1, \ldots, dX_r$.
\end{aexam}

\begin{aprop}
\label{homeless}
Let\/ $A \to B$ and\/ $A \to C$ be continuous 
maps of  adic rings. Set\/ $D \set B\: \wh{\otimes}_{\<\<A}\: C$. 
Then there exists a natural isomorphism
$
\Ohm^1_{B/A}\:\wh\otimes_B\: D \iso \Ohm^1_{D/C}.
$

\end{aprop}

\begin{proof}
(i) Via the universal properties of $\>\Ohm^1$ and of  extension of scalars (from $B$ to $D$),
the Proposition amounts to the statement that for any complete separated $D$-module $L$, ``restriction" induces a bijection $\Der^c_C(D,L)\iso\Der^c_A(B,L).$
But if we give the $D$-algebra $\cD\set D\oplus L$ (where $L^2$=0) the product topology,
under which it is complete and separated, then the standard correspondence between derivations into $L$ and ring homomorphisms into $\cD$ \cite[p.\,118, 20.1.5]{EGAOIV} respects continuity, and so transforms the statement into
bijectivity of the restriction map from continuous $C$-homomorphisms
$D\to\cD$ to continuous $A$-homomorphisms $B\to\cD$, which is easily seen
to hold, by  the universal property of complete tensor products.\vspace{1pt}
\end{proof}

\emph{Remark.} When $A$ and $B$ are noetherian and $A\to B$ is essentially of
pseudo-finite type, one can replace the $\wh{\otimes}$ in \ref{homeless} by $\otimes$,
see  proof of \ref{lem:diff2}. Moreover,\vspace{.6pt} in view of \ref{prop:diff1}\kern.5pt(ii), one can replace $\Ohm^1$ by $\Ohm^m$ for any $m\ge0$.

\subsection{Differentials on formal schemes} 
\label{subsec:difform} 
Let $ f \colon \X \to \Y$ be a morphism of noetherian formal schemes. 
Suppose $\I \subset \OX$ and $\J \subset \OY$ are defining ideals 
such that $\J\OX \subset \I$. There result  morphisms of ordinary\vspace{1pt}
schemes 
$$
X_n \set (\X,\OX/\I^{n+1})\xrightarrow{f_n} (\Y,\OY/\J^{n+1})=:Y_n
\qquad(n>0).
$$
Let $j_n\colon X_n \into \X$ be the canonical closed immersion.
For  ($m\ge 0)$  let $\Omega^m_{X_n/Y_n}$ be the $m$-th exterior power of the sheaf\vspace{-1pt} of relative differentials on~$X_n$,\vspace{.6pt} and  set
$\Ohm^m_f=\Ohm^m_{\X/\Y} \set \inlm{\under{.2}{n}} j_{n*}\Omega^m_{X_n/Y_n}$. 
\index{ WT@$\Ohm^m$ (continuous order-$m$ relative differentials)!for formal-scheme maps}
This $\Ohm^m_f$\vspace{-1.4pt} is independent of the choice of $\>\I$,  $\J$.\vspace{1pt}
For every $n$ there is a natural sheaf homomorphism\vspace{1pt} 
$\cO_{\<\<X_n} \to \Omega^1_{X_n/Y_n}$, and hence applying $j_{n*}$ and\vspace{.6pt} 
taking inverse limits results in a natural sheaf homomorphism
$\hat d_{\X/\Y}=\hat d_f\colon\OX \to \Ohm^1_{\X/\Y}$. \vspace{1.5pt}
\index{ $\Cfr$2@$\hat d$ (universal derivation)! for formal-scheme maps}

Let $\Spf(B) = \U \subset \X$ and $\Spf(A) = \V \subset \Y$ 
be affine open subschemes with $f(\U)\subset \V$. Let $\Spec(A_n)=U_n\subset\U$ and  $\Spec(B_n)=V_n\subset\V$
be  defined as above.
 Then for each $m\ge0$ there are natural isomorphisms 
\begin{equation}\label{eq:diff3d}
\begin{aligned}
\Gamma(\U,\Ohm^m_{\X/\Y}) &= \Gamma(\U,\>\inlm{n}j_{n*}\Omega^m_{X_n/Y_n})
\cong \inlm{n}\Gamma(\U,\>j_{n*}\Omega^m_{X_n/Y_n})\\
&= \inlm{n}\Gamma(U_n,\>\Omega^m_{X_n/Y_n})
\cong \inlm{n}\Gamma(U_n,\>\Omega^m_{U_n/V_n})\\
&\cong \inlm{n}\Omega^m_{B_n/A_n} 
\cong \inlm{n}\Omega^m_{B_n/A}\> 
\underset{\rm\ref{prop:diff1}(ii)}\cong\, \Ohm^m_{B/A}\>.%
\footnotemark\\[-6pt]
\end{aligned}
\end{equation}
\footnotetext{The hypothesis in \ref{prop:diff1} that $B$ be essentially of pseudo-finite type over $A$ is not used in the proof of this last isomorphism.}

\begin{aprop}
\label{prop:diff4}
Let\/ $f \colon \X \to \Y$ be an essentially pseudo-finite type map 
of noetherian formal schemes. Then for\/ 
$m \ge 0$, the natural map is an isomorphism
$$
\bw^m_{\OX}\Ohm^1_{\X/\Y}\iso\inlm{\under{.2}{n}}\bw^m_{\OX} j_{n*}\Omega^1_{X_n/Y_n}=:\Ohm^m_{\X/\Y}\>, 
$$ 
and\/ $\Ohm^m_{\X/\Y}$ is a coherent\/ $\OX$-module. Moreover, if\/ $\U=\Spf(B)$ and\/ $\V=\Spf(A)$ are as above, and if the induced map $A \to B$
is essentially of pseudo-finite type, 
then,  with\/ $\sim\;=\;\sim_B,$ there are natural isomorphisms
$$
\big(\Ohm^m_{B/A}\big)^\sim\iso\Ohm^m_{\U/\V}\iso\Ohm^m_{\X/\Y}\big|_\U^{}.
$$
Furthermore, if\/ $f$ is smooth then\/ $f$ is flat 
and\/  $\Ohm^m_{\X/\Y}$ is locally free of finite rank. 
\end{aprop}
\begin{proof}
All the assertions are local, so we may assume that $\X=\U$ and $\Y=\V$.
For any affine
open subset $\U_b = \Spf(B_{\{b\}})\subset\U$
we have the natural isomorphisms 
\[
\Gamma(\U_b,\Ohm^m_{\X/\Y}) 
\underset{\rm\eqref{eq:diff3d}}\cong \Ohm^m_{B_{\{b\}}/A}
\cong \Ohm^m_{B/A} \otimes_B B_{\{b\}}
\cong \Gamma(\U_b,(\Ohm^m_{B/A})^{\sim}),
\]
where the second isomorphism comes from \ref{exam:diff3a},
and the third  from~\ref{prop:mod1}(iv). It follows that 
the map  corresponding\vspace{-2.5pt} to the isomorphism 
$\Ohm^m_{B/A}\cong \Gamma(\U,\Ohm^m_{\U/\V}) $ of \eqref{eq:diff3d}
is itself an isomorphism\vspace{-1pt}
$
 (\Ohm^m_{B/A})^{\sim} \iso \Ohm^m_{\U/\V}.
$ 
Since $\Ohm^m_{B/A}$\vspace{-1pt} is a finitely-generated $B$-module (see \ref{prop:diff1}\kern.5pt(ii)),  $\Ohm^m_{\U/\V}$ is coherent; and the
isomorphism\vspace{-.5pt} $\bw^m\Ohm^1_{\X/\Y}\iso\Ohm^m_{\X/\Y}$ results 
from \ref{prop:diff1}\kern.5pt(ii) and \ref{lem:mod2a}\kern.5pt(iv).\vspace{1pt}

The assertions concerning smooth maps follow from \ref{prop:morph2}\kern.5pt(ii) and \ref{prop:diff1}\kern.5pt(iii).
\end{proof}

\begin{adefi}
\label{def:diff4a}
Let $f \colon \X \to \Y$ be a smooth morphism of noetherian formal 
schemes. The \emph{relative dimension of~$f\>$}%
\index{relative dimension (of smooth formal-scheme map)} 
is the function---constant on
connected components of~$\X$---taking $x\in\X$ to the
rank of~$\Ohm^1_f$ at~$x$.\vspace{1pt}

When the relative dimension is constant on all of $\X$, we may identify it
with its value (a nonnegative integer).
\end{adefi}

\begin{aexam}
\label{exam:diff5}
An open immersion is smooth of relative dimension 0 (\ref{exam:diff3a}).
The map $\Spf(C) \to \Spf(A)$ obtained from \ref{exam:diff3c} 
is smooth of relative dimension $r$.
\end{aexam}

\begin{adefi}
\label{def:diff6}
Let $f \colon \X \to \Y$ be a smooth morphism of noetherian formal 
schemes. With $d_i$ the
relative dimension of\/ $f$ on the connected component~$\X_i$ of~$\X$,
we denote by $\omega_{\<\<f}$ the 
invertible $\OX$-module whose restriction to~ $\X_i$ is $\Ohm^{d_i}_{\X/\Y}$. 
\end{adefi}

\begin{aprop}
\label{lem:diff7}
Let\/ $f \colon \X \to \Y$ and\/ $g \colon \Y \to \Z$ be smooth maps
of relative dimensions\/ $d$ and\/ $e$ respectively.
Then\/ $gf$ is smooth of relative dimension\/ $d+e$ and
there is a canonical isomorphism 
\[
f^{\ast}\omega_g \otimes_{\OX}\omega_{\<\<f} \iso \omega_{gf}.
\]
\end{aprop}

\begin{proof}
The question is local on $\X$, so we may assume that $\X$, $\Y$, and $\Z$ are connected, and that $d$ and $e$ are integers.

Let $\I,\J,\sK$ be defining ideals in $\OX,\OY,\OZ$ respectively such that
$\J\OX \subset \I$ and $\sK\OY \subset \J$. Reducing modulo the $n$-th
powers of these defining ideals we get maps of schemes
$X_n \xto{\;f_n\;} Y_n \xto{\;g_n\;} Z_n$ and hence 
a sequence
\[ 
0 \to f_n^{\ast}\Omega^1_{Y_n/Z_n} \to \Omega^1_{X_n/Z_n} 
\to \Omega^1_{X_n/Y_n} \to 0. 
\]
Let $j_n\colon X_n \to \X$ be the canonical immersion. Applying
$j_{n*}$ and taking inverse limits we obtain a sequence\vspace{-3pt}
\[ 0 
\to \inlm{n}j_{n*}f_n^{\ast}\Omega^1_{Y_n/Z_n} \to \Ohm^1_{\X/\Z} 
\to \Ohm^1_{\X/\Y} \to 0. \\[-3pt]
\]
Let $i_n\colon Y_n \to \Y$ be the canonical immersion. There are 
natural maps
\[
f^*\inlm{n}i_{n*} \to \inlm{n}f^*i_{n*} \to \inlm{n}j_{n*}f_n^*.\\[-3pt]
\]
Hence there is a natural map
\[
f^*\Ohm^1_{\Y/\Z} = f^*\inlm{n}i_{n*}\Omega^1_{Y_n/Z_n}
\to \inlm{n}j_{n*}f_n^*\Omega^1_{Y_n/Z_n},\\[-3pt]
\]
and so we have a natural sequence
\[ 
0 \to f^{\ast}\Ohm^1_{\Y/\Z} \to \Ohm^1_{\X/\Z} 
\to \Ohm^1_{\X/\Y} \to 0. 
\]
One checks, using \ref{lem:mod4a}, that if $\X=\Spf(C)$, $\Y=\Spf(B)$, and $\Z=\Spf(A)$ are affine and if the induced maps $A \to B \to C$ are 
essentially of pseudo-finite type, then
this sequence is the same as the one obtained by sheafifying, via $\sim_C$, a split exact sequence as in \ref{lem:diff2}.  Exactness being a local property, we have therefore constructed, in the general case, an exact sequence of locally free sheaves (see \ref{prop:diff4}). It follows that the rank 
of~$\Ohm^1_{\X/\Z}$ is~$d+e$, and that there is a canonical isomorphism
\[\postdisplaypenalty 10000
\bw^{d+e}_{\X}\Ohm^1_{\X/\Z} \cong \bw^{e}_{\X}f^{\ast}\Ohm^1_{\Y/\Z}
\otimes \bw^{d}_{\X}\Ohm^1_{\X/\Y}.
\] 
Since $f^*$ commutes with exterior powers, the Lemma results.
\end{proof}

\begin{aprop}
\label{homelessa} Let\/ $f\colon\X\to\Z$ and\/ $g\colon\Y\to\Z$ be formal-scheme maps, with\/ $f$ essentially of pseudo-finite type, so that  the 
projection\/ $q\colon\W\set\X\times_\Z\Y\to \Y$ is also essentially of pseudo-finite type $($cf.~proof of \textup{\ref{lem:form2}):} 
\[
\begin{CD}
\hbox to 0pt{\hss$\X\times_{\Z}\Y=\:$}\W @>q>> \Y \\
@VpVV    @VVgV \\
\X @>>f> \Z
\end{CD}
\]
Then for the
projection\/ $p\colon\W\to \X$  there is a natural isomorphism
\[
p^*\Ohm^1_{\X/\Z} \iso \Ohm^1_{\W/\Y}.
\]
\end{aprop}

\begin{proof} This is just a globalization of the noetherian case of \ref{homeless}, and can be proved similarly (or reduced, via pasting of local maps, to \emph{loc.\,cit.}, see \ref{lem:mod4a}). Indeed, since both $p^*\Ohm^1_{\X/\Z}$ and $\Ohm^1_{\W/\Y}$ are coherent $\OW$-modules (see \ref{prop:diff4}), it suffices to find a natural isomorphism, functorial in the coherent $\OW$-module $\N\<$,
$$
\Hom_{\OW}\bigl(\Ohm^1_{\W/\Y},\>\N\>\bigr)\iso 
\Hom_{\OW}\bigl(p^*\Ohm^1_{\X/\Z}, \>\N\>\bigr),
$$
which can be done, as in the proof of \ref{homeless}, via restriction of (sheafified) derivations. Details are left to the reader.
\end{proof}

\begin{acor}
\label{homelessb}
With hypotheses as in\/ $\ref{homelessa}$, assume further that\/ $f$ 
is smooth of constant relative dimension $d\in\mathbb N$. Then\/ 
$q$ is also smooth of relative dimension\/~$d$. 
Furthermore there exists a natural isomorphism
$p^*\<\omega_{\<\<f} \iso \omega_q$.
\end{acor}

\begin{aprop}
\label{prop:homelessc}
Let $\Z \xto{\;i\;} \X \xto{\;h\;} \Y$ be maps of formal schemes
where $i$ is a closed immersion, $h$ is essentially of
pseudo-finite type and $hi$ is smooth. 
Let $\I$ be the coherent ideal in $\OX$ corresponding to~$i$.
Then there is an exact sequence of $\OZ$-modules
\[
0 \xto{\quad} \I/\I^2\big|_{\Z} \xto{\quad} i^*\Ohm^1_{\X/\Y} 
\xto{\quad} \Ohm^1_{\Z/\Y} \xto{\quad} 0. 
\]
\end{aprop}

\enlargethispage*{5pt}
\begin{proof}
The natural map $i^*\Ohm^1_{\X/\Y} \to \Ohm^1_{\Z/\Y}$ is defined\vspace{-2pt}
in the same manner as the map $f^*\Ohm^1_{\Y/\Z} \to \Ohm^1_{\X/\Z}$
is defined in the proof of~\ref{lem:diff7}. Consider
the natural sheaf homomorphism 
$$
\psi \colon \I \to \OX \xto{\hat d_{\X/\Y}} \Ohm^1_{\X/\Y} 
\to \Ohm^1_{\X/\Y} \otimes_{\OX} \OX/\I. 
$$ 
Let $\U = \Spf(B) \subset \X$,
$\V = \Spf(A) \subset \Y$ be affine open subsets such that 
$h(\U) \subset \V$ and the induced map $A \to B$ is essentially 
of pseudo-finite type. Let $I = \I(\U)$, so that via \ref{prop:mod1}(iii)
(see also  Remark (a)), $I^{\sim_B} = \I$. By \ref{lem:mod2a}(iii)\vspace{.4pt}
 $I^2 = \I^2(\U)$, and by~(i) and~(iii) of~\ref{prop:mod1},
$I/I^2 = (\I/\I^2)(\U)$. 
Using\vspace{.8pt} \ref{prop:diff4} and \ref{lem:mod2a}(ii), 
we see that $\psi$ induces over $\U$, the obvious natural map\vspace{-2pt} 
$\psi(\U) \colon I \to \Ohm^1_{B/A} \onto \Ohm^1_{B/A}\otimes_B B/I$.
Now $\psi(\U)$ sends $I^2$ to~0\vspace{-1pt} and the induced map
$I/I^2 \to \Ohm^1_{B/A}\otimes_B B/I$ is $B/I$-linear.
Since $\U,\V$ can be chosen to be arbitrarily small,\vspace{,4pt}
it follows that $\psi$ sends~$\I^2$ to~0 
and that the induced sheaf homomorphism
$\I/\I^2\big|_{\Z} \to i^*\Ohm^1_{\X/\Y}$ is $\OZ$-linear.
We therefore have a sequence of $\OZ$-modules
\[
0 \xto{\quad} \I/\I^2\big|_{\Z} \xto{\quad} i^*\Ohm^1_{\X/\Y} 
\xto{\quad} \Ohm^1_{\Z/\Y} \xto{\quad} 0. \tag{$*$}
\]
Exactness of $(*)$ is a local property.
If $\Z = \Spf(C)$, $\X = \Spf(B)$, $\Y = \Spf(A)$ are affine,
then the sequence~$(*)$ is the same as the one obtained by 
sheafifying, via $\sim_C$, the exact sequence in \ref{lem:diff3}.   
Therefore $(*)$ is also exact.
\end{proof}

\pagebreak[3]

For the next result we use the following fact. 
If $\W = \Spf(C)$ is an affine formal scheme, then for any
$w \in \W$, 
the complete local ring $D \set (\OWw,m_w)^{\wh{\;}}$ 
is isomorphic to $(C_{\pfr}, \pfr C_{\pfr})^{\wh{\;}}$ where $\pfr$ 
is the open prime in $C$ corresponding to $w$
(see proof of \ref{prop:mod7}(iv)). Therefore the natural
map $C \to D$ is smooth of relative dimension 0 
(see \ref{exam:diff3a}).

\begin{aprop}
\label{prop:app3}
Let\/ $\fXY$ be\vspace{.5pt} an essentially-pseudo-finite-type map of  formal schemes. Let\/ $x\in\X,$\vspace{-1pt}
set\/ $B \set {\OXx},$ $ A \set {\cO_{\Y, f(x)}},$ and let \/  $\wh{B},$ $\wh{A}$ be the completions\vspace{.5pt}  of\/ $B,$ $A$ along their respective\vspace{-2pt} 
maximal ideals. Then, with\/ $\Ohm^m_{f\<,\>x}$ the stalk of\/ $\Ohm^m_{f}$ at\/ $x,$ there are natural isomorphisms
\[
\Ohm^m_{f\<,\>x}\otimes_B \wh{B} \cong \Ohm^m_{\wh{B}/\wh{A}}
\qquad(m\ge 0).
\]
\end{aprop}
\begin{proof}
Let $\U = \Spf(S), \V = \Spf(R)$ be affine\vspace{.5pt} neighborhoods of $x,y$ 
respectively,   such that $f(\U) \subset \V$. 
By \ref{prop:diff4}, 
$\Ohm^m_f\big|_{\U} \cong (\Ohm^m_{S/R})^{\sim_S}\<,$
whence  $\Ohm^m_{f\<,\>x}\cong\Ohm^m_{S/R} \otimes_S B$. 
There are then natural isomorphisms
\[
\Ohm^m_{f\<,\>x}\otimes_B \wh{B} \iso (\Ohm^m_{S/R} \otimes_S B) \otimes_B \wh{B}
\iso \Ohm^m_{S/R} \otimes_S \wh{B} \xto{\;\alpha\;} \Ohm^m_{\wh{B}/R}
\xto{\;\beta\;} \Ohm^m_{\wh{B}/\wh{A}},
\]
where $\alpha$ is the isomorphism which follows from \ref{lem:diff2}, applied
to the sequence $R \to S \to \wh{B}$ (by the remark preceding \ref{prop:app3}, 
$S\to \wh B$ is smooth),  while $\beta$ 
is the isomorphism obtained similarly from the sequence $R \to \wh{A} \to \wh{B}$.
\end{proof}

\begin{acor}
\label{cor:app4}
With notation and assumptions as in\/ \textup{\ref{prop:app3},} assume
further that\/ $f\colon(\X,\Delta_1)\to (\Y,\Delta)$ is a smooth map in\/ $\bbFc$ $($see section\/ \textup{\ref{subsec:formal}).} Let~$d$ be
the relative dimension of\/~$f$ at $x$, $p\set \Delta_1(x),$ and\/ $q\set \Delta(y)$.
Let $m_A$  be the maximal ideal of\/ $A.$ Then
\[
\dim(B/m_AB) = p+d-q.
\]
\end{acor}
\begin{proof} 
By \ref{prop:morph2}\kern.5pt (ii), $B$ is a formally smooth $A$-algebra;\vspace{.6pt} and so by(v) and (vi) of \ref{exam:morph0}, $\wh B$ is a formally smooth $\wh A$-algebra.
By \ref{prop:app3}, $\Ohm^1_{\wh{B}/\wh{A}}$ is free of rank~$d$,  
and by definition of $\bbFc$, 
\[
q-p = \Delta(y) -\Delta_1(x) = \mbox{tr.deg.}_{k(y)}k(x).
\]
So by (\cite[Lemma 3.9]{Hu}),
\[
\dim(B/m_AB) = \dim(\wh{B}/m_A\wh{B}) = d - (q-p). 
\]
\end{proof}

\newpage

\section{Local cohomology and Cousin complexes}
\label{sec:locou}
In \S\ref{subsec:loc} we recall some basic facts about
local cohomology modules and about various relations between 
different definitions of local 
cohomology in the presence of quasi-coherence and torsionness 
conditions. 
In \S\ref{subsec:cousin} we review 
basic definitions and properties of
Cousin complexes, and of the Cousin functor $E$ 
from the derived category to Cousin complexes. 
(For smooth formal-scheme maps $\X\to \Y$, $E$ will play an 
important role in the construction of
a functor from Cousin complexes
over~$\Y$ to Cousin complexes over~$\X$.) 
The final subsection is on Cohen-Macaulay 
complexes, the derived category counterparts of Cousin
complexes. The results there are based on Suominen's\index{Suominen, Kalevi}
work \cite{Suo}, in which it 
is shown that $E$ induces an 
equivalence of categories from the category of 
Cohen-Macaulay complexes to the category of Cousin complexes
(where the latter, and hence the former, is abelian).  

\subsection{Local cohomology}
\label{subsec:loc}

Let $X$ be a topological space, $Z$ a closed subset, and $\F$ an 
abelian sheaf on~$X$. We denote by $\iG{Z}\F$ the 
\index{ $\Hom^{\cont}$2@$\iG{}$ \kern-2pt (torsion functor)!$\iG{z}$@$\iG{Z}$ ($Z$ a closed subset)}
subsheaf of~$\F$ whose sections over any
open $U \subset X$ are those 
elements of~$\F(U)$ whose support lies in~$Z$.
Recall further, from~\S\ref{subsec:tor}, 
the definition of the functor $\iG{\I}$ for a ringed space $(X, \cO_{\<\<X})$ 
and  an $\cO_{\<\<X}$-ideal $\I$, and of the functor $\iG{\afr}$
for an ideal $\afr$ in a ring $A$.

\begin{alem}
\label{lem:loc3}
Let\/ $(X,\cO_{\<\<X})$ be a ringed space, $\F$ an\/ $\cO_{\<\<X}$-module, 
and\/ $x \in X$ a point. If the\/ $\cO_{\<\<X}$-ideal~$\I$ is locally finitely-generated, 
then\/ $(\iG{\I}\<\F\>)_x = \iG{\I_x}\F_x$, 
where each side of the equality has been 
identified naturally with an\/ $\cO_{\<\<X\<,\>x}$-submodule of\/~$\F_x$.
\end{alem}  
\begin{proof}
The equality results from the natural isomorphisms
\begin{align}
\bigl(\dirlm{n}\sHom_{\cO_{\<\<X}}\<\<(\cO_{\<\<X}/\I^n , \>\F\>)\bigr)_x 
&\iso \dirlm{n}\bigl(\sHom_{\cO_{\<\<X}}\<
\<(\cO_{\<\<X}/\I^n , \>\F\>)\bigr)_x \notag \\
&\iso \dirlm{n}\Hom_{\cO_{\<\<X\<,\>x}}(\cO_{\<\<X\<,\>x}/\I^n_x , \>\F_x),\notag
\end{align}
the second isomorphism holding because
locally, $\cO_{\<\<X}/\I^n$ has finite presentation.
\end{proof}

It follows that if $Z$ is the support of~$\cO_{\<\<X}/\I$ then 
for any $\F \in \A(X)$ there is a natural inclusion 
$\iG{\I}\<\F \subset \iG{Z}\F$. In general, this inclusion is not 
an equality. 

\smallskip
\begin{alem}
\label{lem:loc2}
Let\/ $(\X,\cO_{\X})$ be a noetherian formal scheme and\/ $Z \subset \X$ a 
closed subset. Let\/ $\I$ be an open coherent $\OX$-ideal
such that the support of\/ $\cO_{\X}/ \I$ is\/ $Z$. Then for any\/ 
$\F \in \Aqct(\X\>)$ we have\/
\[ \iG{\I}\<\F = \iG{Z}\F. \]
\end{alem}  

\begin{proof}
It suffices to show that for any affine open subset $\U \set \Spf (A)$
we have 
$(\iG{Z}\<\F\>)(\U) \subset (\iG{\I}\F\>)(\U)$, \vspace{.6pt}
the opposite containment having been noted above.
Set $F = \F(\U)$ and $I= \I(\U)$. 
By \ref{prop:mod5}(iii), $\F \in \Avc(\X\>)$ 
and $\F|_{\U}^{} \cong F^{\sim}$, whence by~\ref{lem:mod4c}
and \ref{prop:mod1}(iv), $(\iG{\I}\<\F\>)(\U) = \iG{I}F$. 
In particular, since\vspace{.4pt} 
$\F$ is a torsion $\OX$-module, this equality for an ideal of definition
shows that $F$ is a torsion $A$-module. 
So by~\ref{prop:mod1}(iv) and \ref{lem:mod4b}, 
for any $s \in A$, with $\U_s \set \Spf(A_{\{s\}})$ we have 
$\F(\U_s) = F_s$. If~$s \in I$, then $\U_s \cap Z = \emptyset$,
so the image in~$F_s$ of any $f \in (\iG{Z}\F\>)(\U)$ is zero,
i.e., any such~$f$ is annihilated by a power of~$s$, and hence, 
$I$ being finitely generated, by a power of~$I$.
Thus 
$(\iG{Z}\<\F\>)(\U) \subset \iG{I}F=(\iG{\I}\F\>)(\U)$, as desired.
\end{proof}

For an abelian sheaf $\F$ on a topological space $X,$ and  
$x \in X$, let $\iG{x}\F\subset\F_x$ be the stalk at~$x$ of the 
\index{ $\Hom^{\cont}$2@$\iG{}$ \kern-2pt (torsion functor)!$\iG{y}$@$\iG{x}$ ($x$ a point)}
sheaf $\iG{Z}\F$ with $Z\set\,\ov{\!\{x\}\!}\,$, 
the closure in~$\X$ of the set~$\{x\}$.

\begin{acor}
\label{cor:loc4}
Let\/ $\X$ be a noetherian formal scheme,  and~$\F\in \Aqct(\X\>)$.  
For any\/ $x\in\X$ let\/ $m_x$ be the maximal ideal of the 
local ring~$\cO_{\X\<,\>x}$. 
Then 
$$
\iG{m_x}\F_x = \iG{x}{\F}\subset\F_x.
$$ 
\end{acor}  
\begin{proof}
Let $\I$ be the largest coherent ideal defining $Z = \ov{\{x\}}$. 
Then $\I$ is an open ideal and the stalk of 
$\I$ at $x$ is $m_x$ (since this holds true over the scheme 
obtained by reducing modulo an ideal of definition of $\X$). 
By~\ref{lem:loc3} and~\ref{lem:loc2} we have
$\iG{m_x}\F_x = (\iG{\I}\<\F\>)_x = (\iG{Z}\F\>)_x$. 
\end{proof}

\smallskip

The functors $\iG{?}$ defined above are all left exact. 
We denote by~$\R\iG{?}$
the corresponding right-derived functors  
and by $H^i_? := H^i \R \iG{?}$ the 
corresponding $i$-th right-derived functors. 

{\makebox[0pt]{\raisebox{-1ex}{\hspace{-6em}\;{\Huge $\mathsf{Z}$}}}}
For a complex $\cFb$ on a ringed space\vspace{-1pt} $(X, \cO_{\<\<X})$, 
$H^i_? \cFb \set H^i \R\iG{?}\cFb$ is the \emph{abelian group}
sometimes referred to in the literature as
``local hypercohomology with supports in ?\kern1pt'' and sometimes denoted
``$\>\mathbb H^i$\dots'' 
The symbol 
``$\,H^i_? \cFb\,$'' might also, conceivably, denote the \emph{complex} whose
$n$-th term is $H^i_?\F^n\<$, differentials being the natural
induced ones; in this paper, that will never be so.
\index{ $\Hom^{\cont}$x@$H^i_{?}:= H^i \R \iG{?}$ (hypercohomology supported in ?)}

\smallskip
We now review the construction of local cohomology via 
``direct-limit Koszul complexes". 
\index{Koszul complex}
Let $(X, \cO_{\<\<X})$ be a ringed space.
For an element $t \in \Gamma(X, \cO_{\<\<X})$, let~$\cKb(t)$ be the complex 
which in degrees 0 and 1 is multiplication by~$t$ from 
$\cK^0(t) \set \cO_{\<\<X}$ to $\cK^1(t) \set \cO_{\<\<X}$, 
and is zero elsewhere.
For $0 \le r \le s$, there is a map of complexes 
$\cKb(t^r) \to \cKb(t^s)$ which is identity in degree~0 and
multiplication by~$t^{s-r}$ in degree 1. Thus we get a direct
system of complexes, whose direct limit we denote 
by~$\cKbi(t)$. For any sequence $\bt = (t_1,...,t_n)$ 
of elements in~$\Gamma (X, \cO_{\<\<X})$ we set 
(with $\otimes= \otimes_{\cO_{\<\<X}}$)
\[ 
\cKbi(\bt) \set 
\cKbi(t_1) \otimes \ldots \otimes \cKbi(t_n); 
\]
\index{ $\K$@$\cKbi(\bt)$ ($\subdirlm{}\mkern-3.4mu$ Koszul\vspace{-1.5pt} complex of sheaves)}
\index{ $\K$@$\cKbi(\bt)$ ($\subdirlm{}\mkern-3.4mu$ Koszul\vspace{-1.5pt} complex of sheaves)!$\cKbi(\bt,\cFb) \set  \cKbi(\bt) \otimes \cFb$}%
and for any complex $\cFb$ of $\cO_{\<\<X}$ modules  set
$\cKbi(\bt,\cFb) \set  \cKbi(\bt) \otimes \cFb\<$. 
With $\I$ the $\cO_{\<\<X}$-ideal generated by the sequence 
$\bt$, there are natural identifications 
\[
\iG{\I}\<\F^j = \ker\bigl(\cKi^0(\bt,\F^j) \to \cKi^1(\bt,\F^j)\bigr) \qquad 
   (j\in \mathbb Z),
\]
yielding a map of complexes $\iG{\I}\cFb \to \cKbi(\bt,\cFb)$, whose composition with the natural map $\cKbi(\bt,\cFb) = \cKbi(\bt) \otimes \cFb\to\cKi^0(\bt)\otimes \cFb=\cFb$ is the inclusion $\iG{\I}\cFb \into\cFb$.
Note that the stalk of $\cKbi(t)$ at any point $x \in X$ looks like 
the localization map $\cO_{\<\<X\<,\>x} \to \cO_{\<\<X\<,\>x}[1/t\>]$. 
It follows that 
$\cKbi(\bt)$ is a complex of flat $\cO_{\<\<X}$-modules. In particular,
$\cKbi(\bt,-)$ takes quasi-isomorphisms to quasi-isomorphisms and hence
induces a functor $\D(X) \to \D(X)$, also denoted by
$\cKbi(\bt,-)$. In fact $\cKbi(\bt,-)$ is a $\delta$-functor
(see \S\ref{subsec:conv}, \eqref{conv7}), with 
$\Theta \colon \cKbi(\bt,\cFb[1]) \iso (\cKbi(\bt,\cFb))[1]$ 
given in degree $p+q$ by $(-1)^{p}$ times the identity map of
$\cKi^p(\bt) \otimes \F^{\>q+1}$.

\begin{aprop}
\label{prop:loc5}
Let\/ $(\X, \OX)$ be a noetherian formal scheme, and\/ $\bt$ a 
finite sequence in\/~$\Gamma (\X,\OX)$. Let\/ $\I$ 
be the $\OX$-ideal generated by\/~$\bt$.
Then for any complex of\/ $\A(\X\>)$-injectives\/ $\cLb\<,$ the
natural map\/ $\iG{\I}\cLb \to \cKbi (\bt,\cLb)$ is a 
quasi-isomorphism. Moreover, the induced natural
isomorphism of functors
\[ 
\R\iG{\I} \iso \cKbi (\bt,-) 
\] 
is an isomorphism of\/ $\delta$-functors. Hence there is a\/ $\delta$-functorial isomorphism\/
$$\smash{\R\iG{\I}\cFb \iso\R\iG{\I}\OX \Otimes_{\<\!\OX}\> \cFb}
\qquad\bigl(\cFb\in\D(\X\>)\bigr)$$ 
whose composition with the natural
map\/ \smash{$\R\iG{\I}\OX \Otimes_{\<\!\OX}\> \cFb \to\OX \Otimes_{\<\!\OX}\>\cFb\cong \cFb$} is the natural map\/ $\R\iG{\I}\cFb\to\cFb$.
\end{aprop}  
\begin{proof}
See \cite[Lemma 3.1.1, (1)$\:\Rightarrow\:$(2) ]{AJL1}. Verification of the $\delta$-part is straightforward. (The natural $\delta$-structure on 
$\R\iG{\I}$ can be extracted from \cite[Example 2.2.4]{Li}.) The rest is left to the reader.
\end{proof}

\smallskip

Before stating some consequences of~\ref{prop:loc5} we note that 
the constructions above have obvious analogs over rings.
For a ring $A$ and a sequence~$\bt$ in~$A$, the direct-limit Koszul 
complex $\Kb_{\infty}(\bt)$ of $A$-modules is defined in 
a manner similar to the above.
If $A$ is noetherian and~$\Ib$ a complex of $A$-injectives, 
then the natural map $\iG{\bt A}\Ib \to \Kb_{\infty}(\bt,\Ib)$ 
is a quasi-isomorphism.
\index{ $\K$@$\cKbi(\bt)$ ($\subdirlm{}\mkern-3.5mu$ Koszul\vspace{-1.5pt} complex of modules)}
\index{ $\K$@$\cKbi(\bt)$ ($\subdirlm{}\mkern-3.5mu$ Koszul\vspace{-1.5pt} complex of modules)!$\cKbi(\bt,\Ib) \set  
  \cKbi(\bt) \otimes \Ib$}

\begin{alem}
\label{lem:loc6}
Let\/ $(\X, \OX)$ be a noetherian formal scheme.
\begin{enumerate}
\item For any coherent\/ $\OX$-ideals\/ $\I, \J,$ the natural map\/
$\R\iG{\I+\J} \to \R\iG{\I}\R\iG{\J}$ is an isomorphism.\vspace{1pt}
\item For any sequences\/ $\bt,\bs$ in\/ $G\set\Gamma(\X,\OX)$ such that\/
$\bs G \subset \sqrt{\bt G},$ and for any complex\/ 
$\cFb \in \bC(\X\>),$ the map\/ $($with $\otimes = \otimes_{\OX})$
\[
\cKbi(\bt)\otimes\cKbi(\bs)\otimes\cFb \to \cKbi(\bt)\otimes\cFb
\] 
induced by the natural map\/ $\cKbi(\bs) \to \OX$
is a quasi-isomorphism. \vspace{1pt}
\item For any\/ $x \in \X$ and  $\cFb \in \D(\X\>)$
the natural map\/ $(\R\iG{\I}\cFb)_x \to \R\iG{\I_x}\cFb_x $ 
is a\/ $\D(\OXx)$-isomorphism.
\end{enumerate}
\end{alem}

\begin{proof}
(i). One argument is given in the last four lines on page 25 of \cite{AJL1}.
Another is that the map in question is the canonical isomorphism resulting from
the fact that both source and target are right adjoints of the inclusion into $\D(\X\>)$
of the full subcategory whose objects are complexes whose homology is $(\I+\J)$-torsion, see \cite[p.\,49, Prop.\,5.2.1\kern.5pt(c)]{AJL2}. One can also use~\ref{prop:loc5}, as follows.

The assertion being local, we may assume that 
$\X$ is affine, so that there are sequences $\bt, \bs$ 
in~$\Gamma(\X, \OX)$ generating $\I,\J$ respectively. 
For any $\cFb \in \D(\X\>)$, we claim that the following 
natural diagram commutes:
\[
\begin{CD}
\R\iG{\I+\J}\cFb @>>> \R\iG{\I}\R\iG{\J}\cFb \\
@VVV @VVV \\
\cKbi(\bt,\bs)\otimes\cFb @>>> \cKbi(\bt)\otimes\cKbi(\bs)\otimes\cFb
\end{CD}
\]
The bottom map is an isomorphism, as  
are the vertical maps (by \ref{prop:loc5}), and the assertion results.

Verification of commutativity is straightforward once the diagram is expanded as in
(\ref{lem:loc6}.1) below---where $\cFb$ is assumed, w.l.o.g., to be K-injective,
$\cK_{\bt}$ is an abbreviation for $\cKbi({\bt})$, etc., and $\cK_{\bs}\otimes\cFb\to\cEb$ is a K-injective resolution (see \cite[p.\,19, proof of $(2)'\Rightarrow(2)$]{AJL1}).

\begin{figure}[t]
\centerline{\bf ({\ref{lem:loc6}.1})}
$$
\begin{CD}
 \iG{\I+\J}\cFb @= 
  \underset{\UnderElement{}{\|}{53pt}{}}{\iG{\I}\iG{\J}\cFb}
 @>>> \R\iG{\I}\iG{\J}\cFb  \\[-1.5pt]
 @| @. @| \\
  \underset{\UnderElement{}{\downarrow}{50pt}{}}{\R\iG{\I+\J}\cFb} @. @. \R\iG{\I}\R\iG{\J}\cFb\\[-6pt]
  @. @. @VVV \\
 @. @.  \R\iG{\I}(\cK_{\bs}\otimes \cFb)
    @>>> \R\iG{\I}\cEb\\[-1pt]
  @.  @. @AAA @|  \\
 @. \iG{\I}\iG{\J}\cFb @>>> \iG{\I}(\cK_{\bs}\otimes \cFb)@>>> \iG{\I}\cEb\\[-1.5pt] @. @VVV @VVV @VVV \\[-1.5pt]
 \cK_{\bt,\bs}\otimes \cFb @= \cK_{\bt,\bs}\otimes \cFb @=
    \cK_\bt \otimes \cK_\bs \otimes \cFb @<\Iso<<  \cK_\bt \otimes \cEb
\end{CD}
$$
\end{figure}

\smallskip
(ii). Let $\I,\J$ be the $\OX$-ideals   generated by 
$(\bt),(\bt,\bs)$ respectively. Note that $\iG{\I} = \iG{\J}$.
Let $\cLb$ be  
a K-injective resolution of $\cFb\<$. In the commutative diagram
\[
\begin{CD}
 \cKbi(\bt)\otimes\cKbi(\bs)\otimes\cFb @>>>
        \cKbi(\bt)\otimes\cKbi(\bs)\otimes\cLb @<<< \iG{\J}\cLb \\
 @VVV  @VVV  @| \\
 \cKbi(\bt)\otimes\cFb @>>> \cKbi(\bt)\otimes\cLb @<<< \iG{\I}\cLb
\end{CD}
\]
the horizontal maps are quasi-isomorphisms, and the result follows.

\pagebreak[3]

(iii).  The assertion being local, we may assume that $\X$ is affine,
so that there is a sequence $\bt$ 
in~$\Gamma(\X, \OX)$ generating $\I$.  We may assume further that
$\cFb$ is  K-injective. The following natural 
$\D(\OXx)$-diagram---where $\Kb_{\infty}(\bt)$ is 
the direct-limit Koszul complex on $\bt$
over $\OXx$, $\cK_{\bt}$ is an abbreviation for $\cKbi({\bt})$, etc., and $\cFb_x\to E^\bullet$ is a K-injective $\OXx$-resolution---commutes:

$$
\begin{CD}
\mspace{-15mu}(\R\iG{\I}\cFb)_x =(\iG\I\cFb)_x   
 @>\Iso>\ref{lem:loc3}>\iG{\I_x}\cFb_x 
   @>>>\iG{\I_x}E^\bullet=\R\iG{\I_x}\cFb_x\mspace{-15mu}\\
@V\simeq VV @VVV @VV\simeq V \\
(\cK_\bt\otimes_{\OX}\<\cFb)_x @>\Iso>> K_{\bt}\otimes_{\OXx}\cFb_x @<\Iso<< K_{\bt}\otimes_{\OXx} E^\bullet\\
\end{CD}
$$
The maps in the bottom row are isomorphisms, as  
are the outside vertical maps (by \ref{prop:loc5} and its analog over rings), and the assertion results.
\end{proof}

\begin{alem}
\label{prop:nolab2b}
Let\/ $\fXY$ be a map of noetherian formal schemes, \mbox{$\cGb\in\D(\Y)$} and\/ $\cLb\in\D(\X\>).$ Let\/ $\I,\J$ be 
coherent ideals in~$\OX,\OY$ respectively with\/ 
\mbox{$\J\OX \subset \I$.} 
Then the natural map is an isomorphism
\[
\R\iG{\I}(\bL f^*\R\iG{\J}\cGb \Otimes_{\<\!\Y}\; \cLb) \iso
\R\iG{\I}(\bL f^*\cGb \Otimes_{\<\!\X}\; \cLb). 
\]
\end{alem}

\begin{proof} Using \ref{prop:loc5} and commutativity of $\bL f^*$ and  \smash{$\Otimes$}\vspace{.6pt}  we reduce to showing that the natural map is an isomorphism  
$\R\iG{\I}(\bL f^*\R\iG{\J}\OY) \iso\R\iG{\I}(\OX)$.  By \cite[p.\,53, Prop.\,5.2.8(b)]{AJL2}, the map
$\bL f^*\R\iG{\J}\OY\to\OX$ factors naturally as
$$
\bL f^*\R\iG{\J}\OY\iso\R\iG{\J\OX}\OX\to\OX
$$ (the first map  an isomorphism).
The conclusion follows then from 
\ref{lem:loc6}(i).
\end{proof}

\begin{alem}
\label{lem:loc5a}
In the situation of\/~$\ref{prop:loc5}$, let\/ $\cIb$ be a 
complex of\/~$\Aqct(\X\>)$-injectives. Then the natural map\/ 
$\iG{\I}\cIb \to \cKbi (\bt,\cIb)$ is a quasi-isomorphism. 
In particular, for any quasi-isomorphism\/ $\cIb \to \cLb$
with\/ $\cLb$
 a complex of\/ $\A(\X\>)$-injectives,
the natural map\/ $\iG{\I}\cIb \to \iG{\I}\cLb$ is a quasi-isomorphism.
\end{alem}
\begin{proof}
(Sketch)
Since it suffices to check that over affine open subsets of $\X$ 
the natural map $\varPsi \colon \iG{\I}\cIb \to \cKbi (\bt,\cIb)$ 
is a quasi-isomorphism,
we may assume that $\X = \U =\Spf(A)$ for some noetherian adic 
ring~$A$. Note that $\varPsi$ is obtained by applying the 
functor~$\sim$ of~\ref{prop:mod1} 
to the natural map 
$\psi \colon \iG{\bt A}\Ib \to \Kb_{\infty}(\bt,\Ib)$ where 
$\Ib = \Gamma(\U, \cIb)$.
Since $\cIb$ consists 
of $\Aqct(\U)$-injectives, and so by \ref{prop:mod7}\kern.5pt(ii) is
a direct sum of sheaves of the form $i_xJ_x=J_x^\sim$ (see proof of \ref{lem:mod6}, (i)$\Lra$(iii)\kern1pt),
therefore $\Ib$ is a direct sum of $A$-modules of the form $J_x$, hence consists of $A$-injectives.
Thus, by the ring-theoretic version of \ref{prop:loc5}, $\psi$ is a quasi-isomorphism, whence so is~$\varPsi$ .
\end{proof}

An alternate way of proving the above lemma is as follows. 
By \ref{prop:mod7}(ii), any \emph{indecomposable} 
injective $\cI$ in $\Aqct(\X\>)$ is an $\A(\X\>)$-injective,
and so by \ref{prop:loc5} the natural map 
$\varPsi \colon \iG{\I}\cI \to \cKbi (\bt,\cI)$ is a quasi-isomorphism.
Since $\varPsi$ behaves well with respect to direct sums
($\OX/\I^n$ being coherent, the functor
$\sHom_{\OX}(\OX/\I^n, -)$ commutes with direct sums)
therefore we can extend the result to any $\Aqct(\X\>)$-injective~$\cI$
(see again \ref{prop:mod7}\kern.5pt(ii)).
Finally using way-out type arguments
we extend the result to any complex $\cIb$ of $\Aqct(\X\>)$-injectives (cf.~\cite[p.\,22]{AJL1}).

\begin{aprop}
\label{prop:loc8}
Let\/ $\X$ be a noetherian formal scheme.
\begin{enumerate}
\item For any\/ $\cFb \in \Dqc(\X\>)$, 
we have\/ $\R \iGp{\X} \cFb \in \Dqct(\X\>)$. 
For any\/ $\cFb \in \Dt(\X\>)$, the canonical
map\/ $\R \iGp{\X} \cFb \to \cFb$ is an isomorphism.
\item The natural functors\/
$\D^+(\Aqct(\X\>))  \to \Dqc^+(\At(\X\>)) \to \Dqct^+(\X\>)$ are 
equivalences of categories. In particular, any\/ 
$\cFb \in \Dqct^+(\X\>)$ is isomorphic
to a bounded below complex of\/ $\Aqct(\X\>)$-injectives. 
\end{enumerate}
\end{aprop}  
\begin{proof}
For (i) see  \cite[p.\,49, 5.2.1]{AJL1}.
For (ii), see \cite[Theorem 4.8]{Ye} (or \cite[p.\,57, 5.3.1]{AJL2} for an unbounded version).
\end{proof}

\begin{aprop}
\label{prop:nolab1}
Let\/ $\X$ be a noetherian formal scheme.
Let\/ $Z$ be a closed subset of\/~$\X$ and\/ $\I$ an open 
coherent\/ $\OX$-ideal such that\/ $Z = \Supp(\OX/\I)$. 
\begin{enumerate}
\item For\/ $\cFb\< \in \Dqct^+(\X\>)$ the natural map is an isomorphism\/ 
$\R\iG{\I}\cFb\!\<\< \iso\!\R\iG{Z}\cFb \!$.
\item For any\/ $\cFb \in \Dqc^+(\X\>)$, the  natural
maps are isomorphisms 
\[
\R\iG{\I}\cFb \osi \R\iG{\I}\R\iGp{\X}\cFb
\iso \R\iG{Z}\R\iGp{\X}\cFb.
\]
\end{enumerate}
\end{aprop}
\begin{proof}
(i). By \ref{prop:loc8}(ii), we may assume that $\cFb$ 
is a bounded-below complex of $\Aqct$-injectives. 
Let $\cFb \to \cLb$ be an $\A(\X\>)$-injective resolution
of~$\cFb\<\<$. In the commutative $\D(\X\>)$-diagram
\[
\begin{CD}
\iG{\I}\<\cFb @>a>> \iG{\I}\cLb\cong \R\iG{\I}\<\cFb\!\!\!\!\!  \\
@VbVV @VVdV \\
\iG{Z}\cFb @>>c> \iG{Z}\cLb\cong \R\iG{Z}\cFb\!\!\!\!\!
\end{CD}
\]
the map $a$ is an isomorphism by \ref{lem:loc5a},
$b$ is an isomorphism by \ref{lem:loc2},
and $c$ is an isomorphism because 
$\cFb$ consists of flasque sheaves (see \ref{prop:mod7}(iii)),
which are $\iG{Z}$-acyclic.
Hence $d$ is an isomorphism.

(ii). Since $\I$ is open, 
the first isomorphism is given by \ref{lem:loc6}\kern.5pt(i).  
By \ref{prop:loc8}(i), $\R\iGp{\X}\cFb \in \Dqct^+(\X\>)$, and so the second
isomorphism is given by (i). 
\end{proof}

\pagebreak[3]

Let $\X$ be a noetherian formal scheme\vspace{-1.8pt} and $\cFb\in\Dqc^+(\X\>)$. 
For any $x \in \X$, set $Z \set \ov{\{x\}}$ and let $\I$ be the largest 
coherent ideal defining~$Z\<$, so that $\I$ is open and the stalk $\I_x$
of~$\>\I$ at~$x$ is the maximal ideal $m_x$ of $\OXx$.  Then
\ref{prop:nolab1}(ii) and \ref{lem:loc6}(iii)
give natural isomorphisms
\begin{equation}\label{eq:nolab3}
\R\iG{x}\R\iGp{\X}\cFb = (\R\iG{Z}\R\iGp{\X}\cFb)_x
\cong (\R\iG{\I}\cFb)_x \cong \R\iG{m_x}\cFb_x.\\[3pt]
\end{equation}

\begin{alem}
\label{lem:nolab1}
Let\/ $f \colon \X \to \Y$ be a map of noetherian formal schemes. Then$\>:$
\begin{enumerate}
\item\label{nolab1i1} 
$f^*(\Aqc(\Y)) \subset \Aqc(\X\>)$.
\item\label{nolab1i2} The category 
$\Aqc(\X\>)$ is closed under tensor products.
\item\label{nolab1i3} For any sequence\/ $\bt$ in $\Gamma(\X,\OX),$ and\/
$j\ge0,$ we have $\>\cKi^j(\bt) \in \Avc(\X\>)$.
\end{enumerate}
\end{alem}
\begin{proof} (i) holds for any map of ringed spaces.

(ii) holds over any ringed space.

(iii) holds because $\cKi^j(\bt)$ is a $\dirlm{}\!\!$ of 
finitely-generated free $\OX$-modules. 
\end{proof}

\begin{alem}
\label{lem:psm1}
Let\/ $f \colon \X \to \Y$ be a map of noetherian formal schemes, 
$y \in \Y,$ and $M$ a zero-dimensional\/ $\OYy$-module. Let\/ 
$\M = i_yM$ be the corresponding\/ $\OY$-module \emph{(\ref{lem:mod6}).} 
Suppose\/ $\G\in\Aqc(\X\>)$. Then for any\/ 
$x \in \X$ such that $f(x) \ne y,$  and any\/ $j\ge0,$
\[
H^j_x\R\iGp{\X}(f^{\ast}\<\M \otimes_{\X} \G) = 0  
 = H^j_{m_x}(f^{\ast}\M \otimes_{\X} \G)_x\>.
\]
\end{alem}

\begin{proof}
Set $\F \set f^{\ast}\<\M \otimes_{\X} \G$.
By \ref{lem:mod6}, $\M \in \Aqct(\X\>)$, and hence
by parts (i) and~(ii) of~\ref{lem:nolab1}, $\F \in \Aqc(\X\>)$. 
By \eqref{eq:nolab3},  there is an isomorphism
$H^j_x\R\iGp{\X}\F \cong H^j_{m_x}\F_x$
with~$m_x$the maximal ideal of~$\OXx$.
We consider two cases:\vspace{.6pt}

(a) When $f(x) \not\in\, \overline{\!\{y\}\!}\,$,
the stalk of $\M$ at $f(x)$ is 0, so $(f^{\ast}\<\M)_x = 0$, 
$\F_x = 0$, and  $H^j_{m_x}\F_x = 0$.\vspace{.6pt} 

(b) When  $f(x) \in\overline{\!\{y\}\!}\,$, 
 there exists a non-unit \mbox{$t \in \cO_{\Y,f(x)}$} 
whose image under the natural map $\cO_{\Y,f(x)} \to \OYy$ 
is a unit in~$\OYy$. Then~$t$ acts 
invertibly on~$M$ (since $M$ is an $\OYy$-module) and hence 
on~$\F_x$ and hence on $H^j_{m_x}\F_x$. But
the image of~$t$ under the natural local homomorphism 
$\cO_{\Y,f(x)} \to \OXx$ lies in~$m_x$ and therefore every 
element of the $m_x$-torsion module $H^j_{m_x}\F_x$ 
is annihilated by a power of~$t$. 
Thus $H^j_{m_x}\F_x = 0$.
\end{proof}

\subsection{Cousin complexes}
\label{subsec:cousin}

We will use the notion of  Cousin complex as in~\cite[Chap.\,IV\kern.7pt]{RD}. 
(Additional properties in a more general context may be found in~\cite{Suo}.)
We first review the relevant definitions. 

Throughout this subsection $X$ denotes a noetherian topological space 
in which every irreducible closed subset has a unique generic
point; and $X$ is assumed to be equipped with a filtration 
\index{filtration}
\begin{equation} 
\label{eq:coz-1}
Z^{\<\bullet}: \quad \cdots \supseteq Z^{p-1} \supseteq Z^p 
\supseteq Z^{p+1}\cdots\qquad(Z^p\subset X) 
\end{equation}
satisfying the following conditions (cf.~\cite[p.~240]{RD}):\vspace{1pt}

(a) It is strictly exhaustive, i.e., $Z^{p} = X$ for some 
$p \in \mathbb Z$.

(b) It is separated, i.e., ${\bigcap_p} Z^{p} = \emptyset$.

(c) Each $Z^{p}$ is stable under specialization.

(d) For any $p$, if $x \leadsto x^{\prime}$ is a specialization and 
$x, x^{\prime} \in Z^{p} \setminus Z^{p+1}$ then $x = x^{\prime}$.\vspace{1pt}

\noindent
Corresponding to the filtration $Z^{\<\bullet}$ there is a 
filtration by subfunctors of the identity functor:\vspace{-2pt}
\[ 
\Gamma^{\bullet}: \quad \ldots \supseteq \iG{Z^{p-1}} 
 \supseteq \iG{Z^p} \supseteq \iG{Z^{p+1}} \ldots\\[2pt]  
\]  
Recall that for $Z\subset X$ and $\F$ an abelian sheaf on~$X\<$, $\iG{Z}\>\F$ is the sheaf of 
sections of~$\F$ with support in~$Z$, so that
 $\iG{Z}$ is an idempotent  left exact functor. Recall also that for~$x\in X$, $\iG{x}\F$ is the stalk at $x$ of $\iG{\,\ov{\!\{x\}\!}\,}\F$. 
One sets $\iG{Z^p/Z^{p+1}}\F\set\iG{Z^p}\F/\iG{Z^{p+1}}\F$.
\index{ $\Hom^{\cont}$2@$\iG{}$ \kern-2pt (torsion functor)!$\iG{z}$@$\iG{Z}$ ($Z$ a closed subset)!$\iG{Z^p/Z^{p+1}}\set\iG{Z^p}/\iG{Z^{p+1}}$}

In \cite[p.\,226]{RD}  there is defined, for any flasque abelian sheaf
$\F$ on~$X\<$, an isomorphism 
(with $i_x$ as in~\ref{lem:mod6})
\begin{equation} 
\label{eq:coz0}
\iG{Z^p/Z^{p+1}}\>\F \;\; \iso \bigoplus_{x \in Z^p \setminus Z^{p+1}}
\! \! i_x(\iG{x}\F\>).
\end{equation}
This isomorphism arises thus: for any open $U\subset X$ and $\xi\in\(\iG{Z^p}\F\>)(U)$, the support of $\xi$, a closed subset of $Z^p\cap U$, contains only finitely many 
$x\in Z^p\setminus Z^{p+1}$; hence the stalk of $\iG{Z^p}\>\F$ at any 
such $x$ is $\iG{x}\F$, and  the resulting natural map
\[
\iG{Z^p}\F  \lra
\prod_{x \in  Z^p\setminus Z^{p+1}}i_x(\iG{x}\F\>)
\] 
has kernel $\iG{Z^{p+1}}\F$ and image $\oplus_{\{x \in  Z^p\setminus Z^{p+1}\}}\>i_x(\iG{x}\F\>)$. \vspace{1pt}

As homology commutes with direct sums ($X$ being noetherian), and $i_x$ is an exact functor, we can replace any  homologically bounded-below abelian complex~ $\cFb$ on~$X$ by a flasque resolution%
\footnote{The necessary $\iG{}\<\!$-acyclicity properties of flasque sheaves are given in
\cite[Chap.\,IV, \S1]{RD}.
}
 and deduce  a canonical functorial isomorphism
\begin{equation} \label{H=oplus}
\tag{$\ref{eq:coz0}'$}
H^n_{Z^p/Z^{p+1}}\cFb \set    
H^n\R\iG{Z^p/Z^{p+1}}\cFb \;\; \iso \bigoplus_{x \in Z^p \setminus Z^{p+1}}
 i_x(H^n_x\cFb). 
\end{equation}

An abelian sheaf $\G$ is said to 
\emph{lie on the $Z^p/Z^{p+1}$--skeleton of} $X\>$%
\index{skeleton}
if either of the following two equivalent conditions are satisfied 
(\cite[p.~231]{RD}):
\begin{itemize}
\item There is a family of 
abelian groups $(G_x)\ (x \in Z^p \setminus Z^{p+1})$ and an isomorphism 
$\G \; \cong  \; \oplus \: i_xG_x$.\vspace{1pt}
\item The natural maps $\G \gets \iG{Z^p}\G \to H^0_{Z^p/Z^{p+1}}\G$
are isomorphisms.
\end{itemize}
Thus for  $\cFb$ as above, the sheaf $H^n_{Z^p/Z^{p+1}}\cFb$
lies on the $Z^p/Z^{p+1}$--skeleton of~$X\<$.\vspace{1.5pt}

A complex of abelian sheaves ${\G}^{\bullet}$ on $X$ is called 
{\it a ${Z}^{\bullet}$-Cousin complex\/} if for each~$p \in \mathbb Z$,
\index{Cousin complex}
${\G}^p$~lies on the $Z^p/Z^{p+1}$--skeleton of $X$. The reference to
${Z}^{\bullet}$ is dropped in case of no ambiguity. Since 
$Z^p = X$ for some $p$, a Cousin complex is necessarily 
bounded below. The individual  graded pieces 
${\G}^p$, being direct sums of flasque sheaves, are themselves flasque.

By definition, the underlying graded object
of a Cousin complex~$\cGb$ on $X$ admits a decomposition parametrized by the
points of $X\<$, and this decomposition is {\em finer} than the usual 
$\mathbb Z$-graded decomposition of $\cGb$ as a 
complex. Specifically, one can associate to each $x\in X$ the abelian group $\cGb\<(x)$ underlying ~$\iG{x}\cGb$; and 
since for any group~$G$,  $\iG{x}(i_yG) = 0$ whenever $x \ne y$, the complex $\iG{x}\cGb$ vanishes in all degrees other than that~$p=p(x)$ such that 
$x \in Z^p \setminus Z^{p+1}\<$, so $\cGb\<(x) = \iG{x}\G^{p(x)}$. In 
particular, there is a canonical isomorphism
$$
\phi \colon \G^p\iso \bigoplus_{x \in Z^p \setminus Z^{p+1}}
\! \! i_x\bigl(\cGb\<(x)\bigr),
$$
uniquely determined by the property that for any 
$z \in Z^p \setminus Z^{p+1}\<$, the natural map 
$\cGb\<(z) = \iG{z}\G^{p} \xto{\iG{z}(\phi)} \cGb\<(z)$ is the identity map.

Let $\cGb$  be a Cousin complex of abelian sheaves on~$X$, and let
$\cO_{\<\<X}$ be a sheaf of commutative rings on~$X$.
If each $\cGb\<(x)$ has an $\cO_{\<\<X\<,\>x}$-module structure
then all the component sheaves~$\G^p$ have natural
$\cO_{\<\<X}$-module structures. 
Thus specifying a Cousin complex~$\cGb$ of~$\cO_{\<\<X}$-modules 
with $\cO_{\<\<X}$-linear differentials 
is equivalent to specifying the following data:\looseness=-1
\pagebreak[3]
\begin{itemize}
\item for each $x \in X$, an $\cO_{\<\<X\<,\>x}$-module $\cGb\<(x)$;
\item for each immediate specialization $\>\>x\leadsto y\>\>$ in
$X$, an $\,\cO_{\<\<X,y}$-linear map 
\mbox{$\delta_{x\<,\>y}\colon \cGb\<(x) \xto{\;\;} \cGb\<(y);$}\vspace{2pt}
\end{itemize}
\noindent subject to the following conditions:\vspace{2pt}
\begin{itemize}
\item for each $\xi\in\cGb\<(x),$ $\delta_{x\<,\>y}(\xi)=0$ for all but finitely many $y$;
\item if $x \leadsto z$ is a specialization with $x \in Z^p$
and $z \in Z^{p+2}$, then
$\Sigma_y \;\delta_{y\<,z}\delta_{x\<,\>y} = 0$
where $y$ ranges over the set of all intermediate specializations
$x \leadsto y \leadsto z$ with $x \ne y \ne z$.
\end{itemize}

\smallskip

Let $X$, $\cO_{\<\<X}$ and $Z^{\bullet}$ be as before.  
To every 
$\cFb \in \D^+(X)$ one associates functorially a Cousin
$\cO_{\<\<X}$-complex~$E_{Z^{\<\bullet}}{\F}^{\bullet}\<$, as follows. 
\index{ $\Fgt$@$E$ (Cousin functor associated to a filtration)}
Let $\cLb$ be an injective resolution of~$\cFb\<$.
From the natural exact sequence
\[
0 \to \iG{Z^{p+1}}/\iG{Z^{p+2}}\cLb \to \iG{Z^{p}}/\iG{Z^{p+2}}\cLb
\to \iG{Z^{p}}/\iG{Z^{p+2}} \cLb \to 0
\qquad(p \in \mathbb Z),
\]
one derives connecting homology
homomorphisms (see \S\ref{subsec:conv}, \eqref{conv9}) 
\[
\delta^p \colon H^p_{Z^p/Z^{p+1}}\cFb \lra H^{p+1}_{Z^{p+1}/Z^{p+2}}\cFb,
\]
whence the sequence
\begin{equation}
\label{eq:coz1z}
\cdots\lra H^{p-1}_{Z^{p-1}/Z^{p}}\cFb \xto{\<\delta^{p-\<1}\!}
H^p_{Z^p/Z^{p+1}}\cFb \xto{\;\delta^{p}\;}
H^{p+1}_{Z^{p+1}/Z^{p+2}}\cFb \lra \cdots
\end{equation}
which is in fact a complex, denoted~$E_{Z^{\<\bullet}}{\F}^{\bullet}$, with
degree-$p$ component \smash{$H^p_{Z^p/Z^{p+1}}\cFb\<$}. (The 
filtration~$\Gamma^{\bullet}$ of the identity 
functor gives rise to a spectral sequence 
whose \smash{$E^{p,0}_1$}~terms form the complex 
$E_{Z^{\<\bullet}}{\F}^{\bullet}$. 
The spectral sequence converges to the homology\- of~$\cFb$ when~$X$~has 
finite Krull dimension, see \cite[p.\,227, p.\,241]{RD}.) 
The isomorphism in~$(\ref{eq:coz0}')$ shows that 
$E_{Z^{\<\bullet}}{\F}^{\bullet}$ is a Cousin complex; 
and if $x \in Z^p \setminus Z^{p+1}$ then 
there are natural isomorphisms 
\begin{equation}
\label{eq:coz1}
(E_{Z^{\<\bullet}}{\F}^{\bullet})(x) 
= \iG{x}(E_{Z^{\<\bullet}}{\F}^{\bullet})^p
\underset{(\<\ref{eq:coz0}'\<)\:}\iso 
\iG{x}\<\Bigl(\bigoplus_{y \in Z^p \setminus Z^{p+1}}
\!\!i_y(H^p_y\cFb)\Bigr) \iso H^p_x\cFb\<.
\end{equation}\vspace{2pt}

Here are some basic properties
of Cousin complexes, to be used in the next subsection.

\begin{alem}
\label{lem:coz1a}
Let\/ $X,$ $\cO_{\<\<X}$ and\/ $Z^{\bullet}$ be as before.
Set\/ $E\set E_{Z^{\<\bullet}}$.
Let\/ $\cCb$ be a Cousin\/ $\cO_{\<\<X}$-complex.

\begin{enumerate}
\item
{\rm[a]}  
For any\/ $m\ge n, p \in \mathbb Z$ and\/ $x \in Z^p$
the natural maps are isomorphisms
\begin{align*}
\sigma_{\ge n}\cCb=\iG{Z^{\<n}}\cCb &\iso \R\iG{Z^{\<n}}\cCb,\\
\sigma_{\ge n}\sigma_{< m}\cCb = \iG{Z^n/Z^m}\cCb &\iso 
\R\iG{Z^n/Z^m}\cCb,\\
\cCb(x)[-p] = \iG{x}\cCb &\iso \R\iG{x}\cCb.
\end{align*}
There result  natural isomorphisms$\>:$\vspace{1.5pt}

\noindent{\rm [b]} 
$H^p\R\iG{x}\cCb=H^p_x\cCb \iso \cCb(x)=H^p(\cCb(x)[-p]);$\vspace{1pt} 

\noindent{\rm [c]}
 $H^p_xE\cFb  \underset{[\textup{b}]}\iso 
(E\cFb)(x) \underset{\eqref{eq:coz1}}\iso H^p_x\cFb
\qquad\bigl(\cFb \in \D^+(X)\bigr)$.\vspace{4pt}

\item  The  graded isomorphism\/ $E\cCb \iso \cCb$
made up of the punctual maps
$$
(E\cCb)(x)=\iG{x}H^p\R\iG{Z^p/Z^{p+1}}\cCb \underset{\rm[a]}\iso 
\iG{x}H^p\iG{Z^p/Z^{p+1}}\cCb=\iG{x}\C^p=\cCb(x)\\[-2.6pt]
$$
is an isomorphism of complexes.\vspace{1pt}
\item 
The punctual isomorphism in\/ {\rm(ii)} factors as
 $$
(E\cCb)(x) \underset{\eqref{eq:coz1}}\iso H^p_x
\cCb \underset{\rm(i)[\textup{b}]}\iso \cCb(x);
$$
and hence the following diagram, with\/ $\cFb\<\<,$ $x$ and\/ $p$ as in\/ \emph{(i)[c],} commutes$\>:$ 
\[
\begin{CD}
(EE\cFb)(x) @>{\emph{(ii)}}>{\cCb\;=\;E\cFb}> (E\cFb)(x) \\[-3pt]
@V\eqref{eq:coz1}V\simeq V @V\simeq V\eqref{eq:coz1}V \\
H^p_xE\cFb @>>{\emph{(i)[c]}}> H^p_x\cFb
\end{CD}
\]
\end{enumerate}
\end{alem}
\begin{proof}
For any $Z^{\prime} \subset Z\subset X$, flasque sheaves are acyclic 
for the functors~$\iG{Z}\>$ and~$\iG{Z/Z'}\>$ (\cite[p.~37]{Suo}).
Since a Cousin complex consists of flasque sheaves, (i)  
results.\vspace{1pt} 

For (ii), we note first that in the above definition of $E\cFb$ we can replace $\cLb$ by any  flasque resolution of~$\cFb\<$;  so when $\cFb$ is a Cousin complex $\cCb$, we can take~$\cLb$ to be $\cCb$ itself.

For any integer~$p$,
there is a natural commutative diagram 
with exact rows: 
\[ 
\begin{CD}
0 @>>> \iG{Z^{p+1}/Z^{p+2}}\cCb @>>> \iG{Z^{p}/Z^{p+2}}\cCb 
@>>> \iG{Z^{p}/Z^{p+1}}\cCb @>>> 0 \\
@. @| @| @| \\
0 @>>> \C^{p+1}[-p-1] @>>> \sigma_{\ge p}\sigma_{\le p+1}\cCb
@>>> \C^{p}[-p] @>>> 0   
\end{CD}
\]
An examination of the 
$p$-th connecting homology homomorphisms derived from the rows (see~\eqref{eq:coz1z}) makes (ii) clear.\vspace{1pt}

With $\cCb = E\cFb$, the second assertion in (iii) amounts to commutativity---given by the first assertion---of 
the square in the following diagram. 
\[
\begin{CD}
(E\cCb)(x) @>{\textup{(ii)}}>>  \cCb(x) \\
@VV{\eqref{eq:coz1}}V @| \\
H^p_x\cCb @>{\textup{(i)[b]}}>> \cCb(x)  @>{\eqref{eq:coz1}}>> H^p_x\cFb
\end{CD}
\]

For the first assertion, we can use $\cCb$ instead of $\cLb$ to describe derived functors (see above). Then the definitions of the maps involved tell us that what has to be verified is commutativity of the following diagram (\ref{lem:coz1a}.1) of natural isomorphisms, where the direct sums are all extended over $y \in Z^p \setminus Z^{p+1}$. 
\begin{figure}[t]
$$
\def\1{\smash{\iG{x}H^p\iG{Z^p/Z^{p+1}}\cCb}}
\def\2{\iG{x}H^p\bigl(\!\oplus i_y\iG{y}\cCb \bigr)}
\def\3{\iG{x}\bigl(\!\oplus i_y H^p\iG{y}\cCb \bigr)}
\def\4{H^p\iG{x}\cCb}
\def\5{\iG{x}H^p\bigl(\C^p[-p] \bigr)}
\def\6{\iG{x}H^p\bigl(\>\bigl(\!\oplus i_y\>\cCb(y)\bigr)[-p]\> \bigr)}
\def\7{\iG{x}H^p\bigl(\!\oplus\bigl(i_y\>\cCb(y)[-p]\bigr) \bigr)}
\def\8{\iG{x}\bigl(\!\oplus i_yH^p\bigl(\cCb(y)[-p]\bigr)\bigr)}
\def\9{H^p\bigl(\cCb(x)[-p]\bigr)}
\def\ten{\iG{x}\C^p}
\def\lvn{\iG{x}\bigl(\!\oplus i_y\>\cCb(y)\bigr)}
\def\twv{\cCb(x)}
\begin{CD}
\underset{\UnderElement{}{\downarrow}{32pt}{}{}}\1
    @>>\under{-3}{\displaystyle\square_1}>\5@>>>\ten\\[-1.6pt]
 @. @VVV @VVV \\[-1.6pt]
@. \6 @>>> \underset{\UnderElement{}{\|}{32pt}{}{}}\lvn\\[-1.6pt]
@. @VVV\\[-1.5pt]
\2 @>>> \7\\
@VVV @VVV \\
\3 @>>> \8 @>>> \lvn \\
@VVV @VVV @VVV \\
\4 @>>> \9 @>>> \twv\\[6pt]
\end{CD}
$$
\centerline{\bf(\ref{lem:coz1a}.1)}\vspace{-6pt}
\end{figure}
But in that diagram, commutativity of subdiagram $\square_1$ follows from the description of the map in~\eqref{eq:coz0}; and commutativity of the remaining subdiagrams is straightforward to check.
\end{proof}

Let $ (\X, \Delta) \in \bbFc$ (see \S\ref{subsec:formal}).
Then the codimension function $\Delta$ determines a 
filtration $Z^{\<\bullet}$ of $\X$ with 
\[  
Z^p = \{x \in \X \:|\: \Delta(x) \ge p\}. 
\]
This filtration satisfies the properties (a)-(d) listed near the beginning of \S\ref{subsec:cousin}.
(Property~(a) is satisfied because $\Delta$ is bounded below
by the least of its values at the generic points of 
the finitely many irreducible components of~$\X$.)  
For this $Z^{\<\bullet}$,
a $Z^{\<\bullet}$-Cousin complex is called 
a $\Delta$-Cousin complex, and the Cousin functor~
$E_{Z^{\<\bullet}}$ is denoted by~$E_{\Delta}$.%
\index{Cousin complex!$\Delta$-Cousin complex}%
\index{ $\Delta$!$\Delta$-Cousin complex}%
\index{ $\Fgt$@$E$ (Cousin functor associated to a
filtration)!$E_\Delta$ (\dots filtration given by codimension
function~$\Delta$)} 
The relation of $E_{\Delta}$ to translation 
is immediate from the definition of 
$E_{Z^{\<\bullet}}$:
for any complex~$\cFb$ and integer $n$,
\begin{equation}\label{eq:Cousin+trans}
E_{\Delta-n}\big(\cFb[n]\big)=\big(E_{\Delta}\cFb\big)[n].
\end{equation}

Our main interest is in 
the full subcategory~$\Coz_{\Delta}(\X\>)\subset\bC(\X\>)$ whose objects are those $\Delta$-Cousin complexes whose underlying modules are in~$\Aqct(\X\>)$.
\index{  $\Coz_{\Delta}$ (Cousin $\Aqct$ complexes category)}

\smallskip
Let $(\X, \Delta) \in \bbFc$, $x\in\X$, $\cFb \in \Dqc^+(\X\>)$,
$\cGb \set E_{\Delta}\R\iGp{\X}\cFb$.
Then the composition 
\[
\cGb\<(x) \xto{\; \eqref{eq:coz1} \;} H^p_x\R\iGp{\X}\cFb
\xrightarrow[\;\text{(see \S\ref{subsec:loc})}\;]{\eqref{eq:nolab3}}
H^p_{m_x}\cFb_x
\]
is a natural isomorphism
\begin{equation}
\label{eq:coz3}
\cGb\<(x) \iso H^p_{m_x}\cFb_x \qquad\bigr(\>p \set \Delta(x)\bigl).
\end{equation}
From \ref{lem:mod6} it follows 
that $\cGb \in\Coz_{\Delta}(\X\>)$.

\enlargethispage*{3pt}

 \begin{alem}
\label{lem:psm2}
Let\/ $f \colon (\X,\Delta_1) \to (\Y,\Delta_2)$ be a flat map 
in\/ $\bbFc,$ $\cMb$ a complex\- in\/~$\Coz_{\Delta_2}\!(\Y)$ 
and\/ $\cPb\!$ a flat quasi-coherent\/ $\OX$-module. \kern-.8pt Let\/ $x \in \X$ 
and\/ \mbox{$q\set \Delta_2(f(x))$.} Then for 
any\/ $i \in \mathbb Z,$ with\/  $\sigma$  the truncation functor 
of\/ $\S\ref{subsec:conv}\:\eqref{conv10},$ and\/ $\otimes \set \otimes_{\OX},$  the natural maps are isomorphisms
\begin{alignat*}{3} 
H^i_{m_x}\<(f^{\ast}\sigma_{\ge q}\>\cMb \otimes \cPb)_x\ 
&\<\iso\< 
\ H^i_{m_x}\<(f^{\ast}\<\cMb \otimes \cPb)_x
&&\<\iso\<
\ H^i_{m_x}\<(f^{\ast}\sigma_{\le q}\>\cMb \otimes \cPb)_x\>,  \\
H^i_{\<x}\R\iGp{\X}(f^{\ast}\sigma_{\ge q}\>\cMb \otimes \cPb)
&\<\iso\< 
H^i_{\<x}\R\iGp{\X}(f^{\ast}\<\cMb \otimes \cPb)
&&\<\iso\<
H^i_{\<x}\R\iGp{\X}(f^{\ast}\sigma_{\le q}\>\cMb \otimes \cPb). 
\end{alignat*}
\end{alem}

\begin{proof}
The first two maps are isomorphic, respectively,  to the last two, see
~\eqref{eq:nolab3}, 
so it suffices to show that the last two are isomorphisms. 
Since $f$ and $\cPb$ are flat, the  natural exact sequences
\[ 
0 \to \sigma_{\ge q}\>\cMb \to \cMb \to \sigma_{\le q-1}\>\cMb \to 0, 
\]
\[ 
0 \to \sigma_{\ge q+1}\>\cMb \to \cMb \to \sigma_{\le q}\>\cMb \to 0. 
\]
induce exact sequences
\[ 
0 \to f^{\ast}\sigma_{\ge q}\>\cMb \otimes \cPb \to 
f^{\ast}\<\cMb \otimes \cPb \to 
f^{\ast}\sigma_{\le q-1}\>\cMb \otimes \cPb \to 0, 
\]
\[ 
0 \to f^{\ast}\sigma_{\ge q+1}\>\cMb \otimes \cPb \to 
f^{\ast}\<\cMb \otimes \cPb \to 
f^{\ast}\sigma_{\le q}\>\cMb \otimes \cPb \to 0. 
\]

\pagebreak[3]

\noindent
Applying the derived functors $\R\iG{x}\R\iGp{\X}$ 
to the triangles in~$\Dqc(\X\>)$ arising from
the last two exact sequences, we get triangles in
the derived category of complexes of~$\>\OXx$-modules.
The resulting long exact cohomology sequences  
show that for the  lemma to hold it suffices  that for all $i,$
\begin{equation} 
\label{eq:psm3}
H^i_{\<x}\R\iGp{\X}(f^{\ast}\sigma_{\le q-1}\>\cMb \otimes \cPb) 
\;=\; 0 \;=\; 
H^i_{\<x}\R\iGp{\X}(f^{\ast}\sigma_{\ge q+1}\>\cMb \otimes \cPb). 
\end{equation}

But  (\ref{eq:psm3}) follows by induction from \ref{lem:psm1}.
For example, in the case of~$\sigma_{\ge q+1}\cMb$ 
we have for any $n \ge 1$ an exact sequence 
\begin{equation} 
\label{eq:psm4}
0 \to \sigma_{\ge q+n+1}\>\cMb \to \sigma_{\ge q+n}\>\cMb \to \M^{q+n}[-q-n] \to 0. 
\end{equation}
Since $\M^{q+n} = \oplus_y \; i_y\cMb(y)$ where $y \ne f(x)$ (since $\Delta_2(y) \ne q$), 
\ref{lem:psm1} yields 
\[
H^i_{\<x}\R\iGp{\X}(f^{\ast}\<\M^{q+n}[-q-n] \otimes \cPb) = 0 
\quad \quad \forall \; i.
\]
Therefore from the long exact sequence corresponding to (\ref{eq:psm4})
we conclude that
\[ 
H^i_{\<x}\R\iGp{\X}(f^{\ast}\sigma_{\ge q+n+1}\>\cMb \otimes \cPb) \cong  
H^i_{\<x}\R\iGp{\X}(f^{\ast}\sigma_{\ge q+n}\>\cMb \otimes \cPb). 
\]
Hence by induction on $n$ we see that for all $n \ge 1$ and for all $\: i$, 
\[ 
H^i_{\<x}\R\iGp{\X}(f^{\ast}\sigma_{\ge q+1}\>\cMb \otimes \cPb) \cong  
H^i_{\<x}\R\iGp{\X}(f^{\ast}\sigma_{\ge q+n}\>\cMb \otimes \cPb). 
\]
But for a fixed $i$ and for $n$ sufficiently large, 
the right hand side is 0,
and thus 
\[
H^i_{\<x}\R\iGp{\X}(f^{\ast}\sigma_{\ge q+1}\>\cMb \otimes \cPb) = 0 
\quad \quad \forall \; i.
\]
A similar argument works for $\sigma_{\le q-1}\>\cMb$. 
\end{proof}

Applying \ref{lem:psm2} with $\sigma_{\le q}\cMb$ in place of $\cMb$, and using \ref{lem:psm1}, we conclude:

\begin{acor}
\label{cor:psm3}
In the situation of \textup{\ref{lem:psm2}} there are natural isomorphisms
\begin{small}
$$
\minCDarrowwidth=.23in
\begin{CD}
H^i_{m_x}\<(f^{\ast}\<\cMb \otimes \cPb)_x @>\Iso>>
H^i_{m_x}\<(f^{\ast}\<\M^q[-q] \otimes \cPb)_x @>\Iso>>
H^{i-q}_{m_x}\<(\cMb(y) \otimes_{\OYy} \cPb_x)\\
@V(6)V\simeq V @V(6)V\simeq V @V\simeq V(6)V \\
H^i_{\<x}\R\iGp{\X}(f^{\ast}\<\cMb \otimes \cPb)@>\Iso>>
H^i_{\<x}\R\iGp{\X}(f^{\ast}\<\M^q[-q] \otimes \cPb)\<\< @>\Iso>>
\<H^{i-q}_{\<x}\R\iGp{\X}(f^{\ast}\<i_y\cMb(y\<)\< \otimes \cPb\<)
\end{CD}
$$
\end{small}
\end{acor}

\smallskip
\subsection{Cohen-Macaulay complexes}
\label{subsec:cm}

Let $(\X,\Delta) \in \bbFc$ and let $Z^{\bullet}$ denote the
filtration of~$\X$ induced by $\Delta$, i.e.,
$Z^p = \{x \in \X \:|\: \Delta(x) \ge p\}$.
We say that a complex~$\cFb \in \D^+(\X\>)$ is
{\em Cohen-Macaulay with respect to~$\Delta$} (in short $\Delta$-CM, or
simply CM)
\index{Cohen-Macaulay complex}%
\index{CM|see{Cohen-Macaulay}}%
\index{CM!$\Delta$-CM (Cohen-Macaulay w.r.t. codimension function~$\Delta$)}%
\index{ $\Delta$!$\Delta$-CM (Cohen-Macaulay w.r.t. codimension function~$\Delta$)}%
if it satisfies the following equivalent conditions
(cf.~\cite[p.~242]{RD}, (i)$\:\Leftrightarrow\:$(ii), where
the boundedness assumption on~$\cFb$ is not used):
\begin{enumerate}
\item
\begin{enumerate}
\item For all integers $i < p$, $H^i_{\<\<Z^{\<p}}\cFb = 0$;
\item For integers $i \ge p$ the canonical map
$H^i_{\<\<Z^p}\cFb \to H^i\cFb$ is surjective when
$i=p$ and bijective for $i>p$ (equivalently,
$H^i_{\<X/Z^p}\cFb = 0$ for $i \ge p$);
\end{enumerate}
\item For any
$x \in \X$, $H^i_{\<x}\cFb = 0$ for $i \ne \Delta(x)$
(equivalently, for any integers $i \ne p$,
$H^i_{\<\<Z^p/Z^{p+1}}\cFb = 0$).
\end{enumerate}

Let $\D^+(\X; \Delta)_{\text{CM}}$ be the full subcategory of $\>\D^+(\X\>)$%
\index{ $\D(\A)$ (derived category
of $\A$-complexes)!$\D(X)\set\D(\A(X))$!Z@$\D^+(\X; \Delta)_{\text{CM}}$ ($\Delta$-CM complexes subcat.\  of~$\D^+(\X))$}
whose objects are $\Delta$-CM complexes.
We now recall the relation between $\Delta$-CM complexes
and $\Delta$-Cousin complexes. For the next few results
we set aside considerations of quasi-coherence and torsionness.
For convenience, we suppress the reference to the codimension
function~$\Delta$ in our notation.
Thus $E = E_{\Delta}$ and $\D^+(\X)_{\text{CM}}= \D^+(\X; \Delta)_{\text{CM}}$.

\pagebreak[3]

Let $\Cou(\X\>)\subset\bC(\X\>)$ be the full subcategory of Cousin complexes.
\index{  $\Cou(\X\>)$ (Cousin complex category on $\X\>$)}
A basic property of $\Cou(\X\>)$ is that if two maps $h,k$ from 
$\cFb$ to $\cGb$ in $\Cou(\X\>)$ are 
homotopy equivalent, then $h=k$. (More generally, for integers $p>q$
any map from $\F^p$ to~$\G^q$ is zero since for 
any $x,y \in \X$ and abelian groups $F,G$, 
we have $\Hom(i_x F, i_y G) = 0$ 
if $y \notin \ov{\{x\}}$.) So $\Cou(\X\>)$ can 
be considered as a full subcategory of~$\K^+(\X\>)$.
Then \emph{the~localization functor\/ $Q\colon\K(\X\>) \to \D(\X\>)$  
takes $\Cou(\X\>)$ into} $\D^+(\X\>)_{\text{CM}}$.
Indeed, by~(i) of \ref{lem:coz1a}, for any $\cCb \in \Cou(\X\>)$ 
and $x \in \X$, with $p \set \Delta(x)$ there is an 
isomorphism $H^j_x\cCb \cong H^j(\cCb(x)[-p])$ 
and therefore $H^j_x\cCb = 0$ if $j \ne p$. 
Thus $Q$~induces an additive functor 
\[
Q \colon \Cou(\X\>) \to \D^+(\X\>)_{\text{CM}}.
\]
Let us elaborate. With
\[
\>\>\ov{\<\<\bC\<\<}\>\> \set \Cou(\X\>), \qquad \>\ov{\<\D\!}\, \set \D^+(\X\>)_{\text{CM}},
\]
we shall assume for the rest of this subsection
that $Q$ (resp.~$E$) is a functor from $\>\>\ov{\<\<\bC\<\<}\>\>$ to $\>\ov{\<\D\!}\,$
(resp.~$\>\ov{\<\D\!}\,$ to $\>\>\ov{\<\<\bC\<\<}\>\>$).

\begin{aprop}
\label{prop:coz4}
With preceding notation,
the functor\/ $Q \colon \>\>\ov{\<\<\bC\<\<}\>\> \to \>\ov{\<\D\!}\,$ is an 
\emph{equivalence of categories}
having\/ $E$ as a pseudo-inverse equivalence.
In particular, $E$~is fully faithful. 
\end{aprop}
\begin{proof}
That $Q$ is an equivalence
of categories is one of the main results (Theorem\-\ 3.9) in \cite{Suo}. 
Let $S \colon \>\ov{\<\D\!}\, \to \>\>\ov{\<\<\bC\<\<}\>\>$ be a pseudo-inverse 
equivalence. Since  $EQ \iso 1_{\>\>\ov{\<\<\bC\<\<}\>\>}$ (\ref{lem:coz1a}(ii)),  
there are isomorphisms $E \iso EQS) \iso S$. 
Thus $E$ is a pseudo-inverse of~$Q$.
\end{proof}

\begin{acor}
\label{cor:coz4a}
Let\/ $EQ \iso 1_{\>\>\ov{\<\<\bC\<\<}\>\>}$ be the isomorphism of\/
\textup{\ref{lem:coz1a}(ii).} Then there exists a unique 
isomorphism of functors\/ 
$\Ssf=\Ssf_{\X,\Delta}\colon 1_{\>\ov{\<\D\!}\,} \iso QE$ 
such that the following two induced isomorphisms are
equal$\>:$
\[
E \: = E\/1_{\>\ov{\<\D\!}\,} \:\iso \:EQE, \qquad \qquad  
E \: = 1_{\>\>\ov{\<\<\bC\<\<}\>\>}\/E \: \iso \: EQE.
\]
Furthermore, the following induced isomorphisms are equal$\>:$
\[
QEQ\iso Q1_{\>\>\ov{\<\<\bC\<\<}\>\>} = \:Q, \qquad \qquad  
QEQ\iso 1_{\>\ov{\<\D\!}\,}Q=\:Q.
\]

\end{acor}
\begin{proof}
Since $E$ is fully faithful, for any $\cFb \in \>\ov{\<\D\!}\,$ the map 
$E\cFb \iso EQE\cFb$ given by $E = 1_{\>\>\ov{\<\<\bC\<\<}\>\>}E \iso EQE$
comes via $E$ from a unique functorial isomorphism $\cFb \iso QE\cFb\>$
that  fulfills the first assertion. The second assertion need only be shown after
application of the functor~$E$, at which point the first assertion reduces the problem
to showing that the following isomorphisms are equal
\[
EQEQ\iso EQ1_{\>\>\ov{\<\<\bC\<\<}\>\>} = \:EQ, \qquad \qquad  
EQEQ\iso 1_{\>\>\ov{\<\<\bC\<\<}\>\>}EQ=\:EQ,
\]
which is easy to do after composition with the isomorphism 
$EQ\iso 1_{\>\>\ov{\<\<\bC\<\<}\>\>}.$
\end{proof}

\begin{acor}
\label{cor:coz4b}
For any\/ $\cFb \in \>\ov{\<\D\!}\,$ the isomorphism 
$\Ssf(\cFb) \colon \cFb \iso QE\cFb$ obtained in~\textup{\ref{cor:coz4a}} 
is the unique one satisfying the property that for any\/ $x \in \X,$ 
with\/ $p\set\Delta(x),$ the induced isomorphism\/ 
$H^p_x\,\Ssf(\cFb) \colon H^p_x\cFb \iso H^p_xQE\cFb$ 
is the inverse of the isomorphism in\/ \textup{\ref{lem:coz1a}(i)[c]}.
\end{acor}
\begin{proof}
Let $\phi \colon \cFb \iso QE\cFb$ be an isomorphism.
Consider the following diagram, where the vertical isomorphisms are 
given by the punctual decomposition of~$E(-)$ in~{\eqref{eq:coz1}}: 
\[
\begin{CD}
(E\cFb)(x) @>{(E\phi)(x)}>> 
(EQE\cFb)(x) @>{EQ\;\overset{\under{.4}{\sim}}=
\;1_{\>\>\ov{\<\<\bC\<\<}\>\>}}>> (E\cFb)(x) \\
@V\overset{\under{.4}{\sim}}= VV @VVV  @VV\overset{\under{.4}{\sim}}=V \\
H^p_x\cFb @>>{H^p_x\phi}>
H^p_xQE\cFb @>>{\ref{lem:coz1a}\kern.5pt\textup{(i)[c]}}>  H^p_x\cFb
\end{CD}
\]
The rectangle on the left commutes for functorial reasons,
while the one on the right commutes by \ref{lem:coz1a}\kern.5pt(iii).
If $\phi = \Ssf(\cFb)$, then by \ref{cor:coz4a}, the top row composes 
to an identity map, and hence so does the bottom. Conversely, if the 
bottom row composes to an identity map, then 
so does the top. Since this holds for every $x \in \X$, 
therefore the composite map
$E\cFb \xto{E\phi} EQE\cFb 
\xrightarrow[\textup{\ref{lem:coz1a}(ii)}]{EQ \iso 1_{\>\>\ov{\<\<\bC\<\<}\>\>}} E\cFb$
is identity and hence $\phi = \Ssf(\cFb)$.
\end{proof}

Let $x \in \X$ and $p \set \Delta(x)$. Any $\OX$-complex~$\cCb\<$,
has the subcomplex $\iG{\,\ov{\!\{x\}\!}\,}\cCb$, and for the stalks at~$x$ there is the inclusion map $\gamma=\gamma(x,\cCb)\colon\iG{x}\cCb\into\cCb_x$.\vspace{1pt}
Since ``stalk at $x$" is an exact functor, $Q\gamma$ factors naturally as
$$
\iG{x}\cCb\xto{\ } \R\iG{x}\cCb\xto{\>\bar\gamma\>}\cCb_x
$$
with $\bar\gamma=\bar\gamma(x,\cCb)$  a $\D(\OXx)$-morphism.\vspace{1pt}
So for $\cFb\in\D^+(\X\>)$ there is a natural map
$$
 \eta^{}_1= \eta^{}_1(x,\cFb)\colon\bigl(H^p_x\cFb\bigr)[-p] \underset{\eqref{eq:coz1}}\iso\bigl(E\cFb\bigr)(x)[-p]=
\iG{x} E\cFb\xto{\>\gamma\>}\bigl(E\cFb\bigr)_{\<x}\>.
$$

If $\cFb\in\>\ov{\<\D\!}\,$ then by definition, $H^i_x\cFb=0$ for $i\ne p$. Hence the natural maps (see \S\ref{subsec:conv},\,\eqref{conv10}) are isomorphisms
$$
\R\iG{x}\cFb\iso\tau_{\ge p}\R\iG{x}\cFb\osi\tau_{\le p}\tau_{\ge p}\R\iG{x}\cFb=\bigl(H^p_x\cFb\bigr)[-p].
$$
Let $ \eta^{}_2= \eta^{}_2(x,\cFb)$ be the resulting isomorphism $\R\iG{x}\cFb\iso\bigl(H^p_x\cFb\bigr)[-p]$.

\begin{acor}
\label{cor:coz4c}
For any\/ $\cFb \in  \>\ov{\<\D\!}\,,$ set 
$\cEb\set QE\cFb$ and let\/ $\nu\colon\cFb\to\cEb$ be the isomorphism given by\/
\textup{\ref{cor:coz4a}}.
Then with notation as above the following diagram in\/~$\D(\OXx)$commutes$\>:$
\[
\begin{CD}
\cFb_x @>\Iso>\under{1.2}{\nu_x}> \cEb_x \\[-2pt]
@A\bar\gamma AA @AA{ \eta^{}_1}A \\
\R\iG{x}\cFb @>\Iso>\under{1.2}{ \eta^{}_2}> (H^p_x\cFb)[-p]
\end{CD}
\]
\end{acor}

\begin{proof}
Expand the diagram thus:
$$
\def\1{\cFb_x}
\def\2{\cEb_x}
\def\3{\bigl(E\cFb\bigr)_x}
\def\4{\R\iG x\cFb}
\def\5{\bigl(H^p_x\cFb\bigr)[-p]}
\def\6{\bigl(E\cFb\bigr)(x)[-p]}
\def\7{\iG x \bigl(E\cFb\bigr)}
\def\8{\R\iG x\cEb}
\begin{CD}
\1 @>\Iso>\under{1.2}{\nu_x^{}}>\2 @.\Equals{13em}@. \3 \\
@A\bar\gamma A\mspace{60mu}\displaystyle\square_1A 
    @AA\eta_1^{}A @. 
        @A\displaystyle\square_2\mspace{140mu}A\gamma_x^{} A \\
\4 @>\under{-5}
    {\hbox to 0pt{$\mspace{166mu}\displaystyle\square_3$\hss}}      >\under{1.2}{\eta_2^{}}>\5 @>\Iso>\eqref{eq:coz1}>\6 @= \7 \\[-3pt]
@| @. @. @A\under{1.2}{\rm\ref{lem:coz1a}}A\simeq A\\
\4 @. @.\mspace{-168mu}
\overset
{\under{.5}{\widetilde{\phantom{mmm}}}}
   {\underset
         {\R\iG{x} \nu}{\Rarrow{22.6em}}
    }
 @.\8
\end{CD}
$$
By definition, the rightmost column composes to $\bar\gamma(x,\cEb)$, and so the outer border of this expanded diagram commutes.

Subdiagram $\square_2$ commutes by definition.

According to \ref{cor:coz4b}, application of the homology functor $H^p$
to subdiagram~$\square_3$ produces a commutative diagram. But that means  $\square_3$ itself commutes, because its vertices are complexes which have vanishing homology in all degrees except $p$, and $H^p$ is an equivalence from the full subcategory of such complexes in $\D(\OXx)$ to the category of $\OXx$-modules (with pseudo-inverse $F\mapsto F[-p]$).

It results then that $\square_1$ commutes, as asserted.
\end{proof}

\newpage

\section{Generalized fractions and pseudofunctors }
\label{sec:pseudo}

Below we recall the definition of a pseudofunctor.\index{pseudofunctor (2-functor)}
In \ref{thm:huang1} we give a modified version of 
Huang's pseudofunctor\index{Huang, I-Chiau!pseudofunctor} constructed in~\cite{Hu}.
Various components of this pseudofunctor have an explicit 
description in terms of generalized fractions. Generalized 
fractions are useful in denoting elements of certain 
local cohomology modules. In \S\ref{subsec:gfrac} we 
review the definition of generalized fractions and state 
the relation between two different competing definitions.
In \S\ref{subsec:foris} we review some isomorphisms
that are needed to describe the pseudofunctor 
discussed in~\ref{thm:huang1}. In \S\ref{subsec:proofstwo}
we prove a somewhat technical result which is mainly
used in proving commutativity of the diagram 
in~\eqref{cd:conc0} in~\S\ref{subsec:conc}.

A \emph{contravariant pseudofunctor} `$\,{\#}\,$' (or $(-)^{\#}$)%
\index{pseudofunctor (2-functor)!contravariant}
on a category~$\mathsf{C}$ 
assigns  to each $\mathsf{C}$-object~$X$ a category~$X^{\#}\<$, 
to  each 
$\mathsf{C}$-morphism $f\colon X \to Y$ a contravariant\vspace{-.5pt} functor 
$f^{\#} \colon Y^{\#} \to X^{\#}\<$,
to each $\mathsf{C}$-diagram 
$X \xto{\; f \;} Y \xto{\; g \;} Z$ a functorial
\emph{c\-omparison isomorphism}\vspace{-1pt}
$C^{\#}_{f\<,\>g} \colon f^{\#}g^{\#} \iso (gf)^{\#}\<$, %
\index{ $\Cfr$01@$C^{\#}_{(-),(-)}$ (comparison isomorphism for
generic pseudofunctor ${\boldsymbol{{-}^{\sharp}}}$)}
and to each $\mathsf{C}$-object $Z$  a functorial 
\emph{u\-nit isomorphism}%
\index{ $\Cfr$01@$\delta^{\#}_{(-)}$ (unit isomorphism for
generic pseudofunctor ${\boldsymbol{{-}^{\sharp}}}$)}
$\delta^{\#}_Z \colon (1_Z)^{\#} \iso 1_{Z^{\#}}$,
all subject to the following conditions:
\begin{enumerate}
\renewcommand{\labelenumi}{\roman{enumi})}
\item\label{pf0i5} For every triple of morphisms 
$X \xto{\; f \;} Y \xto{\; g \;} Z \xto{\; h \;} W$ in $\mathsf{C}$
the following associativity diagram commutes.
\[
\begin{CD}
f^{\#}g^{\#}h^{\#}  @>{f^{\#}C^{\#}_{g,h}}>>   f^{\#}(hg)^{\#} \\
@V{C^{\#}_{f\<,g}h^{\#}}VV    @VV{C^{\#}_{f\<,hg}}V \\
(gf)^{\#}h^{\#}  @>>{C^{\#}_{gf\<,h}}>  (hgf)^{\#} 
\end{CD}
\]
\item\label{pf0i4} For any map $f \colon X \to Y$ in $\mathsf{C}$, 
the following two compositions are identity 
\[
(1_X)^{\#}f^{\#} \; \xto{\; C^{\#}_{1_X\<,f} \;} \; f^{\#} 
\; \xto{\; {\delta^{\#}_X}^{-1}f^{\#} \;} \;(1_X)^{\#}f^{\#},
\]
\[
f^{\#}(1_Y)^{\#} \; \xto{\; C^{\#}_{f\<,1_Y} \;} \; f^{\#} 
\; \xto{\; f^{\#}{\delta^{\#}_Y}^{-1} \;} \;f^{\#}(1_Y)^{\#}.
\] 
\end{enumerate}
If necessary, we also use the cumbersome notation of a quadruple 
\[
\Bigl( \;(-)^{\#}, \;(-)^{\#}, \;C^{\#}_{(-),(-)}, \;\delta^{\#}_{(-)}  
\;\Bigr)
\]
to denote the pseudofunctor `${\#}$' where the first entry operates
on the objects of~$\mathsf{C}$, the second on morphisms, 
the third on pairs of composable morphisms and the 
fourth on objects of~$\mathsf{C}$, all having
the obvious meanings as per the definition above. 

\smallskip

A {\it covariant pseudofunctor} is defined by associating%
\index{pseudofunctor (2-functor)!covariant} 
to each morphism in~$\mathsf{C}$ a covariant functor which 
is required to satisfy appropriately modified compatibility 
conditions. We will use the subscript notation 
$(-)_{\#}$ when dealing with covariant pseudofunctors 
(cf.~ \cite[chapter IV]{Hu}).

\medskip
\emph{Remark.} The definition is simpler when
for all~$Z$, $(1_Z)^\#=1_{Z^\#}$ and
$\delta^{\#}_Z$ is the identity map. A pseudofunctor satisfying these
additional conditions is said to be \emph{normalized.} Replacing the functor 
$(1_Z)^\#$ by $1_{Z^\#}$ for every $Z$ and replacing $C_{f\<,\>g}^\#$
by $C_{gf}^\#$ whenever $f$ or $g$ is an identity map  transforms any pseudofunctor $(-)^\#$ into an isomorphic (in the obvious sense) normalized one.\vspace{1pt}

A normalized contravariant pseudofunctor  on $\mathsf{C}$%
\index{pseudofunctor (2-functor)!normalized} is the same thing as a contravariant
2-functor from $\mathsf{C}$ to the 2-category of all categories.

\medbreak

\subsection{Generalized fractions}
\label{subsec:gfrac}\index{generalized fraction}

Let $A$ be a noetherian ring and $M$ an $A$-module. 
Let $\bt = (t_1, \ldots, t_n)$ be a 
sequence in $A$ and let $\afr$ denote an ideal in $A$ 
such that $\sqrt{\afr} = \sqrt{\bt A}$.
Then the elements of the local cohomology 
module~$H^n_{\afr}M$ can be represented in the form of generalized
fractions of the type 
\[
\genfrac{[}{]}{0pt}{}{m}{t_1^{b_1},\ldots, t_n^{b_n}},
\qquad m \in M, \qquad b_i > 0,
\]
however there are two natural ways of doing so which we now 
review. 
We state the relationship between these two ways in 
Lemma~\ref{lem:appsmo5}.

The first method involves using the \v Cech complex. 
Let $X = \Spec(A)$, let~$Z$ denote the closed subset 
of~$X$ defined by the ideal~$\afr$ 
and let $U$ denote the complement $X \setminus Z$. Let 
$\Ufr = (U_1, \ldots, U_n)$
denote the affine open covering of~$U$ associated to the 
sequence~$\bt$ where $U_i = \Spec(A_{t_i})$. 
Let $\Ccb(\Ufr,M)$ denote the \v Cech
complex of $A$-modules for the 
sheaf $M^{\sim}\big|_U$ (see \ref{prop:mod1}), 
corresponding to the cover $\Ufr$ of~$U$. 
Let $M \to \Ib$ denote an 
injective resolution of $M$. 
Then $M^{\sim} \to \Ib{}^{\sim}$ is a flasque resolution,
in fact a resolution by $\Aqc(X)$-injectives. We may therefore
make the identification $H^*(U, M^{\sim}) = H^*\Gamma(U, \Ib{}^{\sim})$. 
The complex $\Ccb(\Ufr,M)$ exists 
between degrees~0 and~$n-1$ and we have a natural map 
\begin{equation}
\label{eq:appsmo3}
M_{t_1 \cdots t_n} = \check{C}^{n-1}(\Ufr,M) 
\onto H^{n-1}\Ccb(\Ufr,M) 
\xto{\;\alpha\;} H^{n-1}(U,M^{\sim}) \xto{\;\eta\;}  H^n_{\afr}M.
\end{equation}
where $\alpha$ is the standard isomorphism relating the \v Cech
cohomology to the usual one (see \cite[III, Lemma~4.4]{Ha}) 
and~$\eta$ is the usual connecting homomorphism in the homology long exact sequence associated to the exact sequence
\[
0 \to \iG{\afr}\Ib \to  \Ib \to \iG{}(U,\Ibsim) \to 0.
\]
The map $\eta$ is 
surjective for $n \ge 1$ and an isomorphism for
$n>1$. Therefore any element $\theta \in H^n_{\afr}M$ is the image of 
an element in $M_{t_1 \cdots t_n}$. This leads to one way of representing
a cohomology class by a generalized fraction, viz., for 
$\theta \in H^n_{\afr}M$, for~$m\in M$ and for positive integers~$b_i$ 
we say that 
\[
\theta = \genfrac{[}{]}{0pt}{}{m}{t_1^{b_1},\ldots, t_n^{b_n}}_{\text{C}}
\] 
(the subscript C stands for \v Cech) if the element 
\[
\frac{m}{t_1^{b_1} \cdots t_n^{b_n}} \in M_{t_1 \cdots t_n}
\]
goes to $\theta$ under the natural map of (\ref{eq:appsmo3}).
We call such a generalized 
fraction a \emph{C-fraction} representing $\theta$.
\index{generalized fraction!C-fraction}%

The second method involves the direct-limit Koszul complexes.
Let $\Kb_{\infty}(\bt)$ denote the direct-limit 
Koszul complex on $\bt$ over $A$. With $\Ib$ denoting an injective 
resolution of $M$, there are natural
quasi-isomorphisms (with $\otimes = \otimes_A$) 
\[
\Kb_{\infty}(\bt) \otimes M \to \Kb_{\infty}(\bt) \otimes \Ib 
\gets \iG{\bt A}\Ib = \iG{\afr}\Ib. \qquad (\text{cf.~\ref{prop:loc5}})
\]
Therefore we have a natural map 
\begin{equation}
\label{eq:appsmo4}
M_{t_1 \cdots t_n} = K_{\infty}^n(\bt) \otimes M \onto 
H^n(\Kb_{\infty}(\bt) \otimes M) \cong H^n_{\afr}M.
\end{equation}
We represent any $\theta \in H^n_{\afr}M$ 
by a generalized fraction  
\[
\theta = \genfrac{[}{]}{0pt}{}{m}{t_1^{b_1},\ldots, t_n^{b_n}}_{\text{K}}
\] 
(the subscript K being for Koszul) if the element 
\[
\frac{m}{t_1^{b_1} \cdots t_n^{b_n}} \in M_{t_1 \cdots t_n}
\]
goes to $\theta$ under the natural map of (\ref{eq:appsmo4}). 
We  call such a generalized fraction 
a \emph{K-fraction} representing $\theta$.\index{generalized fraction!K-fraction}%
\index{K-fraction|see{generalized fraction}}

\medskip
\begin{alem}
\label{lem:appsmo5}
With $A,M,\bt,n$ as above, for any $m \in M$ and integers $b_i >0$ 
with~$i$ between 1 and $n,$ the following relation holds.
\[
\genfrac{[}{]}{0pt}{}{m}{t_1^{b_1},\ldots, t_n^{b_n}}_{\text{K}} 
= (-1)^n
\genfrac{[}{]}{0pt}{}{m}{t_1^{b_1},\ldots, t_n^{b_n}}_{\text{C}}
\]
\end{alem}
\begin{proof}
Before we begin with the proof we will take a short digression 
regarding \v Cech resolutions.
Let $X = \Spec(A)$, $U$, $\Ufr$ be as above. 
Let $\cO_U \to \eCb$ denote the corresponding 
\v Cech resolution of $\cO_U$ on~$U$. 
For an $A$-module $N$, set $N_U^{\sim} \set N^{\sim}\big|_U$.
We now recall some basic facts concerning~$\eCb$.

The complex $\eCb$ consists of flat $\cO_U$-modules.
For any $A$-module $N$, the natural induced map
$N^{\sim}_U \to \eCb \otimes_U N^{\sim}_U$
can be canonically identified with the corresponding \v Cech resolution
and if $N$ is injective, so that $N^{\sim}_U$ consists of 
$\Aqct(U)$-injectives
then $\eCb \otimes_U N^{\sim}_U$
also consists of $\Aqct(U)$-injectives.
Set $\Cb = \Gamma(U, \eCb)$. Then $\Cb$ consists of flat
$A$-modules and there is a natural isomorphism 
\[
\Cb \otimes_A N \iso \Gamma(U, \eCb \otimes_U N^{\sim}_U) 
= \Ccb(\Ufr, N).
\]
It follows that 
if $M \to \Ib$ is an injective resolution then 
in the following diagram having obvious natural maps,
\begin{equation}
\label{eq:appsmo5a}
\begin{CD}
M^{\sim}_U @>>> \eCb \otimes_U M^{\sim}_U \\
@VVV @VV{\beta_1}V \\
\Ibsim_U @>{\beta_2}>> \eCb \otimes_U \Ibsim_U 
\end{CD}
\end{equation}
all the maps are quasi-isomorphisms and $\beta_2$ 
has a homotopy inverse. 
Moreover, there are isomorphisms
\begin{align}
H^i\Ccb(\Ufr, M) = H^i\Gamma(U, \eCb \otimes_U M^{\sim}_U)
 &\xrightarrow[\text{via $\beta_1$}\;\;]{\sim} 
 H^i\Gamma(U, \eCb \otimes_U \Ibsim_U) \notag \\
&\xrightarrow[\text{via $\beta_2^{-1}$}]{\sim} 
 H^i\Gamma(U, \Ibsim_U) = H^i(U, M^{\sim}_U). \notag
\end{align}
It follows that the map $\alpha$ 
in~\eqref{eq:appsmo3} is obtained from this process.
It then follows that the map 
$M_{t_1 \cdots t_n} \onto H^n_{\afr}M$ of \eqref{eq:appsmo3} 
is also obtained from the sequence
\begin{align}
M_{t_1 \cdots t_n} = C^{n-1} \otimes_A M 
 \onto H^{n-1}(\Cb \otimes_A M) 
 &\xrightarrow[\text{via $\beta_1$}\;\;]{\sim}
 H^{n-1}(\Cb \otimes_A \Ib) \notag \\
&\xrightarrow[\text{via $\beta_2^{-1}$}]{\sim} 
 H^{n-1}(\Ibsim(U)) \xto{\;\eta\;} H^{n}\iG{\afr}\Ib. \notag
\end{align}
Finally, we recall that $\Cb$ is a displaced (and truncated) version
of the complex $\Kb \set \Kb_{\infty}(\bt)$. More 
precisely, for all $p$ except $p = -1$, we have $C^p = K^{p+1}$
and $d^p_{\Cb} = d^{p+1}_{\Kb}$ where $d$ stands for the
corresponding differentials. 
In particular, the graded 
maps $C^p \to K^{p+1}$, defined as 
$(-1)^{p+1}$ times the identity map, form a map of complexes 
$\psi \colon \Cb \to \Kb[1]$.

Returning to the proof of \ref{lem:appsmo5}
consider the following diagram with maps described below (we use
$\otimes = \otimes_A$ for the rest of the proof). 
\begin{equation}
\label{cd:appsmo6}
\begin{CD}
H^{n-1}(\Cb \otimes M) @>{\mu_1}>> H^{n-1}(\Kb[1] \otimes M)
  @>{\mu_2}>> H^{n}(\Kb \otimes M) \\
@VV{\mu_3}V @VV{\mu_4}V @VV{\mu_5}V \\
H^{n-1}(\Cb \otimes \Ib) @>{\mu_6}>> H^{n-1}(\Kb[1] \otimes \Ib)
  @>{\mu_7}>> H^{n}(\Kb \otimes \Ib) \\
@AA{\mu_8}A @. @AA{\mu_9}A \\
H^{n-1}(\Ibsim(U)) @>>{\mu_{10}}> H^{n}\iG{\afr}\Ib @= H^{n}\iG{\afr}\Ib
\end{CD}
\end{equation}
The maps $\mu_1,\mu_6$ are induced by $\psi \otimes 1$.
while $\mu_2,\mu_7$ are induced by the natural isomorphism
$\Kb[1] \otimes - \iso (\Kb \otimes -)[1]$ obtained using the
convention in~\S\ref{subsec:conv},~\eqref{conv4}. It follows that 
the underlying maps of complexes of $\mu_2,\mu_7$ are 
identity maps at the graded level.
The maps $\mu_3,\mu_4,\mu_5$ are the obvious 
natural ones. The map~$\mu_8$ is the isomorphism 
induced by $\beta_2$ in \eqref{eq:appsmo5a}.
In particular, the underlying map of complexes of $\mu_8$
is given, at the graded level, by the restriction maps
\[ 
I^{p \sim}(U) \;\xto{\quad}\; \oplus_i \: I^{p}_{t_i} 
\; = \; C^0 \otimes I^{p} \subset (\Cb \otimes \Ib)^p.
\]
The map $\mu_9$ is the usual isomorphism (see \ref{prop:loc5}) 
while $\mu_{10}$ is the connecting homomorphism referred to in the
definition of~$\eta$ in~(\ref{eq:appsmo3}).

It suffices to show that the diagram in (\ref{cd:appsmo6}) commutes. 
Assuming that \eqref{cd:appsmo6} commutes, one proves the 
Lemma as follows. The composition 
$\mu_{10}\mu_8^{-1}\mu_3$ defines the map used in~(\ref{eq:appsmo3}) 
while the composition 
$\mu_9^{-1}\mu_5$ defines the map used in (\ref{eq:appsmo4}). 
From the definitions involved we see that the map $\mu_2\mu_1$ sends the 
cohomology class of the cocycle 
\mbox{$\frac{m}{t_1^{b_1} \cdots t_n^{b_n}} \in C^{n-1} \otimes M$} 
to the cohomology class of 
$(-1)^n \frac{m}{t_1^{b_1} \cdots t_n^{b_n}} \in K^{n} \otimes M$,
thereby proving the Lemma.

In (\ref{cd:appsmo6}) the top two rectangles commute due to 
functorial reasons.
To prove commutativity of the bottom part, consider an element 
$\zeta \in H^{n-1}(\Ibsim(U))$. Since the natural map 
$I^{n-1} \to I^{n-1 \sim}(U)$
is surjective (as $I^{n-1 \sim}$ is flasque) there exists an 
element, say~$z$, in $I^{n-1}$ whose restriction 
to $I^{n-1 \sim}(U)$ is a cocycle mapping
to~$\zeta$. Let~$\delta^{\bullet}$ denote the 
differential of the complex $\Ib$.
Then $\delta^{n-1}z$ maps to zero in~$I^{n \sim}(U)$ and 
so lies in~$\iG{\afr}I^n$ and is an $n$-cocycle 
in~$\iG{\afr}I^n$. From the definition 
of~$\mu_{10}$ we see that $\mu_{10}(\zeta) = [\delta^{n-1}z]$. 
Under the natural map 
$\iG{\afr}I^n \to I^n = K^0 \otimes I^n \subset (\Kb \otimes \Ib)^n$
the element $\delta^{n-1}z$ maps to $1 \otimes \delta^{n-1}z$
and therefore we have 
\[
\mu_9\mu_{10}(\zeta) = [1 \otimes \delta^{n-1}z].
\]
Traveling the other route in the bottom half of~(\ref{cd:appsmo6}), 
under the sequence of natural maps 
$I^n \onto I^{n-1 \sim}(U) \to \oplus_i I^{n-1}_{t_i} 
= C^0 \otimes I^{n-1} \subset (\Cb \otimes \Ib)^{n-1}$,
the element $z$ goes to $y \set (z/1,\ldots,z/1)$. 
Therefore $\mu_8(\zeta) = [y]$. Considering $y$ as an element 
of~$K^1 \otimes I^{n-1}$ via the identification 
$K^1 = C^0$, we see that 
\[
\mu_7\mu_6\mu_8(\zeta) = [-y].
\] 
To finish the proof we must verify that 
$\mu_9\mu_{10}(\zeta) = \mu_7\mu_6\mu_8(\zeta)$, which amounts
verifying that the element 
$\beta = 1 \otimes \delta^{n-1}z \; + \; y$ in $(\Kb \otimes \Ib)^n$ 
is cohomologous to zero, i.e., $\beta$ is a coboundary. 
The element 
$\gamma= 1\otimes z \in K^0 \otimes I^n \subset (\Kb \otimes \Ib)^{n-1}$ 
maps to $\beta$.
\end{proof}

\medskip

From now on, unless specifically mentioned, all generalized
fractions, by default are C-fractions.\index{generalized fraction!C-fraction}%
\index{C-fraction|see{generalized fraction}}

As a final remark note that, with notation as above, if $J$ is an ideal 
in~$A$ such that $M$ is $J$-torsion, then the natural induced
map $H^*_{\afr+J}M \to H^*_{\afr}M$ is an isomorphism.
In particular, we may also represent the elements of 
$H^n_{\afr+J}M$ by generalized fractions.

\medskip

\subsection{Some local isomorphisms}
\label{subsec:foris}
We now describe some local isomorphisms that form the base of 
the construction of a pseudofunctor over local rings. 
By default, any local ring is considered a topological
ring under the topology given by the powers of the 
maximal ideal. If any letter such as $R,S,A,B, \ldots$ denotes 
a local ring then the corresponding maximal ideal is denoted
by $m_R,m_S,m_A,m_B, \ldots$.

Let $R \to S$ be a local homomorphism of noetherian local rings.
Let $M$ be an $m_R$-torsion $R$-module and~$L$ an~$S$-module.
Let $\wh{R},\wh{S}$ denote the respective completions of $R,S$
along the corresponding maximal ideals. Then $M$ can naturally 
be considered as an~$\wh{R}$-module. Set 
$\wh{L} \set L \otimes_S \wh{S}$.
For any integer~$r$, the $S$-modules $H^r_{m_S}(M \otimes_R L)$
and $\Homc_R(S,M)$ (where $\Homc_R$ denotes continuous 
Hom as~$R$-modules) are $m_S$-torsion modules and 
hence naturally inherit the structure of an~$\wh{S}$-module.  
Moreover there are following isomorphisms.

For any integer~$r$, there is a natural isomorphism
\begin{equation}
\label{eq:itloco0}
H^r_{m_S}(M \otimes_R L) \iso 
H^r_{m_{\wh{S}}}(M \otimes_{\wh{R}} \wh{L}_1)
\end{equation}
defined by the following sequence of natural maps 
(with $m = m_S, \wh{m} = \wh{m_S}$)
\[
H^r_{m}(M \otimes_R L) \xto{\; \alpha_1\;}
(H^r_{m}(M \otimes_R L)) \otimes_S \wh{S} \xto{\; \alpha_2\;}
H^r_{m}((M \otimes_R L) \otimes_S \wh{S}) 
\]
\[
\iso H^r_{m}(M \otimes_R \wh{L}_1) 
\xgets{\; \alpha_3\;} H^r_{m}(M \otimes_{\wh{R}} \wh{L}_1) 
\xto{\; \alpha_4\;} H^r_{m\wh{S}}(M \otimes_{\wh{R}} \wh{L}_1) 
= H^r_{\wh{m}}(M \otimes_{\wh{R}} \wh{L}_1),
\]
where $\alpha_i$ are seen to be isomorphisms due to the following
reasons. For $\alpha_1$ we use the fact that 
$H^r_{m_S}(M \otimes_R L)$ is $m_S$-torsion, for $\alpha_2$ 
we use flatness of $\wh{S}$ over $S$, for~$\alpha_3$ we use 
that $M$ is $m_R$-torsion and finally $\alpha_4$ is the 
isomorphism corresponding to extension of scalars.

The following natural map is an isomorphism
\begin{equation}
\label{eq:itloco0a}
\Homc_R(S,M) \xto{\quad} \Homc_{\wh{R}}(\wh{S},M)
\end{equation}
as is seen by taking direct limits over $i$ of the 
following sequence
\[
\Hom_R(S/m_S^i\:, M) \cong \Hom_{R/m_R^i}(\wh{S}/m_{\wh{S}}^i\:, M)
\cong \Hom_{\wh{R}/m_{\wh{R}}^i}(\wh{S}/m_{\wh{S}}^i\:, M).
\]

In the rest of this subsection we shall consider iterated 
versions of the isomorphisms in~\eqref{eq:itloco0} 
and~\eqref{eq:itloco0a}. These isomorphisms 
give the comparison maps of the pseudofunctor in 
\ref{subsec:huang}.

\medskip

\subsubsection{}
Let $A \to B \to C$ be  
local homomorphisms of noetherian local rings.
Let $I$ be an ideal in~$B$ and~$N$ an $I$-torsion $B$-module.
Let $i = \dim (B/I)$ and $j= \dim (C/m_BC)$.
For any $C$-module $L_2$ 
there exists a canonical isomorphism 
\begin{equation}
\label{eq:itloco1}
H^j_{m_C}(H^i_{m_B}(N)\otimes_B L_2) \iso
H^{i+j}_{m_C}(N \otimes_B L_2),
\end{equation}
that is described in terms of generalized fractions as follows.
Let $\bs$ be a system of parameters of length $i$ 
for $m_B/I$ and let $\bt$ be a system of parameters 
of length $j$ for $m_C/m_BC$.
Then (\ref{eq:itloco1}) assigns
\begin{equation}
\tag{\ref{eq:itloco1}$'$}
\genfrac{[}{]}{0pt}{}{\genfrac{[}{]}{0pt}{}{n}{\bs} \otimes l}{\bt}
\xto{\quad} 
\genfrac{[}{]}{0pt}{}{n \otimes l}{\bs,\bt},
\qquad
\qquad
\{n \in N, l \in L_2\}. 
\end{equation}
(By varying $\bs, \bt$ we see that (\ref{eq:itloco1}$'$) 
uniquely determines \eqref{eq:itloco1}. 
The existence and canonicity of \eqref{eq:itloco1} 
follows from \cite[2.5]{Hu}, recalled in $(\ref{eq:itloco2}^*)$ below,
modulo an isomorphism that reverses the order of the 
tensor product and a constant sign factor depending only on~$i,j$.
Also, see \ref{lem:appsmo-1} for another definition.) 
Now suppose $M$ is a zero-dimensional $A$-module and $L_1$ a 
$B$-module. Set $I \set m_AB$. Then the 
module $N = M\otimes_A L_1$
is $I$-torsion and hence from~(\ref{eq:itloco1}) we obtain 
\begin{equation}
\label{eq:itloco2}
H^j_{m_C}(H^i_{m_B}(M\otimes_A L_1)\otimes_B L_2) \iso
H^{i+j}_{m_C}(M\otimes_A (L_1\otimes_B L_2)).
\end{equation}
As in ({\ref{eq:itloco1}$'$}), if $\bs$ and $\bt$ 
denote systems of parameters
for $m_B/m_AB$ and $m_C/m_BC$ respectively, then (\ref{eq:itloco2})
is described in terms of generalized fractions by the rule
\begin{equation}
\tag{\ref{eq:itloco2}$'$}
\genfrac{[}{]}{0pt}{}{\genfrac{[}{]}{0pt}{}{m\otimes l_1}{\bs} 
\otimes l_2}{\bt} \xto{\quad} 
\genfrac{[}{]}{0pt}{}{m \otimes (l_1 \otimes l_2)}{\bs,\bt},
\qquad
\qquad
\{m \in M, l_1 \in L_1, l_2 \in L_2\}. 
\end{equation}
Finally, the isomorphism in \eqref{eq:itloco2} is compatible with 
passage to completion:
Let $\wh{A},\wh{B},\wh{C}$ denote the 
respective completions of $A,B,C$ along the corresponding maximal 
ideals. Set $\wh{L}_1 \set L_1 \otimes_B \wh{B}$,
$\wh{L}_2 \set L_2 \otimes_C \wh{C}$.
The following diagram commutes
\begin{equation}
\tag{$\,\wh{\ref{eq:itloco2}}\,$}
\begin{CD}
H^j_{m_C}(H^i_{m_B}(M\otimes_A L_1)\otimes_B L_2) @>>> 
H^{i+j}_{m_C}(M\otimes_A (L_1\otimes_B L_2)) \\
@VVV  @VVV  \\
H^j_{m_{\wh{C}}}(H^i_{m_{\wh{B}}}(M\otimes_{\wh{A}} \wh{L}_1)
\otimes_{\wh{B}} \wh{L}_2) @>>>
H^{i+j}_{m_{\wh{C}}}(M\otimes_{\wh{A}} (\wh{L}_1\otimes_{\wh{B}} \wh{L}_2))
\end{CD}
\end{equation}
where the rows are isomorphisms obtained using \eqref{eq:itloco2}
and the columns are isomorphisms obtained by applying 
\eqref{eq:itloco0} iteratively. 

For convenience of reference, we also recall the analogue 
of~\eqref{eq:itloco2} in Huang's convention\index{Huang, I-Chiau} for local cohomology
and generalized fractions.
Thus in Huang's convention, where the order of the tensor 
product is reversed, the iterated isomorphism 
of~\eqref{eq:itloco2} is written as (see {\cite[2.5]{Hu}})
\begin{equation}
\tag{$\ref{eq:itloco2}^*$}
H^j_{m_C}(L_2 \otimes_B H^i_{m_B}(L_1 \otimes_A M)) \iso
H^{i+j}_{m_C}((L_2\otimes_B L_1)\otimes_A M),
\end{equation}
and at the level of generalized fractions, is given by
\begin{equation*}
\genfrac{[}{]}{0pt}{}{l_2 \otimes \genfrac{[}{]}{0pt}{}{l_1\otimes m}{\bs}}
{\bt} \xto{\quad} 
\genfrac{[}{]}{0pt}{}{(l_2 \otimes l_1) \otimes m)}{\bt,\bs},
\qquad
\qquad
\{m \in M, l_1 \in L_1, l_2 \in L_2\}. 
\end{equation*}

\medskip

\subsubsection{}
For $A,B,C,M$ as in \eqref{eq:itloco2} there also exists a 
natural isomorphism
\begin{equation}
\label{eq:itloco3}
\Homc_B(C, \Homc_A(B,M)) \iso \Homc_A(C,M) 
\end{equation}
corresponding to ``evaluation at 1''. Furthermore the following
diagram of isomorphisms induced by~\eqref{eq:itloco0a} 
and~\eqref{eq:itloco3} commutes.
\begin{equation}
\tag{$\:\wh{\ref{eq:itloco3}}\:$}
\begin{CD}
\Homc_B(C, \Homc_A(B,M)) @>>> \Homc_A(C,M) \\
@VVV  @VVV  \\
\Homc_{\wh{B}}(\wh{C}, \Homc_{\wh{A}}(\wh{B},M)) 
@>>> \Homc_{\wh{A}}(\wh{C},M) 
\end{CD}
\end{equation}

\medskip

\subsubsection{}
Let $A,B,M$ be as in \eqref{eq:itloco2}, and let $J$ be an ideal
in $A$. Let $L$ be a $B$-module. 
Set $\ov{A} \set A/J$, $\ov{B} \set B/JB$, 
$\ov{L} = L/JL$.  Then for any 
integer $r$ there is a natural map 
\begin{equation}
\label{eq:itloco4}
H^r_{m_{\ov{B}}}(\Hom_A(\ov{A},M) \otimes_{\ov{A}} \ov{L}) \xto{\quad} 
\Hom_B(\ov{B}, H^r_{m_B}(M \otimes_A L))
\end{equation}
defined by the following sequence (where $N \set \Hom_A(\ov{A},M)$,
$m=m_B$, $\ov{m} = m_{\ov{B}}$) 
\begin{align}
H^r_{\ov{m}}(N \otimes_{\ov{A}} \ov{L}) \cong 
H^r_{m}(N \otimes_{\ov{A}} \ov{L}) \cong
H^r_{m}(N \otimes_{A} L) 
&\cong \Hom_B(\ov{B}, H^r_m(N \otimes_A L)) \notag \\
&\xto{\quad} \Hom_B(\ov{B}, H^r_m(M \otimes_A L)). \notag
\end{align}
(Note that for a surjective map $R \to \ov{R}$, 
continuous Hom = usual Hom, i.e.,
$\Hom_R^{\cont}(\ov{R},-) = \Hom_R(\ov{R},-)$.) 
If we further assume that $A \to B$ is formally smooth 
and $L$ is flat over $B$ and if $r = \dim(B/m_AB)$ 
then \eqref{eq:itloco4} is an isomorphism (\cite[3.6]{Hu}).
Finally, the following diagram, induced by \eqref{eq:itloco0},
\eqref{eq:itloco0a} and \eqref{eq:itloco4} commutes.     
\begin{equation}
\tag{$\wh{\,\ref{eq:itloco4}\,}$}
\begin{CD}
H^r_{m_{\ov{B}}}(\Hom_A(\ov{A},M) \otimes_{\ov{A}} \ov{L}) @>>>
\Hom_B(\ov{B}, H^r_{m_B}(M \otimes_A L)) \\
@VVV  @VVV  \\
H^r_{m_{\ov{\wh{B}}}}(\Hom_{\wh{A}}(\ov{\wh{A}},M) \otimes_{\ov{\wh{A}}} 
\ov{\wh{L}}) @>>>
\Hom_{\wh{B}}(\ov{\wh{B}}, H^r_{m_{\wh{B}}}(M \otimes_\wh{A} \wh{L})) 
\end{CD}
\end{equation}

\medskip

\subsubsection{}
Let $A \xto{\;f\;} B \onto A$ be local homomorphisms of noetherian
local rings that factor the identity map on~$A$. Suppose $f$ 
is formally smooth of relative dimension $r$. 
Let $t_1, \ldots,t_r$ be a regular system of parameters 
of $m_B/m_AB$. Let $M$ be 
an $m_A$-torsion $A$-module and $L$ be a rank one free $B$-module
with generator~$g$. Then there exists a natural isomorphism
(depending on the choices $t_1,\ldots,t_r$ and~$g$)
\begin{equation}
\label{eq:itloco5}
\Hom_B(A, H^r_{m_B}(M \otimes_A L)) \iso M
\end{equation}
which we now describe in terms of generalized fractions.
In this case it is convenient to first pass 
to the completions (see \cite[Chp.~5]{Hu} for proofs of statements below).
Upon completing $A,B,L$ we get that $\wh{B}$ can be identified with
a power series~ring over~$\wh{A}$, 
say $\wh{B} \cong \wh{A}[[T_1, \ldots, T_r]]$
where $t_i$ maps to $T_i$.\vspace{-3pt} Further, 
$H^r_{m_{\<\<\wh{B}}}\!(M \otimes_{\wh{A}} \wh{L})$, as an $\wh{A}$-module,
is isomorphic to a direct sum of copies of~$M\<$, i.e., 
there is a natural decomposition
\begin{equation}
\label{eq:itloco6}
H^r_{m_{\wh{B}}}(M \otimes_{\wh{A}} \wh{L}) \cong \oplus M_{i_1,\ldots,i_r}
\, , \qquad M_{i_1,\ldots,i_r} = 
\begin{cases}
M, &\text{if $i_j>0$ for all $j$;} \\
0, &\text{otherwise.}
\end{cases}
\end{equation}
Moreover if $i_j>0$ for all $j$, the canonical inclusion 
$M = M_{i_1,\ldots,i_r} \to 
H^r_{m_{\wh{B}}}(M \otimes_\wh{A} \wh{L})$ 
is given by
\[
m \to \genfrac{[}{]}{0pt}{}{m \otimes g}{T_1^{i_1}, \ldots,T_r^{i_r}},
\qquad m \in M.
\]  
It follows that the map $M_{i_1,\ldots,i_r} \to 
M_{i_1,..,i_j-1,..,i_r}$,
given by multiplication by $T_j$, is an isomorphism for $i_j > 1$
and zero for $i_j = 1$. In particular, the submodule of 
$H^r_{m_{\wh{B}}}(M \otimes_{\wh{A}} \wh{L})$
consisting of elements annihilated by $T_1, \ldots, T_r$ 
is precisely the summand $M_{1,\ldots,1}$. In view 
of~\eqref{eq:itloco0} and~\eqref{eq:itloco0a}, 
we obtain the isomorphism in \eqref{eq:itloco5}. 
Now suppose $\wh{L} = \Ohm^r_{\wh{B}/\wh{A}}$
and $g = dT_1 \wedge \ldots \wedge dT_r$.
Consider the natural maps
\begin{equation}
\tag{$\,\wh{\ref{eq:itloco5}}\,$}
\Hom_{\wh{B}}(\wh{A}, H^r_{m_{\wh{B}}}
(M \otimes_{\wh{A}} \Ohm^r_{\wh{B}/\wh{A}})) \xto{\quad i \quad} 
H^r_{m_{\wh{B}}}(M \otimes_\wh{A} \Ohm^r_{\wh{B}/\wh{A}})
\xto{\; \text{res} \;} M = M_{1,\ldots,1}
\end{equation}  
where $i$ is the canonical inclusion and $\text{res}$ is the
projection map induced by \eqref{eq:itloco6}.

\pagebreak[3] 
\noindent Then $\text{res}$, which \emph{a priori} 
depends on the choice of variables $T_1, \ldots, T_r$, 
is in fact independent of such a 
choice\footnote{The map $\text{res}$ is called the 
residue map.}. Since $\text{res}$ is canonical\index{residue map}%
\index{ OO@res|see{residue map}}
it follows that $\wh{\eqref{eq:itloco5}}$ defines 
a canonical isomorphism.
In terms of generalized fractions 
{$\wh{\eqref{eq:itloco5}}$} is given by
\begin{equation}
\tag{\ref{eq:itloco5}$'$}
m \xto{\quad} \genfrac{[}{]}{0pt}{}
{m \otimes dT_1 \wedge \ldots \wedge dT_r}{T_1, \ldots,T_r}.
\end{equation}  

\medskip

\subsection{Pseudofunctors over local rings}
\label{subsec:huang}

We now describe a canonical pseudo\-functor over local rings, one that
forms the base for the pseudofunctor that we construct over the category
$\mathbb F$ of formal schemes. Let $\Cfr$
denote the category whose objects are noetherian complete 
local rings and whose morphisms are local 
homomorphisms that are essentially of pseudo-finite type 
(see \S\ref{subsec:formal}).  
By \ref{prop:morph2}(i), our definition of $\Cfr$ agrees with 
the one in \cite[p.~28]{Hu}. Note
that $\Cfr$ is anti-isomorphic to the full subcategory of 
connected zero-dimensional formal schemes in $\mathbb F$. 
By a smooth map in $\Cfr$ we mean a map which is formally smooth 
under the topology given by powers of the corresponding maximal ideals.

\medskip

\begin{athm}
\label{thm:huang1}
For any ring $R \in \Cfr$, let $R_{\sharp}$ denote the category of 
zero-dimensional $R$-modules. Then there exists a canonical covariant
pseudofunctor~$(-)_{\sharp}$ on $\Cfr$%
\index{    ${\boldsymbol{{-}_{\sharp}}}$ (covariant\vspace{.6pt} pseudofunctor on $\Cfr$)}
and a choice for isomorphisms
as in \textup{I(i)} and \textup{I(ii)} below, such that 
\textup{II-IV} below are satisfied. 
(The terms $C_{\sharp}^{(-),(-)}$ and~$\delta_{\sharp}^{(-)}$ 
used in \emph{II} refer to the comparison isomorphisms associated 
with~$(-)_{\sharp}$). 
\begin{enumerate}
\renewcommand{\labelenumi}{\Roman{enumi}.}
\renewcommand{\labelenumii}{(\roman{enumii})}
\item \label{huang1:A}
\begin{enumerate}
\item \label{huang1:A:(i)}
If $f \colon A \to B$ is a smooth map of relative dimension $r$ 
in $\Cfr$, then for any $M \in A_{\sharp}$, there is a natural isomorphism
(where $\omega_f = \Ohm^r_{B/A}$)
\[
f_{\sharp}M \iso H^r_{m_B}(M \otimes_A \omega_f).
\]
\item \label{huang1:A:(ii)}
If $f \colon A \to B$ is a surjective map in $\Cfr$, then 
for any $M \in A_{\sharp}$, there is a natural isomorphism
\[
f_{\sharp}M \iso \Hom_A(B, M).
\]
\end{enumerate}
\item \label{huang1:B}
\begin{enumerate}
\item \label{huang1:B:(i)}
Let $A \xto{\; f \;} B \xto{\; g \;} C$ be smooth maps in $\Cfr$
having relative dimensions $r_1,r_2$ respectively. Let $t_1, t_2$ 
denote the transcendence degree of the induced maps of residue fields 
$k_A \to k_B$ and $k_B \to k_C$ respectively. Then for any 
$M \in A_{\sharp}$, the following diagram commutes
\[
\begin{CD}
g_{\sharp}f_{\sharp}M @>{C^{f,g}_{\sharp}}>> (gf)_{\sharp}M \\
@VVV  @VVV \\
H^{r_2}_{m_C}(H^{r_1}_{m_B}(M\otimes_A \omega_f)\otimes_B \omega_g) @>>>
H^{r_2+r_1}_{m_C}(M\otimes_A \omega_{gf})
\end{CD}
\]
where the vertical maps are obtained by using \emph{I.(i)} 
and the bottom row is $(-1)^{t_1r_2}$ times the map given 
by \eqref{eq:itloco2} and the natural
isomorphism $\omega_f \otimes_B \omega_g \iso \omega_{gf}$ 
induced by the exact sequence in \textup{\ref{lem:diff2}}.
\item Let $A \xto{\; f \;} B \xto{\; g \;} C$ be surjective 
homomorphisms in $\Cfr$. Then for any $M \in A_{\sharp}$, 
the following diagram, whose vertical maps are obtained using 
\emph{I.(ii)}, commutes.
\[
\begin{CD}
g_{\sharp}f_{\sharp}M @>{C^{f,g}_{\sharp}}>> (gf)_{\sharp}M \\
@VVV  @VVV  \\
\Hom_B(C, \Hom_A(B,M)) @>{\eqref{eq:itloco3}}>> \Hom_A(C,M) 
\end{CD}
\]
\item Let $f,A,B,M,r$ be as in \emph{I.(i)}. Let $J$ be an ideal in $A$.
Let $\ov{f} \colon \ov{A} \to \ov{B}$ denote the map induced 
by going modulo $J$. Let 
$A\xto{\;i\;}\ov{A}, B\xto{\;j\;}\ov{B}$ denote the canonical 
surjections. Then for any $M \in A_{\sharp}$, the following diagram
commutes
\[
\begin{CD}
\ov{f}_{\sharp}i_{\sharp}M 
@>{(C^{f,j}_{\sharp})^{-1}C^{i,\ov{f}}_{\sharp}}>> 
j_{\sharp}f_{\sharp}M \\
@VVV  @VVV  \\
H^r_{m_{\ov{B}}}(\Hom_A(\ov{A},M) \otimes_{\ov{A}} \omega_{\ov{f}}) @>>> 
\Hom_B(\ov{B}, H^r_{m_B}(M \otimes_A \omega_f))
\end{CD}
\]
where the vertical maps are obtained using \emph{I.(i),(ii)} and the bottom
row is given by \eqref{eq:itloco4} and the natural isomorphism 
$\omega_f \otimes_B \ov{B} \iso \omega_{\ov{f}}$.
\item Let $A \xto{\; f \;} B \xto{\; \pi \;} A$ be maps in $\Cfr$
such that $\pi f = 1_A$ and $f$ is smooth of relative dimension 
$r$. Then for any $M \in A_{\sharp}$, the following diagram, with 
vertical maps induced by \emph{I.(i),(ii),} commutes.
\[
\begin{CD}
\pi_{\sharp}f_{\sharp}M   
 @>{\qquad C^{f,\pi}_{\sharp}\qquad}>> M \\
@VVV  @VV{\delta^A_{\sharp}}V  \\
\Hom_B(A, H^r_{m_B}(M \otimes_A \omega_f)) 
 @>{\wh{\eqref{eq:itloco5}}}>> M 
\end{CD}
\]
\end{enumerate}
\item If $f \colon A \to A$ is the identity map, then for any 
$M \in A_{\sharp}$, the following isomorphisms 
\[  
f_{\sharp}M \xto{\;\emph{by I,(i)}\;} H^0_{m_A}M \xto{\text{canonical}} M, 
\qquad 
f_{\sharp}M \xto{\;\emph{by I,(ii)}\;} \Hom_A(A,M) \xto{\text{canonical}} M, 
\]
agree with the isomorphism induced by 
$\delta^A_{\sharp} \colon f_{\sharp} \iso 1_{A_{\sharp}}$. 
\item For any map $f \colon A \to B$ in $\Cfr$, $f_{\sharp}$ 
takes an injective hull of the residue field~$k_A$ of~$A$
to an injective hull of~$k_B$.
\end{enumerate}
\end{athm}

\medskip

The proof of \ref{thm:huang1} is based on Huang's result
\cite[Theorem 6.12]{Hu}. 
We also need the following lemma for the proof.

\medskip

\begin{alem}
\label{lem:huang2}
Let $\Bigl( (-)_{\alpha}, (-)_{\alpha}, 
C_{\alpha}^{(-),(-)}, \delta_{\alpha}^{(-)} \Bigr)$ 
be a pseudofunctor on $\Cfr$. For any two maps
$R \xto{\;f\;} S \xto{\;g\;} T$ set 
$C_{\beta}^{f,g} \set (-1)^{t_1t_2}C_{\alpha}^{f,g}$ 
where $t_1,t_2$ are the transcendence degrees of the 
induced maps of residue fields $k_R \to k_S$ and
$k_S \to k_T$ respectively.
Then $\Bigl( (-)_{\alpha}, (-)_{\alpha}, 
C_{\beta}^{(-),(-)}, \delta_{\alpha}^{(-)} \Bigr)$
is also a pseudofunctor on $\Cfr$.
\end{alem}
\begin{proof}
This is a straightforward verification.
\end{proof}

\medskip

Proof of \ref{thm:huang1}. 
In \cite{Hu}, Huang\index{Huang, I-Chiau!pseudofunctor} constructs a canonical covariant pseudofunctor 
on~$\Cfr$ which we denote, following Huang, by~`$\#$'. 
For any ring~$R$ in~$\Cfr$, $R_{\#}$ is the category 
of zero-dimensional $R$-modules and so $R_{\#} = R_{\sharp}$.
We define~`$\sharp$' as the pseudofunctor obtained by modifying 
`$\#$' using \ref{lem:huang2}. Then for any map 
$f\colon A \to B$ in $\Cfr$, we have $f_{\sharp} = f_{\#}$.  
In order to prove properties I--IV for~`$\sharp$'
we now recall some properties of~`$\#$'. 

Let $f \colon R \to T$ be a morphism 
in~$\Cfr$ and let $R \xto{\;\eta\;} S \stackrel{\pi}{\onto} T$ 
be a factorization 
of~$f$ where~$\eta$ is a smooth map of relative dimension $r$. 
For any $M \in R_{\#}$, set 
\[
(\eta, \pi)_{\#} M \set \Hom_S(T, H^r_{m_S}(\omega_{\eta} \otimes_R M))
\]
(cf.~\cite[p.~29]{Hu} where Huang uses the somewhat ambiguous 
notation $\eta^0_{\#}$ instead).
The construction of $(-)_{\#}$ provides a natural isomorphism
\begin{equation}
\label{eq:huang2a}
f_{\#}M \iso (\eta, \pi)_{\#} M, \qquad \text{see \cite[p.~33]{Hu}.}
\end{equation}
Note that the order of the tensor product $\omega_{\eta} \otimes_R M$
occurring in~$(\eta, \pi)_{\#} M$ is reverse to the one we work with, 
such as the one in~I.(i) of~\ref{thm:huang1}. Our next task includes 
redescribing Huang's~$(-)_{\#}$ in our convention. 
 
For $f$ as in I.(i), consider 
the factorization $A \xto{\;f\;}B \xto{\;=\;} B$. We 
define the required isomorphism in I.(i), 
to be given by the following sequence
\begin{equation}
\label{eq:huang3}
f_{\sharp}M = f_{\#}M \xto{\;\eqref{eq:huang2a}\;} 
(f, 1_B)_{\#}M \xto{\text{can.}} H^r_{m_B}(\omega_{f} \otimes_A M)
\iso H^r_{m_B}(M \otimes_A \omega_{f})
\end{equation}  
where the last isomorphism is induced by the obvious map
$\omega_{f} \otimes_A M \iso M \otimes_A \omega_{f}$.

For $f$ as in I.(ii), consider
the factorization $A \xto{\;=\;}A \onto B$. We define
the required isomorphism in I.(ii), to be given by the following 
sequence
\begin{equation}
\label{eq:huang4}
f_{\sharp}M = f_{\#}M \xto{\;\eqref{eq:huang2a}\;} 
(1_A,f)_{\#}M \xto{\text{can.}} \Hom_A(B, M).
\end{equation}  

Now under the hypothesis of part II.(i) of the Theorem we claim that
the following diagram commutes
where the vertical maps are obtained by using \eqref{eq:huang3}
and the bottom row is $(-1)^{t_1(r_2 + t_2)}$ times the map given 
by \eqref{eq:itloco2} and the canonical isomorphism
$\omega_f \otimes_B \omega_g \iso \omega_{gf}$ induced by
the exact sequence in \ref{lem:diff2}.
\begin{equation}
\label{eq:huang5}
\begin{CD}
g_{\#}f_{\#}M @>{C^{f,g}_{\#}}>> (gf)_{\#}M \\
@VVV  @VVV \\
H^{r_2}_{m_C}(H^{r_1}_{m_B}(M\otimes_A \omega_f)\otimes_B \omega_g) @>>>
H^{r_2+r_1}_{m_C}(M\otimes_A \omega_{gf})
\end{CD}
\end{equation}
Since $C^{f,g}_{\sharp}=(-1)^{t_1t_2}C^{f,g}_{\#}$,
commutativity of \eqref{eq:huang5} implies that of the
diagram in~II.(i).
To prove that \eqref{eq:huang5} commutes we shall expand it 
vertically using the definition of \eqref{eq:huang3}.
First note that the following diagram commutes due to functorial
reasons.
\[
\begin{CD}
g_{\#}f_{\#}M @>{\eqref{eq:huang3}}>> 
g_{\#}H^{r_1}_{m_B}(M\otimes_A \omega_f) \\
@V{\eqref{eq:huang2a}}V{g_{\#}f_{\#} \,\cong\, (g,1)_{\#}(f,1)_{\#}}V  
@VV{\eqref{eq:huang3}}V \\
H^{r_2}_{m_C}(\omega_g\otimes_B H^{r_1}_{m_B}(\omega_f\otimes_A M)) 
@>{\text{swap}}>> 
H^{r_2}_{m_C}(H^{r_1}_{m_B}(M\otimes_A \omega_f)\otimes_B \omega_g) 
\end{CD}
\]
Therefore we may expand \eqref{eq:huang5} as follows.
\begin{equation}
\label{eq:huang6}
\begin{CD}
g_{\#}f_{\#}M @>{C^{f,g}_{\#}}>> (gf)_{\#}M \\
@V{\eqref{eq:huang2a}}V{g_{\#}f_{\#} \,\cong\, (g,1)_{\#}(f,1)_{\#}}V 
@V{\eqref{eq:huang2a}}V{(gf)_{\#} \,\cong\, (gf,1)_{\#}}V  \\
H^{r_2}_{m_C}(\omega_g\otimes_B H^{r_1}_{m_B}(\omega_f\otimes_A M)) 
@>{\alpha}>> H^{r_2+r_1}_{m_C}(\omega_{gf}\otimes_A M) \\
@V{\beta}VV  @V{\gamma}VV \\
H^{r_2}_{m_C}(H^{r_1}_{m_B}(M\otimes_A \omega_f)\otimes_B \omega_g) 
@>{\delta}>> H^{r_2+r_1}_{m_C}(M\otimes_A \omega_{gf})
\end{CD}
\end{equation}
The vertical maps are the obvious ones. Thus $\beta, \gamma$ are induced 
by switching the order of the tensor products.
The maps $\alpha, \delta$ are defined as 
\[
\alpha \quad = \quad (-1)^{r_1t_2}\Bigl(\quad (\ref{eq:itloco2}^*)
\quad + \quad \omega_g \otimes_B \omega_f \xto{\;\text{canonical}\;} 
\omega_{gf}\Bigr),
\]
\[
\delta \quad = \quad (-1)^{t_1(r_2 + t_2)}
\Bigl(\quad \eqref{eq:itloco2} \quad
+ \quad \omega_f \otimes_B \omega_g \xto{\;\text{canonical}\;} 
\omega_{gf}\Bigr).
\]
Then by the definition of ${C^{f,g}_{\#}}$ in \cite[6.10]{Hu}, 
the top rectangle of \eqref{eq:huang6} commutes. For the 
bottom rectangle we use the
following calculations on generalized fractions with notation 
explained below.
\begin{align}
\genfrac{[}{]}{0pt}{}
{\mu \otimes \genfrac{[}{]}{0pt}{}{\nu\otimes m}{\bs}}{\bt} 
\xto{\;\alpha\;} 
(-1)^{r_1t_2}\genfrac{[}{]}{0pt}{}{(\mu\wedge\nu) \otimes m}{\bt,\bs}
&\xto{\;\gamma\;} (-1)^{r_1t_2}\genfrac{[}{]}{0pt}{}
{m \otimes (\mu\wedge\nu)}{\bt,\bs} \notag \\
&= (-1)^{r_1r_2 + r_1t_2}
\genfrac{[}{]}{0pt}{}{m \otimes (\mu\wedge\nu)}{\bs,\bt}, \notag \\
\genfrac{[}{]}{0pt}{}
{\mu \otimes \genfrac{[}{]}{0pt}{}{\nu\otimes m}{\bs}}{\bt} 
\xto{\;\beta\;} 
\genfrac{[}{]}{0pt}{}
{\genfrac{[}{]}{0pt}{}{m\otimes \nu}{\bs} \otimes \mu}{\bt} 
&\xto{\;\delta\;} (-1)^{t_1(r_2 + t_2)}\genfrac{[}{]}{0pt}{}
{m \otimes (\nu\wedge\mu)}{\bs,\bt} \notag \\
&= (-1)^{r_1r_2 + r_1t_2}
\genfrac{[}{]}{0pt}{}{m \otimes (\mu\wedge\nu)}{\bs,\bt}. \notag
\end{align}
(Here, $m \in M, \mu \in \omega_g, \nu \in \omega_f$ and $\bs,\bt$ 
are systems of parameters $B/m_AB, C/m_BC$ respectively and have 
lengths $r_1,r_2$ respectively. Also, the free modules 
$\Ohm^1_{B/A}$ and~$\Ohm^1_{C/B}$ have ranks $r_1+t_1, r_2+t_2$
respectively.) Thus \eqref{eq:huang5} commutes, thereby 
proving~II.(i).

The remaining cases of part II
of the Theorem are relatively straightforward. Note that 
for a closed immersion, the relative dimension is zero and so is 
the transcendence degree at the residue fields. Thus the 
sign change effected by \ref{lem:huang2} plays no role 
in the remaining cases of~II. Similarly Huang's 
sign factor in~\cite[6.10]{Hu} does not affect the remaining cases.
Thus (ii), (iii) and (iv) of~II follow easily from Huang's\index{Huang, I-Chiau!pseudofunctor}
definition of ${C^{f,g}_{\#}}$ (cf.~4.4, 3.6, 6.5 of \cite{Hu}).

Part III immediately follows from the definition of $\delta^A_{\sharp}$
in~\cite[p.~37]{Hu}. Part IV holds because $(-)_{\#}$ 
takes injective hulls to injective hulls. Finally, 
since $(-)_{\#}$ is canonical, so is $(-)_{\sharp}$.
\qed

\medskip

For convenience of reference we state the following
easily verified remark.

\begin{arem} 
\label{rem:huang7}
Let $f \colon A \to B$ be an isomorphism in $\Cfr$.
Then for $(-)_{\sharp}$ as in~\ref{thm:huang1}, for any 
$M \in A_{\sharp}$, the following two isomorphisms
are equal.
\[
f_{\sharp}M \xto{\text{\ref{thm:huang1}, I.(i)}} 
H^0_{m_B}(M \otimes_A B) = M \otimes_A B
\xto{\text{canonical}} M
\]
\[
f_{\sharp}M \xto{\text{\ref{thm:huang1}, I.(ii)}} \Hom_A(B, M) 
\xto{\text{canonical}} M
\]
\end{arem}

One proves \ref{rem:huang7}
by verifying the analogous statement for Huang's pseudofunctor~$(-)_{\#}$.

\begin{arem}
\label{rem:huang8}
In I.(ii) of \ref{thm:huang1}, 
since $f$ is assumed to be a closed immersion, 
there is no distinction between the functors
$\Hom_A(B,-)$ and $\Homc_A(B,-)$.
However, in general, the former does not preserve the property of
modules being zero-dimensional and so the latter
is preferred. Moreover we have the following.

Let $\Cfr_{\text{rf}}$ be the subcategory of $\Cfr$ consisting of 
residually finite maps, i.e., those maps for which the%
\index{residually finite homomorphism} 
corresponding induced map of residue fields is finite
(e.g., a power series ring
over the base ring). Then there is a canonical choice
for a pseudofunctor on $\Cfr_{\text{rf}}$, viz., one that assigns
to any map $f \colon A \to B$ in $\Cfr_{\text{rf}}$, the functor sending
$M \in A_{\#}$ to $\Homc_A(B,M) \in B_{\#}$, and for 
any pair of composable maps, the comparison map corresponding to 
``evaluation at 1'' (see \eqref{eq:itloco3}).
In \cite[Chp.~7]{Hu}, Huang\index{Huang, I-Chiau!pseudofunctor} shows that the restriction of $(-)_{\#}$ 
to $\Cfr_{\text{rf}}$ is in fact isomorphic to this canonical 
pseudofunctor. This isomorphism generalizes the one resulting 
from~I.(ii) of~\ref{thm:huang1}.
\end{arem}


\subsection{Iterated fractions and translation invariance}
\label{subsec:proofstwo}
The results in this subsection are somewhat technical in nature.
The main result is~\ref{lem:appsmo-1}.

Let $f \colon B\to C$ be a flat local homomorphism of 
noetherian local rings with corresponding 
maximal ideals $m_B,m_C$. Suppose~$\Nb$ is a 
bounded-below complex of $B$-modules having 
the property that there exists an integer~$j$ such that 
$H^l_{m_B}\Nb = 0$ for $l \ne j$. This property
is equivalent to requiring that there exist an integer~$j$
for which the natural truncation maps 
(\S\ref{subsec:conv}, \eqref{conv10})
\begin{equation}
\label{eq:proofstwo0}
\R\iG{m_B}\Nb \gets \tau_{\le j}\R\iG{m_B}\Nb \to
(H^{j}_{m_B}\Nb)[-j]
\end{equation}
are isomorphisms.
This property holds for any shift of $\Nb$ too.
Let $\Lb$ be a bounded-above complex of flat $C$-modules.
From \ref{prop:nolab2b} we know that the natural map 
\[
\R\iG{m_C}(\R\iG{m_B}\Nb \otimes_B \Lb) \to \R\iG{m_C}(\Nb \otimes_B \Lb) 
\]
is an isomorphism. Therefore, applying $H^i_{m_C}$ to the following 
natural map 
\begin{equation}
\label{eq:proofstwo0a}
(H^j_{m_B}\Nb)[-j] \otimes_B \Lb 
\xto{\;\cong\;} \R\iG{m_B}\Nb \otimes_B \Lb
\to \Nb \otimes_B \Lb
\end{equation}
we obtain a natural isomorphism
(for all $i$) 
\begin{equation}
\label{eq:proofstwo1}
H^i_{m_C}((H^j_{m_B}\Nb)[-j] \otimes_B \Lb) 
\iso H^i_{m_C}(\Nb \otimes_B \Lb).
\end{equation}

We elaborate on this isomorphism in the special case 
when $\Nb$ and $\Lb$ consist of one module each. 
Let $N$ be a $B$-module which, as a complex, satisfies 
the hypothesis of (\ref{eq:proofstwo0}). For application 
purposes let $q^{\prime}$ denote the integer for which 
we have the isomorphism
$\R\iG{m_B}N \iso (H^{q^{\prime}}_{m_B}N)[-q^{\prime}]$ 
as in (\ref{eq:proofstwo0}). Let $a,b$ be integers. 
Let $L$ be a flat $C$-module. For any integer $p^{\prime}$, 
we now consider a natural isomorphism  
\begin{equation}
\label{eq:proofstwo1b}
\theta_{a,b} :
H^{p^{\prime}}_{m_C}(H^{q^{\prime}}_{m_B}N \otimes_B L)
\iso H^{p^{\prime} + q^{\prime}}_{m_C} (N\otimes_B L) 
\end{equation}
which is obtained by applying $H^{p^{\prime} + q^{\prime} -a-b}_{m_C}$ to the
following sequence of maps described below. (We are also using the 
convention in \S\ref{subsec:conv}, \eqref{conv8}.)
\begin{align}
(H^{q^{\prime}}_{m_B}N \otimes_B L)[a - q^{\prime} +b] 
 &\to  (H^{q^{\prime}}_{m_B}N)[a - q^{\prime}] \otimes_B L[b] 
 \tag{\ref{eq:proofstwo1b}a} \\
 &\to  (H^{q^{\prime}-a}_{m_B}(N[a]))[a - q^{\prime}] \otimes_B L[b]  
 \tag{\ref{eq:proofstwo1b}b}\\
 &\to N[a]\otimes_B L[b]  \tag{\ref{eq:proofstwo1b}c} \\
 &\to (N\otimes_B L)[a+b]  \tag{\ref{eq:proofstwo1b}d} 
 \end{align}
The maps ({\ref{eq:proofstwo1b}a}) and ({\ref{eq:proofstwo1b}d}) are the 
natural isomorphisms obtained using the conventions 
in~\S\ref{subsec:conv},~\eqref{conv5},
while ({\ref{eq:proofstwo1b}b}) is an equality 
by~\S\ref{subsec:conv},~\eqref{conv8}.
The map ({\ref{eq:proofstwo1b}c}) is obtained from 
(\ref{eq:proofstwo0a}) using $\Nb = N[a]$, $\Lb = L[b]$ 
and~$j = q^{\prime} -a$, and so from 
(\ref{eq:proofstwo1}) we conclude that 
$H^{p^{\prime} + q^{\prime} -a-b}_{m_C}$ applied to ({\ref{eq:proofstwo1b}c})
yields an isomorphism. Thus (\ref{eq:proofstwo1b}) is an isomorphism.

\medskip

\medskip
\begin{alem}
\label{lem:appsmo-1}
In the above situation we have the following.
\begin{enumerate}
\item
The isomorphism $\theta_{a,b}$ in $(\ref{eq:proofstwo1b})$ is independent 
of the choice of $a,b$.
\item
Assume the following.
\begin{enumerate}
\item There exists an ideal $I \subset m_B$ such that $N$ is
$I$-torsion.
\item There exists a sequence $\bs$ of length ${q^{\prime}}$ 
in~$m_B$ such that $\bs$ generates $m_B/I$ up to radicals.
\item There exists a sequence $\bt$ of length ${p^{\prime}}$
in~$m_C$ such that $\bt$ generates $m_C/(m_BC)$ up to radicals.
\end{enumerate}
Then $\theta_{a,b}$
satisfies the following iteration 
formula for generalized fractions 
\[
\theta(\genfrac{[}{]}{0pt}{}{\genfrac{[}{]}{0pt}{}{n}{\bs} \otimes l}{\bt})
\quad = \quad 
\genfrac{[}{]}{0pt}{}{n \otimes l}{\bs,\bt},
\qquad \qquad n\in N, \quad l \in L.
\]
\end{enumerate}
\end{alem}

\begin{proof}
For part (i), it suffices to show that $\theta_{a,b} = \theta_{0,0}$.
Let $\phi_{a,b}$ denote the map 
$(H^{q^{\prime}}_{m_B}N \otimes_B L)[a - q^{\prime} +b] \iso
(N \otimes_B L)[a+b]$ obtained by composing 
(\ref{eq:proofstwo1b}a)--(\ref{eq:proofstwo1b}d).
To prove (i) it therefore suffices to show that 
$\phi_{a,b} = \phi_{0,0}[a+b]$ where the latter map  
denotes shifting $\phi_{0,0}$ by $a+b$.
This reduces to verifying that the outer portion 
of the following diagram commutes 
\[
\begin{CD}
(H^{q^{\prime}}_{m_B}N \otimes_B L)[a - q^{\prime} +b] 
  @>{(\ref{eq:proofstwo1b}\text{a})}>> 
  (H^{q^{\prime}}_{m_B}N)[a-q^{\prime}] \otimes_B L[b] \\
@VVV @VV{(\ref{eq:proofstwo1b}\text{b})}V \\
((H^{q^{\prime}}_{m_B}N)[-q^{\prime}] \otimes_B L)[a+b]  
 @>>> (H^{q^{\prime}-a}_{m_B}(N[a]))[a-q^{\prime}] \otimes_B L[b] \\
@VVV @VV{(\ref{eq:proofstwo1b}\text{c})}V \\
(N \otimes_B L)[a+b] @<{(\ref{eq:proofstwo1b}\text{d})}<< N[a] \otimes L[b]
\end{CD}
\]
where the maps in left column 
describe $\phi_{0,0}[a+b]$, the remaining outer 
edges define~$\phi_{a,b}$ and the 
horizontal map in the middle row is 
the composite of the following isomorphisms
obtained using the conventions in~\S\ref{subsec:conv}, \eqref{conv4}, 
\eqref{conv8}.
\begin{align}
((H^{q^{\prime}}_{m_B}N)[-q^{\prime}] \otimes_B L)[a+b]  
&= ((H^{q'-a}_{m_B}(N[a]))[-q^{\prime}] \otimes_B L)[a+b] \notag \\  
&\iso (H^{q'-a}_{m_B}(N[a]))[-q^{\prime}][a] \otimes_B L[b] \notag  
\end{align}
The upper rectangle commutes by  
Lemma \ref{lem:proofstwo3} below, while the 
lower one commutes due to functoriality 
of the truncation maps. Thus (i) follows.

To prove (ii), by (i) we may assume $a=0=b$.
By \ref{lem:appsmo5} we see that if the required iteration formula
holds for K-fractions then it also holds for C-fractions. 
Therefore it suffices to prove~(ii) 
using K-fractions.

Before we proceed further, we recall some terminology
concerning the local cohomology functors.
Let $R$ be a commutative ring, $\bu$ a finite sequence in $R$
and let~$\bu R$ denote the ideal generated by~$\bu$.
Then we set
\[
H^i_{\bu R}(-) \set H^i\R\iG{\bu R}(-), \qquad
H^i_{\uv{\bu}}(-) \set H^i\Kbi(\bu, -).  
\]
By \ref{prop:loc5}
there is a natural isomorphism 
$H^i_{\bu R}(-) \iso H^i_{\uv{\bu}}(-)$. If $R \to S$ is
a homomorphism, then we distinguish the $\Kbi$-complexes
over~$R$ and~$S$ by using the terms $\Kbi(\bu;R)$ and~$\Kbi(\bu;S)$
respectively, however for any $S$-module $M$, the term 
$H^i_{\uv{\bu}}M$ is unambiguous in view of the canonical
isomorphism 
\[
\Kbi(\bu;R)\otimes_R M \iso \Kbi(\bu;S)\otimes_S M.
\]

Returning to the proof of part (ii), we consider the following   
diagram through which we relate the map~${\theta_{0,0}}$ 
to a map involving the~$\Kbi$-complexes. 
\begin{equation}
\label{cd:appsmo0}
\begin{CD}
H^{p^{\prime}}_{m_C}(H^{q^{\prime}}_{m_B}N \otimes_B L) @<{\alpha}<<  
 \Bigl( \Kbi(\bt) \otimes_C ((\Kbi(\bs) \otimes_B N)[q^{\prime}] 
 \otimes_B L) \Bigr) [p^{\prime}] \\
@VV{\theta_{0,0}}V @VV{\gamma \; = \; ?}V \\
H^{p^{\prime}+q^{\prime}}_{m_C}(N \otimes_B L) @<{\delta}<<
 \Bigl( \Kbi(\bt,\bs)  \otimes_C (N \otimes_B L) \Bigr) 
 [p^{\prime}+q^{\prime}]
\end{CD}
\end{equation}
The objects in the left column are modules thought of as complexes 
existing in degree zero only. The objects in the right column are complexes 
of length $p^{\prime}+q^{\prime}$ since the sequences $\bt$ and $\bs$ 
have length $p^{\prime},q^{\prime}$ respectively by assumption. 
Furthermore the complexes in the right column exist between degrees 
$-p^{\prime}-q^{\prime}$ and 0.
The horizontal maps, $\alpha,\beta$, described below in detail, 
are the natural truncation maps and are surjective at the graded level.
Via the horizontal maps, the elements in the modules 
occurring in the left column can be represented by suitable 
Koszul fractions occurring in  degree zero in the complexes 
in the right column. The proof of~(ii) then 
reduces to first finding a map~${\gamma}$ which 
makes~(\ref{cd:appsmo0}) commute and then chasing the 
image of Koszul fractions through~${\gamma}$.

The horizontal maps in (\ref{cd:appsmo0}) are defined as follows.
The map $\alpha$ is defined by the composition of the maps
\begin{align}
H^{p^{\prime}}_{m_C}(H^{q^{\prime}}_{m_B}N \otimes_B L)
 &\cong H^{p^{\prime}}_{\bt C}(H^{q^{\prime}}_{\bs B}N \otimes_B L) \notag \\
&\cong H^{p^{\prime}}_{\uv{\bt}}
 (H^{q^{\prime}}_{\uv{\bs}}N \otimes_B L) \notag \\
&\gets \Bigl( \Kbi(\bt) \otimes_C (H^{q^{\prime}}_{\uv{\bs}}N \otimes_B L)
 \Bigr) [p^{\prime}] \notag \\
&\gets \Bigl( \Kbi(\bt) \otimes_C ((\Kbi(\bs) \otimes_B N)[q^{\prime}] 
 \otimes_B L) \Bigr) [p^{\prime}] \notag
\end{align}
where the last two maps are obtained by truncation of the appropriate 
$\Kbi$-complex at the highest homology spot. The map $\delta$ 
in (\ref{cd:appsmo0}) is defined by the composition of the maps
\begin{align}
H^{p^{\prime}+q^{\prime}}_{m_C}(N \otimes_B L)
 \cong H^{p^{\prime}+q^{\prime}}_{(\bt,\bs)C}(N \otimes_B L) 
&\cong H^{p^{\prime}+q^{\prime}}_{\uv{\bt,\bs}}(N \otimes_B L) \notag \\
&\gets \Bigl( \Kbi(\bt,\bs)  \otimes_C (N \otimes_B L) \Bigr) 
 [p^{\prime}+q^{\prime}]. \notag
\end{align}
where the last map is obtained by truncating at degree zero.

In order to find $\gamma$ which makes 
(\ref{cd:appsmo0}) commute we 
expand~(\ref{cd:appsmo0}) horizontally using the
definition of~$\alpha$ and vertically using the 
definition of~$\theta_{0,0}$. 
For convenience we do the expansion in two stages 
((\ref{cd:appsmo0}a), (\ref{cd:appsmo0}b) below).

\begin{small}
\begin{figure}\vspace{5pt}
\rotatebox{90}
{\begin{minipage}{8.5in}\vspace{-3mm}
\begin{equation}
\tag{\ref{cd:appsmo0}a}
\begin{CD}
H^{p^{\prime}}_{m_C}(H^{q^{\prime}}_{m_B}N \otimes_B L)
 @>>> H^{p^{\prime}}_{\bt C}(H^{q^{\prime}}_{\bs B}N \otimes_B L) 
 @>>> H^{p^{\prime}}_{\uv{\bt}}(H^{q^{\prime}}_{\uv{\bs}}N \otimes_B L) \\
@|  @| @VV{\gamma_2}V \\
H^{p^{\prime}+q^{\prime}}_{m_C}((H^{q^{\prime}}_{m_B}N 
 \otimes_B L)[-q^{\prime}])  @>>> H^{p^{\prime}+q^{\prime}}_{\bt C}
 ((H^{q^{\prime}}_{\bs B}N \otimes_B L)[-q^{\prime}]) @>>>
 H^{p^{\prime}+q^{\prime}}_{\uv{\bt}}
 ((H^{q^{\prime}}_{\uv{\bs}}N \otimes_B L)[-q^{\prime}]) \\
@V{(\ref{eq:proofstwo1b}a)}VV @VVV @VVV \\
H^{p^{\prime}+q^{\prime}}_{m_C}((H^{q^{\prime}}_{m_B}N)[-q^{\prime}]
 \otimes_B L)  @>>> H^{p^{\prime}+q^{\prime}}_{\bt C}
 ((H^{q^{\prime}}_{\bs B}N)[-q^{\prime}] \otimes_B L) @>>>
 H^{p^{\prime}+q^{\prime}}_{\uv{\bt}}
 ((H^{q^{\prime}}_{\uv{\bs}}N)[-q^{\prime}] \otimes_B L) \\
@A{\eqref{eq:proofstwo0}}AA @AA{\text{truncation}}A 
@AA{\text{truncation}}A \\
H^{p^{\prime}+q^{\prime}}_{m_C}(\R\iG{m_B}N \otimes_B L)  @>>> 
 H^{p^{\prime}+q^{\prime}}_{\bt C}(\R\iG{\bs B}N \otimes_B L) @>>>
 H^{p^{\prime}+q^{\prime}}_{\uv{\bt}}((\Kbi(\bs; B) \otimes_B N) 
 \otimes_B L) \\
@V{\text{canonical}}VV @VV{\gamma_1}V @VV{\text{canonical}}V \\
H^{p^{\prime}+q^{\prime}}_{m_C}(N \otimes_B L)  @>{\delta_1}>> 
 H^{p^{\prime}+q^{\prime}}_{\bt C}\R\iG{\bs C}(N \otimes_B L) @>>>
 H^{p^{\prime}+q^{\prime}}_{\uv{\bt}}(\Kbi(\bs; C) 
 \otimes_C (N \otimes_B L))
\end{CD}
\end{equation}
\medskip 
\qquad\qquad \hrulefill \qquad\qquad
\medskip 
\begin{equation}
\tag{\ref{cd:appsmo0}b}
\begin{CD}
H^{p^{\prime}}_{\uv{\bt}}(H^{q^{\prime}}_{\uv{\bs}}N \otimes_B L) @<<<
 \Bigl( \Kb_1 \otimes_C (H^{q^{\prime}}_{\uv{\bs}}N \otimes_B L) \Bigr)
 [p^{\prime}] @<<<
 \Bigl( \Kb_1 \otimes_C (\Kb_N[q^{\prime}] \otimes_B L) \Bigr)
 [p^{\prime}] \\
@V{\gamma_2}VV  @VV{\gamma_3}V @VV{\gamma_5}V \\
H^{p^{\prime}+q^{\prime}}_{\uv{\bt}}
 ((H^{q^{\prime}}_{\uv{\bs}}N \otimes_B L)[-q^{\prime}]) @<<<
 \Bigl( \Kb_1 \otimes_C (H^{q^{\prime}}_{\uv{\bs}}N \otimes_B L)
 [-q^{\prime}] \Bigr) [p^{\prime}+q^{\prime}] @<<<
 \Bigl( \Kb_1 \otimes_C (\Kb_N[q^{\prime}] \otimes_B L)[-q^{\prime}] \Bigr)
 [p^{\prime}+q^{\prime}] \\
@VVV @VV{\gamma_4}V @VV{\gamma_6}V \\
H^{p^{\prime}+q^{\prime}}_{\uv{\bt}}
 ((H^{q^{\prime}}_{\uv{\bs}}N)[-q^{\prime}] \otimes_B L) @<<<
 \Bigl( \Kb_1 \otimes_C ((H^{q^{\prime}}_{\uv{\bs}}N)[-q^{\prime}] 
 \otimes_B L) \Bigr) [p^{\prime}+q^{\prime}] @<<<
 \Bigl( \Kb_1 \otimes_C (\Kb_N \otimes_B L) \Bigr) 
 [p^{\prime}+q^{\prime}] \\
@A{\text{truncation}}AA @A{\text{truncation}}AA @| \\
H^{p^{\prime}+q^{\prime}}_{\uv{\bt}}(\Kb_N \otimes_B L) @<<<
 \Bigl( \Kb_1 \otimes_C (\Kb_N \otimes_B L) \Bigr) 
 [p^{\prime}+q^{\prime}] @=
 \Bigl( \Kb_1 \otimes_C (\Kb_N \otimes_B L) \Bigr) 
 [p^{\prime}+q^{\prime}] \\
@V{\text{canonical}}VV @V{\text{canonical}}VV @VV{\text{canonical}}V \\
H^{p^{\prime}+q^{\prime}}_{\uv{\bt}}(\Kb_2 
 \otimes_C (N \otimes_B L)) @<<< 
 \Bigl( \Kb_1 \otimes_C (\Kb_2  \otimes_C (N \otimes_B L)) \Bigr) 
 [p^{\prime}+q^{\prime}] @<{\text{canonical}}<<
 \Bigl( \Kbi(\bt,\bs)  \otimes_C (N \otimes_B L) \Bigr) 
 [p^{\prime}+q^{\prime}]
\end{CD}
\end{equation}
\end{minipage}\hspace{-14mm} 
}
\end{figure}
\end{small}

The leftmost column in (\ref{cd:appsmo0}a) corresponds 
to the definition of~$\theta_{0,0}$. 
The horizontal maps on the left side in~(\ref{cd:appsmo0}a) 
are, with the exception of the bottommost map~$\delta_1$,
the isomorphisms induced by the natural maps
$\R\iG{m_C} \to \R\iG{\bt C}$ and $\R\iG{m_B} \to \R\iG{\bs B}$
while~$\delta_1$ is induced by the natural map
\[
\R\iG{m_C} \to \R\iG{(\bt,\bs)C} = \R(\iG{\bt C}\iG{\bs C}) \cong 
\R\iG{\bt C}\R\iG{\bs C}.
\]
The vertical maps in the middle column, with the exception of the 
bottommost one~$\gamma_1$, are defined similar to 
the way the corresponding ones on the left column are defined 
while $\gamma_1$ is induced by the natural isomorphism
\[
\R\iG{\bs B}N \otimes_B L \cong \R\iG{\bs C}(N \otimes_B L). 
\]
The horizontal maps on the right side in~(\ref{cd:appsmo0}a)  
are induced by the natural isomorphisms of the kind
$\R\iG{?}(-) \iso \Kbi(?) \otimes (-)$. The map $\gamma_2$ 
is induced by 
\[
\Theta \colon \Kbi(\bt, (-))[-q'] \iso \Kbi(\bt, (-)[-q'])
\]
where $\Theta$ is the second component of the $\delta$-functor 
pair~$(\Kbi(\bt, (-)), \Theta)$, (cf.~\ref{prop:loc5}).
The remaining vertical maps in the rightmost column 
of~(\ref{cd:appsmo0}a) are the obvious counterparts of the 
maps in the middle column.

The top three rectangles on the left side and 
the bottom three rectangles on the right side 
in~({\ref{cd:appsmo0}a}) commute for obvious reasons. 
The topmost rectangle on the right commutes because the
isomorphism in~\ref{prop:loc5} is one of $\delta$-functors. 
To verify the commutativity of the bottom rectangle 
on the left we expand 
it horizontally using obvious natural maps as follows.
{\small{
\[
\begin{CD}
H^{p^{\prime}+q^{\prime}}_{m_C}(\R\iG{m_B}N \otimes_B L)  @>>> 
 H^{p^{\prime}+q^{\prime}}_{(\bt,\bs)C}(\R\iG{\bs B}N \otimes_B L) @>>>
 H^{p^{\prime}+q^{\prime}}_{\bt C}(\R\iG{\bs B}N \otimes_B L) \\
@VVV @VVV @VV{\gamma_1}V \\
H^{p^{\prime}+q^{\prime}}_{m_C}(N \otimes_B L)  @>>> 
 H^{p^{\prime}+q^{\prime}}_{(\bt,\bs)C}(N \otimes_B L) @>>>
 H^{p^{\prime}+q^{\prime}}_{\bt C}\R\iG{\bs C}(N \otimes_B L) 
\end{CD}
\]
}}
In this diagram the rectangle on the left commutes for functorial
reasons while the one on the right commutes 
by~Lemma \ref{lem:appsmo1} below. Thus ({\ref{cd:appsmo0}a}) commutes.

For ({\ref{cd:appsmo0}b}) we use the notation
\[
\Kb_1 \set \Kbi(\bt;C), \quad \Kb_2 \set \Kbi(\bs;C), \quad
\Kb_N \set  \Kbi(\bs;B) \otimes_B N.
\]
The horizontal maps in ({\ref{cd:appsmo0}b})
are the obvious ones given by truncation. (In each case the 
complex in question has zero homology in higher degrees.)
The leftmost column of~({\ref{cd:appsmo0}b}) is precisely
the rightmost column of~({\ref{cd:appsmo0}a}). 
The remaining columns of~({\ref{cd:appsmo0}b})
are the obvious counterparts of its leftmost one. 
In particular, the convention in~\S\ref{subsec:conv}, \eqref{conv4},
applies for the maps $\gamma_3,\gamma_4,\gamma_5,\gamma_6$.
Commutativity of ({\ref{cd:appsmo0}b}) is obvious from the
functoriality of the maps involved. 

Upon composing the top rows in ({\ref{cd:appsmo0}a}), 
({\ref{cd:appsmo0}b}) we recover the map~$\alpha$ 
of~(\ref{cd:appsmo0}).
In order to verify that the composition of the bottom rows in 
({\ref{cd:appsmo0}a}), ({\ref{cd:appsmo0}b}) is the same 
as~$\delta$ of~(\ref{cd:appsmo0}) 
we reduce to checking commutativity of the 
following diagram whose left column gives $\delta$ 
and whose right column gives the sequence of bottom rows in 
({\ref{cd:appsmo0}a})-({\ref{cd:appsmo0}d}).
(For simplicity, we set $n = p'+q'$, $M = N \otimes_B L$.)
\[
\begin{CD}
H^n_{m_C}M @= H^n_{m_C}M \\
@VVV @VVV \\
H^n_{(\bt,\bs)C}M @>>> H^n_{\bt C}\R\iG{\bs C}M \\
@VVV @VVV \\
H^n_{\uv{\bt,\bs}}M @>>> H^n_{\uv{\bt}}(\Kbi(\bs;C) \otimes_C M) \\
@AAA @AAA \\
\Bigl( \Kbi(\bt,\bs)  \otimes_C M \Bigr) [n] @=
 \Bigl( \Kbi(\bt,\bs)  \otimes_C M \Bigr) [n]
\end{CD}
\]
The horizontal arrows in the middle are the obvious natural ones.
Commutativity of the above diagram follows easily.

Since ({\ref{cd:appsmo0}a}) + ({\ref{cd:appsmo0}b}) = \eqref{cd:appsmo0}
holds along three edges of~\eqref{cd:appsmo0}
we therefore have found a map $\gamma$ which makes \eqref{cd:appsmo0}
commute, viz., $\gamma$ is given by the rightmost 
column in~({\ref{cd:appsmo0}b}). 

To complete the proof of (ii),
let the sequence 
$\bs$ be given by elements $ s_1,\ldots,s_{q^{\prime}}$ and 
$\bt$ by $t_1,\ldots,t_{p^{\prime}}$.
Fix $n \in N, l \in L$. Set 
\[
z_1 \set 
\frac{\frac{n}{s_1 \cdots s_{q^{\prime}}}\otimes l}
{t_1 \cdots t_{p^{\prime}}} \in F^0
\quad \text{where} \quad \Fb \set
\Bigl( \Kbi(\bt) \otimes_C ((\Kbi(\bs) \otimes_B N)[q^{\prime}] 
\otimes_B L) \Bigr) [p^{\prime}]
\]
and 
\[
z_2 \set 
\frac{(-1)^{p^{\prime}q^{\prime}}n \otimes l}{t_1 \cdots t_{p^{\prime}} 
\cdot s_1 \cdots s_{q^{\prime}}} \in G^0
\quad \text{where} \quad \Gb \set
\Bigl( \Kbi(\bt,\bs)  \otimes_C (N \otimes_B L) \Bigr) 
[p^{\prime}+q^{\prime}].
\]
The maps $\alpha$ and $\delta$ of (\ref{cd:appsmo0}) 
induce, in degree zero, the following maps of~$C$-modules 
\[
\alpha^0 \colon F^0 \to H^{p^{\prime}}_{m_C}
(H^{q^{\prime}}_{m_B}N \otimes_B L),
\quad
\delta^0 \colon G^0 \to H^{p^{\prime}+q^{\prime}}_{m_C}(N \otimes_B L).
\]
From the definition of K-fractions in~(\ref{eq:appsmo4}) we see that 
\[
\alpha^0(z_1) = \genfrac{[}{]}{0pt}{}{\genfrac{[}{]}{0pt}{}{n}{\bs} 
\otimes l}{\bt} \quad \text{and} \quad
\delta^0(z_2) = \genfrac{[}{]}{0pt}{}{n \otimes l}{\bs,\bt}. 
\]
It therefore suffices to check that $\gamma^0(z_1) = z_2$ 
where $\gamma^0$ is the degree zero component of~$\gamma$.
Using that~$\gamma$ is given by the rightmost column 
in~({\ref{cd:appsmo0}b}), the verification $\gamma^0(z_1) = z_2$ 
follows easily. (Note that as per~\S\ref{subsec:conv}, \eqref{conv4}, 
the map~$(\gamma_5)^0$ 
of~({\ref{cd:appsmo0}b}) has a sign 
of~$(-1)^{p^{\prime}q^{\prime}}$
while $(\gamma_6)^0$ is the identity map.)
\end{proof}

\medskip

In the proof of \ref{lem:appsmo-1}, we used the following two lemmas.

\medskip

\begin{alem}
\label{lem:appsmo1}
Let $B \to C$ be a homomorphism of noetherian rings. Let $N$ be a $B$-module
and let $L$ be a flat $C$-module. For any ideals 
$\bfr \subset B, \cfr \subset C$, the following diagram,  
with maps described below, commutes.
\[
\begin{CD}
\R\iG{\cfr + \bfr C}(\R\iG{\bfr}N \otimes_B L) @>{\alpha_1}>>
 \R\iG{\cfr}(\R\iG{\bfr}N \otimes_B L) \\
@V{\alpha_2}VV  @V{\alpha_3}VV \\
\R\iG{\cfr + \bfr C}(N \otimes_B L) @>{\alpha_4}>>
\R\iG{\cfr}\R\iG{\bfr C}(N \otimes_B L) 
\end{CD}
\]
The map $\alpha_1$ is induced by the natural map 
$\iG{\cfr + \bfr C} \to \iG{\cfr}$. The map $\alpha_2$ is induced by the 
natural map $\R\iG{\bfr}N \to N$. The map $\alpha_3$ is induced by the 
natural isomorphism $\R\iG{\bfr}N \otimes_B L \iso 
\R\iG{\bfr C}(N \otimes_B L)$. The map $\alpha_4$ is the natural 
isomorphism of derived functors corresponding to the composition
$\iG{\cfr + \bfr C} = \iG{\cfr}\iG{\bfr C}$. (The functor $\iG{\bfr C}$ sends
injective modules to injective ones.)
\end{alem}
\begin{proof}
Set $\bfr' \set \bfr C$. 
We expand the above diagram vertically as follows
where the unlabeled maps are the obvious natural ones.
\begin{equation}
\label{cd:appsmo1a}
\begin{CD}
\R\iG{\cfr + \bfr'}(\R\iG{\bfr}N \otimes_B L) @>{\alpha_1}>>
 \R\iG{\cfr}(\R\iG{\bfr}N \otimes_B L) \\
@VVV  @V{\alpha_3}VV \\
\R\iG{\cfr + \bfr'}\R\iG{\bfr'}(N \otimes_B L) @>>>
 \R\iG{\cfr}\R\iG{\bfr'}(N \otimes_B L) \\
@VVV  @| \\
\R\iG{\cfr + \bfr'}(N \otimes_B L) @>{\alpha_4}>>
\R\iG{\cfr}\R\iG{\bfr'}(N \otimes_B L) 
\end{CD}
\end{equation}
The top retangle obviously commutes. Set $M \set N \otimes_B L$.
Let $M \to \Ib$ be a $C$-injective resolution.
Then $\R\iG{\bfr'}M = \iG{\bfr'}\Ib$. 
Since $\iG{\bfr'}\Ib$ also consists of $C$-injectives,
the lower rectangle of \eqref{cd:appsmo1a} reduces 
to the following diagram, where each map being 
the obvious one, is an equality.
\[
\begin{CD}
\iG{\cfr + \bfr'}\iG{\bfr'}\Ib @>>> \iG{\cfr}\iG{\bfr'}\Ib \\
@VVV  @VVV \\
\iG{\cfr + \bfr'}\Ib @>>> \iG{\cfr}\iG{\bfr'}\Ib
\end{CD}
\]

\end{proof}

\medskip

The following lemma is proven by a direct calculation and 
we omit the proof.

\medskip

\begin{alem}
\label{lem:proofstwo3}
Let $R$ be a commutative ring. Let $\Fb$, $\Gb$ be 
complexes of $R$-modules. For any integers $i,j,k$ 
the following diagram commutes
\[
\begin{CD}
(\Fb \otimes \Gb)[i+j+k] @>{(\theta^{F,G}_{i+j,k})^{-1}}>> 
 \Fb[i+j] \otimes \Gb[k] \\
@V{(\theta^{F,G}_{i,0})^{-1}[j+k]}VV @| \\
(\Fb[i] \otimes \Gb)[j+k] @>{(\theta^{F[i],G}_{j,k})^{-1}}>> 
 (\Fb[i])[j] \otimes \Gb[k]
\end{CD}
\]
where all the maps are obtained using the convention
in~\S\ref{subsec:conv}, \eqref{conv4}.
\end{alem}

\newpage

\section{Pseudofunctorial behavior for smooth maps}
\label{sec:cosmo}

The main result here is Proposition \ref{prop:itloco1}.
In~\S\ref{subsec:gl2lo} we provide the main input that goes 
into defining the isomorphism \eqref{eq:out2} 
of~\S\ref{subsec:outline} (see \eqref{eq:gl2lo1} below). 
In \S\ref{subsec:itcoz}, we consider the situation of 
composition of two smooth maps.
Suppose $f,g$ are smooth maps such that $gf$ exists. Then we
consider the two Cousin-valued functors, one obtained
by using \ref{exam:mainre1}(i)
for $gf$ and the other by using \ref{exam:mainre1}(i)
iteratively for~$g$ and~$f$. 
We provide a comparison isomorphism between
these two in~\eqref{eq:itcoz1} below. In 
Proposition~\ref{prop:itloco1} we show that this global
comparison map is compatible with the local one 
of~\eqref{eq:itloco2} via the isomorphism in~\eqref{eq:gl2lo1} 
below.

In particular, the results of this section allow us 
to define a Coz-valued pseudofunctor for the
category of smooth maps in $\bbFc$ (see \S\ref{subsec:fin-sm}).
The proof of~\ref{prop:itloco1} is somewhat lengthy
and \S\ref{subsec:scheme} and \S\ref{subsec:conc}
are devoted exclusively towards this. 
In~\S\ref{subsec:scheme} we decompose the diagram 
of~\ref{prop:itloco1} into convenient parts,
the non-trivial cases being handled in \S\ref{subsec:conc}.

\subsection{An isomorphism}
\label{subsec:gl2lo}
Let $h \colon (\X,\Delta_1) \to (\Y,\Delta)$ be a smooth map 
in~$\bbFc$ having constant relative dimension, say $d$. 
Let $\cMb$ be a complex in~$\Coz_{\Delta}(\Y)$, i.e., 
$\cMb$ is a $\Delta$-Cousin 
complex of quasi-coherent torsion $\OY$-modules. Let $\cL$ be a 
flat quasi-coherent $\OX$-module. By \ref{lem:nolab1},
the complex $h^*\cMb \otimes_{\X} \cL[d\>]$ consists 
of $\Aqc(\X)$-modules and hence 
$E_{\Delta_1}\R\iGp{\X}(h^*\cMb \otimes_{\X} \cL[d\>]) \in  
\Coz_{\Delta_1}(\X)$. (see \eqref{eq:coz3})

Let $x \in \X, y = h(x), p = {\Delta_1}(x), q = \Delta(y)$.
We now define a natural isomorphism
\begin{equation}
\label{eq:gl2lo1}
 (E_{\Delta_1}\R\iGp{\X}(h^*\cMb \otimes_{\X} 
\cL[d\>]))(x) \iso H^{p+d-q}_{m_x}(\cMb(y) \otimes_{y} \cL_x).
\end{equation}
Set $M = \cMb(y)$. Then (\ref{eq:gl2lo1}) is the composition of 
the following isomorphisms (see explanatory remarks below).
\begin{align}
&(E_{\Delta_1}\R\iGp{\X}(h^*\cMb \otimes_{\X} \cL[d\>]))(x) & &
 \notag \\ 
&\iso H^p_{m_x}(h^*\cMb \otimes_{\X} \cL[d\>])_x \tag{\ref{eq:gl2lo1}a}
 & &\text{(by \eqref{eq:coz3})} \\
&\iso H^p_{m_x}(h^*(i_yM)[-q] \otimes_{\X} \cL[d\>])_x \tag{\ref{eq:gl2lo1}b}
 & &\text{(via truncation)} \\
&\iso H^p_{m_x}(M[-q] \otimes_{y} \cL_x[d\>]) \tag{\ref{eq:gl2lo1}c}
 & &\text{(stalk at $x$)} \\
&\iso H^p_{m_x}((M\otimes_{y} \cL_x)[d-q]) \tag{\ref{eq:gl2lo1}d} 
 & &\text{(see \S\ref{subsec:conv}, \eqref{conv5})} \\
&\iso H^{p+d-q}_{m_x}(M\otimes_{y} \cL_x). \tag{\ref{eq:gl2lo1}e} 
 & &\text{(see \S\ref{subsec:conv}, \eqref{conv8})} 
\end{align}
For (\ref{eq:gl2lo1}b) we use the sequence of natural maps 
\[
\cMb \to \sigma_{\le q}\cMb \gets \sigma_{\ge q}\sigma_{\le q}\cMb 
= \M^q[-q] \hookleftarrow i_yM[-q].
\]
Each one of the complexes in this sequence
is in $\Coz_{\Delta}(\Y)$ and hence by~(\ref{lem:psm2}), the functor 
$H^p_{m_x}(h^*(-) \otimes_{\X} \cL[d\>])$ sends each map in  
this sequence to an isomorphism. Composing these 
isomorphisms appropriately results in ({\ref{eq:gl2lo1}b}).
The map in~({\ref{eq:gl2lo1}c}) is induced by the canonical
identification of the stalk at $x$ (the stalk of
$i_yM$ at~$y$ being~$M$). 
The map in ({\ref{eq:gl2lo1}d}) is induced by the isomorphism
(\S\ref{subsec:conv}, \eqref{conv5})
\[
M[-q] \otimes_y \cL_x[d\>] \iso (M \otimes_y \cL_x)[d-q]
\]
which, as per convention, is $(-1)^{qd}$ times 
the identity map on~$M \otimes_y \cL_x$. 

\medskip

\begin{arem}
\label{rem:gl2lo3}
Another way of obtaining the isomorphism from 
({\ref{eq:gl2lo1}b})-({\ref{eq:gl2lo1}c}) is the following. 
Consider the complex $\cMb_y$. Since $\cMb$ is a Cousin complex,
therefore $(\M^i)_y = 0$ for $i > q$ and furthermore
$(\M^q)_y = \cMb(y) = M$. In particular, there is a natural map 
of complexes $M[-q] \to \cMb_y$, and therefore an induced natural map
$M[-q] \otimes_{y} \cL_x[d\>] \to \cMb_y \otimes_{y} \cL_x[d\>]$.
The following isomorphisms
\begin{align}
H^p_{m_x}(h^*\cMb \otimes_{\X} \cL[d\>])_x &\iso 
 H^p_{m_x}(\cMb_y \otimes_{y} \cL_x[d\>]) \tag{\ref{eq:gl2lo1}b$'$ } \\
&\osi H^p_{m_x}(M[-q] \otimes_{y} \cL_x[d\>]) \tag{\ref{eq:gl2lo1}c$'$}
\end{align}
define the same map as the one resulting from 
({\ref{eq:gl2lo1}b})--({\ref{eq:gl2lo1}c}).
\end{arem}

\medskip

\begin{arem}
\label{rem:gl2lo2}
By using \eqref{eq:nolab3} instead of \eqref{eq:coz3} in ({\ref{eq:gl2lo1}a})
we see that ({\ref{eq:gl2lo1}a})-({\ref{eq:gl2lo1}e}) in fact give an
isomorphism
\[
H^i_x\R\iGp{\X}(h^*\cMb \otimes_{\X} \cL[d\>]) \iso
H^{i+d-q}_{m_x}(\cMb(y) \otimes_{y} \cL_x)
\]
for all $i$. 
\end{arem}

\medskip

\begin{alem}
\label{lem:quasim1}
With notation as above, 
$\R\iGp{\X}(h^*\cMb \otimes_{\X} \cL[d\>])$ is Cohen-Macaulay 
with respect to\/~$\Delta_1$ \textup{(\S\ref{subsec:cm}).} In other words, 
for any\/ $x \in \X$ and any\/ 
$i \in \mathbb Z$ such that\/ $i \not= p = \Delta_1(x)$ we have
\[
H^i_x\R\iGp{\X}(h^*\cMb \otimes_{\X} \cL[d\>]) \; = \; 0 \; = \;
H^{i+d-q}_{m_x}(\cMb(y) \otimes_{y} \cL_x).
\]
\end{alem}
\begin{proof}
Set $M \set \cMb(y)$.
We claim that for any finitely generated submodule~$N$
of~$M$ and any free $\OXx$-module~$L$ of finite rank, 
we have $H^{i+d-q}_{m_x}(N \otimes_{y} L) = 0$ for $i \ne p$.
Since $M \otimes_{y} \cL_x$ is a direct limit 
of modules of the type $N \otimes_{y} L$ and since 
local cohomology commutes with direct limits the 
lemma follows from the claim and~\ref{rem:gl2lo2}.

Set $F \set N \otimes_{y} L$. Since $N$ has finite length
it is a Cohen Macaulay $\OYy$-module. By \ref{prop:morph2}, (ii),
the induced local homomorphism $\OYy \to \OXx$ is flat.
Moreover the fiber ring $\OXx/m_y\OXx$, being formally smooth
over the residue field $\OYy/m_y$, is a regular ring and 
hence a Cohen-Macaulay ring.
Therefore, by \cite[(6.3.3)]{EGAIV} or \cite[p.\,181, Corollary]{Ma}, 
$F$ is a Cohen Macaulay $\OXx$-module. In particular, 
$H^{j}_{m_x}F = 0$ for $j \ne \dim F$.
It suffices to show that 
$\dim F\set \dim \textup{Supp}(F)  = p+d-q$. 

Since $N$ has finite length, we conclude that
\[
\Supp(F) = \Supp((\OYy/m_y) \otimes_y L) = \Supp(\OXx/m_y\OXx).
\] 
Therefore, by \ref{cor:app4}, $\dim(F) = p+d-q$.
\end{proof}

\smallskip

\subsection{Iteration of the Cousin functor}
\label{subsec:itcoz}

\def\Esssg{{\mathbb E}_{\scriptscriptstyle g}}
\def\Esssf{{\mathbb E}_{\mspace{-1mu}\scriptscriptstyle f}}
\def\Esssgf{{\mathbb E}_{\scriptscriptstyle g\!f}}

Let $(\X,\DsssX) \xto{f} (\Y,\DsssY) \xto{g} (\Z,\DsssZ)$ 
be smooth maps in~$\bbFc$ having constant relative dimension 
$d,e$ respectively. 
Let $\cL_1$ be a quasi-coherent flat  $\OY$-module and 
$\cL_2$ a quasi-coherent flat $\OX$-module. 
For any complex $\cFb$ of $\OZ$-modules and 
$\cGb$ of $\OY$-modules set
\begin{align}
\Esssg\cFb \set &E_{\DsssY}\R\iGp{\Y}
  (g^*\cFb \otimes_{\Y} \cL_1[e]), \notag \\
\Esssf\cGb \set &E_{\DsssX}\R\iGp{\X}
  (f^*\cGb \otimes_{\X} \cL_2[d\>]), \notag \\
\Esssgf\cFb \set &E_{\DsssX}\R\iGp{\X}
  ((gf)^*\cFb \otimes_{\X} (f^*\cL_1 \otimes_{\X} \cL_2)[d+e]). \notag
\end{align}
The functors ${\mathbb E}_{(-)}$ map Cousin complexes to 
Cousin complexes:\index{ $\Forget$0@${\mathbb E}_{(-)}$ (Cousin-complex functor associated to smooth map $(-)\>$)} 
\begin{align}
\Esssg(\Coz_{\DsssZ}(\Z)) &\subset \Coz_{\DsssY}(\Y), \notag \\
\Esssf(\Coz_{\DsssY}(\Y)) 
  &\subset \Coz_{\DsssX}(\X), \notag\\
\Esssgf(\Coz_{\DsssZ}(\Z)) &\subset \Coz_{\DsssX}(\X). \notag
\end{align}
Let $\cMb \in \Coz_{\DsssZ}(\Z)$.
There exists a natural isomorphism 
\begin{equation}
\label{eq:itcoz1}
\Esssgf\cMb \iso \Esssf\Esssg\cMb
\end{equation}
which we describe as follows : 
(with $\cGb \set g^*\cMb \otimes_{\Y} \cL_1[e]$) 
\begin{align}
\Esssgf\cMb &\set E_{\DsssX}\R\iGp{\X}((gf)^*\cMb \otimes_{\X} 
  (f^*\cL_1 \otimes_{\X} \cL_2)[d+e]) & \notag \\
&\iso E_{\DsssX}\R\iGp{\X}(f^*g^*\cMb \otimes_{\X} 
 (f^*\cL_1[e] \otimes_{\X} \cL_2[d\>])) 
  & \tag{\ref{eq:itcoz1}a}\\
&\iso \Esssf\cGb
  & \tag{\ref{eq:itcoz1}b}\\
&\osi \Esssf\R\iGp{\Y}\cGb
   & \tag{\ref{eq:itcoz1}c}\\
&\iso \Esssf E_{\DsssY}\R\iGp{\Y}\cGb 
 =: \Esssf\Esssg\cMb 
\tag{\ref{eq:itcoz1}d}
\end{align}
The map ({\ref{eq:itcoz1}a}) is induced by 
$(f^*\cL_1 \otimes_{\X} \cL_2)[d+e] \iso 
f^*\cL_1[e] \otimes_{\X} \cL_2[d\>]$
obtained from our convention in \S\ref{subsec:conv}, \eqref{conv5}.
In particular, ({\ref{eq:itcoz1}a}) is $E\R\iGp{\X}(-)$ of a map 
that is $(-1)^{ed}$ times the identity map at the graded level. 
The maps ({\ref{eq:itcoz1}b}), ({\ref{eq:itcoz1}c}) are the 
obvious natural ones, the latter being an isomorphism 
by~\ref{prop:nolab2b}. By \ref{lem:quasim1}, 
$\R\iGp{\Y}\cGb$ 
is a $\DsssY$-Cohen-Macaulay complex on~$\Y$. We define
({\ref{eq:itcoz1}d}) to be the natural map induced by 
the unique isomorphism $QE_{\DsssY}\R\iGp{\Y}\cGb \cong \R\iGp{\Y}\cGb$ 
obtained from~\ref{cor:coz4a}.

The isomorphism in \eqref{eq:itcoz1} can be described explicitly at 
the punctual level. Fix a point $x$ in~$\X$ and let
$y=f(x)$, $z= g(y)$. Let $\DsssZ(z) = r$,
$\DsssY(y) = q$, $\DsssX(x) = p$.
Let $\OXx$, $\OYy$, $\OZz$ be the corresponding local 
rings and $m_x,m_y,m_z$ the corresponding maximal 
ideals. Set $p^{\prime} \set p+d-q$ and $q^{\prime} \set q+e-r$.

\def\vertloco{\rotatebox{-90}{\makebox[0pt]{$\xrightarrow{\qquad
{\rotatebox{90}{\scriptsize \eqref{eq:gl2lo1}}}\qquad}$}}}

\medskip
\begin{aprop}
\label{prop:itloco1}
The following diagram of $\OXx$-modules commutes.
\[
\begin{CD}
(\Esssf\Esssg\cMb)(x) @<{(\ref{eq:itcoz1})}<<
 (\Esssgf\cMb)(x) \\
@VV{(\ref{eq:gl2lo1})}V    \\
H^{p^{\prime}}_{m_x}((\Esssg\cMb)(y)\otimes_y {\cL_2}_x) @. 
 \vertloco \\
@VV{\text{induced by }(\ref{eq:gl2lo1})}V    \\
H^{p^{\prime}}_{m_x}(H^{q^{\prime}}_{m_y}(\cMb(z)\otimes_z{\cL_1}_y)
 \otimes_y {\cL_2}_x) @>{(\ref{eq:itloco2})}>> 
 H^{p^{\prime}+q^{\prime}}_{m_x}
 (\cMb(z)\otimes_z({\cL_1}_y\otimes_y{\cL_2}_x))
\end{CD}
\]
\end{aprop}

\medskip

The rest of this section is devoted to 
proving Proposition \ref{prop:itloco1}.
In \ref{subsec:scheme} we indicate how to decompose the diagram 
in Proposition \ref{prop:itloco1} into more convenient parts whose 
commutativity is verified in the last subsection.

\subsection{Scheme of proof of \ref{prop:itloco1}}
\label{subsec:scheme}

We now outline the general scheme of our proof 
of Proposition \ref{prop:itloco1}. 
The diagram, whose commutativity is in question, 
can be written as follows.
\def\vertschemea{\rotatebox{-90}{\makebox[0pt]
{$\xrightarrow{\hspace{5em}}\;\,$}}}
\begin{equation}
\label{cd:scheme1}
\begin{CD}
V_1  @<<<  V_2  \\
@VVV   \\
V_3 @. \vertschemea \\
@VVV   \\
V_4  @>>> V_5
\end{CD}
\end{equation}
Without loss of generality, we assume that 
$\cMb = i_zM[-r]$ where~$M =\cMb(z)$. This can be justified using 
the truncation arguments of~(\ref{eq:gl2lo1}b) as
each truncation map is functorial and the corresponding 
induced map on the cohomologies is an isomorphism. 
In order to prove commutativity of~(\ref{cd:scheme1}), 
we expand it in the following way explained below.
\begin{equation}
\label{cd:scheme2}
\begin{CD}
V_1  @= V_1 @<<<  V_2 @=  V_2 \\
@VVV  @VVV @VVV  \\
V_3 @. V_6  @>>> V_7 @.  \vertschemea \\
@VVV  @VVV @VVV  \\
V_4  @= V_4 @>>> V_8 @>>> V_5
\end{CD}
\end{equation}
For the above diagrams and subsequent ones occurring in this section,
we use the notation that any 
map $V_i \to V_j$ occurring in the diagrams is denoted by $e(i,j)$. 
Set $M \set \cMb(z)$ and set $N \set M \otimes_z {\cL_1}_y$.
The vertices occurring in (\ref{cd:scheme2}) are the following.
\begin{align}
V_1 &\set (\Esssf\Esssg\cMb)(x) & V_5 &\set H^{p^{\prime}+q^{\prime}}_{m_x}
  (M \otimes_z({\cL_1}_y\otimes_y{\cL_2}_x)) \notag \\
V_2 &\set (\Esssgf\cMb)(x) & 
V_6 &\set H^p_{m_x}((H^q_{m_y}(N[e-r]))[-q] \otimes_y {\cL_2}_x[d\>]) \notag \\
V_3 &\set H^{p^{\prime}}_{m_x}((\Esssg\cMb)(y)\otimes_y {\cL_2}_x) &
 V_7 &\set H^p_{m_x}(N[e-r] \otimes_y {\cL_2}_x[d\>]) \notag \\
V_4 &\set H^{p^{\prime}}_{m_x}(H^{q^{\prime}}_{m_y}N \otimes_y {\cL_2}_x)
 & V_8 &\set H^{p^{\prime}+q^{\prime}}_{m_x}(N \otimes_y {\cL_2}_x) \notag 
\end{align}
The maps in (\ref{cd:scheme2}) which come from (\ref{cd:scheme1}),
viz., $e(2,1), e(1,3),e(3,4)$ and $e(2,5)$ are, by definition, the ones 
specified in \ref{prop:itloco1}.
The remaining maps in (\ref{cd:scheme2}) are defined as follows.

The map $e(1,6)$ is defined as the composition of the following maps 
\begin{align}
V_1 = (\Esssf\Esssg\cMb)(x) 
 &\iso H^p_{m_x}((\Esssg\cMb)(y)[-q]\otimes_y {\cL_2}_x[d\>]) \notag \\
&\iso H^p_{m_x}((H^q_{m_y}(N[e-r]))[-q] \otimes_y {\cL_2}_x[d\>]) = V_6
 \notag
\end{align}
where the first isomorphism is obtained from the sequence 
({\ref{eq:gl2lo1}a}) to ({\ref{eq:gl2lo1}c}) 
for the functor $\Esssf$ acting on the Cousin complex 
$\Esssg\cMb$ on $\Y$, 
while the second isomorphism is induced from the sequence 
({\ref{eq:gl2lo1}a}) to ({\ref{eq:gl2lo1}d}) for the 
functor $\Esssg$ acting on~$\cMb$.

The map $e(2,7)$ is defined to be the composition of 
\begin{align}
V_2 = (\Esssgf\cMb)(x) &\iso 
   H^p_{m_x}(M[-r] \otimes_z ({\cL_1}_y\otimes_y{\cL_2}_x)[d+e])\notag \\
&\iso H^p_{m_x}(M[-r] \otimes_z ({\cL_1}_y[e]
 \otimes_y{\cL_2}_x[d\>]))\notag \\
&\iso H^p_{m_x}((M\otimes{\cL_1}_y)[e-r] 
 \otimes_y {\cL_2}_x[d\>]) = V_7 \notag 
\end{align}
where the first isomorphism is obtained using the sequence 
({\ref{eq:gl2lo1}a}) to ({\ref{eq:gl2lo1}c})  
while the second and the third maps 
are induced naturally by the isomorphisms
\begin{small}
\[
({\cL_1}_y\otimes_y{\cL_2}_x)[d+e] \iso {\cL_1}_y[e]
 \otimes_y{\cL_2}_x[d\>], \qquad 
M[-r] \otimes_z {\cL_1}_y[e] \iso (M\otimes{\cL_1}_y)[e-r] 
\]
\end{small}%
obtained using \S\ref{subsec:conv}, \eqref{conv5}.

We define $e(6,7)$ using (\ref{eq:proofstwo1}) with the following 
notation. Let $B \to C$ denote the natural local homomorphism  
$\OYy \to \OXx$. Let $\Nb = N[e-r]$, $\Lb = {\cL_2}_x[d\>]$ and~$j =q$. The 
hypothesis of (\ref{eq:proofstwo1}),
which in effect is the hypothesis of~(\ref{eq:proofstwo0}),
is satisfied because for any integer $l$, there are isomorphisms  
\[
H^l_{m_B}\Nb \set H^l_{m_y}(N[e-r]) = 
H^{l+e-r}_{m_y}(M\otimes_z{\cL_1}_y) 
\xrightarrow[\ref{rem:gl2lo2}]{\sim} 
H^{l}_y\R\iGp{\Y}(g^*\cMb \otimes_{\Y} \cL_1[e\>]),
\]
so that by~\ref{lem:quasim1}, these modules vanish for $l \ne q$.

The map $e(6,4)$ is obtained as the composition of the maps
\begin{align}
V_6 = H^p_{m_x}\Bigl((H^q_{m_y}(N[e-r]))[-q] \otimes_y {\cL_2}_x[d\>]\Bigr)
  &\iso H^{p}_{m_x}((H^{q^{\prime}}_{m_y}N)[-q] \otimes_y {\cL_2}_x[d\>])
   \notag \\
&\iso H^{p^{\prime}}_{m_x}(H^{q^{\prime}}_{m_y}N \otimes_y {\cL_2}_x)
   = V_4 \notag 
\end{align}
where the first isomorphism involves \S\ref{subsec:conv}, \eqref{conv8}, 
as in~({\ref{eq:gl2lo1}e}), while the
second one is an imitation of~({\ref{eq:gl2lo1}d})-({\ref{eq:gl2lo1}e}). 

For $e(7,8)$ we follow ({\ref{eq:gl2lo1}d})-({\ref{eq:gl2lo1}e}),
while $e(4,8)$ is gotten using (\ref{eq:itloco1}).
The canonical associativity 
map for tensor products gives $e(8,5)$.

Note that $e(8,5)e(4,8) = e(4,5)$, where the latter is an arrow in 
(\ref{cd:scheme1}). Therefore, to prove that (\ref{cd:scheme1})
commutes it suffices to show that 
the four subdiagrams in~(\ref{cd:scheme2}) commute. 
Of these four, the ones on the extreme left and extreme right may be 
written as follows.
\begin{equation}
\label{cd:scheme3}
\begin{CD}
V_3  @<<< V_1 \\
@VVV  @VVV  \\
V_4 @<<< V_6
\end{CD}
\quad \quad \quad \quad
\begin{CD}
V_2  @>>> V_5 \\
@VVV  @AAA  \\
V_7 @>>> V_8
\end{CD}
\end{equation}
Commutativity of these two is easily verified; we do so 
in this subsection. Commutativity of the remaining two subdiagrams of 
(\ref{cd:scheme2}) is proved in the next subsection.\looseness=-1

The commutativity of the rectangle on the left in~(\ref{cd:scheme3}) 
follows from the
commutativity of the following diagram whose maps are obtained using the 
various isomorphisms ({\ref{eq:gl2lo1}a}) to ({\ref{eq:gl2lo1}e}).
\[
\begin{CD}
V_1 @= V_1 \\[-2pt]
@V{\text{({\ref{eq:gl2lo1}a})-({\ref{eq:gl2lo1}e})}}VV  
  @V{\text{({\ref{eq:gl2lo1}a})-({\ref{eq:gl2lo1}c})}}VV  \\[-3pt]
V_3 @<{\text{({\ref{eq:gl2lo1}d})-({\ref{eq:gl2lo1}e})}}<< 
  H^p_{m_x}\Bigl((\Esssg\cMb)(y)[-q] \otimes_y {\cL_2}_x[d\>]\Bigr) \\[-1pt]
@V{\text{({\ref{eq:gl2lo1}a})-({\ref{eq:gl2lo1}d})}}VV  
   @V{\text{({\ref{eq:gl2lo1}a})-({\ref{eq:gl2lo1}d})}}VV  \\[-2pt]
H^{p^{\prime}}_{m_x}\left((H^q_{m_y}(N[e-r])) \otimes_y {\cL_2}_x\right)
  @<{\text{({\ref{eq:gl2lo1}d})-({\ref{eq:gl2lo1}e})}}<< V_6  \\[-2pt]
@V{\text{({\ref{eq:gl2lo1}e})}}VV  @V{\text{({\ref{eq:gl2lo1}e})}}VV  \\[-2pt]
V_4 @<{\text{({\ref{eq:gl2lo1}d})-({\ref{eq:gl2lo1}e})}}<< 
  H^{p}_{m_x}\left((H^{q^{\prime}}_{m_y}N)[-q] \otimes_y {\cL_2}_x[d\>]\right)
\end{CD}
\]

\enlargethispage*{0pt}

For the rectangle on the right in~(\ref{cd:scheme3}),
we expand  vertically using the definition of the 
vertical maps~$e(2,7)$ and~$e(8,5)$. 
Since both the maps $e(2,7)$ and $e(2,5)$ start from $V_2$ and
factor through the steps {({\ref{eq:gl2lo1}a})-({\ref{eq:gl2lo1}c})},
we may replace $V_2$ by the target of~({\ref{eq:gl2lo1}c}).   
The resulting diagram,
shown below, is seen to be commutative from the 
sign chase indicated along the inner 
sides of the arrows.
\[
\begin{CD}
H^p_{m_x}(M[-r] \otimes_z ({\cL_1}_y\otimes_y{\cL_2}_x)[d+e])
 @>{\text{({\ref{eq:gl2lo1}d})-({\ref{eq:gl2lo1}e})}}>{(-1)^{(d+e)r}}> 
  H^{p^{\prime}+q^{\prime}}_{m_x}
  (M \otimes_z({\cL_1}_y\otimes_y{\cL_2}_x)) \\[-3pt]
@VV{(-1)^{de}}V  @AAA \\[-1pt]
H^p_{m_x}(M[-r] \otimes_z ({\cL_1}_y[e]\otimes_y{\cL_2}_x[d\>]))
@. H^{p^{\prime}+q^{\prime}}_{m_x}(N \otimes_y {\cL_2}_x) \\[-1pt]
@VV{(-1)^{re}}V  @|  \\[-3pt]
H^p_{m_x}(N[e-r] \otimes_y {\cL_2}_x[d\>]) 
  @>{(-1)^{(e-r)d}}>{\text{({\ref{eq:gl2lo1}d})-({\ref{eq:gl2lo1}e})}}> 
  H^{p^{\prime}+q^{\prime}}_{m_x}(N \otimes_y {\cL_2}_x)
\end{CD}
\]



\subsection{Conclusion of proof of \ref{prop:itloco1}}
\label{subsec:conc} 
In the previous subsection, the commutativity of the two rectangles in 
(\ref{cd:scheme3}) was shown. Now we verify the commutativity of 
the remaining two subdiagrams in~(\ref{cd:scheme2}). 
  
First we consider the following rectangle in (\ref{cd:scheme2}).
\begin{equation}
\label{cd:conc0}
\begin{CD}
V_6  @>>> V_7 \\
@VVV  @VVV \\
V_4  @>>> V_8
\end{CD}
\end{equation} 
We claim that the isomorphism $e(7,8)e(6,7)e(6,4)^{-1} \colon V_4 \to V_8$ 
is the same as the isomorphism $\theta_{a,b}$ of (\ref{eq:proofstwo1b}) 
for the choices
$a \set e-r$, $b \set d$ along with the notation $p'= p+d-q$, $q'= q+e-r$. 
Indeed, from the definition of the maps involved 
we see that $e(6,4)^{-1}$ corresponds to 
(\ref{eq:proofstwo1b}a)-(\ref{eq:proofstwo1b}b), 
$e(6,7)$ corresponds to (\ref{eq:proofstwo1b}c) and~$e(7,8)$ 
corresponds to (\ref{eq:proofstwo1b}d). Now the commutativity of 
(\ref{cd:conc0}) follows from part~(ii) of Lemma \ref{lem:appsmo-1}.

Thus the only remaining rectangle in (\ref{cd:scheme2}) whose 
commutativity needs to be checked is the following. 
\begin{equation}
\label{cd:conc1}
\begin{CD}
V_1 @<<<  V_2 \\
@VVV  @VVV \\
V_6  @>>>  V_7 \\
\end{CD}
\end{equation}
Before proceeding further we review our terminology. 
Recall that $\cMb$ is assumed to be based at one point only, i.e.,
$\cMb = i_zM[-r]$ for a suitable $\OZz$-module~$M$. 
Also $N = M \otimes_z {\cL_1}_y$.
Set $\cNb \set (g^*\M\otimes_{\Y}\cL_1)[e-r]$  
where~$\M =i_zM$. We shall henceforth make the canonical 
identification
\begin{equation}
\label{eq:conc1a}
\cNb_y = N[e-r].
\end{equation}
Also note that with $\cGb \set g^*\cMb \otimes_{\Y} \cL_1[e]$ as 
in~\eqref{eq:itcoz1}, using~(\S\ref{subsec:conv}, \eqref{conv5}),
we obtain an isomorphism 
\begin{equation}
\label{eq:conc3}
\cGb \set g^*\cMb \otimes_{\Y} \cL_1[e] \iso 
(g^*\M \otimes_{\Y} \cL_1)[e-r] =: \cNb.
\end{equation}

Now expand (\ref{cd:conc1})
modulo \eqref{eq:conc1a} as follows, with the maps~$\alpha_i$ described below.
\def\alphconc{\rotatebox{90}{\makebox[0pt]{$\scriptstyle \alpha_1$}}}
\def\vertconca{\rotatebox{-90}{\makebox[0pt]
{$\xrightarrow[\alphconc]{\hspace{5em}} \quad$}}}  
\begin{equation}
\label{cd:conc2}
\begin{CD}
(\Esssf\Esssg\cMb)(x)  
 @<{e(2,1)}<< (\Esssgf\cMb)(x) \\
@.  @AA{\alpha_2 \; = \; e(2,7)^{-1}}A  \\
\vertconca @. H^p_{m_x}(\cNb_y\otimes_y{\cL_2}_x[d\>]) \\
@.  @AA{\alpha_3}A  \\
H^{p}_{m_x}((E_{\DsssY}\R\iGp{\Y}\cNb)_y\otimes_y {\cL_2}_x[d\>])
 @<{\alpha_4}<<  H^p_{m_x}((\R\iGp{\Y}\cNb)_y\otimes_y{\cL_2}_x[d\>]) \\
@A{\alpha_5}AA  @AA{\alpha_6}A \\
H^p_{m_x}((H^q_{m_y}\cNb_y)[-q] \otimes_y {\cL_2}_x[d\>])  
 @<{\alpha_7}<<  
 H^p_{m_x}(\R\iG{m_y}\cNb_y \otimes_y {\cL_2}_x[d\>])  
\end{CD}
\end{equation}
The map $\alpha_1$ is the composition of the following 
isomorphisms
\begin{align}
(\Esssf\Esssg\cMb)(x) 
&\xto{\text{(\ref{eq:gl2lo1}a)}} 
H^{p}_{m_x}(f^*\Esssg\cMb \otimes_{\X} {\cL_2}[d\>])_x \notag \\
&\xrightarrow[\ref{rem:gl2lo3}]{\text{({\ref{eq:gl2lo1}b$'$})}}
H^{p}_{m_x}((\Esssg\cMb)_y\otimes_y {\cL_2}_x[d\>]) \notag \\
&\xto{\qquad} H^{p}_{m_x}((E_{\DsssY}\R\iGp{\Y}\cNb)_y\otimes_y 
{\cL_2}_x[d\>]) \notag
\end{align}
where the last map is the obvious one induced by \eqref{eq:conc3}.
The map $\alpha_2$ is the inverse of the isomorphism
$e(2,7)$ via (\ref{eq:conc1a}) while $\alpha_3$ is the obvious 
natural map. The map~$\alpha_4$ is obtained\vspace{1pt} by 
applying the isomorphism of functors\vspace{2pt}
\smash{$QE_{\DsssY} \cong 1_{\D^+(\Y)_{\text{CM}}}$} of~\ref{cor:coz4a} 
on the Cohen-Macaulay complex $\R\iGp{\Y}\cNb$.
The map $\alpha_5$ is induced by the composition
\[
(H^q_{m_y}\cNb_y)[-q] 
\xrightarrow[\eqref{eq:coz3}]{\sim} 
(E_{\DsssY}\R\iGp{\Y}\cNb)(y)[-q]
\xrightarrow[\text{cf.~\ref{rem:gl2lo3}}]{\text{truncation}} 
(E_{\DsssY}\R\iGp{\Y}\cNb)_y.
\]
The map $\alpha_6$ is the composition
\[
\R\iG{m_y}\cNb_y \xto{\; \ref{lem:loc6}\text{(iii)} \;} (\R\iG{\I}\cNb)_y
\xto{\text{canonical}} (\R\iGp{\Y}\cNb)_y
\] 
where $\I$ is an open 
coherent ideal in~$\OY$ defining~$\ov{\{y\}}$.
Finally,~$\alpha_7$ is the isomorphism induced by the truncation maps 
in~(\ref{eq:proofstwo0}).

\enlargethispage*{3pt}

Now we verify that the outer border of~(\ref{cd:conc2}) 
is equivalent to~(\ref{cd:conc1}),
which is the same as verifying the relations, 
$\alpha_5^{-1}\alpha_1 = e(1,6)$ and 
$\alpha_3\alpha_6\alpha_7^{-1} = e(6,7)$. The latter 
follows from the definitions involved while  
the former is a consequence of
the commutativity of the following diagram. Here the
left column represents $\alpha_5^{-1}\alpha_1$
while the right column represents $e(1,6)$ modulo \eqref{eq:conc1a}.
Also, recall that $\cMb$ is assumed to be based only at $z$.\vspace{1pt}
\[
\begin{CD}
(\Esssf\Esssg\cMb)(x)  @= (\Esssf\Esssg\cMb)(x) \\[-1pt]
@V{\text{(\ref{eq:gl2lo1}a) and ({\ref{eq:gl2lo1}b$'$})}}VV  
  @V{\text{(\ref{eq:gl2lo1}a)-(\ref{eq:gl2lo1}c)}}VV \\[-3.5pt]
H^{p}_{m_x}((\Esssg\cMb)_y\otimes_y {\cL_2}_x[d\>]) 
  @<{(\text{\ref{eq:gl2lo1}c}')^{-1}}<< 
  H^{p}_{m_x}((\Esssg\cMb)(y)[-q] \otimes_y {\cL_2}_x[d\>]) \\
@V{\text{induced by \eqref{eq:conc3}}}VV   
 @V{\text{induced by \eqref{eq:conc3}}}VV \\[-1pt]
H^{p}_{m_x}((E_{\DsssY}\R\iGp{\Y}\cNb)_y\otimes_y {\cL_2}_x[d\>]) 
@<{\text{via}}<{\text{truncation}}<
  H^{p}_{m_x}((E_{\DsssY}\R\iGp{\Y}\cNb)(y)[-q] \otimes_y {\cL_2}_x[d\>]) \\[-2pt]
@A{\alpha_5}AA   @A{\text{induced by \eqref{eq:coz3}}}A{\text{ 
 cf.~(\ref{eq:gl2lo1}a)}}A \\
H^p_{m_x}((H^q_{m_y}\cNb_y)[-q] \otimes_y {\cL_2}_x[d\>]) @=
  H^p_{m_x}((H^q_{m_y}\cNb_y)[-q] \otimes_y {\cL_2}_x[d\>])
\end{CD}
\]

To prove commutativity of the upper subdiagram in (\ref{cd:conc2}) 
we expand it horizontally using the definition of~$e(2,1)$
and vertically using the definition of~$\alpha_1$.
For convenience we break the diagram into two parts, see~\ref{cd:conc4},
\ref{cd:conc5} below. 
\begin{figure}
\rotatebox{90}
{
\begin{minipage}{8.5in}\vspace{-3mm}
\begin{equation}
\label{cd:conc4}
\stepcounter{equation}\tag*{(\theequation)}
\quad
\begin{CD}
(\Esssf\Esssg\cMb)(x) @<{(\ref{eq:itcoz1}d\>)}<< 
 (\Esssf\R\iGp{\Y}\cGb)(x) @>{(\ref{eq:itcoz1}c)}>>
 (\Esssf\cGb)(x) \\
@V{(\ref{eq:gl2lo1}a)}V{\text{see \eqref{eq:coz3}}}V @VV{\eqref{eq:coz3}}V 
 @VV{\eqref{eq:coz3}}V \\
H^{p}_{m_x}(f^*\Esssg\cMb \otimes_{\X} {\cL_2}[d\>])_x @<<<
 H^{p}_{m_x}(f^*\R\iGp{\Y}\cGb \otimes_{\X} {\cL_2}[d\>])_x 
 @>>> H^{p}_{m_x}(f^*\cGb \otimes_{\X} {\cL_2}[d\>])_x \\
@V{({\ref{eq:gl2lo1}b'})}V{\text{= localization}}V 
@VV{\text{localization}}V @VV{\text{localization}}V \\
H^{p}_{m_x}((\Esssg\cMb)_y\otimes_y {\cL_2}_x[d\>]) @<<< 
 H^{p}_{m_x}((\R\iGp{\Y}\cGb)_y\otimes_y {\cL_2}_x[d\>]) 
 @>>> H^{p}_{m_x}(\cGb_y\otimes_y {\cL_2}_x[d\>])\\
@VV{\text{via \eqref{eq:conc3}}}V @VV{\text{via \eqref{eq:conc3}}}V 
 @VV{\text{via \eqref{eq:conc3}}}V \\
H^{p}_{m_x}((E_{\DsssY}\R\iGp{\Y}\cNb)_y\otimes_y {\cL_2}_x[d\>]) @<{\alpha_4}<<
 H^{p}_{m_x}((\R\iGp{\Y}\cNb)_y\otimes_y {\cL_2}_x[d\>]) @>{\alpha_3}>>
 H^{p}_{m_x}(\cNb_y\otimes_y {\cL_2}_x[d\>]) 
\end{CD}
\end{equation}
\bigskip \bigskip
\bigskip \bigskip
\bigskip \bigskip
\begin{equation}
\label{cd:conc5}
\stepcounter{equation}\tag*{(\theequation)}
\mkern16mu
\begin{CD}
(\Esssf\cGb)(x) @>{(\ref{eq:itcoz1}b)}>> 
 (E'(f^*g^*\cMb \otimes_{\X} 
 (f^*\cL_1[e] \otimes_{\X} \cL_2[d\>])))(x)  
 @>{(\ref{eq:itcoz1}a)}>> (\Esssgf\cMb)(x)\\
@VV{\eqref{eq:coz3}}V  @VV{\eqref{eq:coz3}}V 
 @V{\text{see (\ref{eq:gl2lo1}a)}}V{\eqref{eq:coz3}}V  \\
H^{p}_{m_x}(f^*\cGb \otimes_{\X} {\cL_2}[d\>])_x @>>>
 H^{p}_{m_x}(f^*g^*\cMb \otimes_{\X} (f^*\cL_1[e]\otimes_{\X} {\cL_2}[d\>]))_x 
 @>>{(-1)^{de}}> 
 H^{p}_{m_x}((gf)^*\cMb \otimes_{\X} (f^*\cL_1\otimes_{\X} {\cL_2})[d+e])_x \\
@V\mkern161muV{\text{localization}}V  @VV{\text{localization}}V 
 @VV{\text{localization}}V  \\
H^{p}_{m_x}(\cGb_y\otimes_y {\cL_2}_x[d\>])  @. 
 H^{p}_{m_x}((g^*\cMb)_y\otimes_y ({\cL_1}_y[e]\otimes_y{\cL_2}_x[d\>]))  
 @. H^{p}_{m_x}(M[-r] \otimes_z ({\cL_1}_y\otimes_y{\cL_2}_x)[d+e]) \\
@V{\text{via \eqref{eq:conc3}}}V{(-1)^{re}}V  
 @VV{(-1)^{re}}V @VV{(-1)^{de}}V  \\
H^{p}_{m_x}(\cNb_y\otimes_y {\cL_2}_x[d\>]) @>{\eqref{eq:conc1a}}>>
 H^{p}_{m_x}((M \otimes_z {\cL_1}_y)[e-r] \otimes_y {\cL_2}_x[d\>]) 
 @>{(-1)^{re}}>>
 H^{p}_{m_x}(M[-r] \otimes_z {\cL_1}_y[e]\otimes_y{\cL_2}_x[d\>]) 
\end{CD}
\end{equation}
\end{minipage}
\hspace{-14.5mm}
}
\end{figure}

The leftmost column of vertical maps in \ref{cd:conc4} 
defines $\alpha_1$ of~\eqref{cd:conc2}. The remaining 
columns follow the same pattern of definition. 
The horizontal maps on the left are induced 
by the isomorphism $QE_{\DsssX} \cong 1$ 
of \ref{cor:coz4a} applied to $\R\iGp{\Y}\cGb, \R\iGp{\Y}\cNb$,  
while the ones on the right are the
canonical ones induced by $\R\iGp{\Y} \to 1$. In particular 
the horizontal maps in the bottommost row are $\alpha_3, \alpha_4$ 
of~\eqref{cd:conc2}. It is obvious that \ref{cd:conc4} commutes.
 
The rightmost column of \ref{cd:conc4} is the same as the 
leftmost column of~\ref{cd:conc5}. We  
use $E' \set E_{\DsssX}\R\iGp{\X}$ for convenience.
The remaining maps in~\ref{cd:conc5} are the obvious ones 
as indicated by the labels. (The minus signs refer to 
the convention in \S\ref{subsec:conv}\eqref{conv5}.)
Commutativity of~\ref{cd:conc5} follows easily. Moreover, 
traveling along the bottom row of~\ref{cd:conc5} and 
then its right column gives the map $\alpha_2$ 
of~\eqref{cd:conc2}. In particular, the upper 
subdiagram in \eqref{cd:conc2} is the same as 
\ref{cd:conc4} + \ref{cd:conc5}, and hence commutes.

For the lower rectangle in (\ref{cd:conc2}), 
upon ``canceling off'' 
$H^{p}_{m_x}$ and~$\otimes_y {\cL_2}_x[d\>]$ from each object, we reduce
to showing that the outer border of the
following diagram (of derived-category maps) commutes.
\[
\begin{CD}
(E_{\DsssY}\R\iGp{\Y}\cNb)_y  @<{QE_{\DsssY} \cong 1}<< (\R\iGp{\Y}\cNb)_y \\
@A{\eqref{eq:coz1} \text{ and truncation}}AA  @AA{\text{canonical}}A \\ 
(H^q_y\R\iGp{\Y}\cNb)[-q]  @<{\text{truncation}}<< 
 \R\iG{y}\R\iGp{\Y}\cNb \\
@AA{\eqref{eq:nolab3}}A  @AA{\eqref{eq:nolab3}}A \\
(H^q_{m_y}\cNb_y)[-q]  @<{\text{truncation}}<< \R\iG{m_y}\cNb_y
\end{CD}
\]
By \ref{cor:coz4c}, the upper rectangle commutes while the lower one 
commutes for functorial reasons. This completes 
the proof that (\ref{cd:conc2}), and hence (\ref{cd:conc1}), commutes. 

The commutativity of the diagrams in  
(\ref{cd:conc1}), (\ref{cd:conc0}) and (\ref{cd:scheme3})
proves that all the four subdiagrams of (\ref{cd:scheme2}) commute.
We have therefore shown that~(\ref{cd:scheme1}) commutes, thereby
proving Proposition \ref{prop:itloco1}.

\newpage

\section{Closed immersions and base change}
\label{sec:cibc}


Having tackled the case of smooth maps in \S\ref{sec:cosmo}, we 
now look at the category of closed immersions. 
For a closed immersion 
$f \colon \X \to \Y$ in $\mathbb F$, the functor~$f^{\flat}$ 
defined in~\S\ref{subsec:climmer} below provides 
the other concrete formula for Cousin complexes.
As in~\S\ref{sec:cosmo}, after relating~$f^{\flat}$ 
to its punctual version, we describe a comparison 
isomorphism for the case of composition of two closed 
immersions, and give the corresponding punctual 
description. The proofs are straightforward in this case.
In~\S\ref{subsec:cartsq} we take up the situation of a 
fibered product of a smooth map and a closed immersion. 
The main result there is Proposition \ref{prop:cartsq9}.  

\subsection{Closed Immersions}
\label{subsec:climmer}

Let $h \colon (\X,\Delta_1) \to (\Y,\Delta)$ be a closed immersion 
in~$\bbFc$. For any complex $\cGb$ on $\Y$ we define 
a complex $h^{\flat}\cGb$ as follows
\[
h^{\flat}\cGb \set h^{-1}\sHom_{\OY}(h_{\ast}\OX, \cGb).
\]
\index{   $(-)^\flat$ (Cousin-complex functor associated to a closed immersion $(-)\>$)}
\smallskip

\begin{aprop}
\label{prop:clim1}
Let notation and conditions be as above.
\begin{enumerate}
\item The functor $h^{\flat}$ takes Cousin complexes to Cousin complexes
and furthermore $h^{\flat}\Coz_{\Delta}(\Y) \subset \Coz_{\Delta_1}(\X)$.
\item The functor $h_*$ maps $\Coz_{\Delta_1}(\X)$ to 
$\Coz_{\Delta}(\Y)$ and furthermore there are natural
isomorphisms $h^{\flat}h_* \iso h^{-1}h_* = 1_{\bC}$ where 
$\bC = \Coz_{\Delta_1}(\X)$.
\end{enumerate}
\end{aprop}
\begin{proof}
Let $\cMb$ be a $\Delta$-Cousin complex on $\Y$. Fix $p \in \mathbb Z$.
By definition, 
there is a natural decomposition 
$\M^p \; \cong \; \oplus_y i_yM_y$ where $y$ ranges over points in~$\Y$ 
such that $\Delta(y)=p$ and $M_y$ is an $\OYy$-module. 
Since $\sHom_{\OY}(h_{\ast}\OX, -)$ 
commutes with direct sums (because $h_{\ast}\OX$ is a coherent 
$\OY$-module) therefore there is a canonical decomposition
\[ 
(h^{\flat}\cMb)^p \; \cong \; \oplus_y 
\; h^{-1}\sHom_{\OY}(h_{\ast}\OX, i_yM_y). 
\]
It suffices to consider only those $y \in \Y$ which have 
a preimage in $\X$. Note that for any sheaf~$\F$ 
on~$\Y$, there is a natural isomorphism
\[
\sHom_{\OY}(\F, i_yM_y) \iso i_y\Hom_{\OYy}(\F_y, M_y)
\]
because, for any open neighborhood $\V$ of~$y$, there 
is a canonical isomorphism
\[
\Hom_{\OV}(\F\big|_{\V}, i_yM_y\big|_{\V}) \iso 
\Hom_{\OYy}(\F_y, M_y).
\]
Therefore, for $x \in \X$ and $y = h(x)$ there are natural isomorphisms
\begin{align} 
h^{-1}\sHom_{\OY}(h_*\OX, i_yM_y)  
&\iso h^{-1}i_y\Hom_{\OYy}((h_*\OX)_y, M_y) \notag \\
&\iso i_x\Hom_{\OYy}(\OXx, M_y). \notag
\end{align}
Thus $h^{\flat}\cMb$ is a $\Delta_1$-Cousin complex on $\X$
with $(h^{\flat}\cMb)(x) \iso \Hom_{\OYy}(\OXx, M_y)$ as 
$\Delta_1$, by definition, is the restriction of~$\Delta$ to~$\X$.
If~$M_y$ is a zero-dimensional $\OYy$-module then
$\Hom_{\OYy}(\OXx, M_y)$ is a zero-dimensional $\OXx$-module.
Therefore, by~\ref{lem:mod6}, 
$h^{\flat}\Coz_{\Delta}(\Y) \subset \Coz_{\Delta_1}(\X)$.

With $x,y$ as before, for any $\OXx$-module
$N$, there is a natural isomorphism $h_*i_xN \cong i_yN$.
Since $h_*$ commutes with direct sums it follows that $h_*$ 
takes Cousin complexes to Cousin complexes. Furthermore,
if $N$ is a zero-dimensional as an $\OXx$-module then 
it is also zero-dimensional as an $\OYy$-module.  
By~\ref{lem:mod6} we see that~$h_*$ maps~$\Coz_{\Delta_1}(\X)$ 
to~$\Coz_{\Delta}(\Y)$. The remaining assertions in~(ii) 
hold more generally for arbitrary complexes on~$\X$.
\end{proof}

\medskip

From the proof of \ref{prop:clim1} we obtain 
a canonical isomorphism
\begin{equation}
\label{eq:clim1a}
(h^{\flat}\cMb)(x) \iso \Hom_{\OYy}(\OXx, \cMb(y))
\end{equation}
which can also be described as follows. Set $p = \Delta_1(x)$.
Then~\eqref{eq:clim1a} equals the map in degree $p$ 
induced by the following map of complexes
\[
(h^{\flat}\cMb)_x = (\sHom_{\OY}(h_{\ast}\OX, \cMb))_y 
\iso \Hom_{{\OYy}}(\OXx, \cMb_y). \quad \text{(cf.~\ref{rem:gl2lo3})}
\]

Let $(\X,\DsssX) \xto{\; f \;} (\Y,\DsssY) \xto{\; g \;} (\Z,\DsssZ)$
be closed immersions in $\bbFc$. Let $\cGb$ be a complex on~$\Z$.
Then there is a natural isomorphism
\begin{equation}
\label{eq:clim2}
f^{\flat}g^{\flat}\cGb \xto{\quad} (gf)^{\flat}\cGb 
\end{equation}
defined as follows. Let $\cI, \cJ$ denote the ideals in $\OZ$ 
that define~$\X$ and~$\Y$ respectively. We make the identification
$g_*\OY = \OZ/\cJ$, $(gf)_*\OX = \OZ/\cI$, $f_*\OX = \OY/\cI\OY$.
We define \eqref{eq:clim2}
by the following sequence of obvious natural maps
\begin{align}
f^{\flat}g^{\flat}\cGb &= 
 f^{-1}\sHom_{\Y}(\OY/\cI\OY, g^{-1}\sHom_{\Z}(\OZ/\cJ, \cGb)) \notag \\
&\iso f^{-1}g^{-1}\sHom_{\Z}(\OZ/\cI,\sHom_{\Z}(\OZ/\cJ, \cGb)) \notag \\
&\iso f^{-1}g^{-1}\sHom_{\Z}(\OZ/\cI, \cGb) = (gf)^{\flat}\cGb. \notag 
\end{align}

Let $x$ be a point in $\X$. Set $y = f(x), z = g(y)$.  Let $A,B,C$ 
denote the local rings $\OZz,\OYy,\OXx$ respectively. 
Let~$\cMb$ be a complex in~$\Coz_{\DsssZ}(\Z)$. Set $M = \cMb(z)$.

\def\vertclim{\rotatebox{-90}{\makebox[0pt]{$\xrightarrow{\qquad
{\rotatebox{90}{\textup{\scriptsize{using \eqref{eq:clim1a}}}}}\qquad}$}}}

\begin{aprop}
\label{prop:clim3}
Under the above conditions, the following diagram commutes. 
\[
\begin{CD}
(f^{\flat}g^{\flat}\cMb)(x) @>{(\ref{eq:clim2})}>> ((gf)^{\flat}\cMb)(x) \\
@V{\textup{using \eqref{eq:clim1a}}}VV @. \\
\Hom_B(C,(g^{\flat}\cMb)(y)) @. \vertclim \\
@V{\textup{induced by \eqref{eq:clim1a}}}VV @. \\
\Hom_B(C, \Hom_A(B,M)) @>{\eqref{eq:itloco3}}>> \Hom_A(C, M)
\end{CD}
\]
\end{aprop}
\begin{proof}
Set $p = \DsssX(x)$. Then the diagram in question is obtained from the
degree $p$ part of the following diagram of complexes whose vertical
maps are given by the usual identification of stalks.
{\small{
\[
\begin{CD}
(f^{-1}\sHom_{\Y}(\OY/\cI\OY, g^{-1}\sHom_{\Z}(\OZ/\cJ, \cMb)))_x
@>{(\ref{eq:clim2})_x}>> (f^{-1}g^{-1}\sHom_{\Z}(\OZ/\cI, \cMb))_x \\
@V{\text{canonical}}VV @V{\text{canonical}}VV \\
\Hom_B(C, \Hom_A(B,\cMb_z)) @>{\eqref{eq:itloco3}}>> \Hom_A(C,\cMb_z)
\end{CD}
\]
}}%
Commutativity of the above diagram follows easily.
\end{proof}

\medskip

\subsection{A Fibered Product}
\label{subsec:cartsq}

In this subsection we are concerned with the situation of
the following diagram in~$\bbFc$ which is cartesian, so that 
$\W= \X\times_{\Y}\Z$,
\begin{equation}
\label{eq:cartsq}
\begin{CD}
(\W,\DsssW) @>g>> (\Z,\DsssZ) \\
@VjVV    @VViV \\
(\X,\DsssX) @>f>> (\Y,\DsssY)
\end{CD}
\end{equation}
and where $f,g$ are smooth maps of constant relative 
dimension $d$ and~$i,j$ are closed immersions. 
Set $E \set E_{\DsssX}$ and $\ov{E} \set E_{\DsssW}$,
the respective Cousin functors on $\X, \W$.
Let $\cLb$ be a complex on~$\X$.
Set 
\[
\Esssg(-) \set \ov{E}\R\iGp{\W}(g^*(-) \otimes_{\W} j^*\cLb),
\qquad
\Esssf(-) \set E\R\iGp{\X}(f^*(-) \otimes_{\X} \cLb).
\]
We use the following convention. Since $j$ is a closed immersion,
therefore the functor~$j_*$ is exact and hence we use~$j_*$ 
to denote the corresponding derived functor~$\R j_*$ too. 
Also note that~$j_*$ maps~$\Aqc(\W)$ to~$\Aqc(\X)$.

\begin{alem}
\label{lem:cartsq1}
Let $\cJ$ be any coherent
ideal in $\OX$. Set $\cI = \cJ\OW$. Let $Z$ be any closed
subset of~$\X$. Set $\ov{Z} = Z \cap \W$. Then for any
$\cGb \in \D^+(\W)$ the following natural maps are isomorphisms.
\[
j_*\R\iG{\cI}\cGb \xto{\quad} \R\iG{\cJ}j_*\cGb, \qquad 
j_*\R\iG{\ov{Z}}\cGb \xto{\quad} \R\iG{Z}j_*\cGb.
\]
\end{alem}
\begin{proof}
For the first isomorphism, we refer to \cite[5.2.8(d)]{AJL2}. Since 
$j_*$ preserves flasqueness and since flasque sheaves are 
acyclic for the functors $\iG{\ov{Z}}$ and $\iG{Z}$, the second 
assertion holds because
$j_*\iG{\ov{Z}} = \iG{Z}j_*$.
\end{proof}

\smallskip

\begin{alem}
\label{lem:cartsq2}
For any $\cFb \in \D^+(\W)$ there is a natural isomorphism
\[
j_{\ast}\overline{E}\cFb \cong Ej_{\ast}\cFb.
\]
\end{alem}
\begin{proof}
Recall that the Cousin functors $E, \ov{E}$ are given by
the $E^{*,0}_1$ terms of the spectral sequences induced
by the corresponding codimension filtrations on $\X, \W$. We prove the
lemma by showing that either side of the required isomorphism
is constructed from the same filtered complex on $\X$. 
For any subset $Z \subset\X$ we set $\ov{Z} = Z \cap \W$. 
Set
\[
Z^p \set \{ x \in \X \, | \, \Delta(x) = p \}
\qquad \text{and} \qquad
\overline{Z}^p \set Z^p \cap \W.
\]
Let $\cIb$ be an injective resolution of $\cFb$.
By definition, $\ov{E}\cFb$ is constructed using the filtration 
$\{\iG{\Zb{p}}\cIb\}_{p \in \mathbb Z}$ of $\cIb$. Since $j_*$ is exact 
$j_*\ov{E}\cFb$ can be constructed using the filtration 
$\{j_*\iG{\Zb{p}}\cIb\}_{p \in \mathbb Z}$ of $j_*\cIb$.
On the other hand note that $j_*\cIb$ consists of flasque sheaves
(which are acyclic for the
functors $\iG{Z^{p}}$). By exactness of $j_*$, the natural map
$j_*\cFb \to j_*\cIb$ is a quasi-isomorphism and hence a flasque resolution.
Therefore $Ej_*\cFb$ can be computed using the filtration 
$\{\iG{Z^{p}}j_*\cIb\}_{p \in \mathbb Z}$ of $j_*\cIb$. Since for 
any subset $Z \subset\X$ we have $j_*\iG{\ov{Z}} = \iG{Z}j_*$, it
follows that $\{j_*\iG{\Zb{p}}\cIb\}_{p \in \mathbb Z}$ and 
$\{\iG{Z^{p}}j_*\cIb\}_{p \in \mathbb Z}$ define the same filtration of
$j_*\cIb$ and hence the result follows.
\end{proof}

\medskip

For any complex $\cGb$ on $\Z$ there are natural isomorphisms
described below.
\begin{align}
j_*(g^*\cGb \otimes_{\W} j^*\cLb) \xto{\;\alpha\;}
 j_*g^*\cGb \otimes_{\X} j_*j^*\cLb 
 &\xto{\;\beta\;} f^*i_*\cGb \otimes_{\X} j_*j^*\cLb \notag \\
&\xgets{\;\gamma\;} f^*i_*\cGb \otimes_{\X} \cLb \label{eq:cartsq3}
\end{align}
Since $j$ is a closed immersion, $j_*$ distributes
over the tensor product and thus we obtain the 
isomorphism~$\alpha$. For $\beta$ we use the 
base-change isomorphism $j_*g^* \cong f^*i_*$.
The canonical map $\gamma$ is an isomorphism 
because $f^*i_*\cGb$ is supported on $\W$.

Let $\cFb$ be a complex in $\Dqc^+(\W)$. Then there are 
natural isomorphisms  
\begin{equation}
\label{eq:cartsq4}
j_*\ov{E}\R\iGp{\W}\cFb \iso Ej_*\R\iGp{\W}\cFb \iso E\R\iGp{\X}j_*\cFb
\end{equation}
where the first one is obtained using the isomorphism constructed 
in~\ref{lem:cartsq2}
and the second one using~\ref{lem:cartsq1}.
Therefore, for any $\cGb \in \Dqc^+(\Z)$, there are 
natural isomorphisms
\begin{align}
j_*\Esssg\cGb = j_*\ov{E}\R\iGp{\W}(g^*\cGb \otimes_{\W} j^*\cLb)
&\xto{\;(\ref{eq:cartsq4}) \;} E\R\iGp{\X}j_*(g^*\cGb \otimes_{\W} j^*\cLb) 
\notag \\
&\xto{\;(\ref{eq:cartsq3}) \;} E\R\iGp{\X}(f^*i_*\cGb \otimes_{\X} \cLb) 
= \Esssf i_*\cGb. \notag
\end{align}
In particular, applying $j^{-1}$ results in an isomorphism
\begin{equation}
\label{eq:cartsq7}
\Esssg\cGb \iso j^{-1}\Esssf i_*\cGb.
\end{equation}

Let $\cMb$ be a complex in $\Coz_{\DsssY}(\Y)$. Then there is a 
sequence of natural maps 
\begin{equation}
\label{eq:cartsq8}
\Esssg i^{\flat}\cMb \iso j^{-1}\Esssf i_*i^{\flat}\cMb
\gets j^{\flat}\Esssf i_*i^{\flat}\cMb 
\to j^{\flat}\Esssf \cMb,
\end{equation}
where the first map is the isomorphism from (\ref{eq:cartsq7}),
the second map is induced by the natural map $j^{\flat} \to j^{-1}$
(and is readily seen to be an isomorphism)
and the third is induced by the natural map $i_*i^{\flat} \to 1$.

Let $w$ be a point in $\W$ and let $x,y,z$ denote the corresponding
images in $\X,\Y,\Z$. Set $p=\Delta_1(x), q=\Delta_2(y)$. Since 
the maps $i,j$ in (\ref{eq:cartsq}) are closed immersions, 
therefore $\Delta_3(z)=q$ and $\Delta_4(w)=p$. Set $p_1 \set p-q+d$.
By \ref{cor:app4} the relative 
dimension of $\OXx/\OYy$ and $\OZz/\OWw$ is~$p_1$.
Set $M \set \cMb(y)$,
$N \set \Hom_{\OYy}(\OZz, M)$.

\medskip

\begin{aprop}
\label{prop:cartsq9}
Let notation be as above.
Assume that $\cMb \in \Coz_{\Delta_2}(\Y)$ and that $\cLb = \cL[d\>]$
where $\cL$ is a quasi-coherent flat $\OX$-module and
$d$ is the relative dimension of $f,g$ in \eqref{eq:cartsq}.  
Then the natural maps in \eqref{eq:cartsq8} are 
isomorphisms and the following diagram commutes.
\[
\begin{CD}
(\Esssg i^{\flat}\cMb)(w) @>{(\ref{eq:cartsq8})}>> 
(j^{\flat}\Esssf \cMb)(w) \\
@VV{(\ref{eq:gl2lo1})}V  @VV{\eqref{eq:clim1a}}V  \\
H_{m_w}^{p_1}((i^{\flat}\cMb)(z) \otimes_z (j^*\cL)_w)  @.
\Hom_{\OXx}(\OWw, (\Esssf \cMb)(x)) \\
@VV{\ref{prop:clim1}, (ii)}V  @VV{(\ref{eq:gl2lo1})}V \\
H_{m_w}^{p_1}(N \otimes_z (j^*\cL)_w) @>{\eqref{eq:itloco4}}>>
\Hom_{\OXx}(\OWw, H_{m_x}^{p_1}(M\otimes_y \cL_x)) 
\end{CD}
\]
\end{aprop}
\begin{proof}
All the maps involved in the diagram in question are functorial
in~$\cMb$ and hence using truncation arguments, 
we may assume without loss of generality, 
that $\cMb = i_yM[-q]$. We now proceed by expanding the diagram 
of the Proposition horizontally according to the definition 
in~\eqref{eq:cartsq8}.
Set ${j_w}_{\flat}(-) \set \Hom_{\OXx}(\OWw, -)$ and set
\begin{align}
\cFb_1 &= g^*i^{\flat}\cMb \otimes_{\W} j^*\cL[d\>],
&\cFb_2 &= f^*i_*i^{\flat}\cMb \otimes_{\X} \cL[d\>], \notag \\
F_1 &= (i^{\flat}\cMb)(z) \otimes_z (j^*\cL)_w,
&F_2 &= (i_*i^{\flat}\cMb)(y) \otimes_y \cL_x. \notag
\end{align}
From \eqref{eq:cartsq3} we obtain a natural 
isomorphism $j_*\cFb_1 \iso \F_2$.
Also, as $\OXx$-modules, there is a natural isomorphism
$F_1 \iso F_2$.

The expanded version of the diagram of the Proposition  
occurs in~\eqref{cd:cartsq10}. The downward arrows in
the upper portion are obtained 
using~\eqref{eq:gl2lo1} or~{\eqref{eq:clim1a}}. 
The horizontal maps on the top row are 
given by the maps in~\eqref{eq:cartsq8}. The 
remaining maps are the obvious natural ones.

\def\vertcartsqa{\rotatebox{-90}{\makebox[0pt]{$\xrightarrow{\qquad
{\rotatebox{90}{\makebox[0pt]{\scriptsize{\hspace{15em} $\Box_1$}}}}
\qquad}$}}}
\def\vertcartsqb{\rotatebox{-90}{\makebox[0pt]{$\xrightarrow{\qquad
\qquad}$}}}

\begin{figure}
\rotatebox{90}
{
\begin{minipage}{8.5in}\vspace{-3mm}
\begin{equation}
\label{cd:cartsq10}
\begin{CD}
(\Esssg i^{\flat}\cMb)(w) @>>> (j^{-1}\Esssf i_*i^{\flat}\cMb)(w)
 @<<< (j^{\flat}\Esssf i_*i^{\flat}\cMb)(w) @>>>
 (j^{\flat}\Esssf \cMb)(w)  \\
@.  @| @VVV @VVV \\
\vertcartsqa @. (\Esssf i_*i^{\flat}\cMb)(x) @<<< 
 {j_w}_{\flat}((\Esssf i_*i^{\flat}\cMb)(x)) 
 @>>>  {j_w}_{\flat}((\Esssf \cMb)(x)) \\
@. @VVV @VVV @VVV \\
H_{m_w}^{p_1}F_1  @>>> H_{m_x}^{p_1}F_2 
 @<<< {j_w}_{\flat}H_{m_x}^{p_1}F_2 
 @>>> {j_w}_{\flat}H_{m_x}^{p_1}(M \otimes_y \cL_x) \\
@VVV  @VVV @VVV @| \\
H_{m_w}^{p_1}(N \otimes_z (j^*\cL)_w) @>>> 
 H_{m_x}^{p_1}(N \otimes_y \cL_x) @<<< 
 {j_w}_{\flat}H_{m_x}^{p_1}(N \otimes_y \cL_x) @>>> 
 {j_w}_{\flat}H_{m_x}^{p_1}(M \otimes_y \cL_x) 
\end{CD}
\end{equation}
\bigskip
\bigskip
\bigskip
\begin{equation}
\label{cd:cartsq11}
\begin{CD}
 (j_*\ov{E}\R\iGp{\W}\cFb_1)(x) @>{\text{\eqref{eq:cartsq4}}}>> 
 (E\R\iGp{\X}j_*\cFb_1)(x) @>{\text{\eqref{eq:cartsq3}}}>>
 (E\R\iGp{\X}\cFb_2)(x) \\ 
@VV{\text{(\ref{eq:gl2lo1}a)}}V   
 @V{\hspace{-7em}\Box_2}V{\text{(\ref{eq:gl2lo1}a)}}V  
 @VV{\text{(\ref{eq:gl2lo1}a)}}V  \\
H_{m_w}^{p}(g^*(i^{\flat}\cMb) \otimes_{\W} j^*\cL[d\>])_{w}  
 @>>>  
 H_{m_x}^{p}(j_*\cFb_1)_{x} @>>>
 H_{m_x}^{p}(f^*(i_*i^{\flat}\cMb) \otimes_{\X} \cL[d\>])_{x} \\
@VV{\text{(\ref{eq:gl2lo1}b)-(\ref{eq:gl2lo1}c)}}V @.  
 @V{\hspace{-17em}\Box_3}V{\text{(\ref{eq:gl2lo1}b)-(\ref{eq:gl2lo1}c)}}V   \\
H_{m_w}^{p}((i^{\flat}\cMb)(z)[-q] \otimes_{z} (j^*\cL[d\>])_{w})  
 @. \makebox[10em]{$\xto{\hspace{12em}}$}  @.
 H_{m_x}^{p}((i_*i^{\flat}\cMb)(y)[-q] \otimes_{y} (\cL[d\>])_{x}) \\
@VV{\text{(\ref{eq:gl2lo1}d)-(\ref{eq:gl2lo1}e)}}V @.  
 @VV{\text{(\ref{eq:gl2lo1}d)-(\ref{eq:gl2lo1}e)}}V   \\
H_{m_w}^{p_1}F_1  @>>> 
 H_{m_x}^{p_1}F_1 @>>> H_{m_x}^{p_1}F_2 
\end{CD}
\end{equation}
\end{minipage}\hspace{-15mm}
}
\end{figure}
\begin{figure}
\rotatebox{-90}
{
\begin{minipage}{8.5in}\vspace{-3mm}
\begin{equation}\minCDarrowwidth=.3 in
\label{cd:cartsq12}
\stepcounter{equation}\mkern40mu\tag*{(\theequation)}
\quad
\begin{CD}
(j_*\ov{E}\R\iGp{\W}\cFb_1)(x) @>>> (Ej_*\R\iGp{\W}\cFb_1)(x)
 @>>> (E\R\iGp{\X}j_*\cFb_1)(x) \\
@VVV  @VVV @VVV \\
(j_*H^p_{\ov{\!\{\!w\!\}\!}}\R\iGp{\W}\cFb_1)_x @>>> 
 (H^p_{\ov{\!\{\!x\!\}\!}}j_*\R\iGp{\W}\cFb_1)_x
 @>>> (H^p_{\ov{\!\{\!x\!\}\!}}\R\iGp{\X}j_*\cFb_1)_x \\
@AAA  @AAA @AAA \\
(j_*H^p_{\cJ\OW}\R\iGp{\W}\cFb_1)_x @>>> 
 (H^p_{\cJ}j_*\R\iGp{\W}\cFb_1)_x
 @>>> (H^p_{\cJ}\R\iGp{\X}j_*\cFb_1)_x \\
@VVV  @. @VVV \\
(j_*H^p_{\cJ\OW}\cFb_1)_x @. \makebox[0pt]{$\xto{\hspace{15em}}$} 
 @. (H^p_{\cJ}j_*\cFb_1)_x \\
@VVV  @. @VVV \\
H_{m_w}^{p}\cFb_{1_w}  @. \makebox[0pt]{$\xto{\hspace{15em}}$} 
 @. H_{m_x}^{p}(j_*\cFb_1)_{x} 
\end{CD}
\end{equation}
\bigskip\bigskip
\bigskip\bigskip
\begin{equation}\minCDarrowwidth=.3 in
\label{cd:cartsq13}
\stepcounter{equation}\mkern30mu\tag*{(\theequation)}
\begin{CD}
H_{m_w}^{p}(g^*i^{\flat}\cMb \otimes_{\W} j^*\cL[d\>])_{w}
 @>>> H_{m_w}^{p}((g^*i^{\flat}\cMb)_w \otimes_{w} (j^*\cL[d\>])_{w}) 
 @>>> H_{m_w}^{p}((i^{\flat}\cMb)(z)[-q] \otimes_{z} (j^*\cL[d\>])_{w}) \\
@V{\text{natural}}VV @. \\
H_{m_x}^{p}(j_*(g^*i^{\flat}\cMb \otimes_{\W} j^*\cL[d\>]))_{x} @. @. 
\rotatebox{-90}{\makebox[0pt]{$\xto{\hspace{5em}}$}} \\
@V{\text{via $\alpha$ of \eqref{eq:cartsq3}}}VV @. \\
H_{m_x}^{p}(j_*g^*i^{\flat}\cMb \otimes_{\X} j_*j^*\cL[d\>])_{x}
 @>>> H_{m_x}^{p}((j_*g^*i^{\flat}\cMb)_x \otimes_{x} (j_*j^*\cL[d\>])_{x}) 
 @>>> H_{m_x}^{p}((i^{\flat}\cMb)(z)[-q] \otimes_{z} (j^*\cL[d\>])_{w}) \\
@V{\text{via $\beta, \gamma$ of \eqref{eq:cartsq3}}}VV 
 @V{\text{via $\beta, \gamma$ of \eqref{eq:cartsq3}}}VV 
 @V{\hspace{-10em}\Box_4}VV \\
H_{m_x}^{p}(f^*i_*i^{\flat}\cMb \otimes_{\X} \cL[d\>])_{x}
@>>> H_{m_x}^{p}((f^*i_*i^{\flat}\cMb)_x \otimes_{x} (\cL[d\>])_{x}) 
@>>> H_{m_x}^{p}((i_*i^{\flat}\cMb)(y)[-q] \otimes_{y} (\cL[d\>])_{x}) 
\end{CD}
\end{equation}
\end{minipage}
\hspace{-15mm}
}
\end{figure}

Commutativity of all the rectangles in \eqref{cd:cartsq10},
except for the top left one 
marked as $\Box_1$, is straightforward to verify. 
We expand $\Box_1$ in~\eqref{cd:cartsq11}. 

The unlabeled horizontal maps in~\eqref{cd:cartsq11}
are the obvious natural isomorphisms.
It remains to verify commutativity of~$\Box_2$ and ~$\Box_3$,
since commutativity of the other rectangles is evident. 

In~\ref{cd:cartsq12}, we expand~$\Box_2$ 
horizontally using~\eqref{eq:cartsq4} and vertically
using~\eqref{eq:coz3} and~\ref{prop:nolab1} 
with the notation that 
$\cJ$ denotes the largest coherent ideal in $\OX$
defining the closed set~$\ov{\{x\}}$.
Commutativity of~\ref{cd:cartsq12} is straightforward to verify.

An expanded and \emph{transposed} version of $\Box_3$ 
occurs in~\ref{cd:cartsq13}.
For the horizontal maps in~\ref{cd:cartsq13} we don't use
{(\ref{eq:gl2lo1}b)-(\ref{eq:gl2lo1}c)}, but instead invoke 
\ref{rem:gl2lo3} (localize first and then truncate). Only
commutativity of $\Box_4$ needs explanation, the rest being obvious.

Set $\cNb \set i^{\flat}\cMb$. Then 
$\cNb \in \Coz_{\DsssZ}(\Z)$ and $i_*\cNb \in \Coz_{\DsssY}(\Y)$. 
Now commutativity of $\Box_4$, 
where the horizontal maps are induced by truncation, follows 
from the commutativity of the following diagram.
\[
\begin{CD}
(j_*g^*\cNb)_x @>>> (j_*g^*\N^q)_x[-q] @>>>\cNb(z)[-q] \otimes_z \OWw \\
@V{\text{base change}}VV @VVV @VVV \\
(f^*i_*\cNb)_x @>>> (f^*i_*\N^q)_x[-q] @>>> (i_*\cNb)(y)[-q] \otimes_y \OXx 
\end{CD}
\]

We have thus shown that \eqref{cd:cartsq10} commutes.
By \eqref{eq:itloco4}, the maps in the bottom row 
of~\eqref{cd:cartsq10} are isomorphisms. Since the vertical 
maps in~\eqref{cd:cartsq10} are isomorphisms, therefore 
the remaining horizontal maps are also isomorphisms.
This shows that~(\ref{eq:cartsq8}) is an isomorphism. 
\end{proof}

\newpage

\section{The retract case}
\label{sec:retract}
The main result of this section, viz., Proposition \ref{prop:ret1a}
is the final key ingredient, along with the main results
of the preceding two sections, needed in proving the main
theorem of this paper.

The title of this section refers to the situation where 
the identity map on an object $(\Y,\DsssY)$ in $\bbFc$ is 
factored as $(\Y,\DsssY) \xto{i} (\X,\DsssX) \xto{h} (\Y,\DsssY)$ 
such that~$i$ is a closed immersion 
and~$h$ a smooth map of constant relative dimension.
In \ref{prop:ret1a} we show that the
Cousin functor obtained by using \ref{exam:mainre1}
for $h,i$ is isomorphic to the identity functor 
on~$\Coz_{\DsssY}(\Y)$.
Residues of differential forms play a role here. 

After covering some homological preliminaries in 
\S\ref{subsec:homo} below, we prove a general result
in \S\ref{subsec:hosp} which provides a 
partial description of the coboundary
map of the Cousin complex of \ref{exam:mainre1}(i). 
Proposition \ref{prop:ret1a} is then an easy consequence 
modulo the local ingredient involving residue maps.

\subsection{Homological preliminaries}
\label{subsec:homo}
Let $\D$ be a triangulated category with translation functor $T$.
Consider a commutative diagram in~$\D$ as below whose 
rows and columns are assumed to be triangles in $\D$.
\begin{equation}
\label{eq:homo1}
\begin{CD}
A_1 @>{\alpha_1}>> A_2 @>{\alpha_2}>> A_3 @>{\alpha_3}>> TA_1 \\
@V{u_1}VV  @V{u_2}VV  @V{u_3}VV  @V{Tu_1}VV \\
B_1 @>{\beta_1}>> B_2 @>{\beta_2}>> B_3 @>{\beta_3}>> TB_1 \\
@V{v_1}VV  @V{v_2}VV  @V{v_3}VV  @V{Tv_1}VV \\
C_1 @>{\gamma_1}>> C_2 @>{\gamma_2}>> C_3 @>{\gamma_3}>> TC_1 \\
@V{w_1}VV  @V{w_2}VV  @V{w_3}VV  \\
TA_1 @>{T\alpha_1}>> TA_2 @>{T\alpha_2}>> TA_3 
\end{CD}
\end{equation}
Suppose $C_1 \cong 0 \cong A_3$. Then we obtain two ways of defining
a map $C_3 \to TA_1$ as follows. Consider the following diagram.  
\begin{equation}
\label{eq:homo1a}
\begin{CD}
B_3 @>{\beta_3}>> TB_1 \\
@V{v_3}VV @A{Tu_1}AA \\
C_3 @. TA_1 \\
@A{\gamma_2}AA @V{T\alpha_1}VV \\
C_2 @>{w_2}>> TA_2
\end{CD}
\end{equation}
Since $C_1 \cong 0$ hence $\gamma_2$ and $Tu_1$ are isomorphisms.
Since $A_3 \cong 0$ hence~$v_3$ and~$T\alpha_1$ are isomorphisms.
Therefore, inverting these isomorphisms, we obtain one map each 
from the upper and lower half of~(\ref{eq:homo1a}) 
viz., $(Tu_1)^{-1}\beta_3v_3^{-1}$ 
and~$(T\alpha_1)^{-1}w_2\gamma_2^{-1}$. 

We now find a condition under which these two maps coincide.
Complete the map 
$A_1 \xto{\:u_2\alpha_1\:} B_2$ to the following triangle 
(uniquely determined up to isomorphism) 
\begin{equation}
\label{eq:homo1b}
A_1 \xto{\quad} B_2 \xto{\quad} X \xto{\quad} TA_1
\end{equation}
in $\D$ for suitable $X$ and maps 
$B_2 \to X \to TA_1$ in~$\D$. 
Consider the following commutative diagram whose rows are
triangles in $\D$. (This may be imagined 
as comparing the triangles in the second row and second column 
of~(\ref{eq:homo1}) via a triangle constructed along the diagonal.)
\begin{equation}
\label{eq:homo2}
\begin{CD}
B_1 @>{\beta_1}>> B_2 @>{\beta_2}>> B_3 @>{\beta_3}>> TB_1 \\
@A{u_1}AA  @|  @.  @A{Tu_1}AA  \\
A_1 @>{u_2\alpha_1}>> B_2 @>>> X @>>>  TA_1 \\
@V{\alpha_1}VV  @|  @.  @V{T\alpha_1}VV  \\
A_2 @>{u_2}>> B_2 @>{v_2}>> C_2 @>{w_2}>>  TA_2
\end{CD}
\end{equation}
From the defining axioms of triangles 
we conclude that there exist maps (not uniquely determined) 
\[
X \xto{\;u^{\prime}\;} B_3 \quad \quad \text{and} \quad \quad 
X \xto{\;{\alpha}^{\prime}\;} C_2
\]
which, when filled in (\ref{eq:homo2}), make it commute.
(Note that $u^{\prime},{\alpha}^{\prime}$ are isomorphisms
since the remaining vertical maps in~(\ref{eq:homo2}) are isomorphisms.)
We therefore obtain the following diagram. 
\begin{equation}
\label{eq:homo2a}
\begin{CD}
X @>{u^{\prime}}>> B_3 \\
@V{{\alpha}^{\prime}}VV  @V{v_3}VV  \\
C_2 @>{\gamma_2}>> C_3
\end{CD}
\end{equation}

\medskip

\begin{alem}
\label{lem:homo4}
Suppose we have a commutative diagram as in $(\ref{eq:homo1})$
with rows and columns as triangles in $\D$ 
and suppose the following conditions hold.
\begin{enumerate}
\item We have $C_1 \cong 0 \cong A_3$ and hence the 
vertical maps in $(\ref{eq:homo1a})$ are invertible. 
\item There is a choice of maps $u^{\prime},{\alpha}^{\prime},$
(satisfying their defining conditions) for which the diagram 
in~$(\ref{eq:homo2a})$ \emph{commutes}.
\end{enumerate}
Then the diagram obtained by replacing the vertical maps 
in~$(\ref{eq:homo1a})$ by their respective inverses, commutes.
\end{alem}
\begin{proof}
We expand (\ref{eq:homo1a}) in the following way 
\[
\begin{CD}
B_3 @= B_3 @>{\beta_3}>> TB_1 \\
@V{v_3}VV  @A{u^{\prime}}AA  @A{Tu_1}AA \\
C_3 @. X @>>> TA_1 \\
@A{\gamma_2}AA  @V{{\alpha}^{\prime}}VV  @V{T\alpha_1}VV \\
C_2 @= C_2 @>{w_2}>> TA_2
\end{CD}
\]
where the two squares on the right side are obtained
from (\ref{eq:homo2}) and hence commute by definition
of~${u^{\prime}}$ and~${{\alpha}^{\prime}}$, 
whereas the rectangle on the left side is the same as the 
diagram in~(\ref{eq:homo2a}).
Condition (ii) therefore implies that the above diagram commutes.  
\end{proof}

\medskip 

Let $\A$ be an abelian category.
Recall that to any exact sequence of complexes 
$0 \to \Lb \to \Mb \to \Nb \to 0$ in~$\bC(\A)$
we can associate a corresponding induced triangle 
$\Lb \to \Mb \to \Nb \to \Lb[1]$ in~$\D(\A)$ (see, e.g., \cite[Example (1.4.4)]{Li}).
This association is natural, in that, any
map of short exact sequences in~$\bC(\A)$ functorially 
gives rise to a corresponding map of triangles. Recall further:

\begin{alem}
\label{lem:homo4aa}
Let $\A$ be an abelian category with enough injectives.
For any triangle $T$ in $\D^+(\A)$ 
there exists in~$\bC^+(\A)$, an exact sequence of complexes 
of injectives given by
$0 \xto{\quad} \Lb \xto{\quad} \Mb \xto{\quad} \Nb \xto{\quad} 0$,
such that the corresponding triangle 
$\Lb \xto{\quad} \Mb \xto{\quad} \Nb \xto{\quad} \Lb[1]$
is isomorphic to~$T$.
\end{alem}

\begin{proof}
Let $\Nb$, $\Lb$ denote complexes of injectives that are 
$\D^+(\A)$-isomorphic to $\Cb$, $\Ab$ respectively. The natural map 
$\Cb[-1] \to \Ab$ induces a map $\Nb[-1] \to \Lb$ which may also 
be thought of as a map in $\K(\A)$ and hence can be represented
by a map, say $\alpha$, in $\C(\A)$. Let $\Mb$ denote the mapping 
cone of $\alpha$. The standard exact sequence involving $\Lb,\Mb,\Nb$
results in a triangle, and from the defining property of triangles we 
see that there is a map of triangles 
\[
\begin{CD}
\Cb[-1] @>>> \Ab @>>> \Bb @>>> \Cb \\
@VVV  @VVV  @VVV  @VVV \\
\Nb[-1] @>>> \Lb @>>> \Mb @>>> \Nb 
\end{CD}
\]
which is in fact an isomorphism of triangles. 
\end{proof}

\medskip

Let $\A$ be an abelian category with enough injectives. Let 
\[
0 \xto{\quad} \Gamma_1 \xto{\quad} \Gamma_2 \xto{\quad} 
\Gamma_3 \xto{\quad} 0
\]
be a sequence of additive functors from $\A$ to~$\A$ which is 
exact on injectives of $\A$. Then any complex $\Gb$ in 
$\D^+(\A)$ gives rise to a triangle 
\begin{equation}
\label{eq:homo4ab}
\R\Gamma_1\Gb \xto{\text{natural}} 
\R\Gamma_2\Gb \xto{\text{natural}} 
\R\Gamma_3\Gb \xto{\quad ? \quad} \R\Gamma_1\Gb[1] 
\end{equation}
in the following manner. Let $\Gb \to \Ib$ be an injective 
resolution so that we may set $\R\Gamma_i\Gb[n]=\Gamma_i\Ib[n]$. 
Then \eqref{eq:homo4ab} is the triangle associated to the exact sequence 
$0 \xto{\quad} \Gamma_1\Ib \xto{\quad} \Gamma_2\Ib \xto{\quad} 
\Gamma_3\Ib \xto{\quad} 0.$

\medskip

\begin{alem}
\label{lem:homo6}
With $\A$ and $\Gamma_i$ as above, let $T$ be a triangle in~$\D^+(\A)$,
say $\Fbp \xto{\quad} \Gb \xto{\quad} \Fb \xto{\quad} \Fbp[1]$,
such that $\R\Gamma_1\Fb \cong 0 \cong \R\Gamma_3\Fbp$.
Then the two rows in the following diagram
(maps being the obvious natural ones) 
give rise to the same map from~$\R\Gamma_3\Fb$ to~$\R\Gamma_1\Fbp[1]$.
\[
\begin{CD}
\R\Gamma_3\Fb @<{\sim}<< \R\Gamma_3\Gb 
 @>{?}>{\text{from \eqref{eq:homo4ab}}}> 
 \R\Gamma_1\Gb[1] @<{\sim}<< \R\Gamma_1\Fbp[1] \\
@| @. @. @| \\
\R\Gamma_3\Fb @<{\sim}<< \R\Gamma_2\Fb @>>> \R\Gamma_2\Fbp[1] @<{\sim}<< 
 \R\Gamma_1\Fbp[1] 
\end{CD}
\]
\end{alem}
\begin{proof}
By \ref{lem:homo4aa} we may assume, without loss of generality, 
that $\Fb,\Fbp,\Gb$ consist $\A$-injectives and fit into an
exact sequence $0 \to \Fbp \to \Gb \to \Fb \to 0$
such that the corresponding induced triangle is isomorphic to $T$.
For $j= 1,2,3,$ set $A_j \set \Gamma_j\Fbp$, $B_j \set \Gamma_j\Gb$,
$C_j \set \Gamma_j\Fb$.
We therefore obtain the following commutative diagram in which 
the rows and columns are exact sequences.
\begin{equation}
\label{cd:homo7}
\begin{CD}
@. 0 @. 0 @. 0 \\
@. @VVV @VVV @VVV \\
0 @>>> A_1 @>>> A_2 @>>> A_3 @>>> 0 \\
@. @VVV  @VVV  @VVV  \\
0 @>>> B_1 @>>> B_2 @>>> B_3 @>>> 0 \\
@. @VVV  @VVV  @VVV  \\
0 @>>> C_1 @>>> C_2 @>>> C_3 @>>> 0 \\
@. @VVV  @VVV  @VVV  \\
@. 0 @. 0 @. 0 
\end{CD}
\end{equation}
The rows and columns give rise to triangles in $\D(\A)$ and so 
we obtain a commutative 
diagram in $\D = \D(\A)$ as in~(\ref{eq:homo1}).
The assumption $\R\Gamma_1\Fb \cong 0 \cong \R\Gamma_3\Gb$ implies 
that $A_3$ and $C_1$ are isomorphic to~0 in~$\D$
and hence condition (i) of Lemma~\ref{lem:homo4} 
is satisfied. In~(\ref{cd:homo7}) 
the natural map $A_1 \to B_2$ is a monomorphism in~$\bC(\A)$.
Let~$X$ denote its cokernel. It follows that the natural
maps $B_2 \to B_3$ and $B_2 \to C_2$ in (\ref{cd:homo7})
factor through the epimorphism $B_2 \to X$ and so we obtain 
the following induced \emph{commutative} diagram in $\bC(\A)$.
\[
\begin{CD}
X  @>>>  B_3\\
@VVV  @VVV \\
C_2  @>>> C_3
\end{CD}
\] 
From the triangle associated to the exact sequence
$0 \to A_1 \to B_2 \to X \to 0$ 
we see that $X$ fits into a triangle as in~(\ref{eq:homo1b}).
Furthermore the maps $X \to B_3$ and $X \to C_2$ of the 
above diagram give a choice for $u^{\prime},{\alpha}^{\prime}$ 
in~\eqref{eq:homo2a}, i.e., make~\eqref{eq:homo2} commute.
Since for this choice of $u^{\prime},{\alpha}^{\prime}$, 
\eqref{eq:homo2a} commutes therefore condition (ii) of 
Lemma \ref{lem:homo4} is satisfied. 
The desired result now follows from Lemma~\ref{lem:homo4}.
\end{proof}

Let $(\X,\Delta)$ be a formal scheme in $\bbFc$ and let 
$\{Z^m\}_{m \in \mathbb Z}$ be the filtration of~$\X$ induced 
by $\Delta$. Let $\F$ be a flasque sheaf on~$\X$. 
Fix an integer $p$. Let $x_1, \ldots, x_k$ be 
points in~$Z^p \setminus Z^{p+1}$.
Let $W = \ov{\{x_1, \ldots, x_k\}}$ 
be the closure of the set $\{x_1, \ldots, x_k\}$.
We claim that for any $i$, the canonical inclusion 
$\iG{x_i}\F \hookrightarrow (\iG{W}\F)_{x_i}$ is surjective.
To that end, let $U$ be an open neighborhood of $x_i$ and let
$s$ be an element of $(\iG{W}\F)(U)$. Let $Y$ be the 
closed set generated by those generic points of $\Supp(s)$
that do not lie in $\ov{\{x_i\}}$. Then the restriction of
$s$ to the open set $V \set U \setminus Y$ lies 
in $(\iG{\ov{\!\{x_i\}\!}}\F)(V)$. Thus our claim 
follows. 

Now assume that there is a point $x' \in Z^{p+1}$ such that $x'$ 
is an immediate specialization of~$x_i$ for each~$i$.
Set $W' \set W \cap Z^{p+1}$. 
Arguing as in the previous paragraph, we see that
the canonical map $\iG{x'}\F \to (\iG{W'}\F)_{x'}$ is an isomorphism.
Next note that there is a surjective map 
\[
\iG{W}\F \onto \iG{W/W'}\F 
\xrightarrow[\text{\cite[p.~227]{RD}, cf.~\eqref{eq:coz0}}]{\sim}
\bigoplus_{j=1}^{j=k}i_{x_j}(\iG{x_j}\F).
\]
Hence, by localizing at $x'$ (resp.~$x_i$ for any $i$), we obtain
a surjection (resp.~an isomorphism)
\[ 
(\iG{W}\F)_{x'} \onto \oplus_j\iG{x_j}\F, \qquad
\text{(resp.~$(\iG{W}\F)_{x_i} \iso \iG{x_i}\F$)},
\]
where the inverse of the last isomorphism is the one described in the
previous paragraph.

For convenience, we shall use $\iG{l:m}$ to denote 
the functor~$\iG{Z^{l}/Z^{m}}$. Consider the following diagram
where all the unlabeled maps are the canonical ones 
and $a$ is defined below.
{\small{
\begin{equation}
\label{cd:coures-1z}
\begin{CD}
0 @>>> \iG{x^{\prime}}\F @>>> (\iG{W}\F)_{x^{\prime}} 
@>>> \oplus_j\iG{x_j}\F @>>> 0 \\
@. @VVV @VVV @VVaV \\
0 @>>> (\iG{{p+1}:{p+2}}\F)_{x^{\prime}} 
@>>> (\iG{{p}:{p+2}}\F)_{x^{\prime}}
@>>> (\iG{{p}:{p+1}}\F)_{x^{\prime}} @>>> 0 
\end{CD}
\end{equation}
}}%
We define $a$ to be the canonical inclusion induced by the 
decomposition in~\eqref{eq:coz0}: 
\[
\bigoplus_j\iG{x_j}\F \subseteq 
\bigoplus_{\{y|\Delta(y)=p, \,\,y \leadsto x^{\prime}\}} \iG{y}\F =
\Bigl(\bigoplus_{\{y|\Delta(y)=p\}} i_y\iG{y}\F\Bigr)_{x^{\prime}}
\xleftarrow[\;\eqref{eq:coz0}\;]{\sim} (\iG{Z^p/Z^{p+1}}\F)_{x^{\prime}}.
\]

\medskip

\begin{alem}
\label{lem:coures0}
The diagram in $(\ref{cd:coures-1z})$ commutes and its rows are exact.
In particular, for any complex $\cGb \in \D(\X)$, 
there is a commutative diagram 
\[
\begin{CD}
\bigoplus_{j=1}^{j=k} H^p_{x_j}\cGb @>>> H^{p+1}_{x^{\prime}}\cGb \\
@VVV @VVV \\ 
(H^p_{Z^{p}/Z^{p+1}}\cGb)_{x^{\prime}} @>>> 
(H^{p+1}_{Z^{p+1}/Z^{p+2}}\cGb)_{x^{\prime}} 
\end{CD}
\]
where the downward arrows are induced from the corresponding ones
in~\eqref{cd:coures-1z}
and the horizontal maps are the usual connecting homomorphisms resulting
from the exact rows. 
\end{alem}
\begin{proof}
Using the natural maps
$\iG{W'} \to \iG{Z^{p+1}}$, $\iG{W} \to \iG{Z^{p}}$,
$\iG{W/W'} \to \iG{Z^{p}/Z^{p+1}}$,
we expand \eqref{cd:coures-1z} as follows.
{\small{
\[
\begin{CD}
0 @>>> \iG{x^{\prime}}\F @>>> (\iG{W}\F)_{x^{\prime}} 
@>>> \oplus_j\iG{x_j}\F @>>> 0 \\
@. @| @| @AA{\wr}A \\
0 @>>> (\iG{W'}\F)_{x^{\prime}} 
@>>> (\iG{W}\F)_{x^{\prime}}
@>>> (\iG{W/W'}\F)_{x^{\prime}} @>>> 0  \\
@. @VVV @VVV @VVV \\
0 @>>> (\iG{Z^{p+1}/Z^{p+2}}\F)_{x^{\prime}} 
@>>> (\iG{Z^{p}/Z^{p+2}}\F)_{x^{\prime}}
@>>> (\iG{Z^{p}/Z^{p+1}}\F)_{x^{\prime}} @>>> 0 
\end{CD}
\]
}}%
It is clear that the above diagram commutes. The middle and the bottom
rows, being localization of exact sequences, are exact and hence the 
top row is also exact.
\end{proof}

Looking at the image of the vertical maps in the commutative
diagram of~\ref{lem:coures0} we in fact see that  
the following diagram of induced natural maps commutes, 
\begin{equation}
\label{cd:coures1}
\begin{CD}
\bigoplus_{j}H^p_{x_j}\cGb @>>> H^{p+1}_{x^{\prime}}\cGb \\
@V{\eqref{eq:coz1}}V{\wr}V @V{\eqref{eq:coz1}}V{\wr}V \\ 
\bigoplus_{j}(E_{\Delta}\cGb)(x_j) @>>> (E_{\Delta}\cGb)(x')
\end{CD}
\end{equation}
where the bottom row is induced 
by the differential of~$E_{\Delta}\cGb$.

\medskip

\subsection{Coboundary for lateral specializations}
\label{subsec:hosp}
\index{specialization!lateral}\index{lateral specialization}
Let $f \colon \X \to \Y$ be a map in $\mathbb F$. We say that 
a specialization $x \leadsto x^{\prime}$ is 
\emph{$f$-lateral} if it is an immediate specialization 
and its image under $f$ is also an immediate specialization.

Let $h \colon (\X,\DsssX) \to (\Y, \DsssY)$ be a smooth map in $\bbFc$
having constant relative dimension~$n$. 
Let $\cMb$ be a complex in $\Coz(\Y,\DsssY)$, and let $\cL$ be a 
quasi-coherent flat $\OX$-module. Set 
$\mathbb E (\cMb) \;\set\; E_{\DsssX}\R\iGp{\X}(h^*\cMb \otimes \cL[n])
\;\in\; \Coz(\X, \DsssX).$
Our aim is to describe the differential of $\mathbb E (\cMb)$
when restricted to the components corresponding to a fixed $h$-lateral
specialization.

Let $y \leadsto y'$ be an immediate specialization in $\Y$.
Let $x' \in h^{-1}\{y'\}$ and $x_1, \ldots, x_k \in h^{-1}\{y\}$
be such that $x'$ is an immediate specialization of $x_i$ 
for each $i$. 
Set 
\[
 q = \DsssY(y), \; p = \DsssX(x_i),
\; M = \cMb(y), \; M^{\prime}  = \cMb(y^{\prime}). 
\]
Then $\DsssY(y^{\prime}) = q+1$ and $\DsssX(x^{\prime}) = p+1$.  
By \ref{cor:app4}, $p_1\set p-q+n$ is the relative dimension 
of~$\cO_{\X,x_i}$ over $\OYy$ and also of $\cO_{\X,x'}$
over $\cO_{\Y,y'}$. Let $\psi$ denote the following map
of $\cO_{\X,x'}$-modules
\[
\bigoplus_j H^{p_1}_{m_{x_j}}(M \otimes_y \cL_{x_j}) 
\iso \bigoplus_j(\mathbb E \cMb)(x_j)
\xto{\quad} (\mathbb E \cMb)(x^{\prime}) \iso 
H^{p_1}_{m_{x^{\prime}}}(M^{\prime} \otimes_{y^{\prime}} \cL_{x^{\prime}}),
\]
where the first and the last isomorphisms are obtained 
using~(\ref{eq:gl2lo1}) and the 
map in the middle is induced by the differential of~$\mathbb{E} \cMb$. 
Let $\partial \colon M \to M^{\prime}$
denote the natural map of $\cO_{\Y,y'}$-modules 
induced by the differential of~$\cMb$. Our aim is to
express $\psi$ in terms of~$\partial$.   

Set $\cO_{x_j} \set \cO_{\X,x_j}$ and $\Oxp \set \cO_{\X,x'}$. 
Let $W = \ov{\{x_1, \ldots, x_k\}}$ be the closure of the set
$\{x_1, \ldots, x_k\}$. 
Let $\I$ be an open coherent ideal
in $\OX$ such that $\Supp(\OX/\I) = W$. Set $I \set \I_{x'}$. 
Note that $I\cO_{\X,x_j} = \I_{x_j}$ is $m_{x_j}$-primary
for each $j$. Indeed, since $\I$ is open, we may first go modulo 
a defining ideal in $\OX$ so that $\X$ may now be assumed to be an 
ordinary scheme. Then the required conclusion follows easily.
In particular, the canonical map 
$\R\iG{m_{x_j}}(M \otimes_y \cL_{x_j}) 
\to \R\iG{I\cO_{x_j}}(M \otimes_y \cL_{x_j})$ is an isomorphism. 

\medskip

\begin{aprop}
\label{lem:hosp0a}
In the above situation, 
consider the following diagram where $\mu_1,\mu_2,\mu_3$ 
are maps of $\cO_{\X,x'}$-modules defined as follows. The map $\mu_1$ is 
defined on each component via the sequence of natural maps
\[
H^{p_1}_{I}(M \otimes_{y^{\prime}} \cL_{x^{\prime}})
\to H^{p_1}_{I}(M \otimes_y \cL_{x_j}) \iso 
H^{p_1}_{I\cO_{x_j}}(M \otimes_y \cL_{x_j}) 
\osi H^{p_1}_{m_{x_j}}(M \otimes_y \cL_{x_j}),
\] 
$\mu_2$ is the one induced by the canonical inclusion 
$\iG{m_{x^{\prime}}} \to \iG{I}$ and $\mu_3$ 
is $(-1)^n$ times 
$H^{p_1}_{I}(-\otimes_{y^{\prime}} \cL_{x^{\prime}})$
applied to $\partial \colon M \to M^{\prime}$.
\[
\begin{CD}
\bigoplus_jH^{p_1}_{m_{x_j}}(M \otimes_y \cL_{x_j}) @>{\psi}>> 
 H^{p_1}_{m_{x^{\prime}}}
 (M^{\prime} \otimes_{y^{\prime}} \cL_{x^{\prime}}) \\
@A{\mu_1}AA  @V{\mu_2}VV \\
H^{p_1}_{I}(M \otimes_{y^{\prime}} \cL_{x^{\prime}}) @>{\mu_3}>> 
 H^{p_1}_{I}(M^{\prime} \otimes_{y^{\prime}} \cL_{x^{\prime}}) 
\end{CD}
\]
Then $\mu_1,\mu_2$ are isomorphisms and the above diagram commutes.
\end{aprop}

\begin{proof}

Using truncation arguments as in (\ref{eq:gl2lo1}b) we 
see that $\psi$ depends only on the modules $M, M'$, i.e., 
we may assume without loss of generality that $\cMb$ satisfies
$\M^j = 0$ for~$j \ne q,q+1$, and  
$\M^q = i_yM$, $\M^{q+1} = i_yM'$. Set 
\begin{align}
\cFb &\set \R\iGp{\X}(f^*\sigma_{\le q}\cMb \otimes_{\X} \cL[n]), \notag \\
\cFbp &\set \R\iGp{\X}(f^*\sigma_{\ge q+1}\cMb \otimes_{\X} \cL[n]), \notag \\
\cGb &\set \R\iGp{\X}(f^*\cMb \otimes_{\X} \cL[n]). \notag
\end{align}
Applying $\R\iGp{\X}(f^*(-) \otimes_{\X} \cL[n])$ to the triangle
\begin{equation}
\label{eq:hosp1}
\sigma_{\ge q+1}\cMb \xto{\quad} \cMb 
\xto{\quad} \sigma_{\le q}\cMb \xto{\quad} (\sigma_{\ge q+1}\cMb)[1]
\end{equation}
associated to the exact sequence 
$0 \to \sigma_{\ge q+1}\cMb \to \cMb \to \sigma_{\le q}\cMb \to 0$
results in the following triangle in~$\D(\X)$
\begin{equation}
\label{eq:hosp0aa}
\cFbp \xto{\quad} \cGb \xto{\quad} \cFb \xto{\;\gamma\;} \cFbp[1].
\end{equation}
The map $\gamma$ is obtained from the $\delta$-functoriality 
of $\R\iGp{\X}(f^*(-) \otimes_{\X} \cL[n])$ 
(\S\ref{subsec:conv}, \eqref{conv6}), which in turn 
involves the three $\delta$-functors, $\R\iGp{\X}$, 
$f^*$ and $\otimes_{\X} \cL[n]$. The first two of these 
commute with translation by construction, and the third, 
as seen from the convention in \S\ref{subsec:conv}, \eqref{conv4}, also 
commutes with translation. 

Set $F_j \set M \otimes_{y} \cL_{x_j}$, 
$F \set M \otimes_{y^{\prime}} \cL_{x^{\prime}}$, 
and $F' \set M' \otimes_{y^{\prime}} \cL_{x^{\prime}}$.
We now expand the diagram in the assertion of the lemma as follows.
\begin{equation}
\label{cd:hosp0aaa}
\begin{CD}
@. \bigoplus_j(\mathbb E \cMb)(x_j) @>{\mu_0}>> (\mathbb E \cMb)(x^{\prime}) \\
@. @V{\mu_4}VV @V{\mu_5}VV \\
\bigoplus_jH^{p_1}_{m_{x_j}}F_j @<{\mu_6}<< 
 \bigoplus_jH^{p}_{x_j}\cFb @>{\mu_7}>>
 H^{p+1}_{x'}\cFbp @>{\mu_8}>> H^{p_1}_{m_{x'}}F' \\
@A{\mu_1}AA @A{\mu_{9}}AA @V{\mu_{10}}VV  @V{\mu_{2}}VV \\
H^{p_1}_{I}F @<{\mu_{11}}<< (H^{p}_W\cFb)_{x'} @>{\mu_{12}}>> 
 (H^{p+1}_W\cFbp)_{x'} @>{\mu_{13}}>> H^{p_1}_{I}F' 
\end{CD}
\end{equation}
The map $\mu_0$ is the obvious one induced by the differential 
of $\mathbb E \cMb$. We define~$\mu_4$ componentwise 
as the composition of the following natural maps 
\[
(\mathbb E \cMb)(x_j) = (E_{\Delta_1}\cGb)(x_j) 
\xrightarrow[\eqref{eq:coz1}]{\sim} H^p_{x_j}\cGb 
\xrightarrow[\eqref{eq:hosp0aa}]{} H^p_{x_j}\cFb.
\]
Note that this amounts to using 
(\ref{eq:gl2lo1}a) and (\ref{eq:gl2lo1}b) without
involving the isomorphism
$H^p_{x_j}\R\iGp{\X}(-) \iso H^p_{m_{x_j}}(-)_{x_j}$
(cf.~\eqref{eq:coz3}). We define~$\mu_6$ by first replacing
$H^p_{x_j}\R\iGp{\X}(-)$ by $H^p_{m_{x_j}}(-)_{x_j}$ and then following the 
remaining steps in~(\ref{eq:gl2lo1}) so that $\mu_6\mu_4$
equals \eqref{eq:gl2lo1}. 
The maps $\mu_5,\mu_8$ are defined in an analogous fashion.
We define $\mu_{11},\mu_{13}$ by following the same steps used
in defining $\mu_6,\mu_8$ respectively; for example, in case 
of~$\mu_{11}$, we first replace 
$(H^p_{W}\R\iGp{\X}(-))_{x'}$ by~$H^p_{I}(-)_{x'}$ 
using the isomorphisms $(\R\iG{W}\R\iGp{\X}(-))_{x'} 
\cong (\R\iG{\I}(-))_{x'} \cong \R\iG{I}(-)_{x'}$
and then follow (\ref{eq:gl2lo1}c)-(\ref{eq:gl2lo1}e). 
In particular, $\mu_4,\mu_5,\mu_6,\mu_8,\mu_{11},\mu_{13}$ 
are all isomorphisms. 
For $\mu_9,\mu_{10}$ we refer to the top 
row of~\eqref{cd:coures-1z}. We define~$\mu_7$ by 
\[
\bigoplus_jH^{p}_{x_j}\cFb \xleftarrow[\eqref{eq:hosp0aa}]{\sim} 
\bigoplus_jH^{p}_{x_j}\cGb 
\xrightarrow[\text{top row}]{\ref{lem:coures0}} 
H^{p+1}_{x'}\cGb \xleftarrow[\eqref{eq:hosp0aa}]{\sim}H^{p+1}_{x'}\cFbp
\]
and define $\mu_{12}$ to be the unique map satisfying 
$\mu_{13}\mu_{12}\mu_{11}^{-1} = \mu_3$.

The rectangles on the bottom left corner and the 
bottom right corner of~\eqref{cd:hosp0aaa}
commute for functorial reasons. 
By \ref{lem:psm1}, we have $\R\iG{x^{\prime}}\cFb \cong 0$ 
and $\R\iG{x_j}\cFbp \cong 0$ and therefore by the exactness
of the top row of \eqref{cd:coures-1z} 
it follows that $\mu_9,\mu_{10}$ are isomorphisms.
In particular, $\mu_1,\mu_2$ are isomorphisms. 

Since $\psi = \mu_8\mu_5\mu_0\mu_4^{-1}\mu_6^{-1}$,
it follows that the commutativity statement of the lemma
is the same as proving that the outer skeleton 
of~(\ref{cd:hosp0aaa}) commutes and so we reduce 
to checking commutativity of the middle two rectangles 
in~(\ref{cd:hosp0aaa}). 

Vertically expanding the topmost rectangle 
in \eqref{cd:hosp0aaa} results
in the following diagram having obvious 
natural isomorphisms as vertical maps.
\[
\begin{CD}
\bigoplus_j(\mathbb E \cMb)(x_j) @>{\text{natural}}>> 
 (\mathbb E \cMb)(x^{\prime}) \\
@| @| \\
\bigoplus_j(H^p_{Z^{p}/Z^{p+1}}\cGb)(x_j) @. 
 (H^{p+1}_{Z^{p+1}/Z^{p+2}}\cGb)(x^{\prime}) \\
@V{\wr}VV @VV{\wr}V \\
\bigoplus_jH^{p}_{x_j}\cGb @>\text{top row}>{\ref{lem:coures0}}> 
 H^{p+1}_{x^{\prime}}\cGb \\
@V{\wr}VV @AA{\wr}A \\
\bigoplus_jH^{p}_{x_j}\cFb @>{\mu_7}>> H^{p+1}_{x^{\prime}}\cFbp
\end{CD}
\]
The left and right columns in the preceding diagram give the maps 
$\mu_4$ and~$\mu_5$ respectively. The commutativity of the
subdiagram on the top was recorded in~({\ref{cd:coures1}})
while the bottom rectangle commutes by definition of $\mu_7$.

So for the lemma it only remains to check that the following 
subrectangle of~(\ref{cd:hosp0aaa}) commutes.
\begin{equation}
\label{cd:hosp0ab}
\begin{CD}
\bigoplus_jH^{p}_{x_j}\cFb @>{\mu_7}>> H^{p+1}_{x^{\prime}}\cFbp \\
@A{\mu_{9}}AA @V{\mu_{10}}VV \\
(H^{p}_W\cFb)_{x'} @>{\mu_{12}}>> (H^{p+1}_W\cFbp)_{x'}
\end{CD}
\end{equation}
To that end we claim that $\mu_{12}$ equals 
$(H^{p}_W(-\gamma))_{x'}$ for $\gamma$
as defined in~(\ref{eq:hosp0aa}).
Assuming the claim, checking that (\ref{cd:hosp0ab}) 
commutes amounts to checking that the following diagram 
``commutes'' where the top row represents $\mu_7$ and the bottom
row represents the rest of~(\ref{cd:hosp0ab}).  
\[
\begin{CD}
\bigoplus_jH^{p}_{x_j}\cFb @<{\sim}<< \bigoplus_jH^{p}_{x_j}\cGb 
@>{\ref{lem:coures0}}>{\text{top row}}> H^{p+1}_{x^{\prime}}\cGb 
@<{\sim}<< H^{p+1}_{x^{\prime}}\cFbp \\
@| @. @. @| \\
\bigoplus_jH^{p}_{x_j}\cFb @<{\sim}<{\mu_9}< (H^{p}_W\cFb)_{x'} 
@>{(H^{p}_W(-\gamma))_{x'}}>> (H^{p+1}_W\cFbp)_{x'}  
@<{\sim}<{\mu_{10}}< H^{p+1}_{x^{\prime}}\cFbp
\end{CD}
\]
We use \ref{lem:homo6} in this situation 
with $\Gamma_1 = \iG{x'}(-)$, $\Gamma_2 = (\iG{W}(-))_{x'}$ 
and $\Gamma_3 = \oplus_j \iG{x_j}$.
Upon applying~$H^p$ to the
diagram of~\ref{lem:homo6} we see that 
the top horizontal map in~{\ref{lem:coures0}}
is $(-1)$ times the map labeled ? in \ref{lem:homo6}
(\S\ref{subsec:conv}, \eqref{conv9}). The remaining maps in the
preceding diagram agree with those in~$H^p$ of the
diagram of~\ref{lem:homo6}.
Therefore, by~\ref{lem:homo6}, (\ref{cd:hosp0ab}) commutes.

Verifying the above claim on $\mu_{12}$ amounts to checking that 
the following diagram commutes.
\begin{equation}
\label{cd:hosp0ac}
\begin{CD}
(H^{p}_W\cFb)_{x'} @>{(H^{p}_W(-\gamma))_{x'}}>> (H^{p+1}_W\cFbp)_{x'} \\
@V{\mu_{11}}VV @V{\mu_{13}}VV \\
H^{p_1}_{I}(M \otimes_{y^{\prime}} \cL_{x^{\prime}}) @>{\mu_3}>>
H^{p_1}_{I}(M^{\prime} \otimes_{y^{\prime}} \cL_{x^{\prime}}) 
\end{CD}
\end{equation}
Let us expand this diagram vertically using the definition 
of $\mu_{11},\mu_{13}, \cFb$ and $\cFbp$. 
\[
\begin{CD}
(H^p_W\R\iGp{\X}(f^*\sigma_{\le q}\cMb \otimes_{\X} \cL[n]))_{x'}
 @>{(H^{p}_W(-\gamma))_{x'}}>> 
 (H^{p+1}_W\R\iGp{\X}(f^*\sigma_{\ge q+1}\cMb \otimes_{\X} \cL[n]))_{x'} \\
@VVV @VVV \\
H^p_{I}(f^*\sigma_{\le q}\cMb \otimes_{\X} \cL[n])_{x'}  @>{\mu_{14}}>> 
 H^{p+1}_{I}(f^*\sigma_{\ge q+1}\cMb \otimes_{\X} \cL[n])_{x'} \\
@| @| \\
H^p_{I}(f^*i_yM[-q] \otimes_{\X} \cL[n])_{x^{\prime}} @>{\mu_{15}}>> 
H^{p+1}_{I}(f^*i_{y'}M'[-q-1] \otimes_{\X} \cL[n])_{x^{\prime}} \\
@V{(\ref{eq:gl2lo1}c)-(\ref{eq:gl2lo1}d)}V{\text{sign}\; = \;(-1)^{qn}}V 
@V{\text{sign}\; = \;(-1)^{qn+n}}V{(\ref{eq:gl2lo1}c)-(\ref{eq:gl2lo1}d)}V \\
H^p_{I}((M\otimes_{y'} \cL_{x^{\prime}})[-q+n]) @>{\mu_{16}}>>
H^{p+1}_{I}((M^{\prime}\otimes_{y'} \cL_{x^{\prime}})[-q-1+n]) \\
@V{\text{(\ref{eq:gl2lo1}e)}}VV @VV{\text{(\ref{eq:gl2lo1}e)}}V \\
H^{p_1}_{I}(M\otimes_{y'} \cL_{x^{\prime}}) @>{\mu_3}>>
H^{p_1}_{I}(M^{\prime}\otimes_{y'} \cL_{x^{\prime}}) 
\end{CD}
\]
Here $\mu_{14}$ is \emph{minus} of $H^{p}_{\pfr}(-)_{x'}$ of the 
composite map (where $\otimes = \otimes_{\X}, \cLb = \cL[n]$)  
\[
f^*\sigma_{\le q}\cMb \otimes \cLb \xto{\;\;\Upsilon\;\;} 
(f^*\sigma_{\ge q+1}\cMb)[1] \otimes \cLb \xto{\;\;\theta\;\;}
(f^*\sigma_{\ge q+1}\cMb \otimes \cLb)[1]
\]
where $\Upsilon$ is induced by 
applying $f^*(-) \otimes \cLb$ to
the last map in~\eqref{eq:hosp1} while 
$\theta$ is obtained using the convention in
\S\ref{subsec:conv}, \eqref{conv4}, and therefore,
is the identity map on the graded level.
It follows from the definition of~$\gamma$ 
in~\eqref{eq:hosp0aa} that the topmost rectangle 
commutes.
Next set $\mu_{15} = \mu_{14}$ so that the second rectangle 
from the top commutes. Keeping in mind that 
$\cMb$ is a two-term complex supported in degrees~$q,q+1$, 
one notes that on the stalks at $y'$ the natural map 
\[
i_yM = \M^q = H^q(\sigma_{\le q}\cMb) 
\xto{\text{see \eqref{eq:hosp1}}} 
H^q((\sigma_{\ge q+1}\cMb)[1]) = \M^{q+1} = i_{y'}M'
\]
is precisely $-\partial$ 
where $\partial$ is the canonical 
$\cO_{\Y,y'}$-linear map $M \to M'$
induced by the differential of~$\cMb$ 
(cf.~\S\ref{subsec:conv}, \eqref{conv9}).
Letting $\mu_{16}$ be the obvious map 
induced by~$(-1)^n\partial$ 
we see that the third rectangle from the top,
after going through (\ref{eq:gl2lo1}c) in the vertical maps,
is obtained as $H^{p}_{\pfr}$ of the following commutative
diagram where $\otimes = \otimes_{y'}, L = \cL_{x'}$.
\[
\begin{CD}
M[-q] \otimes L[n] @>{\partial \otimes 1}>> M'[-q-1][1] \otimes L[n]
@>{\text{no}}>{\text{sign}}> (M'[-q-1] \otimes L[n])[1] \\
@VV{(-1)^{qn}}V @. @V{(-1)^{qn+n}}VV \\
\makebox[0pt]{$(M \otimes L)[-q+n]$} @. \makebox[0pt]{
 $\xto{\hspace{4em} (-1)^n \partial \otimes 1\hspace{4em}}$}  
 @. \makebox[0pt]{$(M' \otimes L)[-q+n]$}
\end{CD}
\]
The rectangle involving $\mu_{16}$ and $\mu_3$ clearly commutes 
and thus \eqref{cd:hosp0ac} commutes.
\end{proof}


\medskip

\begin{acor}
\label{cor:hosp0ad}
In the situation of \textup{\ref{lem:hosp0a}}, 
assume further that there exists a sequence in~$I$ of length $p_1$,
say, $\bt = t_1, \ldots, t_{p_1}$ 
such that the following natural maps are isomorphisms.
\[
H^{p_1}_I(M \otimes_{y'} \cL_{x'}) \to 
H^{p_1}_{\bt\cO_{x'}}(M \otimes_{y'} \cL_{x'})
\qquad
H^{p_1}_I(M' \otimes_{y'} \cL_{x'}) \to 
H^{p_1}_{\bt\cO_{x'}}(M' \otimes_{y'} \cL_{x'})
\]
Then any element of 
$\bigoplus_jH^{p_1}_{m_{x_j}}(M \otimes_y \cL_{x_j})$
(resp.~$H^{p_1}_{m_{x'}}(M' \otimes_{y'} \cL_{x'})$) can be written 
as a sum of generalized fractions of the type
\begin{align}
w = &\genfrac{[}{]}{0pt}{}
{m \otimes l}{t_1^{r_1},\ldots,t_{p_1}^{r_{p_1}}}
, \quad m \in M, l \in \cL_{x'}, r_i > 0 \notag \\
(\text{resp. }&\genfrac{[}{]}{0pt}{}
{m' \otimes l'}{t_1^{r'_1},\ldots,t_{p_1}^{r'_{p_1}}}
, \quad m' \in M', l' \in \cL_{x'}, r'_i > 0) \notag
\end{align}
and $\psi$ of~\textup{\ref{lem:hosp0a}} sends the element represented 
by $w$ to the element represented by the fraction
\[
\genfrac{[}{]}{0pt}{}
{(-1)^n\partial(m) \otimes l}{t_1^{r_1},\ldots,t_{p_1}^{r_{p_1}}}.
\]  
\end{acor}
\begin{proof}
Since $\mu_1,\mu_2$ of \ref{lem:hosp0a} are isomorphisms,
the generalized-fraction representation holds. The description
of $\psi$ in terms of such fractions follows immediately
from the commutative diagram of~\ref{lem:hosp0a}. 
\end{proof}

\medskip

\subsection{Application to the retract case}
\label{subsec:ret}

Let $(\Y,\DsssY) \xto{i} (\X,\DsssX) \xto{h} (\Y,\DsssY)$ 
be a factorization of the identity map on $\Y$ 
in $\bbFc$, where~$i$ is a closed immersion 
and~$h$ a smooth map of constant relative dimension, 
say~$n$. Note that $i$ takes any immediate 
specialization in $\Y$ to an $h$-lateral 
specialization in~$\X$. (\S\ref{subsec:hosp})
\begin{alem}
\label{lem:ret0}
(cf.~\cite[Thm.~2.6]{Ye})
In the above situation, let $y \leadsto y'$ 
be an immediate specialization in~$\Y$ having 
corresponding image $x \leadsto x'$ in~$\X$. 
Then there exists an $\cO_{\X,x'}$-sequence 
$\bt = t_1, \ldots, t_n$ satisfying the following properties.
\begin{enumerate}
\item The sequence $\bt$ 
maps to a regular system of parameters 
in~$\cO_{\X,x'}/m_{y'}\cO_{\X,x'}$ 
and in~$\OXx/m_{y}\OXx$. Moreover, a basis 
of $(\Ohm^1_{\X/\Y})_x$ and $(\Ohm^1_{\X/\Y})_{x'}$
is given by $dt_1,\ldots,dt_n$.
\item Let $\I$ be the largest (open) coherent $\OX$-ideal 
defining $\ov{\{x\}}$ and $\sK$ the largest (open) 
coherent $\OY$-ideal defining $\ov{\{y\}}$. Then 
$\sK_{y'}\cO_{\X,x'} + \bt\cO_{\X,x'} = \I_{x'}$.
\end{enumerate}
\end{alem}

\begin{proof}
The statements are local in nature. 
For any open subset $\U$ in $\X$, with $\V = i^{-1}\U$
we have $h(\U) \subset \V$
and $h$ restricts to a smooth map $\U \to h^{-1}\V$.
Thus we may replace $\X$ by $\U \set \Spf(B)$ 
an affine open neighborhood of~$x'$ in~$\X$,
and $\Y$ by $\V = i^{-1}\U = \Spf(A)$. 
Then the natural induced maps $A \xto{\phi} B \xto{\pi} A$ 
factor the identity map on $A$ and $\phi$ is a smooth map. 

By replacing $\U, \V$ by smaller open subsets if necessary 
such that $x' \in \U$ and $\V = i^{-1}\U$ remain valid, 
we may assume, without loss of generality, 
that~$\Ohm^1_{\U/\V}$ is free so that 
the $B$-module $\Ohm^1_{B/A} \iso \Gamma(\U,\Ohm^1_{\U/\V})$ 
(see \eqref{eq:diff3d} in \S\ref{subsec:difform}) is free
of rank~$n$. Let $\I$ be the ideal in $\OX$ defining the 
closed immersion $i$. Set $J = \ker \pi$. 
Then $\J\big|_{\U} = J^{\sim_B}$.
By~\ref{lem:diff3} and~\ref{prop:homelessc} 
there are natural isomorphisms 
\[ 
J/J^2 \iso \Ohm^1_{B/A} \otimes_B B/J, \qquad
\J/\J^2\big|_{\U} \iso i^*\Ohm^1_{\U/\V}. 
\]
By the Nakayama lemma it follows that $\J_{x'} = J\cO_{\X,x'}$
is generated by~$n$ elements. These elements may be assumed to
be images of sections of $\J$ over a suitably small open affine
neighborhood of $x'$. Again, by shrinking
$\U,\V$ if necessary we may assume that~$J$ is generated 
over~$B$ by~$n$ elements. Let~$\bt = \{t_1, \ldots, t_n\}$ 
be such a generating set. The 
kernel $\J_{x}$ (resp.~$\J_{x'}$) of the natural 
map $\OXx \to \OYy$ (resp.~$\cO_{\X,x'} \to \cO_{\Y,y'}$) 
induced by~$i$ is generated by the natural image of~$\bt$. 
Therefore, by~\cite[Thm.~8.4]{Ma}, the maps 
(completions being along the respective maximal ideals)
\[
\wh{\OYy}[[T_1, \ldots, T_n]] \xto{\;\theta\;} \wh{\OXx},
\qquad
\wh{\cO_{\Y,y'}}[[T_1, \ldots, T_n]] \xto{\;\theta'\;} \wh{\cO_{\X,x'}},
\]
defined by sending $T_j$ to $t_j$, are surjective. Let us 
verify that $\theta, \theta'$ are isomorphisms.
Using \ref{lem:diff3} for the maps
$\OYy \to \wh{\OYy}[[T_1, \ldots, T_n]] \xto{\;\theta\;} \wh{\OXx}$
we see, by comparing ranks of the free modules in \ref{lem:diff3},
that for $\afr = \ker \theta$, we have $\afr/\afr^2 = 0$. By the Nakayama
lemma $\afr = 0$. A similar proof works for $\ker \theta'$. 
Now (i) follows by considering the image of $\bt$ in $\cO_{\X,x'}$.

For (ii) note that $\J \subset \I$ and moreover $\I/\J$ is mapped to 
$\sK$ under the isomorphism $(\OX/\J)\big|_{\Z} \iso \OY$. 
In particular, it holds that $\I = \J + \sK\OX$.
Since $\J_{x'}$ is generated by $\bt$, (ii) results.  
\end{proof}

\medskip

\newcommand{\Esssh}{{\mathbb E}_{\scriptscriptstyle h}}
Let $\cMb \in \Coz(\Y,\DsssY)$. 
Set $\Esssh(-) \set E_{\DsssX}\R\iGp{\X}(h^*(-) \otimes \omega_h[n]).$
We now consider a graded isomorphism 
\begin{equation}
\label{eq:ret0b}
i^{\flat}\Esssh\cMb \xto{\quad} \cMb
\end{equation}
which is defined pointwise, say for $y \in \Y$, 
by the following isomorphism 
(where $M = \cMb(y)$, $x = i(y)$, $q = \DsssY(y) = \DsssX(x)$, 
$\omega = \omega_h$)
\begin{align}
(i^{\flat}\Esssh\cMb)(y) &\xto{\;\eqref{eq:clim1a}\;} 
 \Hom_{\OXx}(\OYy, (\Esssh\cMb)(x)) \notag \\
&\xto{\;\eqref{eq:gl2lo1}\;} \Hom_{\OXx}(\OYy, 
H^n_{m_x}(M \otimes_y \omega_x)) 
\xto{\;(-1)^{qn}\text{ times }\wh{\eqref{eq:itloco5}}\;}  M. \notag
\end{align}
For the last isomorphism we also need \ref{prop:app3} and 
the isomorphisms of~\eqref{eq:itloco0} and~\eqref{eq:itloco0a}.

\medskip

\begin{aprop}
\label{prop:ret1a}
The graded map in \eqref{eq:ret0b} is also a map of complexes.
(cf.~\cite[Lemma 6.13]{Ye})
\end{aprop}
\begin{proof}
Let $y \leadsto y^{\prime}$ be an immediate specialization in~$\Y$ and 
let $x \leadsto x^{\prime}$ denote the corresponding image in~$\X$.
Set $q = \DsssY(y) = \DsssX(x)$, so that $q+1 = \DsssY(y') = \DsssX(x')$.
It suffices to prove that the following diagram, 
whose columns give~\eqref{eq:ret0b}, commutes.
Here~$\psi$ and~$\partial$ have the same meaning as in~\ref{lem:hosp0a}.
\[
\begin{CD}
(i^{\flat}\Esssh\cMb)(y) @>>> (i^{\flat}\Esssh\cMb)(y') \\
@V{\wr}VV  @V{\wr}VV \\
\Hom_{\OXx}(\OYy, H^n_{m_x}(M \otimes_y \omega_x)) @.
 \Hom_{\cO_{\X,x'}}(\cO_{\Y,y'}, H^n_{m_{x'}}
 (M' \otimes_{y'} \omega_{x'})) \\
@VVV  @VVV  \\
H^n_{m_x}(M \otimes_y \omega_x) @>{\psi}>> 
 H^n_{m_{x'}}(M' \otimes_{y'} \omega_{x'}) \\
@V{\text{$(-1)^{qn} \times$ res}}V{\text{of $\wh{\eqref{eq:itloco5}}$}}V  
 @V{\text{$(-1)^{qn+n} \times$ res}}V{\text{of 
 $\wh{\eqref{eq:itloco5}}$}}V \\
M @>{\partial}>> M'
\end{CD}
\]
The commutativity of the top rectangle is straightforward to 
verify. We now verify commutativity of the bottom rectangle by 
a chase involving generalized fractions. 
We use~\ref{cor:hosp0ad} in this regard.
Regarding the notation in \ref{lem:hosp0a}, in 
this situation, we have $W = \ov{\{x\}}$, $p_1 = n$ 
(since $\DsssX(x) = \DsssY(y)$), $\cL = \omega_h$
and $\I$ is the largest coherent ideal defining~$W$. 

Let $\bt$ be the sequence in $\Oxp$
obtained in~\ref{lem:ret0}. Let us verify that $\bt$
satisfies the hypothesis of~\ref{cor:hosp0ad}.
Let $\sK$ be as in (ii) of~\ref{lem:ret0}. 
Since $\sK_{y'}\OYy = \sK_{y} = m_y$, the latter being the maximal
ideal of $\OYy$, and since $M$ is $m_y$-torsion, we see that $M$
is $\sK_{y'}$-torsion as an $\cO_{\Y,y'}$-module. With $I = \I_{x'}$
as in \ref{lem:hosp0a}, by (ii) of~\ref{lem:ret0} we see that the 
canonical map $\R\iG{I}(M \otimes_{y'} \omega_{x'})
\to \R\iG{\bt\Oxp}(M \otimes_{y'} \omega_{x'})$
is an isomorphism. Since $M'$ is $m_{y'}$-torsion 
and by \ref{lem:ret0}(i),
$m_{x'} = m_{y'}\Oxp + \bt\Oxp \subseteq m_{y'}\Oxp + I \subseteq m_{x'}$, 
the other condition also results.
 
By \ref{lem:ret0}(i), $dt_1 \wedge dt_2 \wedge \cdots \wedge dt_n$
is a generator of $\omega_{x'}$ and of~$\omega_x$
(cf.~definition of~\eqref{eq:itloco5}). 
Now chasing the image of the fraction
\[
\genfrac{[}{]}{0pt}{}
{m \otimes dt_1 \wedge dt_2 \wedge \cdots \wedge dt_n}
{t_1^{r_1},\ldots,t_{p_1}^{r_{p_1}}} 
\]
in $H^n_{m_x}(M \otimes_y \omega_x)$ we conclude 
by~\ref{cor:hosp0ad}.
\end{proof}

\newpage

\section{The main theorem}
\label{sec:finale}

Having established the crucial ingredients in \ref{prop:itloco1},
\ref{prop:clim3}, \ref{prop:cartsq9} and~\ref{prop:ret1a}
we are now in a position to prove the Main Theorem stated
in~\S\ref{subsec:outline}. The crucial results 
are brought together in Proposition \ref{prop:fin-fact6} below. This enables
a  fairly quick proof of existence of the basic pseudofunctor $(-)^\sharp$, 
given in \S\ref{subsec:fin-gen}. 

In \S\ref{subsec:trans} we describe the
behavior of $(-)^{\sharp}$ with respect to translations.
\vspace{1pt}


Let us recall some notation discussed in~\S\ref{subsec:outline}.
Let $(\X, \Delta)$ be an object in $\bbFc$. 
Let $\Coz^0_{\Delta}(\X)$ be the category 
consisting of graded maps of 
the underlying graded objects of the complexes
in $\Coz_{\Delta}(\X)$. More formally, an object $\cNb$ 
in $\Coz^0_{\Delta}(\X)$ is a $\mathbb Z$-graded sequence
of $\Aqct(\X)$-modules such that for any $p \in \mathbb Z$,
$\N^p$ lies on the $p$-th skeleton induced by~$\Delta$.
A morphism in $\Coz^0_{\Delta}(\X)$ is simply 
a $\mathbb Z$-graded sequence of maps in~$\Aqct(\X)$.
There is a canonical forgetful functor 
$\Coz_{\Delta}(\X) \to \Coz^0_{\Delta}(\X)$ that 
forgets the differential on a Cousin complex. For
$\cNb \in \Coz_{\Delta}(\X)$ we shall denote
its image in $\Coz^0_{\Delta}(\X)$ also by $\cNb$.

Let $f \colon (\X,\DsssX) \to (\Y,\DsssY)$ be 
a map in $\bbFc$. For any $\cMb \in \Coz^0_{\DsssY}(\Y)$,
we define an object 
$f^{\natural}\cMb \in \Coz^0_{\DsssX}(\X)$ by 
(with notation explained below)
\[
(f^{\natural}\cMb)(x) \set (\wh{f_x})_{\sharp}M, \qquad x \in \X,
\]
where, with $y = f(x)$, we have $M = \cMb(y)$, $\wh{f_x}$ is the 
induced map on the completions of the stalks $\wh{\OYy} \to \wh{\OXx}$
and hence is in $\Cfr$ (\S\ref{subsec:huang}) 
and $(-)_{\sharp}$ is the pseudofunctor
on $\Cfr$ of \ref{thm:huang1}. Via the canonical forgetful functor,
we shall, by abuse of notation, let $f^{\natural}$ take inputs from
$\Coz_{\DsssY}(\Y)$ too.
Since $(-)_{\sharp}$ is a pseudofunctor on~$\Cfr$, we can, 
in an obvious manner, make
$(-)^{\natural}$ into a pseudofunctor on~$\bbFc$.

The key step in defining
$f^{\sharp}\cMb$ of our Main Theorem
is to specify a differential on $f^{\natural}\cMb$. 
We do this first for the subcategory of smooth $\bbFc$-maps and for the 
subcategory of closed immersions, then paste locally via factorizations, and finally use independence of the result of pasting from the choice of factorizations (Proposition \ref{prop:fin-fact6}) to globalize.
   
\subsection{Smooth maps}
\label{subsec:fin-sm}
For any smooth map $f \colon (\X,\DsssX) \to (\Y,\DsssY)$ 
in $\bbFc$ having constant relative dimension $d$ and for any 
$\cMb \in \Coz_{\DsssY}(\Y)$, set \index{ $\Forget$0@${\mathbb E}_{(-)}$ (Cousin-complex functor associated to smooth map $(-)\>$)}
\[
{\mathbb E}_f\cMb \set 
E_{\DsssX}\R\iGp{\X}(f^*\cMb \otimes_{\X} \omega_f[d\>]).
\]
Let $x$ be a point in~$\X$. Set $y= f(x)$, $p=\DsssX(x)$, $q=\DsssY(y)$. 
Let $\wh{f_x}$ be the induced map on the completion of the 
stalks $\wh{\OYy} \to \wh{\OXx}$. Let $p_1$ be the relative dimension
of $\wh{f_x}$, which, by \ref{cor:app4}, equals $p-q+d$. 
Set $M = \cMb(y)$. Consider the following isomorphism
\begin{align}
\label{eq:fin-sm1}
({\mathbb E}_f\cMb)(x) &\xto{\;\eqref{eq:gl2lo1}\;} 
H^{p_1}_{m_x}(M \otimes_y (\omega_f)_x) \\
&\xto{\ref{prop:app3} \text{ and } \eqref{eq:itloco0}} 
 H^{p_1}_{\wh{m_x}}(M \otimes_{\wh{y}} \omega_{\wh{f_x}})
 \xto{\quad\theta\quad} (\wh{f_x})_{\sharp}M = (f^{\natural}\cMb)(x) \notag
\end{align}
where $\theta$ is $(-1)^{(p+d)q+p}$ times the isomorphism in 
\ref{thm:huang1}, I.(i).
As $x$ ranges over $\X$ we therefore obtain a graded isomorphism 
${\mathbb E}_f\cMb \iso f^{\natural}\cMb$.
If $f$ does not have constant relative dimension on $\X$, then we
restrict to connected components of $\X$ and then carry out the 
above procedure. 

In particular, there is now a natural candidate for $f^{\sharp}\cMb$
(with $f$ smooth).

\pagebreak[3]
The isomorphism 
${\mathbb E}_{f} \iso f^{\natural}$ given by~\eqref{eq:fin-sm1} 
behaves well with restriction to open subsets on~$\X$. 
In greater detail, let~$\U$ be an open 
subset of~$\X$ and $(\U, \Delta) \xto{u} (\X,\DsssX)$  
the corresponding open immersion. Set $f_1 = fu$. 
Then the following diagram commutes:
\begin{equation}
\label{eq:fin-sm2}
\begin{CD}
({\mathbb E}_f \cMb)\big|_{\U} @>>> 
 (f^{\natural}\cMb)\big|_{\U} \\
@V{\alpha}VV  @| \\
{\mathbb E}_{f_1} \cMb @>>> f_1^{\natural}\cMb
\end{CD}
\end{equation}
where $\alpha$ is defined through the following sequence of 
obvious natural maps  
(with $E_{\X} \set E_{\DsssX}$, $E_{\U} \set E_{\Delta}$)
\begin{align}
({\mathbb E}_f \cMb)\big|_{\U} &=
\bigl(E_{\X}\R\iGp{\X}(f^*\cMb \otimes_{\X} \omega_f[d\>])\bigr)
 \big|_{\U} \notag \\
&\iso E_{\U}\bigl((\R\iGp{\X}(f^*\cMb \otimes_{\X} \omega_f[d\>]))\big|_{\U}
 \bigr) \notag \\
&\iso E_{\U}\R\iGp{\U}((f^*\cMb \otimes_{\X} \omega_f[d\>])\big|_{\U}) \notag \\
&\iso E_{\U}\R\iGp{\U}(f_1^*\cMb \otimes_{\U} \omega_f\big|_{\U}[d\>]) \notag \\
&\iso E_{\U}\R\iGp{\U}(f_1^*\cMb \otimes_{\U} \omega_{f_1}[d\>]) = 
 {\mathbb E}_{f_1} \cMb. \notag 
\end{align}
The commutativity of \eqref{eq:fin-sm2} is verified punctually as follows.
Let $x, y, p, M$ etc., be as before. Expanding the horizontal maps
we obtain the following diagram where
it suffices to verify that the rectangle on the left commutes.
\[
\begin{CD}
({\mathbb E}_f \cMb)(x) @>(\ref{eq:gl2lo1}a)-(\ref{eq:gl2lo1}c)>> 
 H^{p}_{m_x}(M[-q] \otimes_y \omega_{f,x}[d\>]) 
 @>{\text{Rest of \eqref{eq:fin-sm1}}}>> (\wh{f_x})_{\sharp}M \\
@V{\alpha(x)}VV  @| @| \\
({\mathbb E}_{f_1} \cMb)(x) @>(\ref{eq:gl2lo1}a)-(\ref{eq:gl2lo1}c)>> 
 H^{p}_{m_x}(M[-q] \otimes_y \omega_{f,x}[d\>]) 
 @>{\text{Rest of \eqref{eq:fin-sm1}}}>> (\wh{f_x})_{\sharp}M \\
\end{CD}
\] 

\smallskip\noindent
Expand the rectangle on the left, with the following notation:
\[ 
\cGb \set f^*\cMb \otimes_{\X} \omega_f[d\>], \quad
\cGb_1 \set f_1^*\cMb \otimes_{\U} \omega_{f_1}[d\>], \quad
\Gb \set M[-q] \otimes_y \omega_{f,x}[d\>].
\]

\smallskip
\[
\begin{CD}
(E_{\X}\R\iGp{\X}\cGb)(x) @>>> H^{p}_x\R\iGp{\X}\cGb @>>> H^{p}_{m_x}\cGb_x
 @>>> H^{p}_{m_x}\Gb \\ 
@VVV @VVV \\[1pt]
(E_{\U}(\R\iGp{\X}\cGb)\big|_{\U})(x) @>>> 
 H^{p}_x((\R\iGp{\X}\cGb)\big|_{\U}) @. 
 \rotatebox{90}{\makebox[0pt]{{\raisebox{.5ex}{\makebox[6em]{\hrulefill}}}
 \hspace{-6em}\!\!{\makebox[6em]{\hrulefill}}}} @.
 \rotatebox{90}{\makebox[0pt]{\hspace{-4em}{\raisebox{.5ex}{\makebox[10em]
 {\hrulefill}}}\hspace{-10em}{\makebox[10em]{\hrulefill}}}} \\[1pt]
@VVV @VVV \\[2pt]
(E_{\U}\R\iGp{\U}(\cGb\big|_{\U}))(x) @>>> 
 H^{p}_x\R\iGp{\U}(\cGb\big|_{\U}) @>>> H^{p}_{m_x}\cGb_x\\
@VVV @VVV @VVV \\
(E_{\U}\R\iGp{\U}\cGb_1)(x) @>>> 
 H^{p}_x\R\iGp{\U}\cGb_1 @>>> H^{p}_{m_x}{\cGb_1}_x @>>> H^{p}_{m_x}\Gb\\
\end{CD}
\]
For the  horizontal maps in the first two columns from the left,
we refer to \eqref{eq:coz3} of \S\ref{subsec:cousin}. Thus these two
columns together correspond to (\ref{eq:gl2lo1}a). 
The horizontal maps in the rightmost column 
correspond to (\ref{eq:gl2lo1}b)-(\ref{eq:gl2lo1}c).
The vertical maps are the canonical ones.
Commutativity of the above diagram is straightforward to check and 
thus~\eqref{eq:fin-sm2} commutes.\vspace{1pt}

\pagebreak[3]

The following is an immediate consequence of the 
commutativity of~\eqref{eq:fin-sm2}.

\begin{aprop}
\label{lem:fin-sm3}
Let\/ $f \colon (\X,\DsssX) \to (\Y,\DsssY)$ be a smooth\/ $\bbFc$-map. Let\/~\mbox{$\U\subset \X$} be open  and\/ 
$(\U, \Delta) \xto{\under{.5}{u}} (\X,\DsssX)$  
the corresponding open immersion.\vspace{.5pt} 
If\/ $(-)^{\sharp}$~is defined for smooth maps via\/~\eqref{eq:fin-sm1}$,$
then for\vspace{-1pt} any\/ $\cMb \in \Coz_{\DsssY}(\Y), 
$\mbox{\ $(f^{\sharp}\cMb)\big|_{\U} = (fu)^{\sharp}\cMb$} as complexes.
In particular, 
$u^{\sharp}$ is the restriction of\/ $1_{\X}^{\sharp}$ to\/~$\U$. \looseness=-1
\end{aprop}

\subsection{Closed immersions}
\label{subsec:fin-clim}
Let $f \colon (\X,\DsssX) \to (\Y,\DsssY)$ be a closed immersion 
in~$\bbFc$. Then for any $\cMb \in \Coz_{\DsssY}(\Y)$ and any 
$x \in \X$, with $y=f(x)$, $M = \cMb(y)$, and~$\wh{f_x}$ having the 
obvious meaning as above, we use the following isomorphism
\begin{align}
\label{eq:fin-clim1}
(f^{\flat}\cMb)(x) &\xto{\;\eqref{eq:clim1a}\;} \Hom_{\OYy}(\OXx, M) \\
&\xto{\;\eqref{eq:itloco0a}\;} \Hom_{\wh{\OYy}}(\wh{\OXx}, M)
 \xto{\;\text{\ref{thm:huang1}, I.(ii)}\;}(\wh{f_x})_{\sharp}M = 
 (f^{\natural}\cMb)(x) \notag
\end{align}
to get a graded isomorphism $f^{\flat}\cMb \iso f^{\natural}\cMb$.
In particular, there is now a natural candidate for $f^{\sharp}\cMb$
for $f$ a closed immersion.

The isomorphism $f^{\flat} \iso f^{\natural}$ 
given by \eqref{eq:fin-clim1} commutes with restriction to open subsets
in the following sense. Let~$\U$ be an 
open subset of~$\Y$. Set $\V \set f^{-1}\U$.
Then the following diagram of obvious natural maps commutes.
\[
\begin{CD}
(f^{\flat}\cMb)\big|_{\V} @>>> (f^{\natural}\cMb)\big|_{\V} \\
@VVV @| \\
(f\big|_{\V})^{\flat}(\cMb\big|_{\U}) @>>> 
 (f\big|_{\V})^{\natural}(\cMb\big|_{\U})
\end{CD}
\]
Commutativity of the above diagram is proved by showing it at 
the punctual level where it is straightforward to verify. In summary,
if $(-)^{\sharp}$ is defined for closed immersions
using~\eqref{eq:fin-clim1}, then 
$(f\big|_{\V})^{\sharp}(\cMb\big|_{\U}) = (f^{\sharp}\cMb)\big|_{\V}$. 

\begin{aprop}
\label{rem:fin-clim2}
Let $f \colon (\X, \DsssX) \to (\Y, \DsssY)$ be a map in 
$\bbFc$ that is both a smooth map and a closed immersion. Then 
for any $\cMb \in \Coz_{\DsssY}(\Y)$,
the two natural candidates for a differential 
on~$f^{\natural}\cMb$, viz., the one induced
by~\eqref{eq:fin-sm1} and the one by~\eqref{eq:fin-clim1}, agree.
Moreover, if $\X = \Y$ and $f$ is the identity map $1_{\X}$,
then the graded isomorphism
$\delta^{\natural}_{\X}(\cMb) \colon 1^{\natural}_{\X}\cMb \to \cMb$
also induces the same differential.
\end{aprop}
\begin{proof}
First we show that $f$ is an isomorphism from $\X$ 
onto a union of connected components of~$\Y$. 
Let $\V = \Spf(A)$ be a connected affine open subset of $\Y$ for 
which $f^{-1}\V$ is non-empty. Since $f$ is a closed immersion, 
$f^{-1}\V$ is an affine open subset, say $\U \cong \Spf(A/I)$, 
of $\X$. It suffices to show that $I = 0$. Let $\afr$ 
be a defining ideal in $A$. 
By \ref{prop:morph2}, the map $A \to A/I$ is flat and hence 
the same is true for $A/{\afr} \to A/(\afr + I)$. Hence 
$(\afr + I)/{\afr}$ is idempotent. Since 
$\Spec(A/{\afr})$ ($=\V$ as a topological space)
is connected it follows that $(\afr + I)/{\afr} = 0$ and hence 
$I \subset \afr$. As~$A$ is complete 
w.r.t.~$\afr$, $I$ is in the Jacobson radical of $A$. 
Flatness of $A \to A/I$ implies that $I$ is idempotent and hence 
$I=0$ by the Nakayama lemma. 

Thus we may assume without loss of generality that $f$ 
is an isomorphism. We drop reference to the codimension functions
for the rest of this proof.
Now note that there are canonical isomorphisms 
\[
f^*\cMb \iso {\mathbb E}_f\cMb, 
\qquad  f^*\cMb \iso f^{-1}\cMb \iso f^{\flat}\cMb, 
\]
where the first one is obtained by composing the isomorphisms 
\[
f^*\cMb \iso Ef^*\cMb \iso E\R\iGp{\X}f^*\cMb
\iso E\R\iGp{\X}(f^*\cMb \otimes_{\X} \omega_f[d\>]) = {\mathbb E}_f\cMb
\]
which are based on the following: 
\begin{itemize}
\item The complex $f^*\cMb$ being Cousin, the isomorphism
of \ref{lem:coz1a}(ii) applies.
\item Since $f^*\cMb$ consists of torsion modules, 
the canonical map is an isomorphism
$\R\iGp{\X}f^*\cMb \to f^*\cMb$.
\item In this situation $\omega_f[d\>] = \OX[0]$. (Henceforth we shall identify
$\otimes_{\X} \OX[0]$ with the identity functor.)
\end{itemize}
We claim that the following
two graded isomorphisms are equal, thereby showing 
that~\eqref{eq:fin-sm1} and~\eqref{eq:fin-clim1}
induce the same differential on~$f^{\natural}\cMb$:
\begin{equation}
\label{eq:fin-clim3}
f^*\cMb \iso {\mathbb E}_f\cMb \xto{\;\eqref{eq:fin-sm1}\;} f^{\natural}\cMb,
\qquad
f^*\cMb \iso f^{\flat}\cMb \xto{\;\eqref{eq:fin-clim1}\;} f^{\natural}\cMb.
\end{equation}
The final statement of the proposition then follows from \ref{thm:huang1}, III.

For any $x \in \X$, $\wh{f_x}$ is an isomorphism, and hence 
to prove the above claim it suffices to show that, 
at the punctual level, via the canonical isomorphism
\[
M \iso M \otimes_y \OXx \iso \iG{x}f^*\M^p = (f^*\cMb)(x),  \qquad 
(p = \DsssX(x))
\]
the two graded isomorphisms in~\eqref{eq:fin-clim3} reduce  
to the (inverse of the) corresponding isomorphisms 
in~\ref{rem:huang7}. We now show this reduction 
for the first isomorphism in~\eqref{eq:fin-clim3}, 
leaving the other one to the reader.  

Fix $x \in \X$. Let $y = f(x), M = \cMb(y)$. Consider the following diagram
of natural isomorphisms whose top row is the first isomorphism 
in~\eqref{eq:fin-clim3}. 
\[
\begin{CD}
@. (f^*\cMb)(x) @>>> ({\mathbb E}_f\cMb)(x) @>{\eqref{eq:fin-sm1}}>> 
 (f^{\natural}\cMb)(x) \\
@. @VVV @VV{\eqref{eq:gl2lo1}}V @| \\
M @>>> M \otimes_y \OXx @>>> H^0_{m_x}(M \otimes_y \OXx) 
 @>{\text{rest of \eqref{eq:fin-sm1}}}>> {\wh{f_x}}_{\sharp}M   
\end{CD}
\]
Note that all the signs involved in {\eqref{eq:gl2lo1}} and the 
rest of \eqref{eq:fin-sm1} vanish in this situation. 
The bottom row is inverse of the first isomorphism
in~\ref{rem:huang7} and therefore it suffices to show that 
the above diagram commutes. The rectangle on the right commutes trivially.
We expand the one on the left as follows.
\[
\begin{CD}
(f^*\cMb)(x) @<{\ref{lem:coz1a}\textup{(ii)}}<< (Ef^*\cMb)(x) 
@<<< (E\R\iGp{\X}f^*\cMb)(x) \\
@|  @V{\eqref{eq:coz1}}V{\alpha_1}V  @V{\eqref{eq:coz1}}V{\alpha_2}V  \\
(f^*\cMb)(x) @<{\ref{lem:coz1a}\textup{(i)[b]}}<< H^p_xf^*\cMb 
 @<<< H^p_x\R\iGp{\X}f^*\cMb \\
@VVV  @VV{\beta_1}V  @VV{\beta_2}V  \\
M \otimes_y \OXx @<<< H^0_{m_x}(M \otimes_y \OXx) @= 
 H^0_{m_x}(M \otimes_y \OXx) 
\end{CD}
\]
Only the maps $\beta_i$ require explanation, the other nonlabeled
maps being the obvious ones. Let $\beta_2$ be the unique map 
such that $\beta_2\alpha_2$ gives~{\eqref{eq:gl2lo1}}. (In other words,
$\beta_2$ involves (\ref{eq:gl2lo1}a)-(\ref{eq:gl2lo1}e), except
that~{\eqref{eq:coz1}}, which is part of~(\ref{eq:gl2lo1}a), is 
taken care of by~$\alpha_2$.) We define~$\beta_1$ by following the
same steps as in~$\beta_2$, the only difference being the 
absence of~$\R\iGp{\X}$ (cf.~\eqref{eq:nolab3}, \eqref{eq:coz3},
where~$\R\iGp{\X}$ may be dropped when~$\cFb$ consists of 
torsion modules).

The two rectangles on the right side commute for functorial reasons.
For the one on the top left, see the proof
of~\ref{lem:coz1a}(iii) with $\cCb = f^*\cMb$. The one on the
bottom left expands to the following diagram, whose
commutativity is easily verified.
\[
\begin{CD}
(f^*\cMb)(x) @= H^p\iG{x}f^*\cMb @>>> H^p_xf^*\cMb \\[-2pt]
@| @| @VVV \\[-1.5pt]
\iG{x}f^*\M^p @=H^p(\iG{x}f^*\M^p[-p]) @>>>  H^p_x(f^*\M^p[-p]) \\[-2pt]
@| @. @VVV \\[-1.5pt]
\iG{x}f^*\M^p @.
 {\makebox[0pt]{\quad{\raisebox{.5ex}{\makebox[16em]{\hrulefill}}}
 \hspace{-16em}\!\!{\makebox[16em]{\hrulefill}}}} @.  H^0_xf^*\M^p \\[-3pt]
@VVV @. @VVV \\[-2.5pt]
M \otimes_y \OXx @. 
 {\makebox[0pt]{$\xgets{\hspace{12em}}$}}@. H^0_{m_x}(M \otimes_y \OXx)
\end{CD}
\]   
\end{proof}


\subsection{Factorizable maps}
\label{subsec:fin-fact}
Consider a map $f$ in $\mathbb F$
which admits a factorization of the type $f=f_1f_2$ where 
$f_2$ is a closed immersion and $f_1$ a smooth map. 
The graded isomorphisms in~\eqref{eq:fin-sm1} 
and~\eqref{eq:fin-clim1}, then give a candidate for
$f^{\sharp}$. However since there is no unique 
choice for a factorization of $f$ as a smooth map followed 
by a closed immersion  we must first show that different 
factorizations do lead to the same definition of $f^{\sharp}$. 
We accomplish this in~\ref{prop:fin-fact6}, the main result 
of this section.

\begin{adefi}
\label{def:fin-fact0}
Let $h$ be a map in $\bbFc$ and let $\bh = h_1, \ldots ,h_n$ be a sequence 
of maps in $\bbFc$ such that $h$ factors as 
$h_nh_{n-1}\cdots h_1$ and
each~$h_i$ is a smooth map or a closed immersion:
\[
(\X,\Delta) \xto{\;h_1\;} (\Z_1,\Delta_1) \xto{\;h_2\;} \ldots \xto{\;h_n\;} 
(\Y,\Delta_n).
\]  
For any $\cMb \in \Coz_{\Delta_n}(\Y)$ we define a Cousin complex
in~$\Coz_{\Delta}(\X)$, to be denoted by 
$\bh^{\sharp}\cMb$ or $(h_1| \cdots |h_n)^{\sharp}\cMb\<$, 
as follows. As a graded object we set\vspace{-1pt} 
$\bh^{\sharp}\cMb = h^{\natural}\cMb$. 
%
We set the differential of $\bh^{\sharp}\cMb$ to be the one 
induced by\vspace{-2pt} 
$C_{\bh}^{\natural} \colon 
h_1^{\natural}\cdots h_n^{\natural}\cMb \iso h^{\natural}\cMb  $ 
where $C_{\bh}^{\natural}$ is the natural map derived from the\vspace{-1.5pt} 
pseudofunctoriality of~$(-)^\natural$
and $h_1^{\natural}\cdots h_n^{\natural}\cMb$ is equipped
with a differential by iteratively using~\eqref{eq:fin-sm1} 
and~\eqref{eq:fin-clim1}.
%
\end{adefi}

In the notation of \ref{def:fin-fact0}, 
the underlying object of $\bh^{\sharp}\cMb$
depends only on the composite map $h$ and not on the individual terms
of the sequence~$\bh$. We want to show
that the differential of~$\bh^{\sharp}\cMb$ 
is also something that depends only on~$h$. We first 
show this in four special cases. 
These cases are based on the results in \ref{prop:itloco1},
\ref{prop:clim3}, \ref{prop:cartsq9} and~\ref{prop:ret1a}.
What follows is essentially a rewriting of these results 
in a way that incorporates~$(-)^{\natural}$ via the isomorphisms
in~\eqref{eq:fin-sm1} and~\eqref{eq:fin-clim1}. 
It is important to get signs right, so
we give elaborate proofs whenever necessary. 
 
\smallskip
\enlargethispage*{1pt}
{\textbf{Case 1.}} 
\emph{Suppose $f,g$ are composable smooth maps in $\bbFc$ so
that the composition $gf$ is also smooth. Then we claim that 
$(gf)^{\sharp} = (f|g)^{\sharp}$.} If $f,g$ have constant relative
dimension (and hence so does $gf$ by \ref{lem:diff7}), then 
it suffices to show that in the situation 
of \S\ref{subsec:itcoz}, 
with $\cL_1 = \omega_g$ and $\cL_2 = \omega_f$, 
the following diagram commutes.
\def\verti{\raisebox{2.2pt}{\rotatebox{-90}{\makebox[0pt]
{$\xrightarrow[{\rotatebox{90}{\small using}}]
{\qquad{\rotatebox{90}{\small \eqref{eq:fin-sm1}}}\qquad}$}}}}
\begin{equation}
\label{eq:fin-fact1}
\begin{CD}
({\mathbb E}_{gf}\cMb)(x) 
 @>{\eqref{eq:itcoz1}\text{ and }\ref{lem:diff7}}>> 
 ({\mathbb E}_f{\mathbb E}_g\cMb)(x) \\[-3pt]
@. @V{\text{using \eqref{eq:fin-sm1}}}VV \\[-2pt]
\verti @. (f^{\natural}{\mathbb E}_g\cMb)(x) \\[-2pt]
@.  @V{\text{using \eqref{eq:fin-sm1}}}VV \\[-4pt]
((gf)^{\natural}\cMb)(x) @<{\quad C^{\natural}_{f,g}\quad}<< 
(f^{\natural}g^{\natural}\cMb)(x)
\end{CD}
\end{equation}

Before we verify that \eqref{eq:fin-fact1} commutes we need 
to set up some notation. 
For any map $\phi\colon R \xto{\quad}S$ of 
noetherian local rings that is formally smooth 
of relative dimension~$r$ and for any $R$-module $N$ we set
$\phi_{\ltimes}N \set H^r_{m_R}(N \otimes_R \omega_{\phi})$.
Set $M = \cMb(z)$.
Now we expand \eqref{eq:fin-fact1} as follows with maps described below.
\def\duma{\makebox[0pt]{$\xto{\hspace{12em}}$}}
\def\dumb{\rotatebox{-90}{\makebox[0pt]{$\xto{\qquad\qquad}$}}}
\def\dumbb{\rotatebox{-90}{\raisebox{1.5ex}{\makebox[0pt]
{$\xto{\qquad{\rotatebox{90}{$\scriptstyle \alpha$}}\qquad}$}}}}
\def\dumbbb{\rotatebox{-90}{\raisebox{2.5ex}{\makebox[0pt]
{$\xto{\qquad{\rotatebox{90}{$\scriptstyle \gamma_1$}}\qquad}$}}}}
\def\dumc{\makebox[0pt]{$ {\scriptstyle \Box_5}{\hspace{2em}} 
{\wh{f_x}}_{\sharp}{\wh{g_y}}_{\ltimes}M {\hspace{3em}}$}}
\def\dumd{\makebox[6em]
{$\quad {\scriptstyle \Box_2} \quad\quad 
{\wh{f_x}}_{\ltimes}{g_y}_{\ltimes}M
\qquad {\qquad}$}}
\def\dume{\makebox[0pt]{$\qquad\quad 
{\wh{f_x}}_{\ltimes}{\wh{g_y}}_{\ltimes}M
\qquad {\scriptstyle \Box_4}$}}
\def\dumf{\makebox[0pt]{$\qquad\qquad 
{\wh{f_x}}_{\sharp}{\wh{g_y}}_{\sharp}M
\qquad\quad {\scriptstyle \Box_6}$}}
\[
\begin{CD}
({\mathbb E}_{gf}\cMb)(x) @. \duma @. ({\mathbb E}_f{\mathbb E}_g\cMb)(x) \\
@VV{{\hspace{9em}} \Box_1}V  @. @VVV  \\
{(gf)_x}_{\ltimes}M @<<< {f_x}_{\ltimes}{g_y}_{\ltimes}M 
@<<< {f_x}_{\ltimes}(({\mathbb E}_g\cMb)(y)) \\
@. @VV{{\hspace{4em}} \Box_3}V @VVV \\
\dumb @. \dumd @<<< {\wh{f_x}}_{\ltimes}(({\mathbb E}_g\cMb)(y)) \\
@. @VVV @. \\
{\wh{(gf)_x}}_{\ltimes}M @<{\beta}<< \dume @. \dumbbb \\
@. @V{\gamma_2}VV @. \\
\dumbb @. \dumc @<<< {\wh{f_x}}_{\sharp}(({\mathbb E}_g\cMb)(y)) \\
@. @V{\delta_2}V{\qquad\qquad\quad}V @| \\
{\wh{(gf)_x}}_{\sharp}M @<<< \dumf @. (f^{\natural}{\mathbb E}_g\cMb)(x) \\
@| @V{\Box_7 \quad\qquad}VV @VV{\delta_1}V \\
((gf)^{\natural}\cMb)(x) @<<<  (f^{\natural}g^{\natural}\cMb)(x) 
@= (f^{\natural}g^{\natural}\cMb)(x) 
\end{CD}
\]
The rectangle $\Box_1$ is the same as the diagram in~\ref{prop:itloco1} 
via the isomorphism in~\ref{lem:diff7}. 
Consequently $\Box_1$ commutes. The vertical maps in
$\Box_2,\Box_3$ are the ones induced by going to 
completion (see \eqref{eq:itloco0}) while the remaining horizontal
maps are the obvious induced ones.
Consequently $\Box_2,\Box_3$ commute. 
In~$\Box_4$, $\gamma_1, \gamma_2$ 
are induced by the isomorphism 
${\wh{f_x}}_{\ltimes} \iso {\wh{f_x}}_{\sharp}$
defined exactly as in the map $\theta$ 
in \eqref{eq:fin-sm1}, 
(i.e., $\gamma_i =$ $(-1)^{(p+d)q+p}$ times the 
map induced by \ref{thm:huang1}, I.(i)). 
The horizontal maps being the obvious ones, $\Box_4$ 
commutes. In $\Box_6$, $\delta_1$ is the natural
map induced by the isomorphism~${\mathbb E}_g\cMb \iso g^{\natural}\cMb$
of~\eqref{eq:fin-sm1}.
The remaining maps in~$\Box_6$ spell out the 
definition in~\eqref{eq:fin-sm1}, in particular, $\delta_2$ is
$(-1)^{(q+e)r+q}$ times the isomorphism given 
by~\ref{thm:huang1},~I.(i).
Therefore $\Box_6$ commutes. The vertical maps in $\Box_7$ 
are equalities while the horizontal maps are obtained from the
comparison isomorphism $C^{f,g}_{\sharp}$ and so $\Box_7$
commutes. Finally in $\Box_5$, $\alpha$ is
another instance of the map $\theta$ in \eqref{eq:fin-sm1} 
so that~$\alpha$ is 
$(-1)^{(p+d+e)r+p}$ times the isomorphism of~\ref{thm:huang1}, I.(i).
Note that $\beta$, which is obtained from~\eqref{eq:itloco2},
is $(-1)^{(q-r)(p+d-q)}$ times the 
the bottom row of \ref{thm:huang1},~II.(i), (because
$q-r$ is the transcendence
degree of the induced map of residue fields $k(z) \to k(y)$ and
$p+d-q$ is the relative dimension of~$\OXx$ over~$\OYy$). 
From the commutative diagram in \ref{thm:huang1},~II.(i) and the 
following calculation 
\[
(-1)^{(p+d+e)r+p}(-1)^{(q-r)(p+d-q)} = (-1)^{(p+d)q+p}(-1)^{(q+e)r+q},
\]
we conclude that $\Box_5$ commutes. This proves 
that \eqref{eq:fin-fact1} commutes.

By first restricting to connected components, if necessary, the
the condition on constant relative dimension for $f,g$ is relaxed, 
thus proving Case 1 in general.
As a corollary 
we now have a $\Coz$-valued pseudofunctor $(-)^{\sharp}$ over the  
subcategory of smooth maps in $\bbFc$.

{\textbf{Case 2.}} \emph{Suppose $f,g$ are closed immersions in $\bbFc$ 
such that the composition $gf$ exists and hence is also a closed 
immersion. Then  $(gf)^{\sharp} = (f|g)^{\sharp}$.} 
For this it suffices
to show that in the situation of~\ref{prop:clim3} 
the following diagram commutes. 
\[
\begin{CD}
(f^{\flat}g^{\flat}\cMb)(x) @>>> ((gf)^{\flat}\cMb)(x) \\[2pt]
@VVV  @VVV  \\[-3pt]
(f^{\natural}g^{\natural}\cMb)(x) @>{C^{\natural}_{f,g}}>> 
((gf)^{\natural}\cMb)(x) 
\end{CD}
\]
The commutativity of this diagram follows from~\ref{prop:clim3} 
and the commutative diagram in~$(\wh{\ref{eq:itloco3}})$. 
As a corollary we now have $(-)^{\sharp}$ as a pseudofunctor 
over the subcategory of closed immersions in~$\bbFc$.

{\textbf{Case 3.}} 
\emph{For the diagram in \eqref{eq:cartsq}$,$ we claim that 
$(j|f)^{\sharp} = (g|i)^{\sharp}$.} It suffices to show
that in the situation of \ref{prop:cartsq9}, with
$\cL = \omega_f$, the following diagram commutes. 
\begin{equation}
\label{eq:fin-fact3}
\begin{CD}
({\mathbb E}_gi^{\flat}\cMb)(w) @>>> (j^{\flat}{\mathbb E}_f\cMb)(w) \\[-2pt]
@VVV  @VVV  \\[-2pt]
(g^{\natural}i^{\flat}\cMb)(w) @. (j^{\natural}{\mathbb E}_f\cMb)(w) \\[-2pt]
@VVV  @VVV  \\[-2pt]
(g^{\natural}i^{\natural}\cMb)(w) @>>> (j^{\natural}f^{\natural}\cMb)(w)
\end{CD}
\end{equation}
To prove this we expand \eqref{eq:fin-fact3} as in \eqref{eq:fin-fact4}
with the following notation.
For formally smooth maps of noetherian local rings we use the 
notation of $(-)_{\ltimes}$ as in Case 1.
For any surjective map $\psi \colon R \xto{\quad} S$ of 
noetherian local rings and for any $R$-module~$N$ we set
$\psi_{\rtimes}N \set \Hom_R(S,N)$. In what follows, 
those maps of~\eqref{eq:fin-fact4}
which are not specifically described are assumed 
to have obvious meanings.

\begin{figure}
\rotatebox{90}
{\begin{minipage}{8.5in}\vspace{-10mm}
\begin{equation}\minCDarrowwidth=.35in
\label{eq:fin-fact4}
\begin{CD}
({\mathbb E}_gi^{\flat}\cMb)(w) @. @. 
{\makebox[0pt]{${\hspace{9em}}\xto{\hspace{36em}}$}}@. @. @. 
(j^{\flat}{\mathbb E}_f\cMb)(w) \\
@VV{\hspace{22em}\Box_1}V @. @. @. @. @VVV  \\
{g_w}_{\ltimes}((i^{\flat}\cMb)(z)) @>>> {g_w}_{\ltimes}{i_z}_{\rtimes}M
@>>> {j_w}_{\rtimes}{f_x}_{\ltimes}M @>{\epsilon}>> 
{j_w}_{\rtimes}{\wh{f_x}}_{\ltimes}M @>{\epsilon^{-1}}>> 
{j_w}_{\rtimes}{f_x}_{\ltimes}M
@>>> {j_w}_{\rtimes}(({\mathbb E}_f\cMb)(x)) \\
@VVV @VV{\hspace{8em}\Box_2}V @. @VVV @VVV @VVV  \\
{\wh{g_w}}_{\ltimes}((i^{\flat}\cMb)(z)) @>>> 
{\wh{g_w}}_{\ltimes}{i_z}_{\rtimes}M
@>>> {\wh{g_w}}_{\ltimes}{\wh{i_z}}_{\rtimes}M @>>> 
{\wh{j_w}}_{\rtimes}{\wh{f_x}}_{\ltimes}M @>>> 
{\wh{j_w}}_{\rtimes}{f_x}_{\ltimes}M
@>>> {\wh{j_w}}_{\rtimes}(({\mathbb E}_f\cMb)(x)) \\
@VV{\mu_1}V @VV{\mu_2}V @VV{\mu_3}V @VVV @VVV @VVV  \\
{\wh{g_w}}_{\sharp}((i^{\flat}\cMb)(z)) @>>> 
{\wh{g_w}}_{\sharp}{i_z}_{\rtimes}M
@>>> 
\makebox[0pt]{$\qquad\quad {\wh{g_w}}_{\sharp}{\wh{i_z}}_{\rtimes}M 
\qquad {\scriptstyle\Box_3}$} 
@. {\wh{j_w}}_{\sharp}{\wh{f_x}}_{\ltimes}M @>>> 
{\wh{j_w}}_{\sharp}{f_x}_{\ltimes}M
@>>> {\wh{j_w}}_{\sharp}(({\mathbb E}_f\cMb)(x)) \\
@| @. @VVV @VV{\mu_4}V @. @|  \\
\makebox[0pt]{$\hspace{7em}(g^{\natural}i^{\flat}\cMb)(w)
\hspace{6em}{\scriptstyle\Box_4}$} 
@. @. {\wh{g_w}}_{\sharp}{\wh{i_z}}_{\sharp}M @>{\nu}>> 
\makebox[0pt]{$\hspace{7em}{\wh{j_w}}_{\sharp}{\wh{f_x}}_{\sharp}M 
\hspace{6em}{\scriptstyle\Box_5}$} 
@. @. (j^{\natural}{\mathbb E}_f\cMb)(w) \\ 
@VV{\kappa_1}V @. @| @| @. @VV{\kappa_2}V \\
(g^{\natural}i^{\natural}\cMb)(w) @. 
\makebox[0pt]{{\raisebox{.5ex}{\makebox[10em]{\hrulefill}}}
\hspace{-10em}\!\!{\makebox[10em]{\hrulefill}}}
@. (g^{\natural}i^{\natural}\cMb)(w)
@>>>
(j^{\natural}f^{\natural}\cMb)(w) @. 
\makebox[0pt]{{\raisebox{.5ex}{\makebox[10em]{\hrulefill}}}
\hspace{-10em}\!\!{\makebox[10em]{\hrulefill}}}
@. (j^{\natural}f^{\natural}\cMb)(w) 
\end{CD}
\end{equation}
\end{minipage}\hspace{-12mm}
}
\end{figure}


Modulo the relation $\epsilon\epsilon^{-1} = 1$, 
the diagram in $\Box_1$ is the same as 
that in~\ref{prop:cartsq9} 
and hence $\Box_1$ commutes. 
The maps in $\Box_2$ are obtained from completion and
the isomorphism of~\eqref{eq:itloco4}, so $\Box_2$ commutes.
The maps $\mu_1,\mu_2,\mu_3,\mu_4$ are 
induced by the isomorphisms
${\wh{g_w}}_{\ltimes} \iso {\wh{g_w}}_{\sharp}$
and ${\wh{f_x}}_{\ltimes} \iso {\wh{f_x}}_{\sharp}$
which, in turn, are 
defined exactly as~$\theta$ is defined in~\eqref{eq:fin-sm1} 
(and hence equal $(-1)^{(p+d)q+p}$ times the corresponding map
induced by~\ref{thm:huang1},~I.(i)). 
We obtain~$\nu$ 
from the pseudofunctoriality of~$\sharp$. From the commutative 
diagram in~\ref{thm:huang1},~II.(iii), it follows 
that~$\Box_3$ commutes. 
In~$\Box_4$, $\kappa_1$ is
induced by an application of \eqref{eq:fin-clim1} while $\kappa_2$
in~$\Box_5$ is induced from \eqref{eq:fin-sm1}. 
The remaining maps in~$\Box_4,\Box_5$
spell out the definition of the maps~$\kappa_i$ and hence 
$\Box_4,\Box_5$ commute. The remaining rectangles in~\eqref{eq:fin-fact4}
commute due to functorial reasons. Thus~\eqref{eq:fin-fact4}, and 
hence \eqref{eq:fin-fact3}, commutes.

\enlargethispage*{3pt}
{\textbf{Case 4.}} 
\emph{Suppose $\Y \xto{\;i\;} \X \xto{\;h\;} \Y$ is a factorization
of the identity map on~$\Y$ with~$i$ a closed immersion
and $h$ a smooth map in $\mathbb F$. Then  
$(i|h)^{\sharp} = (1_{\Y})^{\sharp}$.} For this it suffices to show that 
in the situation of \ref{prop:ret1a}, for any $y \in \Y$,
the diagram 
\begin{equation}
\label{eq:fin-fact5}
\begin{CD}
(i^{\flat}\Esssh\cMb)(y) @>{\eqref{eq:ret0b}}>> \cMb(y) \\[-2pt]
@V{\eqref{eq:fin-sm1} = \eqref{eq:fin-clim1}}VV  
 @VV{\delta^{\natural}_y}V  \\[-6pt]
(i^{\natural}h^{\natural}\cMb)(y) @>{C^{\natural}_{i,h}}>> 
 ((1_{\Y})^{\natural}\cMb)(y)\\[1pt]
\end{CD}
\end{equation}
commutes, which can be done along the same lines as the earlier cases.

\pagebreak[3]


Before we put together the four cases above, we will need some 
basic results involving the definition in~\ref{def:fin-fact0}.
For the rest of this subsection, 
we do not write the Cousin complexes that occur
as inputs of the functors we work with. We also drop
reference to the codimension functions. We call a sequence
$\Bf = f_1, \ldots, f_n$ of maps in $\mathbb F$ \emph{composable} 
if, for each~$i$, the target of $f_i$ equals the source 
of~$f_{i+1}$. For any~$i$, $\Bf_{\le i}$ denotes the sequence
$f_1, \ldots, f_i$ and $\Bf_{> i}$ denotes the sequence 
$f_{i+1}, \ldots, f_n$.

\begin{alem}
\label{lem:fin-fact5a}
Let $\Bf = f_1, \ldots, f_n$ be a composable sequence in $\mathbb F$,
such that each~$f_i$ is a smooth map or a closed immersion. 
\begin{enumerate}
\item For any $i$, the canonical graded isomorphism 
\[
(\Bf_{\le i})^{\sharp}(\Bf_{> i})^{\sharp}
= (f_i\cdots f_1)^{\natural}(f_n \cdots f_{i+1})^{\natural} \iso 
(f_n \cdots f_1)^{\natural} 
= \Bf^{\sharp} 
\]
is a map of complexes.
\item Suppose that for some $i$ there exist maps $g_i, g_{i+1}$ 
in $\mathbb F$ such that 
\begin{enumerate}
\item $g_{i+1}g_i = f_{i+1}f_i$ and $g_i, g_{i+1}$ are 
smooth maps or closed immersions;
\item $(g_i| g_{i+1})^{\sharp} = (f_i| f_{i+1})^{\sharp}$.
\end{enumerate}
If $\bg$ denotes the (composable) sequence obtained
by replacing $f_i,f_{i+1}$ in~$\Bf$ by~$g_i,g_{i+1}$ respectively,
then $\Bf^{\sharp} = \bg^{\sharp}$.
\end{enumerate}
\end{alem}
\begin{proof}
The following diagram, whose maps are the canonical graded isomorphisms, 
commutes by pseudofunctoriality. Since the vertical maps 
define the differentials on $\Bf_{\le i},\Bf_{> i}$, (i) follows.
\[
\begin{CD}
(f_i\cdots f_1)^{\natural}(f_n \cdots f_{i+1})^{\natural} @>>>
 (f_n \cdots f_1)^{\natural} \\ 
@VVV @VVV \\
f_1^{\natural}\cdots f_i^{\natural}f_{i+1}^{\natural}\cdots f_n^{\natural} 
 @= f_1^{\natural}\cdots f_n^{\natural} 
\end{CD}
\] 

For (ii), consider the following diagram, whose maps are 
the canonical isomorphisms induced by $(-)^{\natural}$.
\[
\begin{CD}
(f_{i+1}f_i)^{\natural} @= (g_{i+1}g_i)^{\natural} \\
@VVV  @VVV \\
f_i^{\natural}f_{i+1}^{\natural} @>>> g_i^{\natural}g_{i+1}^{\natural} 
\end{CD}
\]
By assumption, if the objects in the bottom row are equipped with a 
differential by using~\eqref{eq:fin-sm1} and~\eqref{eq:fin-clim1},
then the bottom row is a map of complexes. 
Therefore, the analogous statement holds for the
following diagram.
\[
\begin{CD}
(f_n \cdots f_1)^{\natural} @= (f_n \cdots f_1)^{\natural} \\
@VVV  @VVV \\
f_1^{\natural} \cdots f_i^{\natural}f_{i+1}^{\natural} \cdots f_n^{\natural} 
 @>>> f_1^{\natural} \cdots f_{i-1}^{\natural} 
 g_i^{\natural}g_{i+1}^{\natural}f_{i+2}^{\natural}\cdots f_n^{\natural} 
\end{CD}
\]  
Thus (ii) follows.
\end{proof}

\medskip

\begin{aprop}
\label{prop:fin-fact6}
Let $\X \xto{\; i_1 \;} \Z_1 \xto{\; h_1 \;} \Y$ and
$\X \xto{\; i_2 \;} \Z_2 \xto{\; h_2 \;} \Y$ be maps in~$\mathbb F$
such that $h_j$ are smooth maps, $i_j$ closed immersions and 
$h_1i_1 = h_2i_2$. Assume further that $h_1i_1$ is separated.
Then $(i_1|h_1)^{\sharp} = (i_2|h_2)^{\sharp}.$
\end{aprop}
\begin{proof}
Consider the following diagram with obvious choices for the 
maps and notation.
\begin{equation}
\label{eq:fin-fact7}
\begin{CD}
\X @>{d}>> \X \times_{\Y} \X @>{i_1''}>> 
\Z_1 \times_{\Y} \X @>{h_1''}>> \X \\
@. @VV{i_2''}V @VV{i_2'}V @VV{i_2}V \\
@. \X \times_{\Y} \Z_2 @>{i_1'}>> \Z_1 \times_{\Y} \Z_2 @>{h_1'}>> \Z_2 \\
@. @VV{h_2''}V @VV{h_2'}V @VV{h_2}V \\
@. \X @>{i_1}>> \Z_1 @>{h_1}>> \Y
\end{CD}
\end{equation}
Separatedness of $h_1i_1$ and $h_2''i_2''$ implies that 
$d$ is a closed immersion and
by Case~2 and Case~4 above we conclude that 
$(d\,| i_2''|h_2'')^{\sharp} = 1_{\X}^{\sharp} = 
(d\,| i_1''|h_1'')^{\sharp}$. Therefore we obtain
the following isomorphisms of complexes
\begin{align}
1_{\X}^{\sharp}(i_1| h_1)^{\sharp} 
 = (d\,| i_2''|h_2'')^{\sharp}(i_1| h_1)^{\sharp}  
&\xto{\; \alpha\;}(d\,| i_2''|h_2''| i_1| h_1)^{\sharp} \notag \\
&\xto{\; \beta\;} (d\,| i_1''|h_1''| i_2| h_2)^{\sharp} \notag \\
&\xto{\; \gamma\;} (d\,| i_1''|h_1'')^{\sharp}(i_2| h_2)^{\sharp} = 
 1_{\X}^{\sharp}(i_2| h_2)^{\sharp} 
 \notag  
\end{align}
where $\alpha, \gamma$ are obtained using \ref{lem:fin-fact5a}, (i), 
while $\beta$ is defined as follows. We begin from the southwest corner
of~\eqref{eq:fin-fact7}. By Case 3 above applied to the bottom left
rectangle in~\eqref{eq:fin-fact7} and by \ref{lem:fin-fact5a}, (ii),
we may replace the subsequence $h_2'', i_1$ 
by~$i_1',h_2'$. Proceeding in this 
manner for the remaining rectangles in~\eqref{eq:fin-fact7} 
(using Cases~1,~2,~3 above and \ref{lem:fin-fact5a}, (ii)),
we reach the northeast corner of~\eqref{eq:fin-fact7} to complete 
the definition of~$\beta$. 

Since $\alpha$, $\beta$, $\gamma$ are all equalities
and $1_{\X}^{\sharp}$ is isomorphic to the identity
functor on~$\Coz(\X)$ 
we see that $(i_1| h_1)^{\sharp} = (i_2| h_2)^{\sharp}$. 
\end{proof}

\subsection{Constructing $(-)^{\sharp}$ for a general map}
\label{subsec:fin-gen}
In order to extend the preceding results for factorizable
maps as in \ref{prop:fin-fact6} to arbitrary maps 
in~$\mathbb F$ we first need a localization 
result. For the following Lemma, we drop reference to the 
codimension functions and the Cousin complexes that occur 
as inputs to the functors involved.

\begin{alem}
\label{lem:fin-gen1}
Let $\U \xto{\;i\;} \Z \xto{\;h\;} \Y$ be maps in $\mathbb F$
where $i$ is a closed immersion and $h$ is a smooth map.
Let $\U'$ be an open subset of $\U$ so that there exists 
some open subset $\Z'$ of~$\Z$ for which $i^{-1}\Z' = \U'$.
Let $\U' \xto{\;i'\;} \Z' \xto{\;h'\;} \Y$
denote the induced maps. Then 
$(i|h)^{\sharp}\big|_{\U'} = (i'|h')^{\sharp}$.
\end{alem}
\begin{proof}
This follows immediately from the fact that
over smooth maps and closed immersions, 
$(-)^{\sharp}$ behaves well with respect to restriction
to open subsets (cf.~\ref{lem:fin-sm3}).  
\end{proof}

\medskip

We are now in a position to define $(-)^{\sharp}$ of the Main Theorem
over the whole of~$\bbFc$. Let  
$f \colon (\X, \DsssX) \to (\Y, \DsssY)$ be a map in $\bbFc$. 
Let $\sB = \{\U_{\lambda}\}_{\lambda \in \Lambda}$ 
be the collection of all open subsets of~$\X$ such 
that for any $\lambda$, the induced 
map $f_{\lambda} \colon \U_{\lambda} \to \Y$ admits a factorization 
$\U_{\lambda} \xto{\; i_{\lambda} \;} \Z_{\lambda} \xto{\; h_{\lambda} \;}\Y$ 
where~$i_{\lambda}$ is a closed 
immersion and~$h_{\lambda}$ a separated smooth map. By \ref{cor:morph1b}, 
$\sB$ forms a basis for open sets in $\X$.
Fix a complex $\cMb$ in $\Coz_{\DsssY}(\Y)$. 
For any $\U_{\lambda} \in \sB$ and for any factorization 
$f_{\lambda} = h_{\lambda}i_{\lambda}$, set 
$f_{\lambda}^{\sharp}\cMb \set (i_{\lambda}|h_{\lambda})^{\sharp}\cMb$.
By \ref{prop:fin-fact6}, $f_{\lambda}^{\sharp}\cMb$
does not depend on the choice of 
$i_{\lambda},h_{\lambda}$. 
For any $\U_{\lambda} \in \sB$, any open subset 
$\U'$ of $\U_{\lambda}$ is also in $\sB$ and by~\ref{lem:fin-gen1},
with $\U_{\lambda'} \set \U'$ we have
$(f_{\lambda}^{\sharp}\cMb) \big|_{\U_{\lambda'}} 
= (f_{\lambda'}^{\sharp}\cMb)$.
It follows that the differentials\vspace{-1.2pt}
of $(f_{\lambda}^{\sharp}\cMb)$ for $\lambda \in \Lambda$
can be pasted together to yield a differential on~$f^{\natural}\cMb$. 
We take this \emph{canonical} choice as the differential 
of~$f^{\sharp}\cMb$.

It follows that our definition of $f^{\sharp}$ satisfies
properties (ii) and~(iii) of the 
Main Theorem in~\S\ref{subsec:outline}. 
We now address property (i) which amounts to verifying  
that the graded constructs
of~$(-)^{\natural}$ work at the level of complexes too. 

Let us verify that for any map $f \colon (\X, \DsssX) \to (\Y, \DsssY)$ 
in $\bbFc$, $f^{\sharp}$ is a functor of complexes, i.e., 
for any map of complexes 
$\cMb \to \cNb$ in $\Coz_{\DsssY}(\Y)$, the  
canonical graded map $f^{\sharp}\cMb \to f^{\sharp}\cNb$
is also a map of complexes. The verification is local in 
nature and so we reduce to verifying at open sets of the 
type $\U_{\lambda}$ in~$\sB$; here the desired result 
is obvious. 

Under the hypothesis of~(i)(b) (resp.~(i)(c)) of the Theorem, 
we set $C^{\sharp}_{f,g} = C^{\natural}_{f,g}$ 
(resp.~$\delta^{\sharp}_{\X} = \delta^{\natural}_{\X}$). 
From \ref{rem:fin-clim2} it follows that 
for any $(\X, \Delta) \in \mathbb F$ and
$\cMb \in \Coz_{\Delta}(\X)$, the natural
transformation $\delta^{\sharp}_{\X}$ induces 
a map of complexes $1_{\X}^{\sharp}\cMb \iso \cMb$. 
The rest of this subsection is devoted to 
the case of the comparison maps. 

Let $$(\X,\DsssX) \xto{\; f \;} (\Y, \DsssY) 
\xto{\; g \;} (\Z, \DsssZ)$$ be maps in~$\bbFc$ and 
let $\cMb$ be a complex in~$\Coz_{\DsssZ}(\Z)$. We need to
verify that the canonical graded isomorphism
\[
C^{\natural}_{f,g}(\cMb) \colon 
f^{\sharp}g^{\sharp}\cMb \xto{\quad} (gf)^{\sharp}\cMb
\] 
is a map of complexes. Since the definition of $(-)^{\natural}$
is local in nature and since the definition of the differential
for $(-)^{\sharp}$ is based on local constructions,
we may assume without loss of generality that 
$\X,\Y$ and $\Z$ are affine formal schemes,
say $\X= \Spf(C), \Y= \Spf(B)$ and $\Z = \Spf(A)$.
By \ref{lem:morph1} we may factor the natural map $A \to B$ as
\[
A \to P = ((A[X_1,\ldots,X_n])_S, I)^{\wedge} \stackrel{\pi}{\onto} B
\]
and the natural map $B \to C$ as 
\[
B \to Q=((B[Y_1,\ldots,Y_m])_T, J)^{\wedge} \onto C.
\] 
Set $\pfr = \ker\pi$. 
Let $T'\set \pi_Y^{-1}T$ where $\pi_Y\colon P[Y_1,\ldots,Y_m] \onto
B[Y_1,\ldots,Y_m]$ is the map naturally 
induced by $\pi$. Let $J'$ be the inverse image of $J$ under the induced  
surjection $(P[Y_1,\ldots,Y_m])_{T'} \onto (B[Y_1,\ldots,Y_m])_T$;
the kernel of this surjection is also generated by $\pfr$. 
Set $R \set ((P[Y_1,\ldots,Y_m])_{T'}, J')^{\wedge}.$ Then  
\[
Q = ((P/{\pfr}[Y_1,\ldots,Y_m])_T, J)^{\wedge} \cong 
((P[Y_1,\ldots,Y_m])_{T'}/{(\pfr)}, J)^{\wedge}
\cong  R/(\pfr R).
\]
Thus the following statements hold:
\begin{itemize}
\item As a map of noetherian adic rings the natural map $P\to R$ is 
essentially of pseudo-finite type and formally smooth. 
\item $\Spf(Q)$ is the fibered product 
$\Spf(B) \times_{\Spf(P)} \Spf(R)$ in $\mathbb F$.
\end{itemize}
Now consider the following diagram of canonical maps where the 
vertical ones are smooth and the horizontal ones are
closed immersions. 
\[
\begin{CD}
\Spf(C) @>i_2>> \Spf(Q) @>{i_1^{\prime}}>> \Spf(R) \\
@.  @V{h_2^{\prime}}VV  @V{h_2}VV \\
@.   \Spf(B) @>i_1>> \Spf(P) \\
@.  @.  @V{h_1}VV \\
@.  @.  \Spf(A)
\end{CD}
\]
Note that $f={h_2'}i_2$, $g=h_1i_1$.
By construction, we have 
\[
f^{\sharp} = (i_2|h_2')^{\sharp}, \quad
g^{\sharp} = (i_1|h_1)^{\sharp}, \quad
(gf)^{\sharp} = (i_1'i_2|h_1h_2)^{\sharp}. 
\]
To prove that $C^{\natural}_{f,g}$ is a map of complexes
consider the following commutative diagram of canonical
graded isomorphisms.
\[
\begin{CD}
(h_2'i_2)^{\natural}(h_1i_1)^{\natural}  @.
 \makebox[0pt]{$\xto{\hspace{4em}C^{\natural}_{f,g}\hspace{4em}}$}  
 @. (h_1h_2i_1'i_2)^{\natural} \\ 
@VVV @. @VVV \\
i_2^{\natural}h_2'{}^{\natural}i_1^{\natural}h_1^{\natural} @>{\alpha}>>
 i_2^{\natural}i_1'{}^{\natural}h_2^{\natural}h_1^{\natural}
 @>{\beta}>> (i_1'i_2)^{\natural}(h_1h_2)^{\natural} 
\end{CD}
\]
The bottom row, 
when equipped with differentials as in~\S\ref{subsec:fin-sm}
and~\S\ref{subsec:fin-clim}, is seen to be a map of complexes
by using Case~3 of~\S\ref{subsec:fin-fact} for~$\alpha$ 
and Cases~1,~2 for~$\beta$. Therefore~$C^{\natural}_{f,g}$ is
a map of complexes.

\subsection{Compatibility with translations}
\label{subsec:trans}

For any map $f \colon (\X, \Delta) \to (\Y, \Delta')$ 
in $\bbFc$ and for any integer $n$, we denote 
the obvious map $(\X, \Delta-n) \to (\Y, \Delta'-n)$ by~$f^{(n)}$. 
Let $\cMb \in \Coz_{\Delta'}(\Y)$. Then
$\cMb[n] \in \Coz_{\Delta'-n}(\Y)$ and both\vspace{-1.6pt}
$(f^{\sharp}\cMb)[n]$ and ${f^{(n)}}^{\sharp}(\cMb[n])$ are in 
$\Coz_{\Delta-n}(\X)$. The latter two also have the same underlying
graded objects and so we are naturally led to comparing their 
differentials.

Consider the graded isomorphism 
\begin{equation}
\label{eq:trans1}
{f^{(n)}}^{\sharp}(\cMb[n]) \iso  (f^{\sharp}\cMb)[n]
\end{equation}
defined punctually, say at a point $x \in \X$, with $y = f(x)$, by
\begin{equation}
\label{eq:trans2}
({f^{(n)}}^{\sharp}(\cMb[n]))(x) = \wh{f_x}^{\sharp}(\cMb(y)) 
\xto{(-1)^{nt}}
\wh{f_x}^{\sharp}(\cMb(y)) = ((f^{\sharp}\cMb)[n])(x),
\end{equation}
where $t$ is the transcendence degree of the residue field extension
$k(y) \to k(x)$.

\begin{aprop}
\label{prop:trans3}
The graded isomorphism in \eqref{eq:trans1} is an isomorphism of complexes. 
Morever the following hold.
\begin{enumerate}
\item If $m$ is an integer, then the following diagram of isomorphisms
commutes, where the horizontal arrows are obtained using~\eqref{eq:trans1}.
\[
\begin{CD}
{f^{(m+n)}}^{\sharp}(\cMb[m+n]) @. 
\makebox[0pt]{$\xto{\hspace{14em}}$} @. (f^{\sharp}\cMb)[m+n] \\
@| @. @| \\
{f^{(m)(n)}}^{\sharp}(\cMb[m][n]) @>>> ({f^{(m)}}^{\sharp}(\cMb[m]))[n]
@>>> (f^{\sharp}\cMb)[m][n]  
\end{CD} 
\]
\item If $g \colon (\Y, \Delta') \to (\Z, \Delta'')$ is a 
map in $\bbFc$ then the following diagram of isomorphisms
commutes, where the horizontal arrows are obtained using~\eqref{eq:trans1}.
\[
\begin{CD}
{f^{(n)}}^{\sharp}{g^{(n)}}^{\sharp}(\cMb[n]) @>>> 
 {f^{(n)}}^{\sharp}((g^{\sharp}\cMb)[n]) @>>> 
 (f^{\sharp}g^{\sharp}\cMb)[n] \\
@V{\textup{via } C^{\sharp}_{f^{(n)},g^{(n)}}}VV @. 
 @VV{\textup{via }C^{\sharp}_{f,g}}V \\
(g^{(n)}f^{(n)})^{\sharp}(\cMb[n]) @= {(gf)^{(n)}}^{\sharp}(\cMb[n])
@>>> ((gf)^{\sharp}\cMb)[n]  
\end{CD} 
\]
\end{enumerate}
\end{aprop}

\begin{proof}
The horizontal maps in (i) and (ii) are, at the punctual level, 
signed multiples of identity maps, according to the sign $(-1)^{nt}$ in \eqref{eq:trans2}. Hence (i) follows from the additivity of~$nt$ vis-\`{a}-vis
$n$ for fixed $t$, while (ii) follows from the additivity vis-\`{a}-vis~$t$ for fixed $n$.  

To verify that \eqref{eq:trans1} is a map of complexes 
it suffices to do so locally on~$\X$ and~$\Y$ and hence we may assume 
that~$f$ factors as a 
closed immersion into a smooth map. Thus it suffices to assume that $f$
is either a closed immersion or a smooth map. The first case is easily
settled, using, e.g., \eqref{eq:clim1a}. 
(In particular, sign considerations play no role.)
Assume then that $f$ is a smooth map.
 
We use notation as in \S\ref{subsec:fin-sm}. From our convention 
in \S\ref{subsec:conv}\eqref{conv4}, we see that 
\[ f^*\cMb[n] \otimes_{\X} \omega_f[d] = 
(f^*\cMb \otimes_{\X} \omega_f[d])[n]. 
\]

Recall that for any 
complex $\cFb$, $E_{\Delta-n}(\cFb[n]) = (E_{\Delta}\cFb)[n]$
(see ~\eqref{eq:Cousin+trans}). 
It follows that  
\[
{\mathbb E}_{f^{(n)}}(\cMb[n]) = ({\mathbb E}_f\cMb)[n].
\]
To conclude the smooth case, it suffices to show that the
preceding equality is consistent with the isomorphism of \eqref{eq:trans1}, 
i.e., for any $x \in \X$,
the following diagram of punctual isomorphisms commutes. 
\[
\begin{CD}
({\mathbb E}_{f^{(n)}}(\cMb[n]))(x) @= (({\mathbb E}_f\cMb)[n])(x) \\
@V{\text{via \eqref{eq:fin-sm1}}}VV  @VV{\text{via \eqref{eq:fin-sm1}}}V \\
\wh{f_x}^{\sharp}(\cMb(y)) @>{(-1)^{nt}}>> \wh{f_x}^{\sharp}(\cMb(y)) 
\end{CD}
\]
Using $p = \Delta(x), q = \Delta'(y)$ as in \S\ref{subsec:fin-sm},
we have that $t = p-q$. After some diagram chase, using 
the choices of signs involved, one reduces 
the commutativity to the following immediately verifiable calculation.
\[
(-1)^{(p-n)(q-n)+ p-n} (-1)^{n(p-q)} (-1)^{pq +p} = 1
\] 
\end{proof}

\newpage

\section{Residual and dualizing complexes}
\label{sec:residual}

In this section we discuss residual complexes\index{residual complex}
on a formal scheme~$\X$ (Definition \ref{defi:residual}.)%
\footnote
{This definition
agrees with the one for ordinary schemes 
in \cite[VI, \S1]{RD}, and also, for formal schemes of finite 
Krull dimension, with the one in \cite[5.9]{Ye}. 
(Yekutieli's\index{Yekutieli, Amnon} residual complexes are
dualizing, hence cannot exist on formal schemes of infinite
Krull dimension \cite[p.\,283, Cor.\,7.2]{RD}.)
}  
These special Cousin complexes
play an important role in duality theory. They      
can be  characterized as being
the Cousin complexes associated to \emph{pointwise dualizing}
complexes in the derived category~$\Dqct(\X)$ (Proposition 
\ref{prop:residual6}). For any $\bbFc$-map $f$,
$f^\sharp$ takes residual complexes to residual complexes 
(Proposition \ref{prop:residual4}).

On finite-type schemes over finite-dimensional regular noetherian rings, 
Yekutieli\index{Yekutieli, Amnon} and Zhang\index{Zhang, James J.} have found a
duality between the category of coherent sheaves and the category of Cohen-Macaulay complexes with coherent homology \cite[Thm.\,8.9]{YZ}. 
This duality generalizes to all formal schemes having a 
coherent-dualizing complex (Corollary~\ref{cor:CM+coh}). 
(For the formal spectrum of a complete local ring it is essentially
Matlis duality, see Example~\ref{Matlis}.)  There follow quick proofs of a number of facts
about Cousin complexes, for example a
generalization to complexes on formal schemes,
of results of Dibaei, Tousi and Kawasaki\index{Dibaei,
Mohammad Taghi}\index{Tousi, Masoud}\index{Kawaski, Takesi}
about finite generation of the homology of Cousin complexes of finite
type modules  (see Proposition~\ref{E=EHn}).

\enlargethispage*{3pt}

\subsection{$(-)^\sharp$ preserves residual}
\label{subsec:residual to residual}
\begin{adefi}
\label{defi:residual}
We say that a complex~$\cRb$ on a noetherian formal scheme~$\X$ is a residual complex
if it consists of~$\At(\X)$-modules and
if for any defining ideal~$\I$, the complex
$\sHom_{\OX}(\OX/\I, \cRb)$ on the scheme $(\X, \OX/\I)$
is a residual complex. 
\end{adefi}

It is easily shown that~$\cRb$ 
is residual if there exists one defining 
ideal~$\I$ such that for any $n>0$, $\sHom_{\OX}(\OX/\I^n, \cRb)$ 
is residual on $(\X, \OX/\I^n)$. Note that the existence
of a residual complex on a noetherian formal scheme $\X$ 
implies that~$\X$ is in~$\mathbb F$ because the corresponding
statement for ordinary schemes is true.

\begin{alem}[cf.~{\cite[VI, Lemma 5.2]{RD}}]
\label{lem:residual1}
Let $A$ be a noetherian local ring with maximal ideal~$m_A$
and residue field $k$. Let $M$ be an $m_A$-torsion $A$-module. 
Suppose there exists an ideal~$I$ in~$A$ such that for 
any integer~$n>0$, $\Hom_A(A/I^n, M)$ 
is an injective hull of~$k$ over~$A/I^n$. Then $M$ is an 
injective hull of~$k$ over~$A$.
\end{alem}
\begin{proof}
Set $M_n \set \Hom_A(A/I^n, M) \subset M$.
Each $M_n$ is an essential extension of~$k$ and $M = \cup_n M_n$.
Therefore $M$ is also an essential extension of~$k$. In 
particular, there exists an embedding $M \subset E(k)$ 
where $E(k)$ is an injective hull of~$k$ over~$A$. 
Set $E_n \set \Hom_A(A/I^n, E(k)) \subset E(k)$. By hypothesis,
for each $n$, $M_n = E_n$. Therefore, $M=E(k)$.
\end{proof}

\enlargethispage*{3pt}

For the next result we use some basic facts on residual complexes
on ordinary schemes as developed in~\cite[VI, \S1]{RD} 
and~\cite[\S3.1, \S3.2]{BCo}.

\begin{alem}
\label{lem:residual3}
Let $\X$ be a noetherian formal scheme and $\cRb$
a residual complex on~$\X$. Then 
there is a unique codimension 
function~$\Delta$ on~$\X$ such that 
$\cRb \in \Coz_{\Delta}(\X)$. Moreover, for any $x \in \X$,
$\cRb(x)$ is an $\OXx$-injective hull of the residue 
field $k(x)$ at~$x$.
\end{alem}

\begin{proof}
Let $\cI$ be a defining ideal in $\OX$. For $n>0$, set 
\[\postdisplaypenalty=10000
\cRb_n \set \sHom_{\OX}(\OX/\cI^n, \cRb).
\]
Since $\cRb$ consists of $\At(\X)$-modules, it is 
isomorphic to the direct limit of the $\cRb_n$'s. On the ordinary
scheme $X_n \set (\X, \OX/\cI^n)$, the complex $\cRb_n$, 
being a residual complex, 
induces a unique codimension
function $\Delta_n$ for which $\cRb_n \in \Coz_{\Delta_n}(X_n)$.
Moreover, for every $x \in \X$, $\cRb_n(x)$ is an 
$\OXx/\cI_x^n$-injective hull of~$k(x)$. Since 
$\cRb_n \iso \sHom_{\OX}(\OX/\cI^n, \cRb_{n+1})$, 
$\cRb_n$ is also in $\Coz_{\Delta_{n+1}}(X_n)$ 
(cf.~\ref{prop:clim1}).In particular, $\Delta_n = \Delta_{n+1}$, 
which we henceforth denote by~$\Delta$. Moreover, the canonical 
inclusion $\cRb_n \to \cRb_{n+1}$
respects the pointwise decomposition. Taking 
limits we see that $\cRb$ is a $\Delta$-Cousin complex,
where for any $x$, $\cRb(x)$ is a direct limit of 
the~$\cRb_n(x)$'s. Since each $\cRb_n(x)$ is 
an $m_x$-torsion $\OXx$-module, so is~$\cRb(x)$. 
By~\ref{lem:residual1}, $\cRb(x)$ is an 
injective hull of~$k(x)$ over~$\OXx$.
\end{proof}

In light of ~\ref{lem:residual3},
if a codimension function~$\Delta$ is chosen
on~$\X$, then we shall only consider those residual
complexes which induce~$\Delta$. Thus, for $(\X, \Delta) \in \bbFc$,
a residual complex shall always be assumed, by default, 
to lie in~$\Coz_{\Delta}(\X)$. 

\begin{aprop}
\label{prop:residual4}
Let $f \colon (\X, \DsssX) \to (\Y, \DsssY)$ be a map 
in $\bbFc$. If $\cRb$ is a residual complex on~$\Y$, then
$f^{\sharp}\cRb$ is a residual complex on~$\X$. 
\end{aprop}
\begin{proof}
First we prove the result when $X \set \X$ and $Y \set \Y$ 
are ordinary schemes so that $f$ is essentially of finite type.
For any $x \in X$, $(f^{\sharp}\cRb)(x)$ is an 
$\cO_{X,x}$-injective hull of $k(x)$ at $x$ (\ref{thm:huang1}, IV).
So it only remains to check that $f^{\sharp}\cRb$ has coherent 
homology. This being a local property, 
in view of the existence of factorizations as in \ref{cor:morph1b},
it suffices to consider the cases when $f$ is a 
smooth map (of constant relative dimension $n$)
and when $f$ is a closed immersion.
In the smooth situation, by Cohen-Macaulayness 
of~$f^*\cRb \otimes_X \omega_f[n]$ 
(\ref{lem:quasim1}), there are isomorphisms
$f^*\cRb \otimes_X \omega_f[n] \iso E(f^*\cRb \otimes_X \omega_f[n])
\iso f^{\sharp}\cRb$ (see \ref{cor:coz4a} and~\eqref{eq:fin-sm1},
though a-priori, the first isomorphism is only a $\D(\X)$-isomorphism).
Since $\omega_f$ is coherent, $f^{\sharp}\cRb$ has coherent homology.
If $f$ is a closed immersion, in view of the isomorphism 
$f^{\sharp}\iso f^{\flat}$
of~\eqref{eq:fin-clim1}, the proof follows easily from the 
coherence of $f_*\OX$ (cf.~\cite[III, Prop.~6.1]{RD}).

Now we consider the general case, when $\X$, $\Y$
are arbitrary formal schemes. Fix a defining ideal 
$\I$ in~$\OX$ and a defining ideal 
$\J$ in~$\OY$ such that $\J\OX  \subset \I$.
Set $X_n \set (\X, \OX/\I^n), Y_n \set (\Y, \OY/\J^n)$.
Then there is a commutative diagram of natural 
induced maps as follows.
\[
\begin{CD}
X_n @>{f_n}>> Y_n \\
@VV{j_n}V  @VV{i_n}V \\
\X @>f>> \Y
\end{CD}
\]
Using the isomorphisms $i_n^{\sharp} \iso i_n^{\flat}$
and $j_n^{\sharp} \iso j_n^{\flat}$ obtained 
from~\eqref{eq:fin-clim1}, we deduce that
$\cRb$ is residual on~$\Y$ 
$\iff$ $\forall n$, $i_n^{\sharp}\cRb$ is residual on~$Y_n$ 
$\Lra$ $\forall n$, $f_n^{\sharp}i_n^{\sharp}\cRb$ is residual on~$X_n$
$\iff$ $\forall n$, $j_n^{\sharp}f^{\sharp}\cRb$ is residual on~$X_n$
$\Lra$ $f^{\sharp}\cRb$ is residual on~$\X$.
\end{proof}

\subsection{Residual and pointwise dualizing complexes}
\label{sec:ptwise dualizing}
This subsection is devoted to showing a pointwise
dualizing property for residual complexes.
Let $\X$ be a noetherian formal scheme in~$\mathbb F$.
For any point $x \in \X$, let $j_x$ denote the 
canonical map $\Spf(\wh{\OXx}) \to \X$
where $\wh{\OXx}$ is the completion of the local ring~$\OXx$
along the stalk $\I_x$ of a defining ideal $\I$ in $\OX$. 

\begin{alem}
\label{lem:residual5}
Let $\X, x, j \set j_x$ be as above. Set $\W \set \Spf(\wh{\OXx})$. Then$\>:$

{\rm(i)}  The map $j$ is adic, i.e., for any defining ideal $\cI$ in~$\OX$,
$\cI\OW$ is a defining ideal.
In particular, $j^*$ takes~$\At(\X)$ to~$\At(\W)$
and~$\Aqct(\X)$ to~$\Aqct(\W)$.
Moreover, for any $\F \in \At(\X)$, the natural map 
$j^{-1}\F \to j^*\F$ is an isomorphism.

{\rm(ii)} The map $j$ is \'etale\/
$($smooth of relative dimension\/ $0).$
In particular, $j^*$ is exact and 
$\bL j^* = j^*$ takes $\Dt(\X)$ to $\Dt(\W)$
and $\Dqct(\X)$ to $\Dqct(\W)$.

{\rm(iii)} If $w$ denotes the unique closed point of $\W$, 
then for any $\cFb \in \Dqct^+(\X)$, there is a natural 
isomorphism $H^i_wj^*\cFb \iso H^i_x\cFb$.  
\end{alem}

\begin{proof}
First we show that $j$ is adic and \'etale. 
Let $\V = \Spf(A)$ be an affine neighborhood of~$x$. It suffices
to show that the induced map $\W \to \V$ is adic and \'etale,
so we assume $\V = \X$.
Let $I$ be a defining ideal of~$A$.
Then $\I = I^{\sim_A}$ is a defining ideal in~$\OX$
and since $I\OXx = \I_x$ therefore~$I$ generates a defining ideal 
of~$\wh{\OXx}$. By~\ref{lem:mod4a} we see that the image
of~$j^*\I$ in~$\OW$ is the ideal generated by $I\wh{\OXx}$.
Thus $j$ is adic. For \'etaleness we refer to the discussion 
preceding \ref{prop:app3}.

For an adic map, $(-)^*$ sends $\At$-modules 
to $\At$-modules, as can be seen by arguing at the stalks
using \ref{lem:loc3}. Since $(-)^*$ also preserves
quasi-coherence, it sends $\Aqct$-modules to $\Aqct$-ones.
Exactness of $j^*$ follows from flatness of~$j$ which in turn
follows from smoothness of~$j$. Finally, since $j^{-1}$
and~$j^*$ commute with direct limits, to verify that
$j^{-1}\F \iso j^*\F$ for $\F \in \At(\X)$, it suffices
to do so when $\F$ is annihilated by a defining ideal $\I$
in~$\OX$. Now we may descend to the corresponding schemes
and here the result is clear. Thus (i) and~(ii) are proved.

For (iii), by \ref{prop:loc8}(ii), we may assume that
$\cFb$ consists of $\Aqct(\X)$-injectives. 
Then with $m_x, m_w$ as the local rings at $x, w$ 
respectively, there are natural
isomorphisms $H^i_{m_x}\cFb_x \iso H^i_x\cFb
\text{ and } H^i_{m_w}(j^*\cFb)_w \iso H^i_wj^*\cFb$ 
(see~\eqref{eq:nolab3} of \S\ref{subsec:loc}). 
Now (iii) follows by using $(j^*\cFb)_w = (j^{-1}\cFb)_w = \cFb_x$,
and that~$m_w$ is just the completion of~$m_x$.
\end{proof}

For the definition of a \emph{t-dualizing} (=\,\emph{torsion-dualizing}) complex on a noetherian
formal scheme~$\X$ we refer to~\cite[2.5.1]{AJL2}.\index{t-dualizing complex}
(See also \cite[\S5]{Ye}.)
An $\OX$-complex~$\cDb$ 
is \emph{pointwise $t$-dualizing}\index{pointwise t-dualizing complex} 
if $\cDb \in \Dqct^+(\X)$ and for any $x \in \X$, 
with $j_x, \W$ as in \ref{lem:residual5},
$j_x^*\cDb$ is $t$-dualizing on~$\W$. In view
of~\cite[3.1.4]{BCo}, this definition
agrees with the one for ordinary schemes .

We need some notation before proceeding further.
Let $\Y$ be a noetherian formal scheme, and $\I$ a defining ideal
in $\OY$. Let $i$ be the canonical closed immersion
$Y \set (\Y, \OY/\I) \to \Y$. Then for any
$\cFb \in \D(\Y)$ we set 
\[
i^!\cFb \set \R\sHom^{\bullet}_{\Y}(\OY/\I, \cFb) \in \D(Y).
\]
Here $\R\sHom^{\bullet}_{\X}(\OY/\I, \cFb)$ can be considered as 
a complex in $\D(Y)$ in a natural way.
Also, $i^!$ is a right adjoint of 
$i_* = \R i_* \colon \Dqc(Y) \to \D(\Y)$
(see \cite[Examples 6.1.3(4)]{AJL2}).

\begin{aprop}
\label{prop:residual6}
Let $\X$ be a noetherian formal scheme. 

{\rm(i)} Let $\cDb$ be a complex on~$\X$ and let $\I$ be a
defining ideal in $\OX$. For $n>0$, let $i_n$ denote the 
canonical immersion $X_n = (\X, \OX/\I^n) \to \X$. Then $\cDb$ 
is pointwise $t$-dualizing $\iff$ for any $n>0$, $i_n^!\cDb$
is pointwise dualizing on $X_n$.

{\rm(ii)} Let $\cDb$ be a complex on~$\X$. Then $\cDb$ 
is $t$-dualizing $\iff$ $\X$ has finite Krull dimension and
$\cDb$ is pointwise $t$-dualizing.

{\rm(iii)} Any residual complex on $\X$ is pointwise $t$-dualizing.
Conversely, if\/ $\cDb$ is a pointwise $t$-dualizing complex 
on~$\X$, then the following hold:
\begin{enumerate}
\item[(a)] The complex $\cDb$ induces a codimension 
function~$\Delta$ on~$\X$ such that~$\cDb$ is Cohen-Macaulay 
w.r.t.~$\Delta$. 
\item[(b)] The complex $E_{\Delta}\cDb$ is residual, 
so that, via $1 \cong QE_{\Delta}$ 
of~\textup{\ref{cor:coz4a}}, $\cDb$ is isomorphic to
a residual complex. 
\end{enumerate}
\end{aprop}

\pagebreak[3]

\begin{proof} {\rm(i)} The flat-base-change isomorphism
\cite[p.\,8, Thm.\,3]{AJL2} gives:

\begin{alem}
\label{lem:residual7}
Let $\X,x,j_x,\W$ be as in \textup{\ref{lem:residual5}}.
Let $\I$ be a defining ideal in~$\OX$.
Let $X, W$ be the schemes obtained by 
going modulo $\I, \I\OW$ in $\OX, \OW$ 
respectively so that the following diagram 
of induced natural maps is a fiber square. \vspace{-2pt}
\[
\begin{CD}
W @>{j'_x}>> X \\[-2pt]
@V{i'}VV  @VViV \\
\W @>>\under{1.2}{{j_x}}> \X\\[-2pt]
\end{CD}
\]
Then for any $\cFb \in \Dqct^+(\X)$
there is a natural isomorphism
$i'{}^!j^*_x\cFb \iso j'{}^*_xi^!\cFb$.
\end{alem}

Now, in~\ref{lem:residual7} take $i = i_n, X = X_n$.
Then $\cDb$ is pointwise $t$-dualizing $\Lra$ (for any $x \in \X$)
$j_x^*\cDb$ is $t$-dualizing on $\W$ $\Lra$ 
$i'{}^!j_x^*\cDb$ is dualizing on $W$ (\cite[Lemma 2.5.10]{AJL2})
$\Lra$ $j'{}^*_xi^!\cDb$ is dualizing on $W$. Thus $i^!\cDb$
is a pointwise dualizing complex on $X$. The converse is proved by 
reversing the arguments.

(ii). We use the same arguments as in (i) above along with
\cite[V, 7.2 and 8.2]{RD}.

(iii). Suppose $\cRb$ is a residual complex on~$\X$.
First note that if $\X$ has finite Krull dimension,
then $\cRb$ is $t$-dualizing. Indeed, by \ref{lem:residual3} 
and \ref{prop:mod7}(ii), $\cRb$ consists
of $\Aqct$-injectives, therefore, by \cite[2.5.6]{AJL2},
\cite[2.5.10]{AJL2} and \cite[VI, 1.1(a)]{RD}, $\cRb$ is 
$t$-dualizing. Now we drop the hypothesis on Krull dimension
of~$\X$. For $x \in \X$, $j_x$ in \ref{lem:residual5} is 
\'etale and adic. Therefore, by~\eqref{eq:fin-sm1}, we 
obtain an isomorphism $j_x^*\cRb \iso j_x^{\sharp}\cRb$, 
which by \ref{prop:residual4},
is a residual complex on~$\W$. Since $\W$ has finite
Krull dimension, 
$j_x^*\cRb$ is $t$-dualizing on~$\W$. Thus
$\cRb$ is a pointwise $t$-dualizing complex. 

Conversely, suppose $\cDb$ is a pointwise $t$-dualizing complex 
on~$\X$. 
By \ref{prop:loc8}(ii) we may assume that 
$\cDb$ is a bounded below complex of $\Aqct(\X)$-injectives.
For a defining ideal $\I$ in $\OX$ and $n>0$, set 
$\cDb_n \set \sHom_{\OX}(\OX/\I^n, \cDb)$. Then 
$\cDb$ is isomorphic to the direct limit of the $\cDb_n$'s.
By \cite[2.5.6]{AJL2} $\cDb_n \iso i_n^!\cDb$ and hence 
by \cite[2.5.10]{AJL2} $\cDb_n$ is a pointwise dualizing complex on 
$X_n = (\X, \OX/\I^n)$. 
In particular, $\cDb_n$ induces a unique function $\Delta_n$  
for which it is Cohen-Macaulay. Since $\cDb_n \iso
i_{n,n+1}^!\cDb_{n+1}$ where $i_{n,n+1}$ is the canonical 
immersion $X_n \to X_{n+1}$, it follows that 
$\Delta_n = \Delta_{n+1}$, henceforth to be denoted as~$\Delta$.
Since $H^*_x(-)$ commutes with direct limits
we see that the $\Delta$-Cohen-Macaulayness of the $\cDb_n$'s 
carries over to~$\cDb$. By \ref{prop:mod7}(iii),(iv) and
\ref{cor:loc4}, there is an isomorphism  
$H^*_x\cDb \iso H^*_{m_x}\cDb_x$ so that
$E_{\Delta}\cDb \in \Coz_{\Delta}(\X)$. 
Thus, via the isomorphism in \ref{cor:coz4a},
$\cDb$ is $\D(X)$-isomorphic to
a complex $\cRb$ in $\Coz_{\Delta}(\X)$.

If we define $\cRb_n$ in the obvious manner, then
$\cRb_n$ is a pointwise dualizing complex on $X_n$ lying
in $\Coz_{\Delta}(X_n)$. In particular, $\cRb_n$
is a residual complex on $X_n$. Thus~$\cRb$ is a residual complex.
\end{proof}

Finally we prove that Yekutieli's\index{Yekutieli, Amnon}
definition of a residual complex on~$\X$ agrees with ours 
when $\X$ has finite Krull
dimension,\footnote{In this context note that 
Lemma 5.13 in \cite{Ye} requires a   finite-dimensionality \emph{assumption}, see Nagata's  
\index{Nagata, Masayoshi} 
example of an infinite-dimensional noetherian regular ring in \cite[p.\,203]{Nag2}.} 
i.e., we show that a complex $\cRb$ on~$\X$ 
is residual if and only if it is $t$-dualizing and 
there is an isomorphism of $\OX$-modules
$\oplus_p\cR^p \iso \oplus_x i_xJ(x)$ where $J(x)$ is the 
$\OXx$-injective hull of the residue field $k(x)$ at $x$
(see \ref{prop:mod7}(ii)). 
The `only if' part follows from \ref{lem:residual3} 
and \ref{prop:residual6}(ii),(iii). For the `if' part, note
that the corresponding statement holds for ordinary schemes
and hence if we define $\cRb_n$ as in the proof of \ref{lem:residual3} 
then $\cRb_n$ is residual on $X_n$. Thus $\cRb$ is 
residual.

\newcommand{\cDt}{\cD_{\<\mathrm t}}
\newcommand{\Dcs}{{\Dc\<\!\!{}^*}}

\subsection{CM complexes and coherent sheaves}
\label{coherent duality}

For a formal scheme $\X$, let $\Dcs(\X)$ be the essential image of $\R\iGp{\X}|_{\Dc}$, i.e., the full subcategory
of~$\D(\X)$ such that $\E\in\Dcs\Leftrightarrow\E\cong\R\iGp{\X}\F$ with
$\F\in\Dc$. Proposition~\ref{prop:loc8} shows that $\Dcs(\X)\subset\Dqct(\X)$. The functor $\R\iGp{\X}$ is an equivalence from the category $\Dc(\X)$ to the category
$\Dcs(\X)$, with quasi-inverse $\BL\set\R\sHomb(\R\iGp\X\OX,\,-)$  (see
\cite[p.\,24]{AJL2}).%
\index{ $\K(\A)$1@$\BL\set\R\sHomb(\R\iGp\X\OX,\,-)$ (right adjoint of $\R\iGp\X$)} 
In particular, if $\X$ is an ordinary scheme, so that $\R\iGp\X$ is the identity functor on $\D(\X)$, then $\Dcs(\X)=\Dc(\X)$.%
\index{ $\D(\A)$ (derived category
of $\A$-complexes)!$\D(X)\set\D(\A(X))$!$\Dc^*(\X)$
(essential\vspace{.7pt} image of $\D_{\mathrm c}(\X)$)}

By \cite[p.\,26, 2.5.3 and p.\,27, 2.5.5]{AJL2},
$\X$  has a t-dualizing complex $\cRb\in\Dcs(\X)$
iff $\>\X$ has a ``c-dualizing" complex~$\Rc\in\Dc(\X)$, 
which is so, e.g., if $\X$ is locally embeddable in a regular
finite-dimensional formal scheme. In fact $\BL\cRb$ is c-dualizing; and
conversely, if $\Rc$ is c-dualizing then 
$\R\iGp\X\Rc\in\Dcs(\X)$ is t-dualizing.

\begin{aprop}\label{CM+coh}
Let\/ $\X$ be a formal scheme with a\/ t-dualizing complex\/~$\cRb,$ which may be assumed residual, and $\Delta$ the codimension function associated to\/~$\cRb$
$(\<$see \textup{\ref{prop:residual6}, \ref{lem:residual3}).} Assume further that\/ $\cRb\in\Dcs(\X)$.
Then\/ $\cGb\in\Dcs(\X)$ 
is a\/ $\Delta$\textup{-CM} complex if and only if its dual 
$$
\cDt\cGb\set\R\sHomb(\cGb\<,\cRb)
$$
is\/ $(\D(\X)$-isomorphic to$)$ a coherent\/ $\OX$-module.\looseness=-1 
\end{aprop}
\index{ $\Delta$@$\cDt$ (torsion-dualizing functor)}

By \cite[p.\,28, 2.5.8]{AJL2},  the functor $\cDt$ induces, in either direction,  an antiequivalence between $\Dc(\X)$ and $\Dcs(\X)$. Thus:

\begin{acor}\label{cor:CM+coh}
The functor\/ $\cDt$ induces, in either direction,  an antiequivalence between\/ $\Ac(\X)$ and the full subcategory of\/ \textup{$\Delta$-CM} complexes in\/ $\Dcs(\X)$.
\end{acor}

Since for all $\cMb\in\D(\X)$, $\R\sHomb(\cMb\<,\BL\cRb)\cong \R\sHomb(\R\iGp\X\cMb\<,\>\cRb)$ (see \cite[p.\,24, (2.5.0.1)]{AJL2}), therefore:
\begin{acor}
A complex\/ $\cMb\in\Dc(\X)$ is such that\/ $\R\iGp\X\cMb$ is\/ $\Delta$-\textup{CM} 
if and only if, with\/  $\Rc\set\BL\cRb,$ the coherent dual
 $$
 \cD_{\textup{c}}\>\cMb\set  \R\sHomb(\cMb\<,\> \cRb_{\textup{c}})
$$ 
is a coherent\/ $\OX$-module.\vspace{.6pt}

In particular, $\R\iGp\X\OX$ is\/  $\Delta$-\textup{CM} if and only if\/ $\Rc$ is a coherent\/ $\OX$-module.
\end{acor}

\begin{aexam}\label{Matlis}
 Let $A$ be a complete local ring and $\X\set\Spf(A)$ the
formal spectrum of $A$. Then $\A(\X)$ is the category of $A$-modules.
We can take $\cRb$ to be an injective hull of the residue field of $A$,
(considered as a complex vanishing in nonzero degrees) \cite[p.\,25, Example
2.5.2\kern.5pt(3)]{AJL2}. Then $\Delta$ maps the unique point~$x\in\X$ to 0; and a $\Delta$-Cousin complex is simply an $A$-module. Thus the Cousin functor $E$ can be identified with the functor $H^0_x$, and by~\ref{prop:coz4}, $H^0_x$ is an \emph{equivalence} from the category of $\Delta$-CM complexes  to the category of $A$-modules. So if $\cGb$ is $\Delta$-CM then there is a $\D$-isomorphism $\cGb\iso H^0_x\cGb$, whence $\cGb\in\Dcs(\X)\Leftrightarrow H^0_x\cGb\in\Dcs(\X)$. But 
from \cite[p.\,28, Prop.\,2.5.8\kern.5pt(a)]{AJL2}  
and Matlis duality \cite[p.\,148, Thm.\,18.6\kern.5pt(v)]{Ma} it follows  
that $H^0_x\cGb\in\Dcs(\X)$ if and only if $H^0_x\cGb$ is an \emph{artinian} 
$A$-module. Thus: 

\emph{The equivalence of categories\/
$
H^0_x\colon\{\Delta\textup{-CM $\OX$-complexes}\}
\xrightarrow[\ \under{3.3}{\approx}\;]{}
\{A\textup{-modules}\}
$\vspace{-3pt} 
takes the full subcategory of\/
$\Delta$\textup{-CM} complexes in\/ $\Dcs(\X)$ to the full subcategory of \emph{artinian} 
$A$-modules, and transforms the antiequivalence of~\textup{\ref{cor:CM+coh}} into Matlis duality.} 
\end{aexam}

\enlargethispage*{1pt}

{\sc Proof} of~\ref{CM+coh}. 
As above, if $\cGb\in\Dcs$ then $\cDt\cGb\in\Dc$, and $\cGb\cong\cDt\cDt\cGb$.
It will therefore suffice to show that \emph{if\/ $x\in\X$ and $\cFb\in\Dc$ then} 
\[\label{onecoh}
H^i\R\iG{x}\R\sHomb(\cFb\<,\cRb)=0
\iff
\Hr^{\Delta(x)-i}\>\cFb_x=0.\tag{\ref{CM+coh}.1}
\]

Let us prove ~\eqref{onecoh}.
The first assertion in \cite[p.\,33, (5.2.1)]{AJL1} (whose proof applies to any ringed space) gives a canonical isomorphism
$$
\R\iG{\,\ov{\!\{x\}\!}\,}\R\sHomb(\cFb\<,\cRb)\iso
\R\sHomb(\cFb\<,\R\iG{\,\ov{\!\{x\}\!}\,}\cRb).
$$

Recall that the stalk at $x$ of an injective $\OX$-module is an injective $\OXx$-module. (In view of \cite[pp.\:110\kern1pt--111, (5.2.6) and (5.2.8)]{EGAI}, the proof for schemes in \mbox{\cite[p.\,128, 7.12]{RD}} extends to formal schemes.) It follows that $\R\sHomb(\cFb\<,-)$ ``commutes" with the exact functor ``stalk at $x$."
(Since $\cRb$ is a bounded injective complex, one can use \cite[p.\,68, Prop. 7.1]{RD} to reduce the proof to the trivial case $\cFb=\OX^n$.) 
Thus there are canonical isomorphisms
$$
H^i\R\iG{x}\R\sHomb(\cFb\<,\>\cRb)\iso
\Hr^i\R\Homb_{\OXx}\<(\cFb_x,\>\R\iG{x}\cRb).
$$

Since $\R\iG{x}\cRb\cong \cRb(x)[-\Delta(x)]$ and, by~\ref{lem:residual3}, $\cRb(x)$ is an
injective module containing the residue field $k_x$ of~$\OXx$, therefore
$$
\Hr^i\R\Homb_{\OXx}\<(\cFb_x,\>\R\iG{x}\cRb)\cong
\Hom_{\OXx}\<\big(\Hr^{\Delta(x)-i}\>\cFb_x,\cRb(x)\big).
$$

If $\Hr^{\Delta(x)-i}\>\cFb_x\ne 0$   
then any nonzero principal submodule  admits a nonzero $\OXx$-homomorphism into $k_x\subset \cRb(x)$,
and so, $\cRb(x)$ being injective, there exists a nonzero map 
$\Hr^{\Delta(x)-i}\>\cFb_x\to\cRb(x)$, whence
$H^i\R\iG{x}\R\sHomb(\cFb\<,\cRb)\ne 0$. The assertion \eqref{onecoh} results.
\hfill$\square$

\smallskip
\begin{aprop} 
\label{E=EHn}Under the assumptions of\/ \textup{\ref{CM+coh}} there are, for any\/ \mbox{$n\in\mathbb Z$} and\/ ${\cFb}\in\Dc(\X),$ with ${\cFb}'\set\R\sHomb({\cFb}\<,\>\cRb)$ and\/ $H\set H^n\cFb\in\Dc(\X),$  functorial isomorphisms
$$
E_{\Delta-n}\>{\cFb}'\iso E_{\Delta-n}\big((H[-n])'\big)\iso H'[n]\in\Dcs(\X).
$$
Hence for all\/ $n\in\mathbb Z,$ $E_{\Delta-n}\>\Dcs(\X)\subset\Dcs(\X)$. 
\end{aprop}

\begin{proof} Let $\theta\colon\cFb_1\to\cFb_2$ be a $\Dc$-map such that\vspace{.6pt}
$H^n(\theta)$ is an isomorphism $H^n{\cFb}_1\iso H^n{\cFb_2}$. As in the proof of~\ref{CM+coh}, we deduce\vspace{-.6pt} that for each $x\in \X$ the map 
$H^{\Delta(x)-n}_x{\cFb_1}'\to H^{\Delta(x)-n}_x{\cFb_2}'$ induced by $\theta$ is an isomorphism.
Thus $\theta$ induces an isomorphism  
$E_{\Delta-n}{\cFb_1}'\iso E_{\Delta-n}{\cFb_2}'$. 

To get the first isomorphism in the Proposition, 
apply the preceding to the canonical maps 
$\theta_1\colon\cFb\to\tau_{\ge n}\cFb$ and 
$\theta_2\colon(H^n\cFb)[-n]\to\tau_{\ge n}\cFb$, each of which induces a homology isomorphism in degree~$n$. 

For the second isomorphism, \ref{cor:CM+coh} gives that $H'$ is a CM-complex in $\Dcs$, whence 
$$
E_{\Delta-n}\big((H[-n])'\big)
= E_{\Delta-n}\big((H'[n])\big)
\overset{\under{-1.8}{\eqref{eq:Cousin+trans}}}=(E_\Delta H')[n]
\overset{\under{-2.6}{\eqref{cor:coz4a}}}\iso H'[n]\in\Dcs.\\[4pt]
$$

The last assertion holds because by the statement preceding Corollary~\ref{cor:CM+coh}, every $\G\in\Dcs$ is, isomorphic to an ${\cFb}'$.
\end{proof}

\begin{acor}[cf.~{\cite[p.26, Thm.\,3.2]{DT}, 
\cite[Thm.\,4.4]{Kw}}]\label{DTK}
Let\/ $\X,$ $\cRb$ and\/ $\Delta$ be as in~\textup{\ref{CM+coh}.} Let\/ $0\ne\cGb\in\Dcs^-(\X),$ so that\/ $0\ne\cDt\cGb\in\Dc^+(\X)$ and we can  set
$$
m=m(\cGb)\set\min\{\,n\mid H^n\>\cDt\cGb\ne0\,\}.
$$

{\rm (i)} There is a canonical\/ $\D(\X)$-map\/
$\ssf(\cGb)\colon\cGb[-m]\to E_\Delta\big(\cGb[-m]\big)\ne 0$
such that any\/ $\D(\X)$-morphism of\/ $\cGb[-m]$ into a\/ $\Delta$\textup{-CM} complex\/~$\cEb$
in\/ $\Dcs(\X)$ factors uniquely through\/ $\ssf(\cGb)$.\vspace{.6pt}

{\rm (ii)} If\/ $\cGb[-m]$ in\/ \textup{(i)} is\/ $\Delta$\textup{-CM} then\/ $\ssf(\cGb)$  is the isomorphism $\Ssf(\cGb[-m])$ of\/~\textup{\ref{cor:coz4a}.}\vspace{.6pt}

\end{acor}

\begin{proof}
For simplicity, we write $\G$ for $\cGb\<$, \dots, 
$\G'$ for $\cDt\cGb$, \dots. 

Set $\F\set\G'\in\Dc(\X)$ and $H\set H^m(\G)$. By ~\ref{E=EHn} the natural map 
\mbox{$H[-m]\to\F$} induces an isomorphism
$$
E_\Delta\big(\G[-m])
\overset{\under{-1.8}{\eqref{eq:Cousin+trans}}}=
\big(E_{\Delta-m}\>\G\big)[-m]\iso H'\ (\ne 0),\\[4pt]
$$
whose inverse composed with the natural map $\sigma\colon \G[-m]=(\F[m])'\to H'$
is defined to be the map $\ssf(\G)$.

Assertion (i) results now from the sequence of natural isomorphisms
\begin{align*}
\Hom\big(\G[-m],\,\E\big)&\cong
\Hom\big(\E', \>\G'[m]\big)\\
&\cong
\Hom\big(\E'[-m], \>\F\>\big)
\cong
\Hom\big(\E'[-m],\> H[-m]\big)\cong
\Hom\big(H'\mkern-1.5mu,\mkern1.5mu\E\big).
\end{align*}
(For the third isomorphism, recall from~\ref{CM+coh} that $\E'$ is a coherent
$\OX$-module.)

If $\G[-m]$ is $\Delta$-CM then by~\ref{CM+coh}, $\F[m]$ is a coherent $\OX$-module, so that $\sigma$ is an isomorphism, whence so is $\ssf(\G)$; and then functoriality of
$\Ssf$ forces $\ssf(\G)=\Ssf\big(\G[-m]\big)$, proving (ii).
\end{proof}

\begin{aexam}

Suppose that in~\ref{DTK}, $\X$  is an \emph{ordinary} scheme. 
After translating $\cR$ we may assume that $m(\OX)=0$. Then with $K\set H^0(\OX')= H^0(\cR)$, a \emph{canonical\/ $\OX$-module,} \ref{E=EHn} for $\cFb=\cR$ gives $E_\Delta(\OX)\cong K'\<$, whence $m(K)=0$. If $\OX$ satisfies the Serre condition $(\textup{S}_2)$ then the natural map is an isomorphism 
$$
\OX\iso H^0(E_\Delta(O_\X))\cong H^0(K')\  (=\sHom(K,K));
$$ 
so in this case \ref{E=EHn} for $\cFb=K'$ gives 
$$
E_\Delta(K)\cong E_\Delta(K'')\cong \OX'=\cR.
$$

Similar considerations apply with an arbitrary coherent $\OX$-module in place of~$\OX$,
see, e.g., \cite[Thm.\,4.6 (and 1.4)]{Db}.
\end{aexam} 

\begin{small}

Keeping the notation and assumptions of~\ref{DTK} and its proof, we close with some characterizations of the integer $m(\G)$.\vspace{1pt}

Replacing $m$ by  $n<m$ in the proof of~\ref{DTK}\kern.5pt(i) gives
that \mbox{$E_\Delta\big(\G[-n]\big)=E_{\Delta-n\>}\G=0$,} and further that 
any\/ $\D(\X)$-morphism of\/ $\G[-n]$ into a\/ 
$\Delta$\textup{-CM} complex\vspace{-.6pt} in\/ $\Dcs(\X)$ vanishes.  
The~proof of~\ref{CM+coh} shows further that with $\F\set{\G}'$, 
$H_x^{\Delta(x)-m}\G\ne0\iff x\in\Supp(H^m\F)$. 
In particular, \emph{$m\set m(\G)$ is 
the unique integer such that\/}
\vspace{1pt}

(a) \emph{for all\/ $x\in \X$ and\/ $q>\Delta(x)-m$, $H_x^q\G=0$,} \emph{and}

(b) \emph{there is an\/ $x\in\X$ such that\/ $H_x^{\Delta(x)-m}\G\ne0$.}
\vspace{2pt}

\noindent We claim that, furthermore:
\vspace{1pt}

(1) $m=m'\set\min\{\,\Delta(x)-q\mid q\in\mathbb Z,\ x\in\Supp(H^q\G)\,\}$, 
\emph{and}

(2) \emph{for any generic point\/ $x$ of\/ $\Supp(H^m\F),$ $x\in\Supp(H^q\G)$ and\/ $\Delta(x)-q=m.$}
\vspace{2pt}

\noindent In less contorted terms, \emph{$m\set m(\G)$ is 
the unique integer such that}
\vspace{1pt}

(a)$'$ \emph{for all\/ $x\in \X$ and\/ $q>\Delta(x)-m$,
$\Hr^q(\G_x)=0$,}
\emph{and}

(b)$'$ \emph{for any generic point\/ $x$ of\/ $\Supp(H^m\F),$ 
$\Hr^{\Delta(x)-m}(\G_x)\ne 0$.}\vspace{2pt}

\begin{proof}
Let $Z^\bullet$ be the filtration defined by $\Delta$ (i.e., $x\in Z^p\iff \Delta(x)\ge p$.)
In view of~(a), \eqref{H=oplus} gives, for any $j>p-m$,
$$
H^j_{Z^p/Z^{p+1}}\G=\bigoplus_{\Delta(x)=p}\,i_x (H_x^j\>\G)=0.
$$
So if $q+m>d\set\min_{x\in\X}\Delta(x)$, there is a surjection 
followed by isomorphisms
\begin{equation*}
H^q_{Z^{q+m}}\>\G\twoheadrightarrow 
H^q_{Z^{q+m-1}}\>\G\iso
H^q_{Z^{q+m-2}}\>\G\iso
\cdots\iso H^q_{Z^{d}}\>\G=H^q\G.\tag{$*$}
\end{equation*}
Hence $\Supp(H^q\G)\subset \Supp(H^q_{Z^{q+m}}\G)\subset Z^{q+m}$, 
and so $m\le m'$.\vspace{1pt}

Next, among points in $\Supp(H^m\<\F)$ choose one, say $x$,
where the value of $\Delta$ is minimal. We will show that 
$H^{\Delta(x)-m}\G_x\ne 0$, from which follows that $m'\le m$,
proving~(1). 

It also follows that if $z$ is \emph{any} generic point 
of $\Supp(H^m\F$) then, since $z$ has
an open neighborhood $V$ within which $z$ is the \emph{only} such generic point, 
therefore
$$
0\ne H^{\Delta(z)-m(\G|_{V})}(\G|_{V})_z=H^{\Delta(z)-m}\G_z,
$$ 
proving~(2).\vspace{1pt}

Set $\delta\set\Delta(x)$. As before,  $H^{\delta-m}_x\G$ is dual to $\Hr^m\F_x\ne 0$, so 
$H^{\delta-m}_x\G\ne 0$. Similarly, if $\Delta(y)<\delta$, so that $y\notin\Supp(H^m\F)$,
then $H^{\Delta(y)-m}_y\G=0$. Hence, with $q\set\delta-m$,\vspace{1pt} 

(i) $H^q_{Z^{q+m}/Z^{q+m+1}}\G$ has nonzero stalk at $x$; and

(ii) $H^{q-1}_{Z^{q+m-1}/Z^{q+m}}\G=0.$\vspace{1pt}

\noindent From (ii) it follows that the natural map $H^q_{Z^{q+m}}\G\to H^q_{Z^{q+m-1}}\G$ is injective. Then since~$q\ge d-m$, $(*)$ shows that $H^q_{Z^{q+m}}\G\cong H^q\G$ (the case $q=d-m$ being trivial).

 Since for any abelian sheaf $A$, 
the stalk at $x$ of $\iG{Z^{\delta+1}}A$ vanishes,
therefore the natural map is an isomorphism 
$(H^{q}_{Z^{q+m}}\G)_x\iso (H^{q}_{Z^\delta/Z^{\delta+1}}\G)_x$;
so for nonvanishing of $H^{\Delta(x)-m}\G_x$ we need only note that 
by (i), the target of this isomorphism doesn't vanish.
\end{proof}

\end{small}
\newpage

\section{Some explicit descriptions}
\label{sec:sexd}

We conclude this paper by
giving explicit descriptions
of some notions encountered in earlier sections.
In \S\ref{subsec:QE} we concentrate on the Suominen
isomorphism\index{Suominen, Kalevi!isomorphism} %
for Cohen-Macaulay complexes; this is the
isomorphism $1_{\ov{\D}} \iso QE$ defined in~\ref{cor:coz4a}.
There it is defined in an implicit way. Here we
give an explicit description for it in
several cases. (See also Corollary~\ref{DTK}\kern.5pt(ii).)
In \S\ref{subsec:etale} we study
$f^{\sharp}$ for an \'etale map $f$. This includes the case
of completion of a formal scheme along an open ideal
and also the case of open immersions which we covered
in~\ref{lem:fin-sm3} and~\ref{rem:fin-clim2}.

We use the following notation throughout this section.

(a) Let $(\X, \Delta) \in \bbFc$.
Let $\ov{\D}\set \D^+(\X, \Delta)_{\text{CM}}$ be the category of
$\Delta$-CM complexes on $\X$ (\S\ref{subsec:cm}). Then the
\emph{Suominen isomorphism} $\Ssf_{\X, \Delta}$
(or simply $\Ssf$ in case of no ambiguity)%
\index{ OON@$\Ssf$ (Suominen isomorphism)} is the
isomorphism $1_{\ov{\D}} \iso QE$ defined in~\ref{cor:coz4a}.
For any $\cGb \in \ov{\D}$ we denote the induced isomorphism
$\cGb \iso QE\cGb$ by $\Ssf(\cGb)$.

(b) For any smooth map $f \colon (\X, \DsssX) \to (\Y, \DsssY)$
in $\bbFc$ of relative dimension $d$, we use
$f^{\diamond}(-)  \set f^*(-) \otimes_{\X} \omega_f[d\>]$.
\index{ $\Forget$3@$f^{\diamond}$} 

\subsection{The Suominen isomorphism}
\label{subsec:QE}

Let $(\X, \Delta) \in \bbFc$. Let $\cGb$
be a $\Delta$-CM complex on $\X$ with bounded homology.
In \cite[p.~242-246]{RD},
Hartshorne describes a non-canonical
way of constructing an isomorphism $\cGb \iso QE\cGb$.
This construction lacks functorial properties.
Nevertheless it satisfies a property, recalled below,
that makes the isomorphism
functorial when $\cGb$ ranges over \emph{Gorenstein} complexes
(\cite[p.~248]{RD}).
We now show that the Suominen isomorphism $\Ssf(\cGb)$ also
satisfies this property, so that if~$\cGb$ is Gorenstein, then
$\Ssf(\cGb)$ agrees with $\phi(\cGb)$ for $\phi \colon 1 \iso QE$
as constructed in~\cite[p.~249]{RD}.

Let $(\X, \Delta) \in \bbFc$. 
Let $\cFb \in \D^{\text{b}}(\X)$, and let $\cLb$ be an injective
resolution of~$\cFb$. Let $Z^{\bullet}$ denote the
filtration of~$\X$ induced by~$\Delta$.
Let $E^{p,q}_r(\cFb) = E^{p,q}_r(\cLb)$ be the
spectral sequence associated to the
filtered complex $\iG{Z^{\bullet}}\cLb$.
Recall that $E_{\Delta}\cFb$ is given by the $E^{*,0}_1$-terms
of this spectral sequence. Now if~$\cFb$ is also $\Delta$-CM,
then the $E^{p,q}_1$-terms vanish for $q \ne 0$ so that
the spectral sequence degenerates.
Since the spectral sequence converges to the homology
of~$\cLb$, there results, for any $i \in \mathbb Z$,
(with $\cFb$ $\Delta$-CM) a natural isomorphism
\[
\psi^i_{\cFb} \colon H^iE_{\Delta}\cFb = E^{i,0}_2(\cLb) \iso H^i \cLb
\iso H^i \cFb.
\]

Let $\ov{\bC}, \ov{\D}, Q, E$ be as in \S\ref{subsec:cm}. Set
$\ov{\D}{}^{\,\text{b}} \set \ov{\D} \cap \D^{\text{b}}(\X)$.

\begin{aprop}
\label{prop:QE1}
With notation as above,
for any $\cFb \in \ov{\D}{}^{\,\textup{b}}$ and for any integer $i$,
the induced isomorphism $H^i\Ssf(\cFb)$ is inverse to $\psi^i_{\cFb}$.
\end{aprop}
\begin{proof}
Since $H^i\Ssf(-)$ and~$\psi^i_{(-)}$ are functorial in $\cFb$
we may replace $\cFb$ by any isomorphic complex in~$\ov{\D}{}^{\,\text{b}}$.
In particular, since $\cFb$ is isomorphic to a Cousin complex,
therefore, we may assume without loss of generality
that $\cFb = Q\cCb$ where~$\cCb$ is a Cousin complex.

By \ref{cor:coz4a}, the isomorphism
$\Ssf(Q\cCb) \colon Q\cCb \iso QEQ\cCb$
is the same as the $Q$-image of the inverse of the isomorphism
$\phi_{\cCb} \colon EQ\cCb \xto{\textup{\ref{lem:coz1a}(ii)}} \cCb$.
Therefore it suffices to show that
$\psi^i_{\cCb} = H^i\phi_{\cCb}$.

Before proceeding further, we need a definition.
For any complex $\cGb$ on $\X$, let $\Esf\cGb$ denote the
complex given by the $E^{*,0}_1$-terms of the spectral
sequence associated to the filtered complex
$\iG{Z^{\bullet}}\cGb$; the only difference from the definition
of $E\cGb$ being that now we do not replace $\cGb$ by an injective
resolution. If $\cGb$ consists of flasque sheaves, then
there is a canonical isomorphism $\Esf\cGb \iso E\cGb$.
If $\cGb$ is $\Delta$-Cousin, then $\Esf\cGb = \cGb$ (\ref{lem:coz1a}(ii))
and moreover the filtration $\{\iG{Z^n}\cGb\}_{n \in \mathbb Z}$
is now given by truncations (\ref{lem:coz1a}(i)[a]). Therefore the
$E^{p,q}_1$-terms of the associated spectral sequence vanish for $q \ne 0$.
It follows that if $\cGb$ is a $\Delta$-Cousin complex, then the canonical
maps
$$
H^i\cGb = H^i\Esf\cGb \xto{\text{deg.~sp.~seq.~for $\iG{Z^{\bullet}}\cGb$}}
H^i\cGb
$$
compose to the identity map.

Returning to the proof, let $\cCb \to \cLb$ be an injective resolution.
Consider the following commutative diagram of natural isomorphisms.
\[
\begin{CD}
H^i\cCb @>{\alpha}>{\cCb = \Esf\,\cCb}> H^i\Esf\,\cCb
 @>>{\text{deg.~sp.~seq.~for $\iG{Z^{\bullet}}\cCb$}}> H^i\cCb \\
@. @V{\beta}VV @VV{\gamma}V \\
@. H^i\Esf\,\cLb
 @>{\delta}>{\text{deg.~sp.~seq.~for $\iG{Z^{\bullet}}\cLb$}}> H^i\cLb
\end{CD}
\]
By definition, $\Esf\,\cLb = E\cCb$ and moreover
$(\beta\alpha)^{-1} = H^i\phi_{\cCb}$ and
$\gamma^{-1}\delta = \psi^i_{\cCb}$.
Since the top row composes to the identity map, it follows that
$H^i\phi_{\cCb} = \psi^i_{\cCb}$.
\end{proof}

Our next goal concerns adic smooth maps in $\bbFc$. We begin with a
preliminary result on fibers.

For any map of formal schemes $f \colon \X \to \Y$, the fiber
of $f$ over a point $y \in \Y$ is the
formal scheme $\X_y \set \X \times_{\Y} \Spec(k(y))$ where
$k(y)$ is the residue field of~$\OYy$.
As is the usual practice, we shall identify
the underlying topological space of $\X_y$
with its canonical image in $\X$.

As usual, for points $x \in \X$ and $y \in \Y$, we shall denote the
maximal ideals at $\OXx,\OYy$ by $m_x, m_y$ respectively.
The following basic facts on fiber spaces follow easily from the
definitions.
For any $x \in \X_y$, the stalk of the structure sheaf of $\X_y$ at
$x$ is canonically isomorphic to $\OXx/m_y\OXx$.
If~$f$ is adic, then $\X_y$ is an ordinary scheme.

\begin{alem}
\label{lem:fiber1}
Let $f \colon (\X, \Delta_1) \to (\Y, \Delta_2)$ be a smooth map
in $\bbFc$ of constant relative dimension $d$. Let $y \in \Y$.
Then for any $x \in \X_y$, $0 \le \Delta_2(y) - \Delta_1(x) \le d$.
Moreover, if $f$ is adic, then the following conditions are equivalent.
\begin{enumerate}
\item $x$ is a generic point of $\X_y$.
\item $m_y\OXx = m_{x}$.
\item $\Delta_1(x) = \Delta_2(y) - d$.
\end{enumerate}
\end{alem}
\begin{proof}
By definition, $\Delta_1 = f^{\sharp}\Delta_2$ (\ref{exam:formal1}) so that
$\Delta_2(y) - \Delta_1(x)$ equals the transcendence degree of
the residue field extension $k(y) \to k(x)$ and hence is nonnegative.
The other inequality follows from \ref{cor:app4} which also
immediately implies (ii)~$\iff$~(iii). If (iii) holds,
then for any nontrivial specialization $x' \leadsto x$
we get a contradiction because
$\Delta_1(x') < \Delta_1(x) = \Delta_2(y) - d$. Thus (iii) $\Lra$ (i).

If $f$ is adic and smooth, then $\X_y$ is a smooth ordinary
scheme over $\Spec(k(y))$ and hence is a disjoint union of integral schemes.
Therefore if $x$ is a generic point of~$\X_y$ then the local ring $R$ of
$\X_y$ at $x$ is a field. In view of the canonical isomorphism
$R \iso \OXx/m_y\OXx$, (i) $\Lra$ (ii) follows.
\end{proof}

Let $f \colon (\X, \DsssX) \to (\Y, \DsssY)$ be an \emph{adic} smooth map
in $\bbFc$ of constant relative dimension $d$.
Let $\cMb$ be a complex in $\Coz_{\DsssY}(\Y)$. Since $\cMb$
consists of torsion $\OY$-modules,
$f$ being adic implies that $f^*\cMb$ consists of torsion $\OX$-modules.
In particular, $f^{\diamond}\cMb = f^*(-) \otimes_{\X} \omega_f[d\>]$
also
consists of torsion $\OX$-modules and hence the canonical $\D(\X)$-map
$\R\iGp{\X}f^{\diamond}\cMb \to f^{\diamond}\cMb$ is an isomorphism.
Therefore, by \ref{lem:quasim1}, $f^{\diamond}\cMb$ is a $\DsssX$-CM complex.
Set $E \set E_{\DsssX}$. Our aim is twofold :
\begin{itemize}
\item To explicitly construct a quasi-isomorphism
$\eta_f(\cMb) \colon f^{\diamond}\cMb \to Ef^{\diamond}\cMb$ in~$\bC(\X)$.
\item To show that the natural image of $\eta_f(\cMb)$ in $\D(\X)$
is the same as the Suominen isomorphism~$\Ssf(f^{\diamond}\cMb)$.
\end{itemize}

Before proceeding further we need to set up some notation.
Much of the notation concerns the fact that both
$f^{\diamond}\cMb$ and $Ef^{\diamond}\cMb$ also have
a graded decomposition arising from the punctual
grading on~$\cMb$, i.e., one parametrized by the points
of~$\Y$. For the rest of this subsection, we shall
assume that $f,\X, \DsssX,\Y,\DsssY,d$ are fixed.
Set $\omega \set \omega_f$.
Here then is some notation that we shall use for the
rest of this subsection.

\begin{enumerate}
\item For any $y \in \Y$ with $M = \cMb(y)$, set
$f^{\diamond}(\cMb,y) \set f^*i_yM \otimes_{\Y} \omega$.
For any integer $a$, there results natural isomorphism
\[
(f^{\diamond}\cMb)^a \iso \bigoplus_{\{ y \in \Y|\, \DsssY(y) = a+d \}}
f^{\diamond}(\cMb,y).
\]
\item For any $x \in \X$, with
$y=f(x), p = \DsssX(x)-\DsssY(y)+d, M = \cMb(y)$, set
\[
G_{x,\cMb} \set H^{p}_{m_x}(M \otimes_{y} \omega_{x}).
\]
\item For any $y \in \Y$ and any integer~$a$, set
$\E_{\cMb}^a(y) \set \bigoplus_x i_x(G_{x,\cMb})$ where~$x$ varies
over points in~$\X_y$ such that $\DsssX(x)= a$.
\item For any integers $a,b$ set
\[
\E_{\cMb}^{a,b} \set \bigoplus_{\{ \,y \in \Y|\, \DsssY(y) = b \}}
\E_{\cMb}^a(y) \quad
= \bigoplus_{\{\,x \in \X|\,\DsssX(x) = a,\,\DsssY(f(x)) = b\}}
i_x(G_{x,\cMb}).
\]\index{ $\E_{\cMb}^{a,b}$\vspace{1pt}}%
By \ref{lem:fiber1}, if $b<a$ or $b>a+d$, then $\E_{\cMb}^{a,b}$
(= $\E^{a,b}$ for simplicity) has an empty sum and
hence equals~0. The points contributing to
$\E^{a,a+d}$ are ones that are generic in their fiber while those
contributing to $\E^{a,a}$ are ones that are closed in their fiber.
\item For any integer~$a$ set
\[
\E_{\cMb}^a \set \;\; \bigoplus_{b\in \mathbb{Z}}\, \E_{\cMb}^{a,b}
\;\;= \;\; \E^{a,a} \oplus \E^{a,a+1} \oplus \cdots \oplus
\E^{a,a+d}.
\]
\end{enumerate}
\index{ $\E_{\cMb}^a$}%

Let us first record that for any integer $a$ there is a natural isomorphism
\begin{equation}
\label{eq:QE1aa}
(Ef^{\diamond}\cMb)^a \iso \E^a_{\cMb},
\end{equation}
defined as follows.
Let $x \in \X$ be such that $\DsssX(x) = a$. For
$y = f(x), M = \cMb(y)$, $q = \DsssY(y)$ and $p_1 = a - q + d$,
consider the natural isomorphism
\[
(Ef^{\diamond}\cMb)(x) \xto{(-1)^{qd} \text{ times \eqref{eq:gl2lo1}}}
H^{p_1}_{m_x}(M\otimes_{y}\omega_{x}) = G_{x,\cMb}.
\]
By summing up over all $x\in\X$ such that $\DsssX(x) = a$,
we obtain~\eqref{eq:QE1aa}.

Our next step is to define, for any integer $a$, a canonical
map $(f^{\diamond}\cMb)^a \to \E^a_{\cMb}$.
Let $y \in \Y$ be such that $\DsssY(y) = a+d$.
Let $x \in \X_y$ be such that $\DsssX(x) = a$. By~\ref{lem:fiber1},
$x$ is generic in the fiber and $m_y\OXx = m_x$. Hence,
with $M = \cMb(y)$ we obtain
\begin{equation}
\label{eq:QE1aaa}
(f^{\diamond}(\cMb,y))_x \cong M \otimes_{y} \omega_{x}
= \iG{m_x}(M \otimes_{y} \omega_{x})
= H^0_{m_x}(M \otimes_{y} \omega_{x}) = G_{x,\cMb}
\end{equation}
where the first of the three
equalities holds because $m_y$-torsionness of~$M$
implies that $M \otimes_{y} \omega_{x}$ is $m_y\OXx$-torsion.
For a fixed $y$, by varying~$x$ over all the generic points
of~$\X_y$ (equivalently, over points in~$\X_y$ such that
$\DsssX(x) = a$) we therefore obtain a map
\[
f^{\diamond}(\cMb,y) \xto{\text{canonical}} \bigoplus_x
i_x(f^{\diamond}(\cMb,y))_x \xto{\;\eqref{eq:QE1aaa}\;}
\bigoplus_x i_xG_{x,\cMb} = \E_{\cMb}^a(y).
\]
Finally, taking direct sums over all $y \in \Y$ such that $\DsssY(y) = a+d$,
we obtain a map
\begin{equation}
\label{eq:QE1ab}
(f^{\diamond}\cMb)^a \xto{\quad} \E^{a,a+d}_{\cMb}
\hookrightarrow \E^a_{\cMb}.
\end{equation}
Combining \eqref{eq:QE1aa} with \eqref{eq:QE1ab} results in a
natural graded map
\begin{equation}
\label{eq:QE1ac}
\eta_f(\cMb) \colon f^{\diamond}\cMb \xto{\quad} Ef^{\diamond}\cMb.
\end{equation}

The rest of this subsection is devoted to showing that \eqref{eq:QE1ac}
is a map of complexes, in fact a quasi-isomorphism whose
derived-category image is the Suominen isomorphism
$\Ssf(f^{\diamond}\cMb)$. We first show this in the special case when
$\cMb$ is concentrated in one spot (\ref{prop:QE6}, \ref{prop:QE6a})
and in \ref{prop:QE9}, \ref{prop:QE13}, the general case is considered.

\medskip

\begin{alem}
\label{lem:QE1a}
Let $f \colon (\X, \DsssX) \to (\Y, \DsssY)$ be an \emph{adic} smooth map
in $\bbFc$ of constant relative dimension $d$. Let $Z^{\bullet}$
be the filtration on~$\X$ induced by~$\DsssX$. Let $\cMb$ be a
complex in $\Coz_{\DsssY}(\Y)$ such that $\cMb$ is concentrated
in only one degree, so that $\cNb \set f^{\diamond}\cMb$ is
also concentrated in one degree, say $a$ (in other words,
$\cNb = \N^a[-a]$ and $\cMb = \M^{a+d}[-a-d\>]$).
For any $y \in \Y$, set $\N(y) \set f^{\diamond}(\cMb,y)$.
Then the following hold.
\begin{enumerate}
\item $H^a_{Z^a}\cNb = \iG{Z^a}\N^a = \N^a = H^a\cNb$.
\item $\E^n_{\cMb} = 0$ for $n<a$ or $n>a+d$.
\item $\E^{a,a+d}_{\cMb} = \E^a_{\cMb}$ or, in other words,
$\E^{a,b}_{\cMb} = 0$ for $b \ne a+d$.
Moreover, for any $x \in \X$ such that $\DsssX(x) = a$
and $\DsssY(f(x)) = b \ne a+d$, we also have
\[
H^a_x\cNb = H^0_x\N^a = \iG{x}\N^a = 0.
\]
\item Let $x \in \X$, $y = f(x)$ and $M = \cMb(y)$. If $\DsssX(x) = a$ and
$\DsssY(y) = a+d$ then
\[
M \otimes_{y} \omega_{x} \cong
(\N(y))_x = \iG{x}(\N(y)) = \iG{x}\N^a = \N^a_x.
\]
\end{enumerate}
\end{alem}

\begin{proof}
(i). Since $\cNb = \N^a[-a]$, hence $H^a_{Z^a}\cNb = \iG{Z^a}\N^a$
and $\N^a = H^a\cNb$. Hence the natural map
$\gamma \colon H^a_{Z^a}\cNb \to H^a\cNb$
is injective. Since $\cNb$ is $\DsssX$-CM (\ref{lem:quasim1}), therefore
by definition of Cohen-Macaulayness, $\gamma$ is also surjective.

(ii). If $n<a$ or $n>a+d$, then for any $x \in \X$ such that
$\DsssX(x) = n$, we have $\DsssY(f(x)) \ne a+d$ by \ref{lem:fiber1}.
Since $\cMb$ is concentrated in degree $a+d$ only, hence
$\cMb(f(x)) = 0$. Thus $G_{x,\cMb} = 0$ and the result follows.

(iii). Let $b \ne a+d$. Fix $x \in \X$ such that $\DsssX(x) = a$ and
$\DsssY(f(x)) =b$. In the canonical
decomposition $\N^a = \oplus_y \N(y)$, we have $y \ne f(x)$ for all $y$
contributing in the summation, since $y$ only ranges over points
with codimension $a+d$.
Hence, by \ref{lem:psm1}, $\R\iG{x}\N^a = 0$.
The desired conclusion follows.

(iv). Since $\N(y)$ is an $\Aqct(\X)$-module, the
first two isomorphisms follow from~\eqref{eq:QE1aaa}
and \ref{cor:loc4}. For any other point $y' \in \Y$ such
that $\DsssY(y') = a+d$ there is no specialization $y' \leadsto y$
and hence $(i_{y'}M')_y = 0$ for any $\cO_{\Y,y'}$-module $M'$.
In particular,
$(\N(y'))_x = (f^*i_{y'}(\cMb(y')) \otimes_{\X} \omega)_x = 0$.
In view of the canonical decomposition of $\N^a$,
the last two equalities of (iv) result.
\end{proof}

\medskip

\begin{alem}
\label{lem:QE1b}
Let notation and assumptions be as in \textup{\ref{lem:QE1a}}.
Then the map
\[
\psi \colon
(f^{\diamond}\cMb)^a = \N^a = H^a_{Z^a}(\cNb) \xto{\text{canonical}}
H^a_{Z^a/Z^{a+1}}(\cNb) = (E\cNb)^a = (Ef^{\diamond}\cMb)^a
\]
equals the map of~\eqref{eq:QE1ac} in degree $a$.
\end{alem}

\begin{proof}
Consider the following diagram wherein $\oplus_y, \oplus_x$ and $\oplus_z$
are indexed over the sets $\{y\in\Y | \DsssY(y) = a+d \}$,
$\{x\in\X_y | \DsssX(x) = a \}$, and
$\{z\in\X | \DsssX(z) = a \}$ respectively.
\[
\qquad \qquad
\begin{CD}
@. (Ef^{\diamond}\cMb)^a \\
@. @| \\
H^a_{Z^a}\cNb @>>> H^a_{Z^a/Z^{a+1}}\cNb  \\
@| @| \\
H^0_{Z^a}\N^a @>>> H^0_{Z^a/Z^{a+1}}\N^a @>{\sim}>> \oplus_zi_zH^0_{z}\N^a \\
@| @. @| \\
\makebox[0pt]{\hspace{-6em}$(f^{\diamond}\cMb)^a = $} \N^a
 @>>> \oplus_y \oplus_x i_x\N^a_x
 {\makebox[0pt]{\hspace{-4em}\smash{\raisebox{5ex}{$\scriptstyle{\Box_1}$}} }}
 @= \oplus_z i_z\iG{z}\N^a\\
@VV{\wr}V @AA{\wr {\makebox[0pt]{\hspace{12em}{$\scriptstyle{\Box_2}$}}}}A
 @AA{\wr}A \\
\oplus_y \N(y) @>>> \oplus_y \oplus_x i_x(\N(y))_x
 @= \oplus_y \oplus_x i_x\iG{x}(\N(y))  \\
@. @VV{\wr {\makebox[0pt]{\hspace{12em}{$\scriptstyle{\Box_3}$}}}}V
 @VV{\wr}V \\
@.  \oplus_y \oplus_xi_x(M_y \otimes_y \omega_{x}) @=
 \oplus_y \oplus_xi_x\iG{m_x}(M_y \otimes_y \omega_{x}) \\
\end{CD}
\]
The identifications in $\Box_2$ and $\Box_3$ follow from
\ref{lem:QE1a}(iii), (iv). Rest of the maps are the obvious
canonical ones. The commutativity of $\Box_1$ follows
from the description of the punctual decomposition in~\eqref{eq:coz0} while
that of the remaining rectangles is straightforward.
To prove the Lemma, it suffices to check that both, $\psi$
and the degree~$a$ component of~\eqref{eq:QE1ac},
occur as the composition over the
two different paths along the outer border
from $(f^{\diamond}\cMb)^a$ to
$(Ef^{\diamond}\cMb)^a$ in the above diagram.
The case of $\psi$ is clear, and the other one follows
upon examining the definitions involved.
\end{proof}

\medskip

\begin{aprop}
\label{prop:QE6}
Let notation and assumptions be as in \textup{\ref{lem:QE1a}}.
Then the graded map $\eta_f(\cMb) \colon \cNb \to E\cNb$
of \eqref{eq:QE1ac} is a quasi-isomorphism of complexes.
\end{aprop}
\begin{proof}
By assumption, $\cNb$ is concentrated only in degree $a$ and in
view of the isomorphism \eqref{eq:QE1aa},
by \ref{lem:QE1a}(ii) it follows that $(E\cNb)^n = 0$ for $n<a$.

Since $\cNb$ is $\DsssX$-CM, therefore $\cNb$ and $E\cNb$
are isomorphic in $\D(\X)$ say, via $\Ssf(\cNb)$,
and hence to prove the Proposition
it suffices to show that in~\eqref{eq:QE1ac},
$\N^a$ is mapped bijectively to the kernel of
$\delta = \delta^a_{E\cNb} \colon H^a_{Z^{a}/Z^{a+1}}\cNb \to
H^{a+1}_{Z^{a+1}/Z^{a+2}}\cNb$. By \ref{lem:QE1b}, it suffices to
show that the sequence of canonical maps
\[
H^a_{Z^{a}}\cNb \to H^a_{Z^{a}/Z^{a+1}}\cNb \xto{\;\delta\;}
H^{a+1}_{Z^{a+1}/Z^{a+2}}\cNb
\]
is exact.

By Cohen-Macaulayness of $\cNb$ we have
$H^{a}_{Z^{a+1}/Z^{a+2}}\cNb = 0$. Since $\delta$ corresponds
to the $a$th connecting homomorphism in the cohomology
long exact sequence associated to the exact sequence
\[
0 \xto{\quad} \iG{Z^{a+1}/Z^{a+2}} \xto{\quad} \iG{Z^{a}/Z^{a+2}}
\xto{\quad} \iG{Z^{a}/Z^{a+1}} \xto{\quad} 0,
\]
we see that $\ker(\delta) = H^a_{Z^{a}/Z^{a+2}}\cNb$. Therefore
it suffices to show that the natural map
$H^a_{Z^{a}}\cNb \to H^a_{Z^{a}/Z^{a+2}}\cNb$
is an isomorphism. Since $\cNb$ is CM, therefore
$H^a_{Z^{a+2}}\cNb = H^{a+1}_{Z^{a+2}}\cNb = 0$. Therefore from the
cohomology long exact sequence associated to the exact sequence
\[
0 \xto{\quad} \iG{Z^{a+2}} \xto{\quad} \iG{Z^{a}} \xto{\quad}
\iG{Z^{a}/Z^{a+2}} \xto{\quad} 0
\]
we obtain the desired conclusion.
\end{proof}

\smallskip

\begin{aprop}
\label{prop:QE6a}
Let notation and assumptions be as in \textup{\ref{prop:QE6}}. Then
the derived-category image of $\eta_f(\cMb)$ is $\mathsf{S}(\cNb)$.
\end{aprop}
\begin{proof}
Since $\cNb$ is concentrated in degree $a$ and
$(E\cNb)^n = 0$ for $n < a$ (see proof of~{\ref{prop:QE6}}),
therefore any derived-category isomorphism
$\cNb \iso E\cNb$ is represented by a quasi-isomorphism
$\cNb \to E\cNb$ which is moreover completely
described by the corresponding (injective) map in degree~$a$,
$\N^a \to (E\cNb)^a$. Let $S$ be the quasi-isomorphism
representing $\Ssf(\cNb)$ and $S^a$ the induced map in degree~$a$.
Set $\eta \set \eta_f(\cMb)$ and let $\eta^a$ be the
induced map in degree~$a$. The Proposition amounts to saying
that $S^a = \eta^a$.

For any sheaf $\F$ on $\X$ and for any family of $\OXx$-modules
$G_x$ where $x$ varies over points in $\X$, if $\Hom_{\OXx}(\F_x, G_x) = 0$
for all $x$, then $\Hom_{\X}(\F, \oplus_xi_xG_x)= 0$.
In particular, since $(E\cNb)^a$ lies on the $Z^a/Z^{a+1}$-skeleton,
to show that $S^a = \eta^a$, it suffices to show that for any
$x \in \X$ with $\DsssX(x) = a$, the natural induced maps
$S^a_x$ and~$\eta^a_x$ from~$\N^a_x$ to~$(E\cNb)^a_x$ are equal.
Moreover, by \ref{lem:QE1a}(iii), it suffices to
confine our attention to $x \in \X$ such that
$\DsssX(x) = a$ and $\DsssY(f(x)) = a+d$.
In particular, the identifications of \ref{lem:QE1a}(iv) apply.

It suffices to show that for $\phi= \eta^a_x$ or $\phi = S^a_x$,
the following maps compose to the identity map
\[
H^a_x\cNb = H^0_x\N^a = \N^a_x \xto{\; \phi \;} (E\cNb)^a_x =
(E\cNb)(x) \xto{\eqref{eq:coz1}} H^a_x\cNb.
\]
Indeed, since the remaining maps are isomorphisms, $\phi$ is uniquely
determined by the condition that the composition be identity.

\pagebreak

In the case of $\phi= \eta^a_x$ we conclude using the outer
border of the following diagram of natural isomorphisms
whose top row gives $\eta^a_x$ and whose commutativity is
straightforward to verify.
{\small{
\[
\hspace{-5em}
\begin{CD}
\N^a_x @>>> M\otimes_y\omega_x @>>> H^0_{m_x}(M\otimes_y\omega_x)
 @>{(-1)^{qd}\times\eqref{eq:gl2lo1}^{-1}}>> (E\cNb)(x)
{\makebox[0pt]{\hspace{5em}$ = (E\cNb)^a_x$}} \\
@| @. @V{\hspace{-10em}\text{see \ref{lem:QE1a}(iv)}}VV
 @VV{\eqref{eq:coz1}}V \\
\iG{x}\N^a @.
 \makebox[0pt]{\raisebox{.6ex}{\makebox[10em]{\hrulefill}}
 \hspace{-10em}\!\!\makebox[10em]{\hrulefill}}
 @. H^0_x\N^a @= H^a_x\cNb
\end{CD}
\]
}}

\smallskip
In the case of $\phi = S^a_x$ we use the following commutative diagram
of natural isomorphisms, where the top row composes to the identity map
by~\ref{cor:coz4b}.
\[
\begin{CD}
H^a_x\cNb @>{H^a_xS}>> H^a_xE\cNb @>{\ref{lem:coz1a}\textup{(i)[c]}}>> H^a_x\cNb \\
@| @V{\hspace{-12em}\;\;H^a_x\,\sigma_{\le a}}V{H^a_x\,\sigma_{\le a}}V \\
H^a_x(\N^a[-a]) @>{H^a_x(S^a[-a])}>> H^a_x((E\cNb)^a[-a]) @.
 \rotatebox{90}{\makebox[0pt]{\hspace{-4em}{\raisebox{.5ex}{\makebox[10em]
 {\hrulefill}}}\hspace{-10em}{\makebox[10em]{\hrulefill}}}} \\
@| @| \\
\iG{x}\N^a @>{\iG{x}S^a}>> \iG{x}(E\cNb)^a \\
@| @| \\
\N^a_x @>{S^a_x}>> (E\cNb)^a_x @>{\eqref{eq:coz1}}>> H^a_x\cNb
\end{CD}
\]
\end{proof}

Now we tackle the general case when $\cMb$ is no longer assumed to be
concentrated in one spot.

For integers $r,s,t$ let $\delta^{r,s,t}_{\cMb}$ denote the
map $\E^{r,s}_{\cMb} \to \E^{r+1,t}_{\cMb}$ defined as follows.
\begin{align}
\E^{r,s}_{\cMb} \xto{\text{inclusion}} \E^{r}_{\cMb}
 &\xto{\,\quad\eqref{eq:QE1aa}\quad\,} (Ef^{\diamond}\cMb)^r \notag \\
 &\xto{\text{differential}} (Ef^{\diamond}\cMb)^{r+1}
 \xto{{\eqref{eq:QE1aa}}^{-1}} \E^{r+1}_{\cMb} \xto{\text{projection}}
 \E^{r+1,t}_{\cMb} \notag
\end{align}

\medskip

\begin{alem}
\label{lem:QE7}
Let $f, \X, \DsssX, \Y, \DsssY, d,$ be as in \textup{\ref{lem:QE1a}}. Let
$\cMb \in \Coz_{\DsssY}(\Y)$, no longer assumed to be concentrated
in one degree. For any integers $r,s,t$ the following hold.
\begin{enumerate}
\item If $s>t$, then $\delta^{r,s,t}_{\cMb} = 0$.
\item Assume $s \le t$. For any integer $n$, if  $n \le s$ so that
$\E^{r,s}_{\cMb} = \E^{r,s}_{\sigma_{\ge n}\cMb}$ and
$\E^{r+1,t}_{\cMb} = \E^{r+1,t}_{\sigma_{\ge n}\cMb}$,
then $\delta^{r,s,t}_{\cMb} = \delta^{r,s,t}_{\sigma_{\ge n}\cMb}$.
Similarly, if $n \ge t$, then
$\delta^{r,s,t}_{\cMb} = \delta^{r,s,t}_{\sigma_{\le n}\cMb}$.
\item For any integer $a$,
$\delta^{a,a+d,a+d}_{\cMb} = \delta^{a,a+d,a+d}_{\M^{a+d}[-a-d\>]}$.
\end{enumerate}
\end{alem}

\begin{proof}
(i) We claim that more generally, any map from $\E^{r,s}_{\cMb}$
to $\E^{r+1,t}_{\cMb}$ is zero when $s>t$. Indeed, let $i_xG_x$
be a component of $\E^{r,s}_{\cMb}$ and $i_{x'}G_{x'}$ a component
of~$\E^{r+1,t}_{\cMb}$. It suffices to show that
$\Hom_{\X}(i_xG_x, i_{x'}G_{x'}) = 0$ or, more generally, that
$x' \notin \ov{\{x\}}$. By definition, $\DsssY(f(x))=s$ and
$\DsssY(f(x')) = t$. Therefore $s>t \Lra f(x') \notin \ov{\{f(x)\}}
\Lra x' \notin \ov{\{x\}}$.

(ii) Suppose $n \le s$. The canonical map $\sigma_{\ge n}\cMb \to \cMb$
induces a map of complexes
$Ef^{\diamond}\sigma_{\ge n}\cMb \to Ef^{\diamond}\cMb$.
Consider the following diagrams where the bottom row of the
diagram on the left is the same as the top row of the diagram
on the right.
\[
\begin{CD}
\E^{r,s}_{\sigma_{\ge n}\cMb} @= \E^{r,s}_{\cMb} \\
@VVV @VVV \\
\E^{r}_{\sigma_{\ge n}\cMb}  @>>> \E^{r}_{\cMb}  \\
@V{\wr}V{\hspace{5em} \Box_1}V @V{\wr}VV \\
(Ef^{\diamond}\sigma_{\ge n}\cMb)^r @>>> (Ef^{\diamond}\cMb)^r \\
@VVV @VVV \\
(Ef^{\diamond}\sigma_{\ge n}\cMb)^{r+1} @>>> (Ef^{\diamond}\cMb)^{r+1}
\end{CD}
\begin{CD}
(Ef^{\diamond}\sigma_{\ge n}\cMb)^{r+1} @>>> (Ef^{\diamond}\cMb)^{r+1} \\
@V{\wr}V{\hspace{5em} \Box_2}V @V{\wr}VV \\
\E^{r+1}_{\sigma_{\ge n}\cMb}  @>>> \E^{r+1}_{\cMb}  \\
@VVV @VVV \\
\E^{r+1,t}_{\sigma_{\ge n}\cMb} @= \E^{r+1,t}_{\cMb} \\
\end{CD}
\]
The left columns of the two diagrams together give
$\delta^{r,s,t}_{\sigma_{\ge n}\cMb}$ while the right ones
right give $\delta^{r,s,t}_{\cMb}$.
The rectangles $\Box_1,\Box_2$ commute because the isomorphism
in~\eqref{eq:QE1aa}, which uses \eqref{eq:gl2lo1},
factors through truncations.
Commutativity of the remaining rectangles is obvious
and so the desired conclusion follows.
A similar argument works for the corresponding statement when $n \ge t$.

(iii) Using $r=a$ and $s = t = n = a+d$ in (ii) we conclude
by truncating on the left and right.
\end{proof}

\begin{aprop}
\label{prop:QE9}
Let notation and assumptions be as in \ref{lem:QE7}. Then the
graded map $\eta_f(\cMb) \colon f^{\diamond}\cMb \to Ef^{\diamond}\cMb$
defined in \textup{\eqref{eq:QE1ac}} is a map of complexes.
\end{aprop}

\begin{proof}
We have to verify that for any integer $a$,
the following diagram commutes.
\[
\begin{CD}
(f^{\diamond}\cMb)^a @>>> \E^{a,a+d}_{\cMb} @>>> \E^{a}_{\cMb}
 @>>> (Ef^{\diamond}\cMb)^a \\
@VVV @. @. @VVV \\
(f^{\diamond}\cMb)^{a+1} @>>> \E^{a+1,a+1+d}_{\cMb} @>>> \E^{a+1}_{\cMb}
 @>>> (Ef^{\diamond}\cMb)^{a+1}
\end{CD}
\]
Let us first verify that for any integer $n \ge 0$, the composite map
\[
(f^{\diamond}\cMb)^a \xto{\quad} \E^{a,a+d}_{\cMb}
\xto{\delta^{a,a+d, a+d-n}_{\cMb}} \E^{a+1, a+d-n}_{\cMb}
\]
is zero. If $n>0$ then we conclude by \ref{lem:QE7}(i). If $n=0$,
then by \ref{lem:QE7}(iii) we may assume without loss of generality
that $\cMb$ is concentrated in degree $a+d$ only and then we conclude
by \ref{prop:QE6}.

We have thus reduced the Proposition to proving that the
following diagram commutes.
\[
\begin{CD}
(f^{\diamond}\cMb)^a @>{\eqref{eq:QE1ab}}>> \E^{a,a+d}_{\cMb} \\
@VVV @VV{\delta^{a,a+d, a+d+1}_{\cMb}}V \\
(f^{\diamond}\cMb)^{a+1} @>>{\eqref{eq:QE1ab}}> \E^{a+1,a+1+d}_{\cMb}
\end{CD}
\]

\pagebreak[3]

By definition, \eqref{eq:QE1ab} respects the grading
on each side indexed by the points of~$\Y$. In particular,
it suffices to verify that for any immediate specialization
$y \leadsto y'$ in $\Y$ such that $\DsssY(y)= a+d$,
the following diagram commutes, with the columns being,
as above, induced by the differential
of $\cMb$ and $Ef^{\diamond}\cMb$.
\[
\begin{CD}
f^{\diamond}(\cMb, y) @>>> \E^{a}_{\cMb}(y) \\
@VVV @VVV \\
f^{\diamond}(\cMb, y') @>>> \E^{a+1}_{\cMb}(y')
\end{CD}
\]
Before simplifying this diagram further, we need some notation.
Set $M \set \cMb(y)$, $M' \set \cMb(y')$.
By \ref{lem:fiber1}, the definition of
$\E^{a+1}_{\cMb}(y') \subset \E^{a+1,a+1+d}_{\cMb}$
involves a finite direct sum of modules of the type $i_{x'}G_{x'}$
where $x'$ ranges over the set of generic points of the fiber
over~$y'$ and $G_{x'} = H^0_{m_x'}(M'\otimes_{y'}\omega_{x'})$.
Fix one such point $x'$.
Let $x_1, \ldots, x_k, \ldots, x_l$ be all the points that
are generic in the fiber over $y$ where $x_1, \ldots, x_k$ are exactly
those that specialize to $x'$. Any such specialization
$x_j \leadsto x'$ is necessarily immediate since
$\DsssX(x_j) = a$ and $\DsssX(x') = a+1$.
Note that the induced map from the summand
$\bigoplus_{j>k}i_{x_j}G_{x_j}$ of~$\E^{a}_{\cMb}(y)$
(where $G_{x_j} = H^0_{m_{x_j}}(M\otimes_{y}\omega_{x_j})$ for any~$j$)
to~$\E^{a+1}_{\cMb}(y')$ is zero
because, $x' \notin \ov{\{x_j\}}$ for $j>k$ implies more generally that
$\Hom_{\OX}(i_{x_j}G_{x_j}, i_{x'}G_{x'}) = 0$.

Coming back to the previous diagram, we may restrict the target
$\E^{a+1}_{\cMb}(y')$ of the diagram to the component $i_{x'}G_{x'}$
and moreover, since for any sheaf~$\F$, we have
$\Hom_{\OX}(\F, i_{x'}G_{x'}) \iso \Hom_{\cO_{\X,x'}}(\F_{x'}, G_{x'})$,
we may localize the previous diagram at $x'$.
Then verifying its commutativity reduces to verifying
that of the rectangle on the right in the following diagram,
where, as expected, $\alpha$ is induced by the differential of
$f^{\diamond}\cMb$ and $\beta$ is induced by the differential of
$Ef^{\diamond}\cMb$ via \eqref{eq:QE1aa}.
\begin{equation}
\label{cd:QE10}
\begin{CD}
M \otimes_{y'} \omega_{x'} @>{\sim}>> (f^{\diamond}(\cMb, y))_{x'}
 @>{\text{to stalk at $x_j$}}>{\text{and then } \eqref{eq:QE1aaa}}>
 \bigoplus_{j \le k} H^0_{m_{x_j}}(M\otimes_{y}\omega_{x_j}) \\
@V{\partial \otimes 1}VV @V{\alpha}VV @V{\beta}VV \\
M' \otimes_{y'} \omega_{x'} @>{\sim}>>
(f^{\diamond}(\cMb, y'))_{x'} @>{\eqref{eq:QE1aaa}}>>
H^0_{m_x'}(M'\otimes_{y'}\omega_{x'})
\end{CD}
\end{equation}
Let $\partial$ denote the $\cO_{\Y,y'}$-linear map $M \to M'$
induced by the differential of~$\cMb$. The rectangle on the left clearly
commutes. It remains to prove that the outer border of~\eqref{cd:QE10}
commutes. We apply \ref{lem:hosp0a} in this situation with essentially
the same notation as used there. (The map $h$ and its relative dimension
$n$ are now denoted by $f, d$ respectively.) Note that in this case
we have $p_1 = 0$ and $\cL = \omega$. Let~$\I$ be the largest coherent
ideal defining the closed set $\ov{\{x_1,\ldots,x_k\}}$ and
let $I= \I_{x'}$ as in~\ref{lem:hosp0a}.
We make the following claims.
\begin{enumerate}
\item The $\cO_{\X,x'}$-modules $M \otimes_{y'} \omega_{x'}$
and $M' \otimes_{y'} \omega_{x'}$ are $I$-torsion.
\item The top and bottom rows of \eqref{cd:QE10}
give the maps $\mu_1,\mu_2$ respectively of~\ref{lem:hosp0a}.
\item The map $\beta$ of \eqref{cd:QE10} equals $(-1)^{d}\psi$
for $\psi$ as in~\ref{lem:hosp0a}.
(Note that $d$ corresponds to $n$ used in~\ref{lem:hosp0a}.)
\end{enumerate}
Assuming these claims, \eqref{cd:QE10} commutes by \ref{lem:hosp0a},
thus proving the Proposition.

Both (ii) and (iii) are straightforward to verify. The rest of this proof
is devoted to showing~(i).

Let $\sK$ be the largest open coherent ideal in~$\OY$ defining
the closed set $\ov{\{y\}}$. Since~$f$ is adic, therefore
$\J \set \sK\OX$ is an open ideal in~$\OX$. Since the ordinary
scheme $X \set (\X, \OX/\J)$ is smooth over the ordinary
scheme $Y \set (\Y, \OY/\sK)$,
therefore the generic points of (the irreducible components of) $X$
all map to the generic point of~$Y$, viz.,~$y$. In particular,
the generic points of the fiber over~$y$ correspond exactly to the
generic points of~$X$. Note that $\J \subset \I$.
Set $\ov{\I} \set \I/\J$. Then the closed set defined by~$\ov{\I}$
contains $\Spec(\cO_{X,x'})$, and hence the stalk~$\ov{I}$
of~$\ov{\I}$ at~$x'$ is a nilpotent ideal. But~$\ov{I}$ can be canonically
identified with $\I_{x'}/\J_{x'}$. It follows that proving~(i)
is equivalent to proving that $M \otimes_{y'} \omega_{x'}$
and $M' \otimes_{y'} \omega_{x'}$ are $\J_{x'}$-torsion.
Since $\J_{x'} = \sK_{y'}\cO_{\X,x'}$ we therefore reduce to
showing that $M,M'$ are $\sK_{y'}$-torsion as $\cO_{\Y,y'}$-modules.
In the case of~$M$ we conclude from the fact that
$\sK_{y'}\OYy = \sK_y = m_y$ and~$M$ is an $m_y$-torsion
$\OYy$-module. In the case of~$M'$ we conclude from the fact that
$\sK_{y'} = m_{y'}$.
\end{proof}

\begin{aprop}
\label{prop:QE13}
With notation and assumptions as in\/ \textup{\ref{prop:QE9},} the
natural map $\eta_f(\cMb) \colon f^{\diamond}\cMb \to Ef^{\diamond}\cMb$
is a quasi-isomorphism whose derived-category
image $Q\eta_f(\cMb)$ is $\Ssf(f^{\diamond}\cMb)$.
\end{aprop}

\begin{proof}
Set $\eta \set \eta_f$.
It suffices to verify that $Q\eta(\cMb) = \Ssf(f^{\diamond}\cMb)$,
because that implies that $\eta$ is a quasi-isomorphism.

By \ref{cor:coz4b},
the validity of $Q\eta(\cMb) = \Ssf(f^{\diamond}\cMb)$
amounts to that of the punctual statement that,
for any $x \in \X$, with $p = \DsssX(x)$, the composite isomorphism
\[
H^p_xf^{\diamond}\cMb \xto{H^p_x\eta(\cMb)} H^p_xEf^{\diamond}\cMb
\xrightarrow[\text{for $\cFb = f^{\diamond}\cMb$}]
{\textup{\ref{lem:coz1a}(i)[c]}}
H^p_xf^{\diamond}\cMb
\]
is the identity. Let us verify this punctual statement.

Fix $x \in \X$ with $p = \DsssX(x)$ and $q = \DsssY(f(x))$.
Since ${H^p_x\eta(-)}$ and the
isomorphism obtained from \ref{lem:coz1a}(i)[c] are functorial in $\cMb$,
in view of the natural isomorphisms (\ref{lem:psm2})
\[
H^p_xf^{\diamond}\cMb \iso H^p_xf^{\diamond}\sigma_{\le q}\cMb,
\qquad
H^p_xEf^{\diamond}\cMb \iso H^p_xEf^{\diamond}\sigma_{\le q}\cMb,
\]
we may replace $\cMb$ by $\sigma_{\le q}\cMb$. By a similar argument
we may truncate in the other direction and thus we reduce to the case
when $\cMb = \M^q[-q]$ is concentrated in only one degree.

By \ref{prop:QE6a} for $\cMb = \M^q[-q]$, we already have
$Q\eta(\cMb) = \Ssf(f^{\diamond}\cMb)$ and hence the punctual
statement corresponding to $x$ now follows by referring back
to~\ref{cor:coz4b}.
\end{proof}

\begin{arem}
Suppose $f\colon (\X,\,\DsssX) \to (\Y,\,\DsssY)$ is as in 
Lemma\,\ref{lem:QE1a}. Let
\[
L_f({\mathcal{M}}^\bullet) (= L_f) \colon 
{{\mathbb{E}}}_f({\mathcal{M}}^\bullet)=E(f^\diamond{\mathcal{M}}^\bullet)
\iso f^\sharp({\mathcal{M}}^\bullet)
\]
be the isomorphism in the construction of $f^\sharp({\mathcal{M}}^\bullet)$
in \S\,\ref{subsec:fin-sm} (see the second paragraph of {\it loc.cit.}).
We would like to summarize the relationship between the global maps
$L_f$, $\eta_f$ and various natural local isomorphisms (e.g. \eqref{eq:coz1},
\eqref{eq:itloco0} or the isomorphisms in Theorem\,\ref{thm:huang1}) in
a special case, viz.~when ${\mathcal{M}}^\bullet$ is concentrated at a
point. The purpose is to gather together the many sign twists we have used
into one place.

Suppose $x$ is a point on $\X$. Let $p=\DsssX(x)$, $y=f(x)$, $q=\DsssY(y)$,
$S={\widehat{\OXx}}$, $R=\widehat{\OYy}$,
$\varphi\colon R\to S$ the natural map induced by $f$,
$t=\dim(S/m_RS)=p-q+d$, $M$ a zero dimensional $R$--module, ${\mathcal F}
=i_yM$, ${\mathcal{M}}^\bullet={\mathcal{F}}[-q]$ and set
\[
\varepsilon(x) := \DsssX(x)\DsssY(y) + \DsssX(x).
\]
The natural quasi-isomorphism
\[
{\mathbb{E}}_f({\mathcal{M}}^\bullet)(x)[-p] = 
\iG{x}{{\mathbb{E}}_f({\mathcal{M}}^\bullet)} \longrightarrow
\R\iG{x}{{\mathbb{E}}_f{{\mathcal{M}}^\bullet}}
\]
obtained by applying $\iG{x}$ to an injective resolution 
${\mathbb{E}}_f({\mathcal{M}}^\bullet)\to {\mathcal{I}^\bullet}$
of the flasque complex ${\mathbb{E}}_f({\mathcal{M}}^\bullet)$ gives
on applying $H^p$ the isomorphism
\[
\rho\colon {\mathbb{E}}_f({\mathcal{M}}^\bullet)(x) \xrightarrow{\sim} 
H^p_x{\mathbb{E}}_f({\mathcal{M}}^\bullet),
\]
which is the inverse of the map in \ref{lem:coz1a}(i)[b]. In view of
Propositions \ref{cor:coz4b} and \ref{prop:QE13} and the definitions
in \eqref{eq:nolab3}, \eqref{eq:coz3}, \eqref{eq:itloco0}, \eqref{eq:gl2lo1},
\eqref{eq:fin-sm1}, we see that the following diagram commutes up to
a sign of $(-1)^{\varepsilon(x)}$:
\[
\xymatrix{
H^p_x({\mathbb{E}}_f({\mathcal{M}}^\bullet)) & & &
{\mathbb{E}}_f({\mathcal{M}}^\bullet)(x) \ar[lll]_{\Iso}^{\rho}
\ar[d]^{L(x)} \\
H^p_x(f^\diamond{\mathcal{M}}^\bullet) \ar[u]^{H^p_x(Q_\X\eta_f)} \ar@{=}[d]
& & & (f^\sharp{\mathcal{M}}^\bullet)(x) \ar@{=}[d] \\
H^p_x(f^*\F[-q]\otimes\omega_f[d]) \ar@{=}[d] & & & \varphi_\sharp{M}
\ar[d]^{\ref{thm:huang1},I(i)} \\
H^t_x(f^*\F\otimes\omega_f) \ar[rrr]_{\eqref{eq:nolab3},\eqref{eq:itloco0}}
& & & H^t_{m_S}(M\otimes\omega_{S/R}) \ar@{}[uuulll]|{(-1)^{\varepsilon(x)}}
}
\]
\end{arem}
\medskip


\subsection{Truncations} In this subsection we gather together
some properties of $f^{\sharp}$ with a view to future applications.
These results are not used in the next subsection.

\newcommand{\fp}{{\mathfrak{p}}}
\newcommand{\De}{{\Delta}}
\newcommand{\Cozs}[1]{{\mathrm {Coz}}_{#1}}
\newcommand{\sh}[1]{{#1}^{\sharp}}
\newcommand{\shr}[1]{{#1}^{\boldsymbol{!}}}
\newcommand{\sha}[1]{{#1}^{\boldsymbol{\sharp}}}
\newcommand{\fm}{{\mathfrak{m}}}
\newcommand{\Tr}[1]{{\mathrm {Tr}}_{{#1}}}
\newcommand{\ttr}[1]{{\mathrm {\tau}}_{#1}}
\newcommand{\bbF}{{\mathbb F}}

\newtheorem{thm}{Theorem}[subsection]

\newcommand{\Cref}[1]{Corollary~\textup{\ref{#1}}}
\newcommand{\Dref}[1]{Definition~\textup{\ref{#1}}}
\newcommand{\Eref}[1]{Example~\textup{\ref{#1}}}
\newcommand{\Lref}[1]{Lemma~\textup{\ref{#1}}}
\newcommand{\Pref}[1]{Proposition~\textup{\ref{#1}}}
\newcommand{\Rref}[1]{Remark~\textup{\ref{#1}}}
\newcommand{\Sref}[1]{Section~\textup{\ref{#1}}}
\newcommand{\Ssref}[1]{Subsection~\textup{\ref{#1}}}
\newcommand{\Tref}[1]{Theorem~\textup{\ref{#1}}}

Recall that if $\F$ is an abelian sheaf on a topological space $X$,
$x$ a point on~$X$ and~$G$ an abelian group, then
$\Hom_{\text{Sheaves}}(\F, i_xG) \iso \Hom_{\text{Groups}}(\F_x, G)$.
In particular, if $\F_x = 0$, then any map
from~$\F$ to~$i_xG$ is zero.

\begin{alem}
\label{lem:zeroC}
Let $f:(\X,\,\De')\to (\Y,\,\De)$ be a map in $\bbFc$ and let
$x\in\X$, $y\in\Y$ be points such that $\De'(x)=\De(y)$. If~$x$
is not closed in its fiber or if $y \ne f(x)$, then
for any $\OXx$-module $C$ and any $\OYy$-module~$M$,
$\Hom_{\OY}(f_*i_xC, i_yM) = 0$.
\end{alem}
\begin{proof}
By definition, $\De' = f^{\sharp}\De$ and hence $\De'(x) \le \De(f(x))$
with equality holding only if~$x$ is closed in its fiber.
Now note that $f_*i_xC = i_{f(x)}C$. Thus if we assume
$\Hom_{\OY}(i_{f(x)}C, i_yM) \ne 0$, then $(i_{f(x)}C)_y \ne 0$
and hence $y \in \ov{\{f(x)\}}$. But this implies that
$\De(f(x)) \le \De(y)$ with equality holding only if $y = f(x)$.
Since $\De'(x) = \De(y)$ we therefore obtain that~$x$ is
closed in its fiber and $y = f(x)$, which is a contradiction.
\end{proof}

\begin{adefi}
\label{def:del-x-y}
Let $(\Y,\,\De) \in \bbFc$ and
let $\cMb\in \Cozs{\De}(\Y)$. For any $\OY$-endomorphism of the
total $\OY$-module $T(\cMb)$ of the graded module $\cMb$, we
get, for every pair of points $y$, $y'$ in $\Y$ an $\OY$-linear map
$i_y(\cMb(y))\to i_{y'}(\cMb(y'))$.
We denote by $\delta_{\cMb}(y,y'):i_y(\cMb(y))
\to i_{y'}(\cMb(y'))$ the map induced by the coboundary map on~$\cMb$.
\end{adefi}

\begin{alem}
\label{lem:prefilter}
Let $f: (\X,\,\De')\to (\Y,\,\De)$
be a map in $\bbFc$ and let $\cCb$ be a Cousin complex on $(\X,\,\De')$.
Then $\delta_{\cCb}(x,x')=0$ if~$f(x')$ is not a specialization of~$f(x)$.
\end{alem}

\begin{proof}
If $f(x')$ is not a specialization of $f(x)$ then $x'$ is not a
specialization of $x$ and hence any map from $i_x(\cCb(x))$ to
$i_{x'}(\cCb(x'))$ is necessarily zero.
\end{proof}


\begin{alem}
\label{lem:exact}
Let $f:(\X,\,\De')\to (\Y,\,\De)$ be
a smooth morphism in $\bbFc$. Let
\[
0 \longrightarrow \cFb_1 \longrightarrow \cFb_2 \longrightarrow \cFb_3
\longrightarrow 0
\]
be an exact sequence of complexes in $\Cozs{\De}(\Y)$. Then the
induced sequence
\[
0 \longrightarrow \sh{f}\cFb_1 \longrightarrow \sh{f}\cFb_2 \longrightarrow
\sh{f}\cFb_3 \longrightarrow 0
\]
of complexes in $\Cozs{\De'}(\X)$ is exact.
\end{alem}

\begin{proof}
Let $x \in \X$. For any integer $i$ set
$L_x^i(-)  \set H^i_x\R\iGp{\X}(f^*(-)\otimes_{\X}\omega_f[d\>])$.
Set $p \set \De'(x)$. Since $f$ is a flat map
and $\omega_f$ a flat $\OX$-module, we get
a long exact sequence
\[
\ldots \xto{\quad} L^{p-1}_x(\cFb_3) \xto{\quad} L^{p}_x(\cFb_1)
\xto{\quad} L^{p}_x(\cFb_2) \xto{\quad} L^{p}_x(\cFb_3)
\xto{\quad} L^{p+1}_x(\cFb_1) \xto{\quad} \ldots
\]
By Cohen-Macaulayness (\ref{lem:quasim1}) the $L^i_x$-terms
vanish for $i \ne p$ so that we obtain a short exact sequence
\[
0 \xto{\quad} L^{p}_x(\cFb_1) \xto{\quad} L^{p}_x(\cFb_2)
\xto{\quad} L^{p}_x(\cFb_3) \xto{\quad} 0.
\]
Since this holds for every $x \in \X$, in view of the isomorphism
$L^p_x(-) \iso (f^{\sharp}(-))(x)$ (see \eqref{eq:fin-sm1}, \eqref{eq:gl2lo1})
the desired conclusion follows.
\end{proof}

Let $(\X,\,\De')\xrightarrow{f} (\Y,\,\De)$ be a map in $\bbFc$.
Then any Cousin complex in $\Coz_{\De'}(\X)$ admits natural filtrations
arising from the codimension function on $\Y$. The filtrations
are obtained as follows.
For $p,n \in {\mathbb Z}$ set
\[
\Sigma^p_n\set \{x\in\X\,|\,\De'(x)=n, \De(f(x))\ge p\}
\]
and
\[
S^n_p\set \{x\in\X\,|\,\De'(x)=n, \De(f(x))\le p\}.
\]

Suppose $\cCb\in \Cozs{\De'}(\X)$.
Then $\cCb$ has a natural decreasing filtration
$\{F^p(\cCb)\}_p$ by subcomplexes in $\Cozs{\De'}(\X)$
defined as follows. For $p\in {\mathbb Z}$
define the $n$th graded piece of $F^p(\cCb)$ by
\[
F^p(\cCb)^n \set \bigoplus_{x\in\Sigma^p_n}i_x\cCb(x).
\]
By \ref{lem:prefilter}, $F^p(\cCb)$ is stable under the action of the
coboundary map $\delta_{\cCb}$ on $\cCb$ so that the
restriction of $\delta_{\cCb}$ to $F^p(\cCb)$ gives a differential
$\delta_{F^p(\cCb)}$ on $F^p(\cCb)$ that makes it a subcomplex of $\cCb$.
In a similar vein we can define, for $p\in {\mathbb Z}$,
a complex $G_p(\cCb)$ whose $n$th graded piece is
\[
G_p(\cCb)^n \set \bigoplus_{x\in S^n_p} i_x\cCb(x).
\]
Viewing $G_p(\cCb)$ as the cokernel of the inclusion $F^{p+1}(\cCb)\to \cCb$,
we obtain a differential $\delta_{G_p(\cCb)}$ on $G_p(\cCb)$.
Note that for any $x,x' \in \X$ if
$\De(f(x)), \De(f(x'))\le p$, then
$\delta_{G_p(\cCb)}(x,x')=\delta_\cCb(x,x')$.
There results a short exact sequence of Cousin complexes on $(\X,\,\De')$:
\begin{equation}
\label{eq:F-G-exact}
0\lra F^{p+1}(\cCb) \lra \cCb \lra G_p(\cCb) \lra 0.
\end{equation}
Note that $F^p$ and $G_p$ are functorial on $\Cozs{\De'}(\X)$.

Let $\cFb \in \Coz_{\De}(\Y)$.
With $\sigma_{\ge p}$ and $\sigma_{\le p}$ the truncation
functors (\S\ref{subsec:conv}, \eqref{conv10})
we claim that the following hold.
\begin{equation}\label{eq:brutal}
\begin{split}
F^p(\sh{f}\cFb) & = \sh{f}(\sigma_{\ge p}\cFb) \\
G_p(\sh{f}\cFb) & = \sh{f}(\sigma_{\le p}\cFb) \\
F^p(\sh{f}\cFb)/F^{p+1}(\sh{f}\cFb) & = \sh{f}(\F^p[-p]).
\end{split}
\end{equation}
All these are straightforward to verify using the arguments
in~(ii) and~(iii) of~\ref{lem:QE7}.

\subsection{\'Etale maps}
\label{subsec:etale}\index{e@{\'etale map}}
By an \emph{\'etale} map of formal schemes we mean a smooth map
of relative dimension~0.
In the case of an \'etale map $f$, the formula for $f^{\sharp}$
can be simplified further.
We begin with a few preliminaries concerning a criterion for
a complex to be Cousin.

\begin{alem}
\label{lem:etale1}
Let $(\X, \Delta) \in \bbFc$. A complex $\cFb$ on~$\X$ is $\Delta$-Cousin
if and only if it satisfies the following condition:
For any $x \in \X$ and integers $i,j$ such that $i \ne 0$ or
$j \ne \Delta(x)$, we have $H^i_x\F^j = 0$.
\end{alem}
\begin{proof}
If $\cFb$ is $\Delta$-Cousin, then the condition stated in the Lemma
follows from \ref{lem:coz1a}(i)[a]. Let $Z^{\bullet}$ be the filtration
of~$\X$ induced by~$\Delta$. For the `if' part, we must show that
for any integer $j$, $\F^j$ lies on the $Z^j/Z^{j+1}$-skeleton,
i.e.,
\[
\F^j = \iG{Z^j}\F^j \qquad \text{and} \qquad
H^0_{Z^{j+1}}\F^j = 0 = H^1_{Z^{j+1}}\F^j.
\]
Note that the condition of the Lemma can be rephrased as saying that
if $i \ne 0$ or $j \ne n$ then $H^i_{Z^n/Z^{n+1}}\F^j = 0$ .

Recall that for $p \ll 0$, $Z^p = \X$. Therefore to show that
$\F^j = \iG{Z^j}\F^j$, it suffices to show that the natural maps
\[
\iG{Z^j}\F^j \to \iG{Z^{j-1}}\F^j \to \iG{Z^{j-2}}\F^j \to \ldots
\to \iG{\X}\F^j = \F^j
\]
are isomorphisms. For every $n$ there are exact sequences
\[
0 \to \iG{Z^{n+1}} \to \iG{Z^{n}} \to \iG{Z^{n}/Z^{n+1}} \to 0.
\]
By hypothesis, $H^0_{Z^{n}/Z^{n+1}}\F^j = 0$ for all $j \ne n$. Thus
$\iG{Z^{n+1}}\F^j = \iG{Z^{n}}\F^j$ for $n<j$.

For verifying the other condition of Cousinness,
first assume that $\X$ has finite Krull dimension so that
$\Delta$ is bounded above and hence for $p \gg 0$, $Z^p = \emptyset$.
To show that $H^i_{Z^{j+1}}\F^j =0$ it suffices to show that the
natural maps
\[
H^i_{Z^{j+1}}\F^j \to H^i_{Z^{j+2}}\F^j \to \ldots \to H^i_{\emptyset}\F^j = 0
\]
are isomorphisms. By hypothesis,
$H^i_{Z^{n}/Z^{n+1}}\F^j = 0$ for all $i$ and for $n>j$, whence
the desired conclusion follows.

If $\X$ is not finite-dimensional, we use a localization argument,
cf.~\cite[p.~242]{RD}. For any $x \in \X$ let $\X_{(x)}$ denote
the space of all the generizations of~$x$. In fact
$\X_{(x)} = \Spf(\wh{\OXx})$ where the completion is with
respect to the stalk~$\I_x$ of any defining ideal~$\I$ in~$\OX$.
Then $\{Z^p_x \set Z^p \cap \X_{(x)}\}_{p \in \mathbb{Z}}$
is the induced filtration on~$\X_{(x)}$. Let $f \colon \X_{(x)} \to \X$ be
the canonical inclusion. By~\ref{cor:etale1c} below,
the condition of the Lemma also holds for the complex $f^{-1}\cFb$
on~$\X_{(x)}$. Moreover we have
$(H^i_{Z^{j+1}}\F^j)_x \iso (H^i_{Z^{j+1}_x}f^{-1}\F^j)_x$. But the latter
is~0 since $\X_{(x)}$ has finite Krull dimension. Since
$(H^i_{Z^{j+1}}\F^j)_x = 0$ for every $x \in \X$, therefore
$H^i_{Z^{j+1}}\F^j = 0$.
\end{proof}

The next few results concern the localization argument used in the
proof of~\ref{lem:etale1}. Lemma~\ref{lem:etale1a} below is used
in~\cite[p.~242]{RD} without proof.

\begin{alem}
\label{lem:etale1a}
Let $X$ be a noetherian topological space such that any
nonempty irreducible closed subset $Z$ contains a unique generic
point. Let $Y$ be a subset of~$X$ that is stable under generization.
Let $f\colon Y \to X$ denote the canonical inclusion. Then, for
any flasque sheaf $\F$ on~$X$, its restriction $f^{-1}\F$ on $Y$
is also a flasque sheaf.
\end{alem}
\begin{proof}
First note that the \emph{presheaf} restriction of $\F$ is flasque, i.e.,
if $\G$ denotes the presheaf on~$Y$ that assigns to any open subset
$V$ of $Y$, the set $\dirlm{U}\F(U)$ where~$U$ ranges over
open subsets of~$X$ containing $V$, then for any open subset
$V_1$ of~$Y$, the canonical map $\G(Y) \to \G(V_1)$ is surjective.
Indeed, if $U_1$ is any open subset of~$X$ containing $V_1$, then
in the following commutative diagram of canonical maps
\[
\begin{CD}
\F(X) @>{\rho_1}>> \F(U_1) \\
@V{\rho_2}VV @V{\rho_3}VV \\
\G(Y) @>{\rho_4}>> \G(V_1)
\end{CD}
\]
$\rho_1, \rho_2, \rho_3$ are all surjective and hence $\rho_4$
is also surjective. Therefore it suffices to show that $\G$ is
already a sheaf. For the rest of this proof the flasqueness
assumption on~$\F$ is not needed.

By quasi-compactness, it suffices to verify the sheaf property
for open coverings consisting of finitely many open sets.
Let $V$ be an open subset of $Y$ and
let $\{V_i\}_{i=1}^{n}$ be an
open covering of $V$. Set $V_{ij} \set V_i \cap V_j$. Let $s_i$
be sections of $\G$ over $V_i$ such that $s_i$ and $s_j$
agree on $V_{ij}$. We must show that there is a unique section $s$
of $\G$ over $V$ whose restriction to $V_i$ is $s_i$.
Before proceeding further we remark that for any open subset $U$
of $X$ and any $\alpha \in \F(U)$, if $\beta$ is the induced element
in $\G(U \cap Y)$ then $\Supp(\beta) = \Supp(\alpha) \cap Y$.

By definition, for each $i$, there exists an open subset $U_i$ of~$X$
and an element $t_i \in \F(U_i)$ such that $V_i = U_i \cap V$
and $t_i$ maps to $s_i$ under the canonical map $\F(U_i) \to \G(V_i)$.
Set $U_{ij} \set U_i \cap U_j $ and set
$t_{ij} \set t_i - t_j \in \F(U_{ij})$. Since $t_{ij}$
vanishes on $V_{ij} = U_{ij} \cap V$, therefore
its support $Z_{ij}$ is a closed subset of $U_{ij}$ avoiding~$Y$. Let
$\ov{Z_{ij}}$ denote the closure of~$Z_{ij}$ in~$X$.
We claim that $\ov{Z_{ij}} \cap Y = \emptyset$.

Since $Z_{ij}$ is an open dense subset of $\ov{Z_{ij}}$, therefore $Z_{ij}$
contains all the generic points of the finitely many irreducible
components of~$\ov{Z_{ij}}$. Moreover, $\ov{Z_{ij}}$ is
the set of all the specializations of these generic points. Since
each generic point is outside~$Y$ and since~$Y$ is stable under
generization, it follows that $\ov{Z_{ij}} \cap Y = \emptyset$.

Set $Z \set \cup_{i,j} \ov{Z_{ij}}$ and $U_i' \set U_i \setminus Z$.
Let $t_i'$ be the restriction of~$t_i$ to~$U_i'$.
Then $t_i'$ and $t_j'$ agree on $U_i'\cap U_j'$ and hence there exists
a section $t'$ of $\F$ over $U' \set \cup_i U_i'$ whose restriction
to $U_i'$ is $t_i'$. Let $s$ be the image of $t'$ in $\G(V)$.
Then the restriction of $s$ to $V_i$ equals $s_i$. It remains to
demonstrate the uniqueness of $s$.

Let $\hat{s} \in \G(V)$ be such that its restriction to each $V_i$
is~$s_i$. There is an open subset $U''$ of~$U'$ and an element
$\hat{t} \in \F(U'')$ such that $\hat{t}$ maps to~$\hat{s}$.
Let $t''$ be the restriction of~$t'$ to~$U''$.
Then the support of $t''-\hat{t}$
is a closed subset~$W$ of~$U''$ that avoids~$Y$.
Arguing as above, we see that the closure $\ov{W}$ of~$W$
in~$X$ also avoids~$Y$. Then the restriction of~$\hat{t}$
to the open set $U''\setminus \ov{W}$ agrees with that
of~$t''$ and hence $\hat{s} = s$.
\end{proof}

\begin{alem}
\label{lem:etale1b}
Let $f:Y \to X$ be as in \textup{\ref{lem:etale1a}}.
Let $Z \subset X$ be a subset that is stable under specialization.
Then for any sheaf $\F$ on~$X$,
$f^{-1}\iG{Z}\F = \iG{Z \cap Y}f^{-1}\F$.
\end{alem}

\begin{proof}
The canonical map $f^{-1}\iG{Z}\F \to f^{-1}\F$ is injective and factors
through $\iG{Z \cap Y}f^{-1}\F \hookrightarrow f^{-1}\F$. Therefore it suffices
to check that the natural induced map
$\phi \colon f^{-1}\iG{Z}\F \to \iG{Z \cap Y}f^{-1}\F$ is surjective.

Let $V$ be an open subset of $Y$. Let $s \in (\iG{Z \cap Y}f^{-1}\F)(V)$.
Then there exists an open subset $U$ of~$X$ and an element
$t \in \F(U)$ such that $U \cap Y = V$ and~$t$ maps to~$s$ under the
canonical map $\F(U) \to (f^{-1}\F)(V)$. Then
$\Supp(t) \cap Y = \Supp(s) \subset Z$. Let~$W$ be the set of all
specializations of those generic points of (the irreducible components of)
$\Supp(t)$ that do not lie in~$Y$. Then~$W$ is a closed set avoiding~$Y$.
The restriction of~$t$ to the open set $U \setminus W$ then has support
in~$Z$. Therefore the image of~$t$ in $f^{-1}\iG{Z}\F(V)$ is an element
whose image under~$\phi$ is~$s$.
\end{proof}

\medskip

\begin{acor}
\label{cor:etale1c}
Let $X$ be as in \textup{\ref{lem:etale1a}}. Let $x$ be a point in $X$
and let $Y \set X_{(x)}$ be the set of all generizations of $x$.
Let $f \colon Y \to X$ be the canonical inclusion map. Let $Z_2\subset Z_1$
be subsets of $X$ that are stable under specialization.
Set $Z_i' \set Z_i \cap Y$. Then for any complex $\cFb$
and any integer $i$, there is a canonical isomorphism
$(H^i_{Z_1'/{Z_2'}}f^{-1}\cFb)_x \iso  (H^i_{Z_1/{Z_2}}\cFb)_x$.
\end{acor}

\begin{proof}
Since $f^{-1}$ is exact, it preserves quasi-isomorphisms
and hence we may replace $\cFb$ by a flasque resolution.
By \ref{lem:etale1a}, $f^{-1}\cFb$ consists of flasque sheaves
on~$Y$. Thus we reduce to showing that for any flasque sheaf $\F$ on $\X$,
$(H^i\iG{Z_1'/{Z_2'}}f^{-1}\F)_x \iso  (H^i\iG{Z_1/{Z_2}}\F)_x$.
Since $H^i$ commutes with localization and $f^{-1}$ preserves
stalks, we reduce to showing that
$\iG{Z_1'/{Z_2'}}f^{-1}\F = f^{-1}\iG{Z_1/{Z_2}}\F$.
Therefore we conclude from the following commutative diagram
where the rows are exact and the two leftmost columns are equalities
by \ref{lem:etale1b}.
\[
\begin{CD}
0 @>>> \iG{Z_2'}f^{-1}\F @>>> \iG{Z_1'}f^{-1}\F @>>>
 \iG{Z_1'/{Z_2'}}f^{-1}\F @>>> 0 \\
@. @| @| @VVV \\
0 @>>> f^{-1}\iG{Z_2}\F @>>> f^{-1}\iG{Z_1}\F @>>>
 f^{-1}\iG{Z_1/Z_2}\F @>>> 0 \\
\end{CD}
\]
\end{proof}

\medskip

\begin{aprop}
\label{prop:etale2}
Let $f \colon (\X,\DsssX) \to (\Y,\DsssY)$ be a map in $\bbFc$
such that for any $x \in \X$, with $y \set f(x)$, $\OXx/m_y\OXx$
is a finite-dimensional vector space over the residue field $k(y)$
at $y$ (in particular, $\DsssX(x) = \DsssY(y)$).
Then $f$ has finite fibers. Moreover, for any
$\cMb \in \Coz_{\DsssY}(\Y)$, the following hold.
\begin{enumerate}
\item $\iGp{\X}f^*\cMb \in \Coz_{\DsssX}(\X)$.
\item Let $y \in \Y$. Set $M \set \cMb(y)$. If $f^{-1}\{y\} = \emptyset$, then
$\iGp{\X}f^*i_yM = 0$. Otherwise, if $f^{-1}\{y\} = \{x_1, \ldots, x_n\}$,
then, 
with $p = \DsssX(x_j) = \DsssY(y)$, the $\DsssX$-Cousin complex
$\iGp{\X}f^*i_yM[-p]$ is only concentrated at the points $x_j$
($1 \le j \le n$), and for any $j$
there is a natural isomorphism
\[
(\iGp{\X}f^*\cMb)(x_j) = (\iGp{\X}f^*i_yM[-p])(x_j)
\iso M \otimes_y \cO_{\X,x_j}
\]
uniquely determined by the natural map
\[
\iG{x_j}\iGp{\X}f^*\M^p = \iG{x_j}\iGp{\X}f^*i_yM \xto{\quad}
(f^*i_yM)_{x_j} \iso M \otimes_y \cO_{\X,x_j}.
\]
\end{enumerate}
\end{aprop}
\begin{proof}
Let $x \in \X, y \in \Y$ with $y =f(x)$. Then the local ring
of the fiber space~$\X_y$ at~$x$ is isomorphic to~$\OXx/m_y\OXx$
and hence, by hypothesis, has Krull dimension zero.
Therefore~$x$ does not have any nontrivial
generization in~$\X_y$. Since this holds for any element $x$
in~$\X_y$, it follows that~$\X_y$ is zero-dimensional and
hence consists of finitely many points.

Note that $\iGp{\X}f^*\cMb$ consists
of $\Aqct(\X)$-modules. For any integer $j$
in view of the canonical decomposition
$\M^j \iso \oplus_yi_y\cMb(y)$ (where $\DsssY(y) = j$),
we deduce, using the same arguments as in~\ref{lem:psm1} that
for any $x \in \X$, if $j \ne \DsssY(f(x))= \DsssX(x)$, then
$H^i_x(\iGp{\X}f^*\M^j) = 0$ for all $i$.

For this paragraph, fix $x \in \X$.
Set $y \set f(x)$, $p \set \DsssX(x) = \DsssY(y)$, $M = \cMb(y)$.
Then there are natural isomorphisms
\[
H^i_x(\iGp{\X}f^*\M^p) \iso H^i_{m_x}(f^*\M^p)_x
\iso H^i_{m_x}(\M^p_y \otimes_y \OXx) \iso H^i_{m_x}(M \otimes_y \OXx).
\]
By hypothesis, $m_y\OXx$ is $m_x$-primary and hence
$m_y$-torsionness of $M$ implies that $M \otimes_y \OXx$ is $m_x$-torsion.
Thus $H^0_{m_x}(M \otimes_y \OXx) = M \otimes_y \OXx$ and
moreover $H^i_{m_x}(M \otimes_y \OXx) = 0$ for $i \ne 0$. In particular,
$H^i_x(\iGp{\X}f^*\M^p) = 0$ for $i \ne 0$.

By \ref{lem:etale1} it now follows that
$\iGp{\X}f^*\cMb \in \Coz_{\DsssX}(\X)$ and moreover for any $x \in \X$
we obtain a canonical isomorphism
$(\iGp{\X}f^*\cMb)(x) \iso \cMb(f(x)) \otimes_{f(x)} \OXx$.

For (ii), with $y \in \Y$ fixed, we may assume that $\cMb$
is concentrated only at~$y$.
It follows from the above description that if $f^{-1}\{y\} = \emptyset$,
then for all $x \in \X$, $(\iGp{\X}f^*\cMb)(x) = 0$ and hence
$\iGp{\X}f^*\cMb = 0$. Also, if $f^{-1}\{y\} = \{x_1, \ldots, x_n\}$
then $\iGp{\X}f^*\cMb$ is only concentrated at the points $x_j$.
By chasing through the isomorphisms in the previous paragraphs
for $i = 0$, we obtain the desired isomorphism in~(ii).
\end{proof}

Let us note that an \'etale map satisfies the hypothesis on~$f$
in \ref{prop:etale2}. Indeed, if $f$ is a smooth map of relative
dimension $d$, then $d$ is the sum of the
Krull dimension of~$(\OXx/m_y\OXx)$ and the transcendence degree
of the residue field extension $k(y) \to k(x)$ (cf.~\ref{cor:app4}).
Therefore $d=0$ implies that $m_y\OXx$ is $m_x$-primary and
that $k(x)$ is a finite extension of $k(y)$. Thus $(\OXx/m_y\OXx)$
has finite length over $k(y)$.

We are now in a position to describe $f^{\sharp}$ for $f$ an \'etale map.
Note that if $f$ is an \'etale map, then $\omega_f = \OX$. Moreover, as
in the proof of \ref{rem:fin-clim2},
we shall identify the functor $(-)\otimes \OX[0]$ with the
identity functor.

\begin{aprop}
\label{prop:etale3}
Let $f \colon (\X,\DsssX) \to (\Y,\DsssY)$ be an \'etale map in $\bbFc$.
Let $\cMb \in \Coz_{\DsssY}(\Y)$. Then $\iGp{\X}f^*\cMb$ is canonically
isomorphic to $f^{\sharp}\cMb$. More precisely, with $E = E_{\DsssX}$,
the following canonical maps are isomorphisms
\begin{equation}
\label{eq:etale4}
\iGp{\X}f^*\cMb \xto{\textup{\ref{lem:coz1a}(ii)}}
E\iGp{\X}f^*\cMb \to E\R\iGp{\X}f^*\cMb
\xto{\eqref{eq:fin-sm1}} f^{\sharp}\cMb,
\end{equation}
and for any $x \in \X$, with $y \set f(x), M \set \cMb(y)$,
the corresponding punctual isomorphism is given by
\begin{equation}
\label{eq:etale5}
(\iGp{\X}f^*\cMb)(x) \xto{\textup{\ref{prop:etale2}(ii)}} M \otimes_y \OXx
\iso M \otimes_{\hat{y}} \wh{\OXx} \xto{\textup{\ref{thm:huang1}, I.(i)}}
\wh{f_x}_{\sharp}M.
\end{equation}
\end{aprop}

\begin{proof}
Let $p = \DsssX(x)$.
We relate \eqref{eq:etale4} with~\eqref{eq:etale5} through a diagram, which
for convenience, is broken into two parts, viz., \eqref{cd:etale6}
and~\eqref{cd:etale7} below. The rightmost column of~\eqref{cd:etale6}
is the same as the leftmost one in~\eqref{cd:etale7}. The top row of
\eqref{cd:etale6} + \eqref{cd:etale7} gives~\eqref{eq:etale5} while
the leftmost column of~\eqref{cd:etale6} gives the first two maps
in~\eqref{eq:etale4}. The remaining portion of the outer border of
\eqref{cd:etale6} + \eqref{cd:etale7} gives~{\eqref{eq:fin-sm1}},
which is the third map used in~\eqref{eq:etale4}. Note that in this
situation the sign obtained from the map $\theta$ occurring
in~{\eqref{eq:fin-sm1}} is trivial because now $d=0$ and $p = q$.

\begin{figure}\vspace{5pt}
\rotatebox{90}
{\begin{minipage}{8.5in}\vspace{-3mm}
\begin{equation}
\label{cd:etale6}
\begin{CD}
\iG{x}(\iGp{\X}f^*\cMb)^p @.
 \makebox[0pt]{\hspace{12em}{\raisebox{.5ex}{\makebox[20em]{\hrulefill}}}
 \hspace{-20em}\!\!{\makebox[20em]{\hrulefill}}}
 @. @. \iG{x}\iGp{\X}f^*\M^p \\
@V{\text{\ref{lem:coz1a}(ii)}}V{\makebox[0pt]{\hspace{32em}$\Box_1$}}V
 @. @. @|  \\
\iG{x}(E\iGp{\X}f^*\cMb)^p @>>> H^p_x\iGp{\X}f^*\cMb
 @>{\text{truncation}}>{\text{on $\cMb$}}>
 H^p_x(\iGp{\X}f^*\M^p[-p]) @= H^0_x\iGp{\X}f^*\M^p \\
@VVV @VVV @VVV @VVV \\
\iG{x}(E\R\iGp{\X}f^*\cMb)^p @>>> H^p_x\R\iGp{\X}f^*\cMb
 @>{\text{truncation}}>{\text{on $\cMb$}}> H^p_x(\R\iGp{\X}f^*\M^p[-p])
 @= H^0_x\R\iGp{\X}f^*\M^p
\end{CD}
\end{equation}
\bigskip \bigskip \bigskip\bigskip\bigskip
\begin{equation}
\label{cd:etale7}
\begin{CD}
\iG{x}\iGp{\X}f^*\M^p @>>> (\iGp{\X}f^*\M^p)_x @>>>
 (f^*\M^p)_x @>>> M \otimes_y \OXx @>>> \wh{f_x}_{\sharp}M \\
@| @AAA @AAA @AAA \\
H^0_x\iGp{\X}f^*\M^p @=
 (H^0_{\ov{\{x\}}}\iGp{\X}f^*\M^p)_x  @>{\ref{cor:loc4}}>> \iG{m_x}(f^*\M^p)_x
@>>> \iG{m_x}(M \otimes_y \OXx) \\
@VVV @VVV @| \\
H^0_x\R\iGp{\X}f^*\M^p @= (H^0_{\ov{\{x\}}}\R\iGp{\X}f^*\M^p)_x
 @>{\eqref{eq:nolab3}}>>H^0_{m_x}(f^*\M^p)_x
\end{CD}
\end{equation}
\end{minipage}\hspace{-14mm}
}
\end{figure}

The commutativity of all the subrectangles in \eqref{cd:etale6} and
\eqref{cd:etale7} except $\Box_1$ are verified easily. For $\Box_1$, set
$\cCb \set \iGp{\X}f^*\cMb$. We may then rewrite $\Box_1$ as the following.
\[
\begin{CD}
\iG{x}\C^p @. {\makebox[0pt]{{\raisebox{.5ex}{\makebox[10em]{\hrulefill}}}
 \hspace{-10em}\!\!{\makebox[10em]{\hrulefill}}}} @. \iG{x}\C^p\\
@VVV @. @| \\
\iG{x}(E\cCb)^p @>>> H^p_x\cCb @>>> H^p_x(\C^p[-p])\\
\end{CD}
\]
For commutativity of $\Box_1$ we now refer to the proof
of \ref{lem:coz1a}(iii). Note that every map in $\Box_1$
is an isomorphism.

All the horizontal maps in \eqref{cd:etale6} are isomorphisms
and all the maps (horizontal as well as vertical) in \eqref{cd:etale7}
are isomorphisms. Therefore the vertical maps in \eqref{cd:etale6},
in particular the ones in the leftmost column are also isomorphisms.
In particular the middle maps in \eqref{eq:etale4} is also an isomorphism.
\end{proof}

\begin{aprop}
\label{prop:etale8}
Let $f \colon (\X,\DsssX) \to (\Y,\DsssY)$ be an etale map in $\bbFc$.
Let $\cMb \in \Coz_{\DsssY}(\Y)$. Set $\cFb \set \R\iGp{\X}f^*\cMb$.
Then the following diagram commutes in~$\D(\X)$ and moreover
all the maps in it are $\D(\X)$-isomorphisms.
\[
\begin{CD}
\iGp{\X}f^*\cMb @>>> \R\iGp{\X}f^*\cMb \\
@VV{\textup{\ref{lem:coz1a}(ii)}}V @VV{\textup{\ref{cor:coz4a}}}V \\
E\iGp{\X}f^*\cMb @>>> E\R\iGp{\X}f^*\cMb
\end{CD}
\]
Thus, the Suominen isomorphism
$\R\iGp{\X}f^*\cMb \iso E\R\iGp{\X}f^*\cMb$ is the composition
of the following (explicit) isomorphisms
\[
\R\iGp{\X}f^*\cMb \xto{(\text{canonical})^{-1}} \iGp{\X}f^*\cMb
\xto{\textup{\ref{lem:coz1a}(ii)}}
E\iGp{\X}f^*\cMb \xto{\text{canonical}} E\R\iGp{\X}f^*\cMb.
\]
\end{aprop}

\begin{proof}
The vertical maps in the given diagram are $\D(\X)$-isomorphisms
and by \ref{prop:etale3}, the bottom row is an isomorphism.
Therefore only commutativity of the above diagram needs to be shown.
The complexes involved in the above diagram are all $\DsssX$-CM in
$\D(\X)$. Let $\ov{\bC}, \ov{\D}$ be as in \S\ref{subsec:cm}.
Since $E = E_{\DsssX}$ is fully faithful as a functor from
$\ov{\D}$ to $\ov{\bC}$ (\ref{prop:coz4}), it suffices
to verify that $E$ applied to the above diagram gives a commutative
diagram.

Set $\cCb \set \iGp{\X}f^*\cMb \in \ov{\bC}$ and
$\cNb \set \R\iGp{\X}f^*\cMb \in \ov{\D}$.
There is a canonical map $Q\cCb \to \cNb$ and $E$ of the above
diagram may be written as the following one.
\[
\begin{CD}
EQ\cCb @>>> E\cNb \\
@V{EQ(\cCb \iso EQ\cCb)}VV @VV{E(\cNb \iso QE\cNb)}V \\
EQEQ\cCb @>>> EQE\cNb
\end{CD}
\]
The vertical map on the left is also the same as
$(EQ\cCb) \iso EQ(EQ\cCb)$ (see proof of \ref{cor:coz4a})
while the vertical map on the right is, by the defining property
of the Suominen isomorphism, the same as
$(E\cNb) \iso EQ(E\cNb)$. Thus the preceding
diagram commutes for functorial reasons.
\end{proof}

Completion of a formal scheme along an open ideal is an etale
map. In this case we obtain further concrete
descriptions of $(-)^{\sharp}$. We begin with a preparatory Lemma.

\pagebreak[3]

\begin{alem}
\label{lem:etale9}
Let $(\Z, \DsssZ) \in \bbFc$. Let $\I$ be an open coherent ideal in $\OZ$.
Let $\kappa \colon (\X, \DsssX) \to (\Z, \DsssZ)$ be the completion map
corresponding to the completion of $\Z$ along $\I$.
Let $\cFb$ be a complex of $\Aqct(\Z)$-modules such that for any
integer $j$, $\F^j$ is a direct sum of modules of the type
$i_zF_z$ for $z \in \Z$, where $F_z$ is a zero-dimensional $\OZz$-module.
Set $\kappa^{!} \set \kappa^{-1}\iG{\I}$.
Then the following hold.
\begin{enumerate}
\item The canonical maps $\kappa^{-1}\iG{\I}\cFb \to
\kappa^{*}\iG{\I}\cFb \to \iGp{\X}\kappa^{*}\cFb$ are isomorphisms
and for any integer $j$, $\kappa^{-1}\iG{\I}\F^j$ is a direct sum
of modules of the type $i_xF_z$ for $x \in \X$, where $z = \kappa(x)$
and $F_z$ the corresponding $\OZz$-module.
\item Suppose $\cFb\in \Coz_{\DsssZ}(\Z)$. Then
$\kappa^{!}\cFb \in \Coz_{\DsssX}(\X)$ and
for any $x \in \X$, with $z = \kappa(x), F = \cFb(z)$, we have
$(\kappa^{!}\cFb)(x) = (\kappa^{!}i_zF)(x) = F$.
\item Suppose $\cFb$ consists of $\Aqct(\Z)$-injectives. Then
$\kappa^{!}\cFb$ consists of $\Aqct(\X)$-injectives
and the canonical map
$\iGp{\X}\kappa^{*}\cFb \to \R\iGp{\X}\kappa^{*}\cFb$
is an isomorphism.
\end{enumerate}
\end{alem}
\begin{proof}
That $\kappa^{-1}\iG{\I}\cFb \to \kappa^{*}\iG{\I}\cFb$
is an isomorphism follows immediately from the fact that
$\iG{\I}\cFb$ consists of $\I$-torsion modules.
At any rate, for (i), since
$\kappa^{-1}, \kappa^{*}, \iG{\I}, \iGp{\X}$ commute
with direct sums, we may replace $\cFb$ by
a single sheaf $\F \set i_zF$.
Since $\F \in \Aqct(\Z)$ therefore
$\iG{\I}\F = \iG{|\X|}\F$ where we identify
the underlying topological space $|\X|$ of $\X$ by its image
in $\Z$, viz., $\Supp(\OZ/\I)$.
If $z \in |\X|$, then $\iG{|\X|}(i_zF) = i_zF$, and
if $z\notin|\X|$, then $\iG{|\X|}(i_zF) = 0$.
In the former case, if~$x$ is the unique point mapping to~$z$,
then $\kappa^{-1}i_zF = i_xF$. Moreover,
since $\iG{\I}i_zF =i_zF$ is an $\I$-torsion module, therefore
the natural maps
$\kappa^{-1}i_zF \to \kappa^{*}i_zF \gets \iGp{\X}\kappa^{*}i_zF$
are isomorphisms. Thus~(i) and~(ii) follow.

For (iii), by (i) above and \ref{prop:mod7},
it only remains to be verified that
the natural map $\iGp{\X}\kappa^{*}\cFb \to \R\iGp{\X}\kappa^{*}\cFb$
is an isomorphism. Since it suffices to prove this locally,
we may assume that $\I\OX$ is generated by global sections.
Then we conclude by~\ref{lem:loc5a}.
\end{proof}

\begin{aprop}
\label{prop:etale10}
Let $\kappa, \X, \DsssX, \Z, \DsssZ$ and $\I$ be as in
\textup{\ref{lem:etale9}}.
Then for any $\cMb \in \Coz_{\DsssZ}(\Z)$,
the canonical graded isomorphism
$\kappa^{-1}\iG{\I}\cMb \iso \kappa^{\sharp}\cMb$, given at
the punctual level, with $x \in \X, z = \kappa(x), M=\cMb(z)$, by
\[
(\kappa^{-1}\iG{\I}\cMb)(x) = M \iso M \otimes_z \OXx
\xrightarrow[\textup{and \ref{thm:huang1}, I.(i)}]{\textup{completions} }
\wh{\kappa_x}_{\sharp}M,
\]
is an isomorphism of complexes.
\end{aprop}
\begin{proof}
It suffices to prove that the outer border of
the following diagram commutes
\[
\begin{CD}
(\kappa^{-1}\iG{\I}\cMb)(x) @>{\text{\ref{lem:etale9}(i)}}>>
 (\iGp{\X}\kappa^{*}\cMb)(x)
 @= (\iGp{\X}\kappa^{*}\cMb)(x)  \\
@| @VV{\text{\ref{prop:etale2}(ii)}}V  @VV{\eqref{eq:etale5}}V\\
M @>>> M \otimes_z \OXx @>>> \wh{\kappa_x}_{\sharp}M
\end{CD}
\]
for then the Proposition follows from \ref{prop:etale3}
because the isomorphisms in~{\ref{lem:etale9}(i)} and~\eqref{eq:etale4}
are maps of complexes.
The square on the right commutes by definition of~{\eqref{eq:etale5}}.
For the square on the left we expand as follows
and the required commutativity is verified easily.
\[
\begin{CD}
(\kappa^{-1}\iG{\I}\cMb)(x) @>>> (\kappa^{*}\iG{\I}\cMb)(x)
 @>>> (\iGp{\X}\kappa^{*}\cMb)(x) \\
@VVV @VVV @VVV \\
(\kappa^{-1}\iG{\I}i_zM)_x @>>> (\kappa^{*}\iG{\I}i_zM)_x
 @>>> (\iGp{\X}\kappa^{*}i_zM)_x \\
@VVV @VVV @VVV \\
(\kappa^{-1}i_zM)_x @>>> (\kappa^{*}i_zM)_x @= (\kappa^{*}i_zM)_x \\
@VVV @VVV @VVV \\
M @>>> M \otimes_z \OXx @= M \otimes_z \OXx
\end{CD}
\]
\end{proof}

\begin{arem}
\label{rem:etale11}
In the situation of \ref{prop:etale10}, for any $z \in \Z$, with
$x = \kappa(z)$, the completions of the local rings $\OZz$ and $\OXx$
along their maximal ideals are canonically isomorphic to each other.
Therefore, if we make the identification $\wh{\OZz} = \wh{\OXx}$,
then we find that
\[
1_{\Z}^{\sharp}\kappa^{-1}\iG{\I}\cMb =
\kappa^{-1}\iG{\I}1_{\X}^{\sharp}\cMb =
\kappa^{\sharp}\cMb,
\qquad
\kappa_*\kappa^{\sharp}\cMb = \iG{\I}1_{\X}^{\sharp}\cMb.
\]
\end{arem}

We conclude this subsection with some remarks concerning the Suominen
isomorphism $\Ssf$. We have already seen two instances, viz.,
\ref{prop:QE13} and \ref{prop:etale8} where an explicit form
for~$\Ssf$ is obtained. In fact we can generalize further.
We shall use the principle,
already used in the proof of \ref{prop:etale8}, which is
that for any CM complex $\cNb$, to find an explicit candidate
for $\Ssf(\cNb)$, it suffices to find an explicit isomorphism
from $\cNb$ to a Cousin complex $\cCb$, because, $\Ssf(\cCb)$
is represented by the isomorphism of \ref{lem:coz1a}(ii).
We elaborate below.

Let $(\X, \DsssX) \xto{\;f_1\:} (\Z, \DsssZ) \xto{\;f_2\:} (\Y, \DsssY)$
be maps in $\bbFc$ where $f_1$ is an etale map and $f_2$ is an
adic smooth map of relative dimension $d$.
Let $\cMb \in \Coz_{\DsssY}(\Y)$. Set $f \set f_2f_1$.
Let us obtain an explicit description of
$\Ssf(\R\iGp{\X}f^{\diamond}\cMb)$.

There are canonical isomorphisms
$f^{\diamond}\cMb \iso f_1^{\diamond}f_2^{\diamond}\cMb
= f_1^*f_2^{\diamond}\cMb$ and hence there are
natural explicit $\D(\X)$-isomorphisms
\[
\R\iGp{\X}f^{\diamond}\cMb \iso \R\iGp{\X}f_1^*f_2^{\diamond}\cMb
\xrightarrow[\eqref{eq:QE1ac}]{\sim}
\R\iGp{\X}f_1^*E_{\DsssZ}f_2^{\diamond}\cMb
\xleftarrow[\ref{prop:etale8}]{\sim}
\iGp{\X}f_1^*E_{\DsssZ}f_2^{\diamond}\cMb.
\]
Since $\cCb \set \iGp{\X}f_1^*E_{\DsssZ}f_2^{\diamond}$ is a
Cousin complex, therefore, via the above isomorphism,
$\Ssf(\R\iGp{\X}f^{\diamond}\cMb)$
is now given by the isomorphism of \ref{lem:coz1a}(ii).

Now assume further that $f_1 = \kappa$ is the completion map as
in~\ref{prop:etale10}. Choose an $\Aqct(\Z)$-injective resolution
$f_2^{\diamond}\cMb \to \cIb$. Since $\kappa$ is flat the induced
map $\kappa^* f_2^{\diamond}\cMb \to \kappa^*\cIb$ is also a quasi-isomorphism.
There are natural isomorphisms
\[
\R\iGp{\X}f^{\diamond}\cMb \iso \R\iGp{\X}\kappa^*\cIb
\xrightarrow[\text{and \ref{lem:etale9}(i)}]{\text{\ref{lem:etale9}(iii)}}
\kappa^{-1}\iG{\I}\cIb \osi \kappa^{-1}\iG{\I}E_{\DsssZ}f_2^{\diamond}\cMb
\]
where the last map is obtained by applying $\kappa^{-1}\iG{\I}$
to the homotopy-unique quasi-isomorphism
$Ef_2^{\diamond}\cMb \to \cIb$ induced by the quasi-isomorphism
$\eta_{f_2}(\cMb)$ (\ref{prop:QE13}). Now with
$\cCb \set \kappa^{-1}\iG{\I}E_{\DsssZ}f_2^{\diamond}\cMb$,
(see \ref{lem:etale9}(ii))
we see that the isomorphism in~\ref{lem:coz1a}(ii)
represents $\Ssf(\R\iGp{\X}f^{\diamond}\cMb)$.


The above results can in fact be applied to any smooth map,
albeit locally. Indeed, any smooth map $f \colon \X \to \Y$
in $\mathbb F$ can be factored locally on $\X$
as $f\big|_{\U} = f_2f_1$ where $f_2\colon \Z \to \Y$ is an
adic smooth map and $f_1\colon \U \to \Z$ is an \'etale map.
This can be shown by an argument
involving locally choosing a differential basis of~$\Ohm^1_f$.
A recent result of Alonso Tarr\'{\i}o, \index{Alonso Tarr\'{\i}o, Leovigildo} 
Jerem\'{\i}as L\'{o}pez\index{Jerem\'{\i}as L\'{o}pez, Ana} and
P\'{e}rez Rodr\'{\i}guez \cite{AJP}\index{P\'{e}rez Rodr\'{\i}guez, Marta}
says that~$f_1$ can in fact
be chosen to be a completion map such as $\kappa$ above.
(However their result is stated only for formally smooth maps
which are of pseudo-finite type; it is possible though, that
their result holds more generally for essentially pseudo-finite type
too.) Thus all the descriptions given above may be applied locally.

\vspace {10pt}
\enlargethispage*{5pt}

\providecommand{\bysame}{\leavevmode\hbox to3em{\hrulefill}\thinspace}
\providecommand{\MR}{\relax\ifhmode\unskip\space\fi MR }
\providecommand{\MRhref}[2]{%
  \href{http://www.ams.org/mathscinet-getitem?mr=#1}{#2}
}
\providecommand{\href}[2]{#2}

\setcounter{footnote}{0}
\renewcommand{\thefootnote}{\fnsymbol{footnote}}

\newpage

\centerline{\textsc{Numbered displays}}
\begin{small}
\begin{table}[h]

\begin{tabular}{lcc|clcc|ccc}
\underbar{Display} &\underbar{Page}&\quad\qquad  &\quad\qquad &\underbar{Display} &\underbar{Page} &\quad\qquad  &\quad\qquad &\underbar{Display}&\underbar{Page}\\[1ex]
\ \ (1)&    9 &&  &\ \ (35a)&  57 &&	&(59)&     77\\
\ \ (2)&    9 &&  &\ \ (35b)&  57 &&	&(60)&     77\\
\ \ (3)&    9 &&  &\ \ (35c)&  57 &&	&(61)&     79\\
\ \ (4)&    24 && &\ \ (35d)&  57 &&	&(62)&     79\\
\ \ (5)&    27 && &\ \ (36)&   58 &&	&(63)&     79\\
\ \ (6)&    36 && &\ \ (36a)&  60 &&	&(64)&     80\\
\ \ (7)&    36 && &\ \ (36b)&  60 &&	&(65)&     80\\
\ \ (8)&    37 && &\ \ (37)&   62 &&	&(66)&     81\\
\ \ (8$'$)& 37 && &\ \ (38)&   63 &&	&(67)&     82\\
\ \ (9)&    38 && &\ \ (38a)&  63 &&	&(68)&     83\\
\ \ (10)&   38 && &\ \ (38b)&  63 &&	&(69)&     84\\
\ \ (11)&   40 && &\ \ (38c)&  63 &&	&(70)&     85\\
\ \ (12)&   40 && &\ \ (38d)&  63 &&	&(71)&     85\\
\ \ (13)&   41 && &\ \ (38e)&  63 &&	&(72)&     85\\
\ \ (14)&   41 && &\ \ (38b$'$)&   64 &&	&(73)&     86\\
\ \ (15)&   46 && &\ \ (38c$'$)&   64 &&	&(74)&     87\\
\ \ (16)&   46 && &\ \ (39)&   65 &&	&(75)&     89\\
\ \ (17)&   47 && &\ \ (39a)&  65 &&	&(76)&     91\\
\ \ (18)&   48 && &\ \ (39b)&  65 &&	&(77)&     92\\
\ \ (19)&   49 && &\ \ (39c)&  65 &&	&(78)&     93\\
\ \ (20)&   49 && &\ \ (39d)&  65 &&	&(79)&     94\\
\ \ (21)&   49 && &\ \ (40)&   65 &&	&(80)&     95\\
\ \ (22)&   50 && &\ \ (41)&   66 &&	&(81)&     97\\
\ \ (23)&   50 && &\ \ (42)&   67 &&	&(82)&     98\\
\ \ ($\wh{23}$)&   50 && &\ \ (43)&   68 &&	&(83)&     97\\
\ \ (24)&   51 && &\ \ (44)&   68 &&	&(84)&     100\\
\ \ ($\wh{24}$)&   51 && &\ \ (45)&   68 &&	&(85)&     102\\
\ \ (25)&   51 && &\ \ (46)&   68 &&	&(86)&     102\\
\ \ ($\wh{25}$)&   51 && &\ \ (47)&   68 &&	&(87)&     114\\
\ \ (25$'$)& 52 && &\ \ (48)&  70 &&	&(88)&     115\\
\ \ (26)&   51 && &\ \ (49)&   70 &&	&(89)&     115\\
\ \ (27)&   54 && &\ \ (50)&   73 &&	&(90)&     115\\
\ \ (28)&   54 && &\ \ (51)&   73 &&	&(91)&     120\\
\ \ (29)&   54 && &\ \ (52)&   73 &&	&(92)&     123\\
\ \ (30)&   54 && &\ \ (53)&   74 &&	&(93)&     124\\
\ \ (31)&   55 && &\ \ (54)&   74 &&    &(94)&     127\\
\ \ (32)&   56 && &\ \ (55)&   75 &&    &(95)&     127\\
\ \ (33)&   56 && &\ \ (56)&   75 &&    &(96)&     128\\
\ \ (34)&   56 && &\ \ (57)&   76 &&    &(97)&     128\\
\ \ (35)&   56 && &\ \ (58)&   76 &&                &&\\
\end{tabular}

\end{table}
\end{small}



\begin{thebibliography}{10}

\bibitem{AJL1}
L.\,Alonso Tarr\'{\i}o, A.\,Jerem\'{\i}as L\'{o}pez, J.\,Lipman, 
\emph{Local homology and cohomology on schemes,} 
Ann.~Sci.~\'{E}cole~Norm.~Sup.~
  \textbf{{\bf 30}} (1997), no.\,1, 1--39.

\bibitem{AJL2}
\bysame, \emph{Duality and flat base change on formal schemes,} 
Contemporary Mathematics {\bf 244}, Amer.~Math.~Soc., Providence, RI,
1999, 3--90.

\bibitem{AJL3}
\bysame, \emph{\ Greenlees-May duality on formal schemes,}\kern1em
Contemporary Mathematics {\bf 244}, Amer.~Math.~Soc., Providence, RI,
1999, 93--112.  



\bibitem{AJP}
L.\,Alonso Tarr\'{\i}o, A.\,Jerem\'{\i}as L\'{o}pez, 
M.\,P\'{e}rez Rodr\'{\i}guez,  \emph{Infinitesimal geometry of formal schemes,}
preprint.


\bibitem{BCo}
B.\,Conrad, \emph{Grothendieck duality and base change,} Lecture Notes in
  Mathematics, {\bf 1750}, Springer, Berlin, 2000.

\bibitem{Db}
M.\,T.\,Dibaei, \emph{A study of Cousin complexes through the dualizing complexes,}
preprint, math.AC/0407187.


\bibitem{DT}
\bysame, M.\,Tousi,
\emph{A generalization of the dualizing complex structure and its 
applications,} J.\ Pure and Applied Algebra {\bf 155} (2001), 17--28.
 
\bibitem{EGAOIV}
A.\,Grothendieck, J.\,Dieudonn\'{e}, 
\emph{\'{E}lements de g\'{e}om\'{e}trie alg\'{e}brique,} 
{{\bf 0$_{\textup{IV}}$}},  
Publications Math. I.H.E.S {\bf 20}, (1964).%
\footnote[2]{Downloadable from
{\tt <http://www.numdam.org/numdam-bin/feuilleter?j=PMIHES\&sl=0>}.}


\bibitem{EGAIV} \bysame,
 \emph{\'El\'ements de  G\'eom\'etrie Alg\'ebrique IV,
 \'Etude locale des sch\'emas et des morphismes de sch\'emas.}
 Publications Math. I.H.E.S. {\bf 24} (1965), {\bf 28} (1966),
 {\bf 32} (1967).\footnotemark[2]


\bibitem{EGAI}
\bysame, \emph{\'{E}lements de g\'{e}om\'{e}trie alg\'{e}brique, {{\bf I}}},
  Springer, Berlin, 1971.

\bibitem{RD}
R.\,Hartshorne, \emph{Residues and duality,} Lecture Notes in Math., {\bf
  20}, Springer, Berlin, 1966.

\bibitem{Ha}
\bysame, \emph{Algebraic geometry,} Graduate Texts in Mathematics, 
No.\,{\bf 52}, Springer, 
Berlin, 1977.


\bibitem{Hu}
I.\,C. Huang, \emph{Pseudofunctors on modules with zero-dimensional support,}
  Memoirs Amer.\ Math.~Soc. {\bf 114} (1995), no.\,{\bf 548}.

\bibitem{Hu2}
\bysame, \emph{An explicit construction of residual complexes,} J.~Algebra
  \textbf{{\bf 225}} (2000), no.\,2, 698--739.

\bibitem{Hubl}
R.\,H\"{u}bl, \emph{Graded duality and generalized fractions,} J.~Pure Appl.~
  Algebra \textbf{{\bf 141}} (1999), no.\,3, 225--247.
\bibitem{Kw}
T.\,Kawasaki, \emph{Finiteness of Cousin cohomologies,} preprint,
\hfill\newline
{\tt <\ http://www.comp.metro-u.ac.jp/~kawasaki/articles.html\#preprint\ >}
\bibitem{Li}
J.\,Lipman, \emph{Notes on derived categories,}  
  {\tt<\;http://www.math.purdue.edu/\~{}lipman/\;>}.
  

\bibitem{Lu}
W.\,L\"{u}tkebohmert, \emph{On compactification of schemes,} Manuscr.~ Math.~
  \textbf{{\bf 80}} (1993), no.\,1, 95--111.

\bibitem{Ma}
H.\,Matsumura, \emph{Commutative ring theory,} Cambridge University Press, 
Cambridge, 1986.

\bibitem{Nag}
M.\,Nagata, \emph{Imbedding of an abstract variety in a complete variety,} J.~
  Math.~Kyoto Univ.~\textbf{{\bf 2}} (1962), 1--10.

\bibitem{Nag2}
\bysame, \emph{Local Rings,} John Wiley \& Sons, New York, 1962.

\bibitem{Nee}
A.\,Neeman, \emph{The Grothendieck duality theorem via 
Bousfield's techniques and Brown representability,}
J.~Amer.~Math.~Soc.~\textbf{{\bf 9}} (1996), 205--236.


\bibitem{Sa}
P.\,Sastry, \emph{Residues and duality for algebraic schemes,} 
Compositio Math.~\textbf{{\bf 101}} (1996), 133--178.

\bibitem{Sa2}
\bysame \emph{Duality for Cousin complexes,} this volume.

  
\bibitem{Suo}
K.\,Suominen, \emph{Localization of sheaves and {C}ousin complexes,} Acta
Math.~\textbf{{\bf 131}} (1973), 27--41.

\bibitem{Ve}
J.\,L. Verdier, \emph{Base change for twisted inverse image of coherent
  sheaves,} Algebraic Geometry, 
Oxford Univ.~Press, London, 1969, pp.\,393--408.

\bibitem{Ye}
A.\,Yekutieli, \emph{Smooth formal embeddings and the residue
  complex,} Canadian
  J.~Math.~\textbf{{\bf 50}} (1998), no.\,4, 863--896.

\bibitem{YZ}
\bysame,  J.\,J.\,Zhang,
\emph{Rigid dualizing complexes and perverse sheaves on schemes,}
preprint, math.AG/0405570.
\end{thebibliography}
\end{document}